\documentclass[12pt,reqno]{amsart}
\usepackage{amssymb}
\usepackage[all]{xy}
\title[Compactification of moduli schemes]
{Compactifications of smooth families 
and of moduli spaces of polarized manifolds}
\author[Eckart Viehweg]{Eckart Viehweg}
\address{Universit\"at Duisburg-Essen, Mathematik, 45117 Essen, Germany}
\email{viehweg@uni-essen.de}
\thanks{This work has been supported by the DFG-SP
``Globale Methoden in der Komplexen Geometrie'', by the DFG-Leibniz program, 
by the SFB/TR 45 ``Periods, moduli spaces and arithmetic of algebraic varieties'', 
and by the European Commission through its 6th Framework Programme
"Structuring the European Research Area" and the contract Nr.
RITA-CT-2004-505493 for the provision of Transnational Access
implemented as Specific Support Action.}
 \fussy
\begin{document}
%%%%%%%%%%%%%%%%%%%% Text italic %%%%%%%%%%%%%%%%%%%%%%%%%%%%
\theoremstyle{plain}
\newtheorem{thm}{Theorem}[section]
\newtheorem{theorem}[thm]{Theorem}
\newtheorem{lemma}[thm]{Lemma}
\newtheorem{corollary}[thm]{Corollary}
\newtheorem{proposition}[thm]{Proposition}
\newtheorem{construction}[thm]{Construction}
\newtheorem{addendum}[thm]{Addendum}
\newtheorem{variant}[thm]{Variant}
\newtheorem{lemmadef}[thm]{Lemma and Definition}
%%%%%%%%%%%%%%%%%%%% Text roman %%%%%%%%%%%%%%%%%%%%%%%%%%%%%
\theoremstyle{definition}
\newtheorem{notations}[thm]{Notations}
\newtheorem{notation}[thm]{Notation}
\newtheorem{ass-not}[thm]{Assumptions and Notations}
\newtheorem{con-not}[thm]{Conclusion and Notations}
\newtheorem{lem-not}[thm]{Lemma and Notations}
\newtheorem{problem}[thm]{Problem}
\newtheorem{remark}[thm]{Remark}
\newtheorem{remarks}[thm]{Remarks}
\newtheorem{definition}[thm]{Definition}
\newtheorem{claim}[thm]{Claim}
\newtheorem{assumption}[thm]{Assumption}
\newtheorem{assumptions}[thm]{Assumptions}
\newtheorem{properties}[thm]{Properties}
\newtheorem{example}[thm]{Example}
\newtheorem{examples}[thm]{Examples}
\newtheorem{conjecture}[thm]{Conjecture}
\newtheorem{constr}[thm]{Construction}
\newtheorem{aconstr}[thm]{Allowed Constructions}
\newtheorem{condition}[thm]{Condition}
\newtheorem{conditions}[thm]{Conditions}
\newtheorem{setup}[thm]{Set-up}
\newtheorem{Variant}[thm]{Variant}
\newtheorem{comments}[thm]{Comments}
\newtheorem{conventions}[thm]{Conventions}
\numberwithin{equation}{thm}
\renewcommand{\theequation}{\arabic{section}.\arabic{thm}.\arabic{equation}}
% Skriptbuchstaben
\newcommand{\sA}{{\mathcal A}}
\newcommand{\sB}{{\mathcal B}}
\newcommand{\sC}{{\mathcal C}}
\newcommand{\sD}{{\mathcal D}}
\newcommand{\sE}{{\mathcal E}}
\newcommand{\sF}{{\mathcal F}}
\newcommand{\sG}{{\mathcal G}}
\newcommand{\sH}{{\mathcal H}}
\newcommand{\sI}{{\mathcal I}}
\newcommand{\sJ}{{\mathcal J}}
\newcommand{\sK}{{\mathcal K}}
\newcommand{\sL}{{\mathcal L}}
\newcommand{\sM}{{\mathcal M}}
\newcommand{\sN}{{\mathcal N}}
\newcommand{\sO}{{\mathcal O}}
\newcommand{\sP}{{\mathcal P}}
\newcommand{\sQ}{{\mathcal Q}}
\newcommand{\sR}{{\mathcal R}}
\newcommand{\sS}{{\mathcal S}}
\newcommand{\sT}{{\mathcal T}}
\newcommand{\sU}{{\mathcal U}}
\newcommand{\sV}{{\mathcal V}}
\newcommand{\sW}{{\mathcal W}}
\newcommand{\sX}{{\mathcal X}}
\newcommand{\sY}{{\mathcal Y}}
\newcommand{\sZ}{{\mathcal Z}}
% Sonderbuchstaben mit Doppellinie
\newcommand{\A}{{\mathbb A}}
\newcommand{\B}{{\mathbb B}}
\newcommand{\C}{{\mathbb C}}
\newcommand{\D}{{\mathbb D}}
\newcommand{\E}{{\mathbb E}}
\newcommand{\F}{{\mathbb F}}
\newcommand{\G}{{\mathbb G}}
\newcommand{\BH}{{\mathbb H}}
\newcommand{\I}{{\mathbb I}}
\newcommand{\J}{{\mathbb J}}
\newcommand{\BL}{{\mathbb L}}
\newcommand{\M}{{\mathbb M}}
\newcommand{\N}{{\mathbb N}}
\newcommand{\BP}{{\mathbb P}}
\newcommand{\Q}{{\mathbb Q}}
\newcommand{\R}{{\mathbb R}}
\newcommand{\BS}{{\mathbb S}}
\newcommand{\T}{{\mathbb T}}
\newcommand{\U}{{\mathbb U}}
\newcommand{\V}{{\mathbb V}}
\newcommand{\W}{{\mathbb W}}
\newcommand{\X}{{\mathbb X}}
\newcommand{\Y}{{\mathbb Y}}
\newcommand{\Z}{{\mathbb Z}}
\newcommand{\fC}{{\mathfrak C}}
\newcommand{\rk}{{\rm rk}}
\newcommand{\rounddown}[1]{\llcorner{#1}\lrcorner}
\newcommand{\fM}{{\mathfrak M}}
\def\bigtimes{\mathop{\mbox{\Huge$\times$}}}
%%%%%%%% Arrows %%%%%%%%%%%%%%%%%%%%%%%%%%%%%%%
\catcode`\@=11
\let\amp@rs@nd@\relax
\newdimen\bigaw@
\newdimen\minaw@
\minaw@16.08739\ex@
\newdimen\minCDaw@
%\minCDaw@2.5pc
\newif\ifCD@
\def\minCDarrowwidth#1{\minCDaw@#1}
\def\>#1>#2>{\amp@rs@nd@\setbox\z@\hbox{$\scriptstyle
 \;{#1}\;\;$}\setbox\@ne\hbox{$\scriptstyle\;{#2}\;\;$}\setbox\tw@
 \hbox{$#2$}\global\bigaw@\minaw@
 \ifdim\wd\z@>\bigaw@\global\bigaw@\wd\z@\fi
 \ifdim\wd\@ne>\bigaw@\global\bigaw@\wd\@ne\fi
 \ifdim\wd\tw@>\z@
 \mathrel{\mathop{\hbox to\bigaw@{\rightarrowfill}}\limits^{#1}_{#2}}\else
 \mathrel{\mathop{\hbox to\bigaw@{\rightarrowfill}}\limits^{#1}}\fi
\amp@rs@nd@}
\def\<#1<#2<{\amp@rs@nd@\setbox\z@\hbox{$\scriptstyle
 \;\;{#1}\;$}\setbox\@ne\hbox{$\scriptstyle\;\;{#2}\;$}\setbox\tw@
 \hbox{$#2$}\global\bigaw@\minaw@
 \ifdim\wd\z@>\bigaw@\global\bigaw@\wd\z@\fi
 \ifdim\wd\@ne>\bigaw@\global\bigaw@\wd\@ne\fi
 \ifdim\wd\tw@>\z@
 \mathrel{\mathop{\hbox to\bigaw@{\leftarrowfill}}\limits^{#1}_{#2}}\else
 \mathrel{\mathop{\hbox to\bigaw@{\leftarrowfill}}\limits^{#1}}\fi
\amp@rs@nd@}
%%%%%%%%%%%%%%%%%%%%%%%%%%%%%%%%%%%%%%%%%%%%%%%%%%%%%%%%%%%%%%%%%%%
\begin{abstract}
Let $M_h$ be the moduli scheme of canonically polarized manifolds with Hilbert polynomial $h$. We construct for $\nu\geq 2$ with $h(\nu)>0$ a projective compactification $\overline{M}_h$ of the reduced moduli scheme $(M_h)_{\rm red}$ such that the ample invertible sheaf $\lambda_\nu$, corresponding to $\det(f_*\omega_{X_0/Y_0}^\nu)$ on the moduli stack, has a natural extension $\overline{\lambda}_\nu\in {\rm Pic}(\overline{M}_h)_\Q$. A similar result is shown for moduli of polarized minimal models of Kodaira dimension zero. In both cases ``natural''
means that the pullback of $\overline{\lambda}_\nu$ to a curve $\varphi:C\to \overline{M}_h$, induced by
a family $f_0:X_0\to C_0=\varphi^{-1}(M_h)$, is isomorphic to $\det(f_*\omega_{X/C}^\nu)$ whenever $f_0$ extends to 
a semistable model $f:X\to C$. 

Besides of the weak semistable reduction of Abramovich-Karu and the
extension theorem of Gabber there are new tools, hopefully of interest by themselves. In particular we will need a theorem on the flattening of multiplier sheaves in families, on their compatibility with pullbacks and on base change 
for their direct images, twisted by certain semiample sheaves. 
\end{abstract}
\maketitle
\tableofcontents
%\section*{Introduction}\label{in}
\renewcommand{\thethm}{\arabic{thm}} 
%\renewcommand{\theequation}{\arabic{equation}}     
%%%%%%%%%%% Introduction %%%%%%%%%%%%%%%%%%%%%%%%%%%%%%%%%%%
Let $h_0:S_0\to C_0$ be a smooth family of complex projective manifolds over a non-singular curve
$C_0$. Replacing $C_0$ by a finite covering $\hat{C}_0$ one can extend the family 
$\hat{h}_0:\hat{S}_0=S_0\times_{C_0}\hat{C}_0\to \hat{C}_0$ to a semistable family $\hat{h}:\hat{S}\to \hat{C}$.
The model $\hat{S}$ is not unique, but the sheaves $\sF^{(\nu)}_{\hat{C}}=\hat{h}_*\omega_{\hat{S}/\hat{C}}^\nu$ are
independent of $\hat{S}$ and compatible with further pullback.
For a smooth family $f_0:X_0\to Y_0$ of $n$-folds over a higher dimensional base
the existence of flat semistable extension over a compactification $X$ of $X_0$ is not known,
not even the existence of a flat Cohen-Macaulay family, except if the fibres are curves or 
surfaces of general type.

It is the aim of this article to perform such constructions on the sheaf level.
So we fix a finite set $I$ of positive integers, and construct 
a finite covering $W_0$ of $Y_0$, and a compactification $W$ of $W_0$
such that for $\nu\in I$ the pullbacks of $f_{0*}\omega_{X_0/Y_0}^\nu$ extend to natural locally free and numerically effective (nef) sheaves $\sF_{\hat{Y}}^{(\nu)}$. The word ``natural'' means, that one has compatibility with pullback for certain morphisms $\hat{Y} \to W$. The precise statements are:
\begin{theorem}\label{in.1}
Let $f_0:X_0\to Y_0$ be a smooth projective morphism of quasi-projective reduced schemes
such that $\omega_F$ is semiample for all fibres $F$ of $f_0$. Let $I$ be a finite
set of positive integers. Then there exists a projective compactification $Y$ of $Y_0$,
a finite covering $\phi:W\to Y$ with a splitting trace map, and for $\nu\in I$ a locally free sheaf $\sF^{(\nu)}_{W}$ on $W$ with:
\begin{enumerate}
\item[i.] For $W_0=\phi^{-1}(Y_0)$ and $\phi_0=\phi|_{W_0}$ one has
$\phi_0^*f_{0*}\omega_{X_0/Y_0}^\nu=\sF^{(\nu)}_{W}|_{W_0}$.
\item[ii.] Let $\xi:\hat{Y}\to W$ be a morphism from a non-singular projective variety $\hat{Y}$ with
$\hat{Y}_0=\xi^{-1}(W_0)$ dense in $\hat{Y}$. Assume either that $\hat{Y}$ is a curve, or that
$\hat{Y}\to W$ is dominant. For some $r\geq 1$ let $X^{(r)}$ be a non-singular projective model of
the $r$-fold product family $\hat{X}_0^r=(X_0\times_{Y_0}\cdots \times_{Y_0}X_0)\times_{Y_0}\hat{Y}_0$
which admits a morphism $f^{(r)}:X^{(r)}\to \hat{Y}$. Then 
$\displaystyle f^{(r)}_*\omega_{X^{(r)}/\hat{Y}}^\nu=\bigotimes^r\xi^*\sF^{(\nu)}_W.$
\end{enumerate}
\end{theorem}
The formulation of Theorem~\ref{in.1} is motivated by what is needed to prove positivity
properties of direct image sheaves. 
\begin{theorem}\label{in.2} The conditions i) and ii) in Theorem~\ref{in.1} imply:
\begin{enumerate}
\item[iii.] The sheaf $\sF^{(\nu)}_{W}$ is nef.
\item[iv.] Assume that for some $\eta_1,\ldots,\eta_s \in I$ and for some $a_1,\ldots,a_s\in \N$
the sheaf $\bigotimes_{i=1}^s\det(\sF^{(\eta_i)}_{W})^{a_i}$ is ample with respect to $W_0$. Then, if
$\nu\geq 2$ and if $\sF^{(\nu)}_{W}$ is non-zero, it is ample with respect to $W_0$.
\end{enumerate}
\end{theorem}
In Section~\ref{ne} we recall the definition of the positivity properties  ``nef'',  ``ample with respect to $W_0$'' and of ``weakly positive over $W_0$'' for locally free sheaves $\sF$ on $W$. Obviously $\sF$ is nef if and only if its pullback under a surjective morphism $\varphi:W'\to W$ is nef. For ``nef and ample with respect to $W_0$'' the same holds if $\varphi$ is finite over $W_0$. ``Weakly positive over $W_0$'' is compatible with finite coverings with a splitting trace map, i.e. if $\sO_{W}$ is a direct factor of $\varphi_*\sO_{W'}$. 

In Section~\ref{pd1} we show that part iii) of Theorem~\ref{in.2} follows from Theorem~\ref{in.1}. 
Unfortunately, as we will explain in~\ref{pd1.5}, the verification of the property iv) is much harder.
Here multiplier ideals will enter the scene. Whereas in a neighborhood of a smooth fibre $F$ one can bound the threshold, introduced in~\ref{pd1.1}, in terms of invariants of $F$, a similar result fails close to the boundary. So we need a variant of parts i) and ii), allowing certain multiplier sheaves, introduced in Section~\ref{de}. as well as the Flattening Theorem ~\ref{fm.5} for multiplier ideal sheaves on total spaces of morphisms, and their compatibility with alterations of the base and fibre products. So the proof of part iv) will only be given at the end of Section~\ref{pd2}.

There are two main ingredients which will allow the construction of $W$ and $\sF_W^{(\nu)}$ in Section~\ref{ex}.
The first one is the Weak Semistable Reduction Theorem \cite{AK} recalled in Section~\ref{ws}. Roughly speaking it says that a given morphism $f:X\to Y$ between projective varieties, with a smooth general fibre can be flattened over some non-singular alteration of $Y$ without allowing horrible singularities of the total space. However one pays a price. One has to modify the smooth fibres as well. As explained in Section~\ref{di} this Theorem has some strong consequences for the compatibility of certain sheaves on the total space of a family with base change and products, similar to those stated in part ii) of Theorem~\ref{in.1}. 

The second ingredient is Gabber's Extension Theorem, stated (and proved) in \cite[Section 5.1]{Vie}, which we will recall in Section~\ref{ex}. 

The comments made in~\ref{pd1.5} and in~\ref{di.9} could serve as a ``Leitfaden'' for the second part of the article. Here we try to indicate why certain constructions contained in Sections~\ref{fm}--\ref{va} are needed for the proof of the Theorem~\ref{in.2} iv).\vspace{.1cm} 

From Theorem~\ref{in.2} one finds by Lemma~\ref{ne.6} that the restriction of $\sF_W^{(\nu)}$ to $W_0$ is weakly positive over $W_0$. We will show in Section~\ref{pd2} that part iv), in a slightly modified version, also restricts to $W_0$. Since $W_0\to Y_0$ has a splitting trace map,
one obtains by Lemma~\ref{ne.7} the ``weak positivity'' and ``weak stability'' for the direct images of powers of the dualizing sheaf, already shown in \cite[Section 6.4]{Vie}.
\begin{corollary}\label{in.3} Under the assumptions made in Theorem~\ref{in.1} one has:
\begin{enumerate}
\item[a.] The sheaves $\sF_{Y_0}^{(\nu)}=f_{0*}\omega_{X_0/Y_0}^\nu$
are weakly positive over $Y_0$.
\item[b.] Assume that for some positive integers $\eta_1,\ldots,\eta_s$ and $a_1,\ldots,a_s\in \N$
the sheaf $\bigotimes_{i=1}^s\det(\sF^{(\eta_i)}_{Y_0})^{a_i}$ is ample. Then for all
$\nu\geq 2$ the sheaf $\sF^{(\nu)}_{Y_0}$ is either ample or zero.
\end{enumerate}
\end{corollary}
As explained in \cite{Vie} this is just what is needed for the construction of quasi-projective moduli schemes $M_h$ for families of canonically polarized manifolds with Hilbert polynomial $h$. At the time \cite{Vie} was written, the Weak Semistable Reduction Theorem of Abramovich and Karu was not known. So we were only able to use Gabber's Extension Theorem to construct $W$ and $\sF_{W}^{(\nu)}$ for $\nu=1$, and correspondingly to prove the weak positivity just for $\sF_{Y_0}^{(1)}$. A large part of \cite{Vie} is needed to reduce the proof of Corollary~\ref{in.3} to this case.
Having $W$ and $\sF_{W}^{(\nu)}$ for all $\nu$ clarifies this part considerably. We could not resist to recall in Section~\ref{co} how to apply Corollary~\ref{in.3} to construct $M_h$ together with an ample invertible sheaf.

There are several ways. One can first construct the moduli scheme as an algebraic space, and then show the existence of an ample sheaf. Or one can use geometric invariant theory, and stability criteria. 
Guided by personal taste, we restrict ourselves to the second method in Section~\ref{co}, applying the Stability Criterion \cite[Theorem 4.25]{Vie}. 

If one uses instead the first method, starting from the existence of $M_h$ as an algebraic space, 
it has been shown in \cite{Vie} how to deduce from Corollary~\ref{in.3} the quasi-projectivity of the normalization of $M_h$. The starting point is Seshadri's Theorem on the elimination of finite isotropies (see \cite[Theorem 3.49]{Vie}) or the direct construction in \cite{Kol}. Both allow to get a universal family $f_0:X_0\to Y_0$ over some reduced covering $\gamma_0:Y_0\to M_h$. Then one can try to apply arguments similar to those used in the proof of Lemma~\ref{ne.9} and in Section~\ref{pr} to get the quasi-projectivity of $(M_h)_{\rm red}$ hence of $M_h$ itself. 

As stated in the report on \cite{ST} in the ``Mathematical Reviews'', J. Koll\'ar pointed out 
that the proof of the quasi-projectivity of the algebraic moduli space $M_h$ seems to contain a gap, even in the canonically polarized case. The authors claim without any justification that for a certain line bundle, which descends to a quotient of the Hilbert scheme, the curvature current descends as well. In a more recent attempt to handle moduli of canonically polarized manifolds Tsuji avoids this point by claiming that a certain determinant sheaf extends to some compactification in a natural way, again without giving an argument. Suitable variants of Theorems~\ref{in.1} and~\ref{in.2} could allow to fill those gaps, and to get another proof of the quasi-projectivity of $M_h$, replacing the 
GIT-approach in Section~\ref{co} by the analytic methods presented in the second part of \cite{ST}.

Either one of the constructions of moduli schemes mentioned above gives an explicit ample sheaf on $M_h$. 
\begin{notations}\label{in.4}
For $\nu, \ p \in \N$ we write $\lambda_{0,\nu}^{(p)}$ for an invertible sheaf satisfying 
\begin{enumerate}
\item[($\divideontimes$)]
If a morphism $\varphi:Y_0\to M_h$ factors through the moduli stack, hence if it is induced
by a family $f_0:X_0\to Y_0$, one has $\varphi^*\lambda_{0,\nu}^{(p)}=\det(f_{0*}\omega_{X_0/Y_0}^\nu)^p$.
\end{enumerate}
\end{notations}
Of course, $\lambda_{0,\nu}^{(p)}$ can only exist if $H^0(F,\omega_F^\nu)\neq 0$ for all manifolds $F$ parametrized
by $M_h$, for example in the canonically polarized case if $\nu\geq 2$ and $h(\nu)\neq 0$.
As we will recall in Addendum~\ref{co.1} for those values of $\nu$ the sheaf $\lambda_{0,\nu}^{(p)}$ is ample.\vspace{.1cm}

As indicated in the title of this article we want to construct compactification of moduli schemes $M_h$.
Assume for a moment that $M_h$ is reduced and a fine moduli scheme, hence that there is a universal family
$\sX_0\to M_h$ with $\sX_0$ reduced. Here one may choose $p=1$ and applying Theorems~\ref{in.1} and~\ref{in.2}, and Lemma~\ref{ne.9} it is easy to see that $\lambda_{0,\nu}^{(1)}=\det(g_{0*}\omega^\nu_{\sX/M_h})$ extends to an invertible sheaf $\lambda_{\nu}^{(1)}$ on $\overline{M}_h$, which is nef and ample with respect to $M_h$ for $\nu\geq 2$. It is compatible with the restriction to curves, provided the induced family has a smooth general fibre and a semistable model.
In Section~\ref{pr} we will use a variant of Theorems~\ref{in.1} and~\ref{in.2} to obtain a similar result for coarse moduli schemes, using the Seshadri's-Koll\'ar construction mentioned above.

\begin{theorem}\label{in.5} Let $M_h$ be the coarse moduli scheme of canonically polarized manifolds
with Hilbert polynomial $h$.
Given a finite set $I$ of integers $\nu\geq 2$ with $h(\nu)>0$, one finds
a projective compactification $\overline{M}_h$ of $(M_h)_{\rm red}$ and for $\nu\in I$ and some $p>0$
invertible sheaves $\lambda_\nu^{(p)}$ on $\overline{M}_h$ with:
\begin{enumerate}
\item $\lambda_\nu^{(p)}$ is nef and ample with respect to $(M_h)_{\rm red}$.
\item The restrictions of $\lambda_\nu^{(p)}$ and of $\lambda_{0,\nu}^{(p)}$ to $(M_h)_{\rm red}$ coincide.
\item Let $\varsigma:C\to \overline{M}_h$ be a morphism from a non-singular curve
with $C_0=\varsigma^{-1}(M_h)$ dense in $C$ and such that $C_0\to M_h$ is induced by a family
$h_0:S_0\to C_0$. If $h_0$ extends to a semistable family $h:S\to C$, then
$\varsigma^*\lambda_\nu^{(p)}=\det(h_*\omega_{S/C}^\nu)^p$.
\end{enumerate}
\end{theorem}
It would be nicer to have an extension of $\lambda_{0,\nu}^{(p)}$ to an invertible sheaf $\lambda_\nu^{(p)}$ on a compactification of $M_h$ itself, but we were not able to get hold of it. 
On the other hand, since the compatibility condition in part (3) only sees the reduced structure of
$M_h$, such an extension would not really be of help for possible applications of Theorem~\ref{in.5}.

The compactification $\overline{M}_h$ depends on the set $I$ and the points in $\overline{M}_h\setminus M_h$ have no 
interpretation as moduli of geometric objects. Shortly after a first version of this article was
submitted there was a ``quantum leap'' in the minimal model program due to \cite{BCHM}
(see also \cite{Siu}). By \cite{Kar} the existence of a minimal model in dimension $\deg(h)+1$ 
should allow the construction of a compactification $\overline{M}_h$ which has an interpretation as a moduli 
scheme. Unfortunately at the present time a proper explanation of this implication is
not in the literature and there is not even a conjectural picture explaining how to construct geometrically meaningful compactifications of moduli of polarized manifolds. 

Only part of what was described up to now carries over to families or moduli of smooth minimal models with an arbitrary polarization. Theorems~\ref{in.1} and~\ref{in.2} apply, but even if $f_0:X_0\to Y_0$ is the universal family over a fine moduli scheme, the sheaf $\det(\sF_{Y_0}^{(\nu)})$ might not be ample. Theorem~\ref{ex.12} is a generalization of Theorem~\ref{in.1} for direct images of the form $f_{0*}(\omega_{X_0/Y_0}^\nu \otimes \sL_0^\mu)$ with $\sL_0$ semiample over $Y_0$. The corresponding variant of Corollary~\ref{in.3} is stated in~\ref{co.2} and we will sketch how to use it to show the existence of quasi-projective moduli schemes in the second half of Section~\ref{co}. 
However we are not able to generalize Theorem~\ref{in.2}. So we are not able to apply Lemma~\ref{ne.9} which will be essential for the proof of Theorem~\ref{in.5} in Section~\ref{pr}
and we are not able to extend the natural ample sheaf to some compactification.\vspace{.1cm} 

The situation is better for the moduli functor $\fM_h$ of polarized minimal manifolds $(F,\sH)$ of Kodaira dimension zero and with Hilbert polynomial $h$. Replacing the corresponding moduli scheme $M_h$
by a connected component, we may assume that for some $\nu > 0$ and for all $(F,\sH)\in \fM_h$ one has 
$\omega_F^\nu=\sO_F$. The notation in~\ref{in.4} carries over and there exists a sheaf $\lambda^{(p)}_{0,\nu}$ with the property ($\divideontimes$) or equivalently with $f_0^*\varphi^*\lambda_{0,\nu}^{(p)}=\omega_{X_0/Y_0}^{p\cdot \nu}$. As we will recall in Addendum~\ref{co.4} the sheaf $\lambda^{(p)}_{0,\nu}$ is again ample. In fact, the natural
ample invertible sheaf first looks quite different, but an easy calculation identifies it with
some power of $\lambda^{(p)}_{0,\nu}$. This calculation only extends to boundary points if the 
polarization is {\em saturated}, as explained in Remark~\ref{sa.1}. This together with the need to consider multiplier ideals makes the notations even more unpleasant, but the general line of arguments remain as in the canonically polarized case. 

\begin{theorem}\label{in.6} Let $M_h$ be the coarse moduli scheme of polarized manifolds $(F,\sH)$
with $\omega_F^\nu=\sO_F$, for some $\nu>0$ and with Hilbert polynomial $h(\mu)=\chi(\sH^\mu)$.
Then there exists a projective compactification $\overline{M}_h$ of $(M_h)_{\rm red}$ and for some $p>0$
an invertible sheaf $\lambda_\nu^{(p)}$ on $\overline{M}_h$ with:
\begin{enumerate}
\item $\lambda_\nu^{(p)}$ is nef and ample with respect to $(M_h)_{\rm red}$.
\item Let $Y_0$ be reduced and $\varphi:Y_0\to M_h$ induced by a family $f_0:X_0\to Y_0$ in $\fM_h(Y_0)$.
Then $\varphi^*\lambda_\nu^{(p)}=f_{0*}\omega_{X_0/Y_0}^{p\cdot\nu}$.
\item Let $\varsigma:C\to \overline{M}_h$ be a morphism from a non-singular curve $C$,
with $C_0=\varsigma^{-1}(M_h)$ dense in $C$ and such that $C_0\to M_h$ is induced by a family
$h_0:S_0\to C_0$. If $h_0$ extends to a semistable family $h:S\to C$, then
$\varsigma^*\lambda_\nu^{(p)}=h_*\omega_{S/C}^{p\cdot\nu}$.
\end{enumerate}
\end{theorem}
Consider, for example, the moduli schemes $\sA_g$ of $g$-dimensional polarized Abelian varieties. To stay close to the usual notation we write $\overline{\sA}_g$ for the compactification in Theorem~\ref{in.6} in this case. One may assume that there is a morphism $\overline{\sA}_g\to \sA^*_g$ to the Baily-Borel compactification $\sA^*_g$. Choosing $\nu=1$ the ample sheaf $\lambda_{0,1}^{(p)}$ extends to an ample sheaf
on $\sA^*_g$. So the sheaf $\lambda_1^{(p)}$ in Theorem~\ref{in.6} is semiample and $\sA^*_g$
is the image under the morphism defined by a high power of $\lambda_1^{(p)}$.

In general we are not able to verify in Theorem~\ref{in.6} the semiampleness of $\lambda_\nu^{(p)}$.
One of the obstacles is the missing geometric interpretation of the boundary points as moduli
of certain varieties. So the theorem can only be seen as a very weak substitute for the Baily-Borel compactification.
Nevertheless, since $\lambda_\nu^{(p)}$ is nef and ample with respect to $M_h$ the degree of
$\varsigma^*\lambda_\nu^{(p)}=h_*\omega_{S/C}^{p\cdot\nu}$ can serve as a height function for
curves in the moduli stack. An upper bound for this height function in terms of the genus of $C$ and $\#(C\setminus C_0)$ was given in \cite{VZ2}, and 
it played its role in the proof of the Brody hyperbolicity of the moduli stack of canonically polarized manifolds in \cite{VZ1}. In both articles we had to use unpleasant ad hoc arguments to control the positivity along the boundary of the moduli schemes and some of those arguments were precursors of methods used here. 

A second motivation for this article was the hope that compactifications could help to generalize the 
uniform boundedness, obtained in \cite{Cap} for families of curves, to families of higher dimensional manifolds. The missing point was the construction of moduli of morphisms from curves to the corresponding
moduli stacks, as it was done in \cite{AV} for compact moduli problems. In between this has been achieved
in \cite{KL} for families of canonically polarized manifolds, using Theorem~\ref{in.5}. 
It is likely that referring to Theorem~\ref{in.6} instead, their methods allow to handle polarized manifolds of Kodaira dimension zero, as well.\vspace{.1cm}

I was invited to lecture on the construction of moduli schemes at the workshop "Compact moduli spaces and birational geometry" (American Institute of Mathematics, 2004), an occasion to reconsider some of the constructions in \cite{Vie} in view of the Weak Semistable Reduction Theorem. A preliminary version of this article, handling just the canonically polarized case, was written during a visit to the I.H.E.S., Bures sur Yvette September and October 2005. I like to thank the members of the Institute for their
hospitality. 

I am grateful for the referee's suggestions how to improve the presentation of the results and
of the methods leading to their proofs.
\begin{conventions} All schemes and varieties will be defined over the field $\C$ of complex numbers (or over an algebraically closed field $K$ of characteristic zero). A quasi-projective variety $Y$ is a reduced quasi-projective scheme. In particular we do not require $Y$ to be irreducible or connected. A locally free sheaf on $Y$ will always be locally free of constant finite rank and a finite covering will denote a finite surjective morphism.\\[.1cm]
$\bullet$ If $\Pi$ is an effective divisor and $\iota: Y\setminus \Pi_{\rm red} \to Y$, then $\sO(*\cdot \Pi)=\iota_*\sO_{Y\setminus \Pi_{\rm red}}$.\\[.1cm]
$\bullet$ An {\em alteration} $\Psi:\hat{Y}\to Y$ is a proper, surjective, generically finite morphism
between quasi-projective varieties.\\[.1cm]
$\bullet$ An alteration $\Psi$ is called a {\em modification} if it is birational. If $U\subset Y$ is an open subscheme with $\Psi|_{\Psi^{-1}(U)}$ an isomorphism, we say that the center of $\Psi$ lies in $Y\setminus U$.\\[.1cm]
$\bullet$ For a non-singular (or normal) alteration or modification we require in addition that $\hat{Y}$ is non-singular (or normal).\\[.1cm]
$\bullet$ A modification $\Psi$ will be called a {\em desingularization (or resolution of singularities)}, if $\hat{Y}$ is non-singular and if the center of $\Psi$ lies in the singular locus of $Y$.\\[.1cm]
$\bullet$ Given a Cartier divisor $D$ on $Y$ we call $\Psi$ a {\em log-resolution} (for $D$) if it is a non-singular modification and if $\Psi^*D$ is a normal crossing divisor.\\[.1cm]
$\bullet$ If $f_0:X_0\to Y_0$ is a projective morphism, we call $f:X\to Y$ a {\em projective model} of $f_0$ if
$X$ and $Y$ are projective, $Y_0$ open in $Y$ and $X_0 \cong f^{-1}(Y_0)$ over $Y_0$.\\[.1cm]
$\bullet$ If $f:X\to Y$ is a projective morphism, we call $f_0:X_0\to Y_0$ {\em the smooth part of $f$} if
$Y_0\subset Y$ is the largest open subscheme with $X_0= f^{-1}(Y_0)\to Y_0$ smooth. In particular, if $X$ and $Y$ are non-singular, or if $f$ is a mild morphism, as defined in~\ref{ws.1}, $Y_0$ is dense in $Y$.\vspace{.1cm}

Finally the numbering of displayed formulas follows the one of the theorems, lemmas etc. Hence 
\eqref{eqpd1.1} is after Definition~\ref{pd1.1} and before Lemma~\ref{pd1.2}.
\end{conventions}
\section{Numerically effective and weakly positive sheaves}\label{ne}
\renewcommand{\thethm}{\arabic{section}.\arabic{thm}} 
\begin{definition}\label{ne.1}
Let $\sG$ be a locally free sheaf on a projective reduced variety $W$. Then
$\sG$ is {\em numerically effective (nef)} if for all morphisms $\tau:C\to W$
from a projective curve $C$ and for all invertible quotients $\tau^*\sG \to \sL$
one has $\deg(\sL)\geq 0$.
\end{definition}

\begin{definition} \label{ne.2}
Let $\sG$ be a locally free sheaf on a quasi-projective reduced variety $W$ and let $W_0 \subset W$ be an open dense subvariety. Let $\sH$ be an ample invertible sheaf on $W$.
\begin{enumerate}
\item[a)] $\sG$ is {\em globally generated over $W_0$} if the natural
morphism $H^0 (W, \sG) \otimes \sO_W \to \sG$ is surjective over $W_0$.
\item[b)] $\sG$ is {\em weakly positive over $W_0$} if for all $\alpha >0$ there exists some $\beta > 0$
such that $S^{\alpha \cdot \beta} (\sG) \otimes \sH^{\beta}$ is globally generated over $W_0$.
\item[c)] $\sG$ is {\em ample with respect to $W_0$} if for some $\eta>0$ the sheaf
$S^\eta(\sG)\otimes \sH^{-1}$ is weakly positive over $W_0$, or equivalently, if
for some $\eta'>0$ one has a morphism
$\bigoplus \sH \to S^{\eta'}(\sG)$, which is surjective over $W_0$.
\end{enumerate}
\end{definition}
It is obvious, that ``nef'' is related to ``weakly positive'' and that it is
compatible with pullbacks.
\begin{lemma}\label{ne.3} For a locally free sheaf $\sG$ on a projective variety $W$ the following
conditions are equivalent:
\begin{enumerate}
\item $\sG$ is nef.
\item $\sG$ is weakly positive over $W$.
\item There exists a projective surjective morphism $\xi: \hat{Y}\to W$ with $\xi^*\sG$ nef.
\item The sheaf $\sO_{\BP(\sG)}(1)$ on $\BP(\sG)$ is nef.
\item There exists some integer $\mu>0$ such that for all projective surjective morphisms $\xi:\hat{Y}\to W$
and for all ample invertible sheaves $\hat{\sH}$ on $\hat{Y}$ the sheaf $\hat{\sH}^\mu\otimes\xi^*\sG$ is nef. 
\end{enumerate}
\end{lemma}
\begin{remark}\label{ne.4}
As we will see in the proof it is sufficient in~\ref{ne.3} (5) to require the existence a tower of finite maps $\xi:\hat{Y}\to W$, such that for each $N>0$ there is some $\xi:\hat{Y}\to W$ with $\xi^*\sH$ the $N$-th power of an invertible sheaf. Such coverings exist by \cite[Lemma 2.1]{Vie}, and one even may assume that they have splitting trace maps. 
\end{remark}
\begin{proof}
The equivalence of the first four conditions has been shown in \cite[Proposition 2.9]{Vie}, and of course
they imply (5). The equivalence of (5) and (2) is a special case of \cite[Lemma 2.15, 3)]{Vie}.
Nevertheless let us give the argument. Let $\sH$ be ample and invertible on $W$.
Let $\pi:C\to W$ be a curve and $\sN$ an invertible quotient
of $\pi^*\sG$ of degree $d$. By \cite[Lemma 2.1]{Vie} for all $N$ there exist a finite covering
$\xi:\hat{Y}\to W$ and an invertible sheaf $\hat{\sH}$ with $\xi^*\sH=\hat{\sH}^N$.  
By assumption $\hat{\sH}^\mu\otimes\xi^*\sG$ is nef, hence if $\tau:\hat{C}\to C$ is a finite covering
such that $\pi$ lifts to $\pi':\hat{C}\to \hat{Y}$ one has
$$
0\leq\deg(\tau)\cdot d + \mu\cdot\deg(\pi'^*\hat{\sH})=\deg(\tau)\cdot (d + \frac{\mu}{N}\cdot\deg(\pi^*\sH)).
$$
This being true for all $N$, the degree $d$ can not be negative.  
\end{proof}
Obviously the notion ``nef'' is compatible with tensor products, direct sums, symmetric products
and wedge products. For the corresponding properties for weakly positive, one has to work a bit 
more, or to refer to \cite[Section 2.3]{Vie}. 
\begin{lemma}\label{ne.5}
Let $\sF$ and $\sG$ be locally free sheaves on $W$.
\begin{enumerate}
\item Let $\sL$ be an invertible sheaf. Assume that for all $\alpha >0$ there exists some $\beta > 0$
such that $S^{\alpha \cdot \beta} (\sG) \otimes \sL^{\beta}$ is globally generated over $W_0$.
Then $\sG$ is weakly positive over $W_0$. In particular Definition~\ref{ne.2} b), is independent
of $\sH$.
\item If $\sG$ is weakly positive over $W_0$ and if $\xi:\hat{Y}\to W$ is a dominant morphism, then
$\xi^*\sG$ is weakly positive over $\xi^{-1}(W_0)$.
\item If $\sG$ is weakly positive over $W_0$ and if $\sG \to \sF$ is a morphism, surjective over
$W_0$, then $\sF$ is weakly positive over $W_0$.
\item If $\sF$ and $\sG$ are weakly positive over $W_0$, the the same holds for
$\sF\oplus\sG$, for $\sF\otimes\sG$, for $S^\nu(\sG)$ and for $\bigwedge^\mu(\sG)$, where
$\nu$ and $\mu\leq \rk(\sG)$ are natural numbers.
\end{enumerate}
\end{lemma}
The equivalence of (1) and (3) in~\ref{ne.3} does not carry over to ``weakly positive over $W_0$''; one needs in addition that the morphism is finite with a splitting trace map.
\begin{lemma}\label{ne.6} 
For a locally free sheaf $\sG$ on $W$ and an open and dense subscheme $W_0\subset W$
the following conditions are equivalent:
\begin{enumerate}
\item $\sG$ is weakly positive over $W_0$.
\item $\bigotimes^r \sG$ is weakly positive over $W_0$ for some $r>0$.
\item $S^r \sG$ is weakly positive over $W_0$ for some $r>0$.
\item There exists an invertible sheaf $\sA$ on $W$ such that $\sA\otimes S^r(\sG)$
is weakly positive over $W_0$, for all $r >0$.
\item For or all (or some) ample invertible sheaves $\sA$ on $W$ and for all $r>0$ the sheaf $\sA\otimes S^r(\sG)$
is ample with respect to $W_0$. 
\item There exists an alteration $\phi: \widetilde{W} \to W$ such that 
$\phi^*\sG$ is weakly positive over $\phi^{-1}(W_0)$, and such that
for $\widetilde{W}_0=\phi^{-1}(W_0)$ the restriction $\phi_0:\widetilde{W}_0\to W_0$
is finite with a splitting trace map (i.e. with a splitting of $\sO_{W_0}\to \phi_{0*}\sO_{\widetilde{W}_0}$). 
\item There exists a constant $\mu>0$ such that for all $\xi:\hat{Y} \to W$ and for all
ample invertible sheaves $\sH'$ on $\hat{Y}$ the sheaf $\sH'^\mu\otimes \xi^*\sG$
is weakly positive over $\xi^{-1}(W_0)$.
\end{enumerate}
\end{lemma}
The remark~\ref{ne.4} applies to~\ref{ne.6} (7) as well, if one assumes that
for all the $\xi:\hat{Y}\to W$ the trace map splits.
\begin{proof}
The equivalence of the first three conditions has been shown in \cite[Lemma 2.16]{Vie}.
The equivalence of (1), (4) and (5) follows directly from the definition,
and the equivalence of (1), (6) and (7) is in \cite[Lemma 2.15]{Vie}.
\end{proof}
Let us consider next the condition ``ample with respect to $W_0$''. 
\begin{lemma}\label{ne.7} Let $\sG$ and $\sF$ be locally free sheaves on $W$ and let $W_0\subset W$ be open and dense.
\begin{enumerate}
\item $\sG$ is ample with respect to $W_0$ if and only if there exists an ample invertible sheaf $\sH$
on $W$ and a finite morphism $\sigma:W'\to W$ with a splitting trace map, and with
$\sigma^*\sH=\sH'^\eta$, for some positive integer $\eta$, such that $\sigma^*(\sG)\otimes \sH'^{-1}$ is weakly positive over $\sigma^{-1}(W_0)$.
\item If $\sF$ is ample with respect to $W_0$ and if $\sG$ is weakly positive over
$W_0$, then $\sF\otimes \sG$ is ample with respect to $W_0$. In particular, the Definition~\ref{ne.2} c), is independent of the ample invertible sheaf $\sH$.
\item If $\sF$ is invertible and ample with respect to $W_0$, and if $S^\eta(\sG)\otimes\sF^{-1}$
is weakly positive over $W_0$, then $\sG$ is ample over $W_0$.
\item The following conditions are equivalent:
\begin{enumerate}
\item $\sG$ is ample with respect to $W_0$.
\item There exists an alteration $\phi: \widetilde{W} \to W$ with
$\widetilde{W}_0=\phi^{-1}(W_0)\to W_0$ finite and with a splitting trace map, such that 
$\phi^*\sG$ is ample with respect to $\widetilde{W}_0$. 
\end{enumerate}
\item If $\sG$ is ample with respect to $W_0$ and if $\sG \to \sF$ is a morphism, surjective over $W_0$, then $\sF$ is ample with respect to $W_0$.
\item If $\sF$ and $\sG$ are both ample with respect to $W_0$, then the same holds for
$\sF\oplus\sG$, for $S^\nu(\sG)$ and for $\bigwedge^\mu(\sG)$, where $\nu$ and $\mu\leq \rk(\sG)$ are natural numbers.
\item If $\sF$ is an invertible sheaf, then $\sF$ is ample with respect to
$W_0$, if and only if for some $\beta>0$ the sheaf $\sF^\beta$ is globally generated over $W_0$
and the induced morphism $\tau: W_0\to \BP(H^0(W,\sF^\beta))$ finite over its image.
\end{enumerate}
\end{lemma}
\begin{proof}
If in (1) the sheaf $\sG$ is ample with respect to $W_0$ there is some $\eta$ such that $S^\eta(\sG)\otimes\sH^{-1}$
is weakly positive. By \cite[Lemma 2.1]{Vie} there is a covering $\sigma:W'\to W$ with a splitting trace map, such that $\sigma^*\sH$ is the $\eta$-th power of an invertible sheaf $\sH'$, necessarily ample.
Then $\sigma^*(S^\eta(\sG)\otimes \sH^{-1})$ is weakly positive over $\sigma^{-1}(W_0)$, hence
by~\ref{ne.6} the sheaf $\sigma^*\sG\otimes \sH'^{-1}$ as well. On the other hand, the weak positivity
of $\sigma^*(\sG)\otimes \sH'^{-1}$ in (1) implies that $\sigma^*S^\eta(\sG) \otimes \sigma^*\sH^{-1}$ is weakly positive over $\sigma^{-1}(W_0)$, hence $S^\eta(\sG) \otimes \sH^{-1}$ is weakly positive over $W_0$, using again~\ref{ne.6}.

For (2) one can use (1), assume that $\sG\otimes \sH^{-1}$ is weakly positive, and then apply
\ref{ne.5} (4). In the same way one obtains (6). Part (3) is a special case of (2) and (5)
follows from~\ref{ne.5} (3).

Let us next verify (7). If $\sF$ is ample with respect to $W_0$, one has for a very ample invertible sheaf $\sH$ on $W$ and for some $\eta'$ a morphism $\bigoplus^s \sH \to \sF^{\eta'}$, surjective over $W_0$. 
Let $V$ denote the image of $H^0(W,\bigoplus^s \sH)$ in $H^0(W,\sF^{\eta'})$. 
Then $\sF^{\eta'}$ is generated by $V$ over $W_0$ and one has embeddings 
$$
W \>>> \bigtimes^s \BP(H^0(W,\sH)) \>>> \BP(\bigotimes^sH^0(W,\sH)).
$$ 
The restriction of the composite to $W_0$ factors through 
$$
W_0\to \BP(V)\subset \BP(\bigotimes^sH^0(W,\sH)),
$$ 
and $W_0\to \BP(V)$, hence $W_0\to \BP(H^0(W,\sF^{\eta'}))$ are embeddings. 

If on the other hand $\sF^\beta$ is globally generated over $W_0$ and if
$$
\tau:W_0\to \BP=\BP(H^0(W,\sF^{\beta}))
$$ 
is finite over its image, consider a blowing up $\phi:\widetilde{W}\to W$ with centers outside of $W_0$ such that $\tau$ extends to a morphism
$\tau':W\to \BP$. We may choose $\phi$ such that for some effective exceptional divisor $E$ the sheaf $\sO_{\widetilde{W}}(-E)$ is $\tau'$-ample. For $\alpha$ sufficiently large $\sA=\sO_{\widetilde{W}}(-E)\otimes \tau'^*\sO_{\BP}(\alpha)$ will be ample. Replacing $E$ and $\alpha$ by some multiple, one may assume that for a given ample sheaf $\sH$ on $W$ the sheaf $\phi^*\sH^{-1}\otimes\sA$ is globally generated, hence nef. Since one has an inclusion $\sA\to \phi^*\sF^{\eta'\cdot\alpha}$, which is an isomorphism over $\phi^{-1}(W_0)$, the sheaf $\phi^*\sF^{\eta'\cdot\alpha}\otimes \sH^{-1}$ is weakly positive over
$\phi^{-1}(W_0)$, and by~\ref{ne.6} one obtains the weak positivity of
$\sF^{\eta'\cdot\alpha}\otimes \sH^{-1}$.

For (4) we use (7). Consider in (4), a), an ample invertible sheaf $\sF$ on $W$. Obviously the condition (7) holds for $\phi^*\sF$, hence this sheaf is again ample with respect to
$\phi^{-1}(W_0)$. If $\sG$ is ample with respect to $W_0$, by definition $S^\nu(\sG)\otimes \sF^{-1}$ is weakly positive over $W_0$. Then by~\ref{ne.7} (6) the sheaf $\phi^*S^\nu(\sG)\otimes \phi^*\sF^{-1}$
is weakly positive over $\phi^{-1}(W_0)$ and (4), b), follows from (3).

So assume that the condition b) in (4) holds.
Let $\sH$ and $\sA$ be ample invertible sheaves on $W$ and $Y$. Then $\sA\otimes \phi^*\sH$ is ample.
By definition we find some $\beta$ such that $S^\beta(\phi^*\sG)\otimes \sA^{-1}\otimes \phi^*\sH^{-1}$ is weakly positive over $Y_0$. So $S^\beta(\phi^*\sG)\otimes \phi^*\sH^{-1}$ has the same property, 
and by Lemma~\ref{ne.6} $S^\beta(\sG)\otimes \sH^{-1}$ is weakly positive over $W_0$.
\end{proof}

\begin{lemma}\label{ne.8}
A locally free sheaf $\sG$ on $W$ is ample with respect to $W_0$ if and only if
on the projective bundle $\pi:\BP(\sG)\to W$ the sheaf $\sO_{\BP(\sG)}(1)$ is ample with respect to
$\BP_0=\pi^{-1}(W_0)$.
\end{lemma}
\begin{proof}
If $\sG$ is ample with respect to $W_0$ choose a very ample invertible sheaf $\sH$ on
$W$ and for some $\eta'>0$ the morphism 
$$
\bigoplus^s \sH \>>> S^{\eta'}(\sG)=\pi_*\sO_{\BP(\sG)}(\eta'),
$$
surjective over $W_0$. The composite
$$
\bigoplus^s \pi^*\sH \>>> S^{\eta'}(\pi^*\sG)\>>> \sO_{\BP(\sG)}(\eta')
$$
induces a rational map $\iota:\BP(\sG)\to \BP^{s-1}$, 
whose restriction to $\BP_0=\pi^{-1}(W_0)$ is an embedding, and $\sO_{\BP(\sG)}(\eta')$ is globally generated over $\BP_0$. So by~\ref{ne.7} (7) $\sO_{\BP(\sG)}(1)$ is ample with respect to
$\BP_0$.

Assume now that $\sO_{\BP(\sG)}(1)$ is ample with respect to $\BP_0$. Choose ample invertible sheaves
$\sH$ on $W$ and $\sA$ on $\BP(\sG)$ such that $\pi^*\sH^{-1}\otimes\sA$ is globally generated. Then for some $\eta'$ and for all $\alpha >0$ one has  morphisms
$$
\bigoplus \pi^*\sH^\alpha \>\Psi >> \bigoplus \sA^\alpha \>\Phi >> \sO_{\BP(\sG)}(\eta'\cdot \alpha)
$$
with $\Psi$ surjective and $\Phi$ surjective over $\BP_0$. For $\alpha$ sufficiently large, this defines a rational map $\BP(\sG)\to \BP^M\times W$ whose restriction to $\BP_0$
is an embedding. For $\beta\gg 1$ the multiplication map
$$
S^\beta(\bigoplus \sH^\alpha) \>>> \pi_*\sO_{\BP(\sG)}(\eta'\cdot\beta\cdot \alpha)=
S^{\eta'\cdot\beta\cdot \alpha}(\sG)
$$
will be surjective over $W_0$, hence $\sG$ ample with respect to $W_0$.
\end{proof}
For the compatibility of ``ample with respect to $W_0$'' under arbitrary finite morphisms
one either needs that the non-normal locus of $W_0$ is proper (see \cite[Proposition 2.22]{Vie}
and the references given there) or one has to add the condition ``nef'':
\begin{lemma}\label{ne.9}
For a locally free sheaf $\sG$ on a projective variety $W$, and for an open dense subscheme $W_0\subset W$
the following conditions are equivalent:
\begin{enumerate}
\item $\sG$ is nef and ample with respect to $W_0$.
\item There exists a finite morphism $\sigma:W'\to W$ such that $\sG'=\sigma^*\sG$ is nef and ample with respect to $W'_0=\sigma^{-1}(W_0)$.
\item There exists an alteration $\phi:\widetilde{W}\to W$ with $\phi^{-1}(W_0)\to W_0$ finite,
such that $\phi^*\sG$ is nef and ample with respect to $\widetilde{W}_0=\phi^{-1}(W_0)$.
\end{enumerate}
\end{lemma}
\begin{proof} Of course (1) implies (2) and (2) implies (3). In order to see that (3) implies (2)
choose for $\sigma:W'\to W$ the Stein factorization of $\phi:\widetilde{W}\to W$.
Since $\widetilde{W}\to W'$ is an isomorphism over $W'_0$~\ref{ne.7} (4) says that $\sG'=\sigma^*\sG$ is ample with respect to $W'_0$ if and only if $\phi^*\sG$ is ample with respect to $\widetilde{W}_0$. Since by~\ref{ne.3} the same holds for nef one obtains (2).

Remark that (2) implies that the sheaf $\sG$ is nef, as well as the sheaf $\sO_{\BP(\sG)}(1)$
on $\BP(\sG)$. Consider the induced morphism $\sigma':\BP(\sG')\to \BP(\sG)$. Lemma~\ref{ne.8} implies that $\sO_{\BP(\sG')}(1)=\sigma'^*\sO_{\BP(\sG)}(1)$ is ample with respect to the preimage of $W_0'$ if and only if
$\sG'$ is ample with respect to $W_0'$, and that the same holds for $\sG$ instead of $\sG'$.

So it will be sufficient to consider an invertible nef sheaf $\sG$ on $W$, and a finite covering
$\sigma:W'\to W$, such that $\sG'=\sigma^*\sG$ is ample with respect to $W'_0$, and we have to show that
$\sG$ is ample with respect to $W_0$. 

As we have already seen, that (1) implies (2), we may replace $W'$ by any dominating finite covering. In particular we may assume $W'$ to be normal. By \cite[Lemma 2.2]{Vie} the morphism
$\sigma:W'\to W$ factors like
$$
W' \> \gamma >> W'' \> \rho >> W,
$$
where $\rho$ has a splitting trace map, and where $\gamma$ is birational.

At this point we could also apply Lemma~\ref{ex.3}, replacing $W'$ by a larger normal covering.
In any case $\gamma^*\rho^*\sG$ is again ample with respect to $\gamma^{-1}\rho^{-1}(W_0)$ and
by~\ref{ne.7} (4) one knows  the equivalence of (1) and (2) with $W'$ replaced by $W''$. 
Hence it is sufficient to study $V\to W''$, and by abuse of notations we
may assume that $W'$ is normal and $\sigma$ birational.

Let $\xi:\hat{Y}\to W'$ be a desingularization, $\delta=\sigma\circ\xi:\hat{Y}\to W$
and let $U\subset W$ be the complement of the center of $\delta$.
Choose a sheaf of ideals $\sJ$ on $W$ with $\sO_W/\sJ$ supported in $W\setminus U$
and such that $\sigma_*\sigma^*\sJ$ maps to $\sO_W$. One can assume that
$\delta^*\sJ/{\rm torsion}$ is invertible hence of the form $\sO_{\hat{Y}}(-E)$ for an effective divisor supported in $\hat{Y}\setminus \delta^{-1}(U)$. Then $\delta_*\sO_{\hat{Y}}(-E)$ is contained in $\sO_W$. One may assume in addition that $\sO_{\hat{Y}}(-E)$ is $\delta$-ample. Finally choose an ample invertible sheaf $\sH$ on $W$, such that $\delta^*\sH\otimes\sO_{\hat{Y}}(-E)$ is ample and such that
$\sH\otimes \sigma_*\sO_{W'}$ is generated by global sections.

By assumption, for some $\eta$ there are morphisms
\begin{equation}\label{eqne.1}
\bigoplus \sigma^*\sH \>>> \sigma^*\sG^{\eta} \mbox{ \ \ and hence \ \ }
\bigoplus\delta^* \sH \>>> \delta^*\sG^\eta,
\end{equation}
surjective over $W'_0$ and $\delta^{-1}(W_0)$, respectively. Blowing up a bit more, we can assume that the image of the second map is of the form $\delta^*\sG^\eta\otimes \sO_{\hat{Y}}(-\Delta)$ for a divisor $\Delta$.
Then $\delta^*\sG^\eta\otimes \sO_{\hat{Y}}(-\Delta-E)$ as a quotient of an ample sheaf will be ample. Replacing
$\eta$, $\Delta$ and $E$ by some multiple, one may also assume that
$$
\delta^*\sG^\eta\otimes \sO_{\hat{Y}}(-\Delta-E)\otimes \omega_{\hat{Y}}^{-1} \otimes \delta^*\sH^{-1}
$$
is ample. Define $\sI'=\xi_*\sO_{\hat{Y}}(-\Delta-E))$ on $W'$ and $\sI=\sigma_*\sI'$. 

Since $\sG$ is nef, for all $\alpha \geq \eta$ and for all $\beta\geq -1$ the sheaf
$$
\delta^*\sG^\alpha\otimes \sO_{\hat{Y}}(-\Delta-E)\otimes \delta^*\sH^\beta
$$
has no higher cohomology. For $\beta \gg 1$ this can only hold if for all $i>0$
$$
R^i\delta_*(\delta^*\sG^\alpha\otimes \sO_{\hat{Y}}(-\Delta-E))=0.
$$
For $\beta=-1$ one finds that $H^i(W,\sG^\alpha\otimes \sH^{-1} \otimes \sI)=0$.

For some $\beta\gg 1$ the sheaf $\sigma^*\sH^{\beta-2}\otimes \sI'$ is generated by global sections.
Using the left hand side of (\ref{eqne.1}) one obtains a morphism 
$$
\bigoplus \sigma^*\sH^{\beta}\otimes \sI' \>>> \sigma^*\sG^{\eta\cdot\beta}\otimes \sI',
$$
surjective over $W'_0$. Therefore the sheaf $\sigma^*(\sG^{\eta\cdot \beta}\otimes \sH^{-2}) \otimes \sI'$ will be globally generated over $W'_0$ hence there are morphisms
$\bigoplus\sigma^*\sH \to \sigma^*\sH^{\beta-1}\otimes \sI'$ and
\begin{equation}\label{eqne.2}
\bigoplus \sH \otimes \sigma_*\sO_{W'} \>>> \sG^{\eta\cdot \beta}\otimes \sH^{-1} \otimes \sI,
\end{equation}
surjective over $W'_0$ and $W_0$. By the choice of $\sH$ the left hand side of \eqref{eqne.2} is globally generated over $W_0$, hence the right hand side as well. For all positive multiples $\alpha$ of $\eta \cdot \beta$, one has an exact sequence
$$
0\to H^0(W,\sG^\alpha\otimes \sH^{-1}\otimes \sI) \>>>
H^0(W,\sG^\alpha\otimes \sH^{-1}) \>>> H^0(W,\sG^\alpha\otimes \sH^{-1}|_{T'})\to 0,
$$
where $T'$ denotes the subscheme of $W$ defined by $\sI$. If $T'\cap W_0=\emptyset$ we are done. 
Otherwise let $T$ be the closure of $T'_{\rm red}\cup W_0$ in $W$. So there is a coherent sheaf $\sF$, supported on $T$ and an inclusion $\sF\to \sO_{T'}$ which is an isomorphism on $W_0\cap T'=W_0\cap T$. 

By induction on the dimension of $W$ we may assume that $\sG|_T$ is ample with respect to $T\cap W_0$. Then for each $\beta'>0$ one finds $\eta'$ and morphisms
$$
\bigoplus \sH^{\beta'-1}|_T \>>> (\sG^{\eta'\cdot\beta'}\otimes \sH^{-1})|_T,
$$
surjective over $Z\cap W_0$. Choose $\beta'$, such that $\sF\otimes \sH|_T^{\beta'-1}$ is globally generated, and $\alpha=\eta'\cdot\beta'$ a multiple of $\eta\cdot\beta$. Then the sheaf
$(\sG^{\alpha}\otimes \sH^{-1})|_T \otimes \sF$
is globally generated over $T\cap W_0$, as well as $\sG^\alpha\otimes \sH^{-1}|_{T'}$.
Since all global sections of this sheaf lift to $H^0(W,\sG^\alpha\otimes \sH^{-1})$
we find that $\sG^{\alpha}\otimes \sH^{-1}$ is globally generated over $W_0$.
\end{proof}
\section{Positivity of direct images I}\label{pd1}
Examples of nef sheaves are direct images of powers of
dualizing sheaves. We will see in this section that Theorem~\ref{in.1} is just what is needed to verify this property in Theorem~\ref{in.2} for $\sF^{(\nu)}_{\hat{Y}}$. Then the compatibility $\sF_{\hat{Y}}^{(\nu)}=\xi^*\sF_W^{(\nu)}$ 
allows to deduce the nefness of $\sF_W^{(\nu)}$. At the end of the section we will make a first step towards
Theorem~\ref{in.2} iv), assuming that $W_0$ is non-singular. This will allow in~\ref{pd1.5}
to explain why we have to include the study of multiplier ideals in Sections~\ref{fm}--\ref{va}.

Let us recall the definition of multiplier ideal sheaves and of the corresponding threshold. 
\begin{definition}\label{pd1.1}
Let $Z$ be a normal projective variety with at most rational Gorenstein singularities, let $\sN$ be an invertible sheaf on $Z$ and let $D$ be the zero divisor of a section of $\sN$. 
\begin{enumerate}
\item[i.] For $a\in \Q$ the {\em multiplier ideal} is defined as 
$\sJ(-a\cdot D)=\tau_*\omega_{\widetilde{Z}/Z}\otimes\sO_{\widetilde{Z}}(-[a\cdot \widetilde{D}])$,
where $\tau:\widetilde{Z}\to Z$ is a log-resolution, $\widetilde{D}=\tau^*D$, and 
where $[a\cdot \widetilde{D}]=\rounddown{a\cdot \widetilde{D}}$ denotes the integral part of the $\Q$-divisor $a\cdot \widetilde{D}$. 
\item[ii.] For $b>0$ one defines the {\em threshold}
$$
e(b\cdot D):={\rm Min}\big\{a\in \Z_{>0}; \ \sJ\big(-\frac{b}{a}\cdot D\big)=\sO_Z\big\}
$$
\item[iii.] Finally $e(\sN):={\rm Max}\big\{e(D); \ D \mbox{ the zero divisor of a section of }\sN\big\}$.
\end{enumerate}
\end{definition}
Most authors write $\sJ(a\cdot D)$ instead of $\sJ(-a\cdot D)$. We prefer the second notion, since
for a smooth divisor $\sO_Z(-D)=\sJ(-D)$, and since it is closer to the classical notion
$\omega_Z\{-a\cdot D\} = \omega_Z\otimes \sJ(-a\cdot D)$
used in \cite{Vie} and \cite{EV}.

One easily shows that the multiplier ideal is independent of the log resolution. 
In \cite{EV} and \cite{Vie} one finds a long list of properties of multiplier ideals and of $e(b\cdot D)$ and $e(\sN)$. In particular, if $\sN$ is a globally generated invertible sheaf on $Z$, then 
\begin{equation}\label{eqpd1.1}
\sJ(-a\cdot D)=\sJ(-a\cdot (D+H))
\end{equation} for the divisor $H$ of a general section of $\sN$ and for $0\leq a <1$. In fact, using the notation introduced above, $\widetilde{H}=\tau^*H$ will be non-singular and it intersects $\widetilde{D}$ transversely.
Then $[a\cdot \widetilde{D}]=[a\cdot (\widetilde{D}+\widetilde{H})]$. 

Let $f_0:X_0\to Y_0$ be a flat morphism over a non-singular variety $Y_0$, with irreducible normal fibres with at most rational singularities. Then for a $\Q$-divisor $\Delta_0$ on $X_0$ not containing fibres, the threshold $e(\Delta|_{f_0^{-1}(y)})$ is upper semi-continuous for the Zariski topology (see \cite[Proposition 5.17]{Vie}). This implies in \cite[Corollary 5.21]{Vie}
that for $Z$ and $\sN$ as in Definition \ref{pd1.1} one has:
\begin{equation}\label{eqpd1.2}
e(\sN)=e({\rm pr}_1^*\sN\otimes \cdots \otimes {\rm pr}_r^*\sN) \mbox{ \ \ for \ \ } Z^r=Z\times \cdots \times Z. 
\end{equation}
If one replaces $\Z_{>0}$ in Definition \ref{pd1.1} ii) by $\Q_{>0}$ one obtains 
the inverse of the logarithmic threshold. 

The multiplier ideals occur in a natural way as direct images of relative dualizing sheaves for certain alterations:
\begin{lemma}\label{pd1.2}
Let $Z'$ and $Z$ be normal with rational Gorenstein singularities and let $\phi:Z'\to Z$ be an alteration. 
If $\sO_Z(D)=\sL^N$ for an invertible sheaf $\sL$ and if $\phi^*D$ is divisible by $N$, then $\sJ\big(-\frac{1}{N}\cdot D\big)$ is a direct factor of $\sL^{-1}\otimes \phi_*\omega_{Z'/Z}$.    
\end{lemma}
\begin{proof}
The sheaf $\phi_*\omega_{Z'/Z}$ does not change, if we replace $Z'$ by a non-singular modification.
So we may assume that $Z'$ is non-singular and that it dominates a log resolution
$\tau:\widetilde{Z}\to Z$ for $D$. Writing $\pi:Z'\to \widetilde{Z}$ for the induced morphism,
$\pi^*(\tau^* D)$ is still divisible by $N$. So $\pi$ factors through the 
cyclic covering $\widetilde{\pi}:\widetilde{Z}'\to \widetilde{Z}$, obtained by taking the $N$-th root out of $\tau^* D$.
By \cite[Section 3]{EV} the sheaf 
$$
\tau^*\sL\otimes \omega_{\widetilde{Z}/Z}\otimes \sO_{\widetilde{Z}}\big(-\big[\frac{1}{N}\cdot\tau^*D\big]\big)
$$
is a direct factor of $\widetilde{\pi}_*\omega_{\widetilde{Z}'/Z}$. The latter is a direct factor
of $\pi_*\omega_{Z'/Z}$. Applying $\tau_*$ one obtains $\sL\otimes \sJ\big(-\frac{1}{N}\cdot D\big)$
as a direct factor of $\phi_*\omega_{Z'/Z}$.
\end{proof}
One starting point for the study of positivity of direct image sheaves is the following 
corollary of Koll\'ar's Vanishing Theorem.
\begin{lemma}\label{pd1.3} Let $X$ be a projective normal variety with at most rational Gorenstein singularities, let
$f:X \to Y$ be a surjection to a projective $m$-dimensional variety $Y$, and let
$U\subset Y$ be open and dense. Let $\sA$ be a very ample invertible sheaf on $Y$, let $\sM$ be an invertible sheaf on $X$ let $\Gamma$ be an effective divisor, and let
$\sE$ be a locally free sheaf on $Y$, weakly positive over $U$. Assume that for some $N>0$ there is a morphism $\sE\to f_*\sM^N(-\Gamma)$ for which the composite
$$
f^*\sE\>>> f^*f_*\sM^N(-\Gamma)\>>> \sM^N(-\Gamma)
$$
is surjective over $V=f^{-1}(U)$. Then for all $\beta$ the sheaf
$$
\sA^{m+2}\otimes f_*\Big(\sM^\beta \otimes \omega_{X}\otimes \sJ\big(-\frac{\beta}{N}\Gamma\big)\Big)
$$
is globally generated over $U$.
\end{lemma}
\begin{proof} We can replace $X$ by a desingularization.
The sheaf $\sA^N \otimes \sE$ is ample with respect to $U$, hence for some $M>0$
the sheaf $\sA^{N\cdot M} \otimes S^M(\sE)$ is globally generated over $U$. Blowing up $X$ with centers outside of $V$ we may assume that the image $\sB$ of the evaluation map
$f^*S^M(\sE) \to \sM^{N\cdot M}(-M\cdot\Gamma)$ is invertible. Let $D$ be the divisor, supported in $X\setminus V$
with $\sB\otimes \sO_{X}(D)=\sM^{N\cdot M}(-M\cdot \Gamma)$. Then 
$$
\sB\otimes f^*(\sA^{N\cdot M})=\sM^{N\cdot M}(-M\cdot \Gamma-D)\otimes f^*(\sA^{N\cdot M})
$$ 
is generated by global sections over $V$. Blowing up again, we find a divisor $\Delta$ supported 
in $X\setminus V$ such that $\sM^{N\cdot M}(-M\cdot \Gamma-D-\Delta)\otimes f^*(\sA^{N\cdot M})$
is generated by global sections, and such that $\Gamma+D+\Delta$ is a normal crossing divisor.

$\sM^{N\cdot M}(-M\cdot \Gamma-D-\Delta)\otimes f^*\sA^{N\cdot M}$ is semiample.
As in \cite[Corollary 2.37, 2)]{Vie} Koll\'ar's Vanishing Theorem implies that the sheaf
$$
\sA^\iota\otimes f_*\Big(\sM^\beta\big(-\big[\frac{\beta}{N\cdot M}(M\cdot\Gamma-D-\Delta)\big]\big)\otimes \omega_{X} \otimes \hat{f}^*\sA\Big)
$$
has no higher cohomology for $\iota>1$. Then by an argument due to N. Nakayama (see \cite[Lemma 2.11]{Kaw})
$$
\sP=\sA^{m+1}\otimes f_*\Big(\sM^\beta\big(-\big[\frac{\beta}{N\cdot M}(M\cdot\Gamma-D-\Delta)\big]\big)\otimes \omega_{X} \otimes \hat{f}^*\sA\Big)
$$
is generated by global sections. On the other hand, $\sP$ is contained in 
$$
\sA^{m+2}\otimes f_*\Big(\sM^\beta\big(-\big[\frac{\beta}{N}\Gamma\big]\big)\otimes \omega_{X}\Big),
$$ 
and since $(D+\Delta)\cap V=\emptyset$, both coincide over $U$. 
\end{proof}
Kawamata's Semipositivity Theorem, saying that the direct image sheaves $\sF_{\hat{Y}}^{(1)}$ in Theorem~\ref{in.2} iii) are nef, can be shown by using Lemma~\ref{pd1.3}. This in turn implies part iii) for all $\nu$. We will give a slightly different argument:
\begin{proof}[Proof of Theorem~\ref{in.2} iii)] 
Let $\hat{Y}\to W$ be any non-singular alteration of $W$.
In Theorem~\ref{in.1} ii) the sheaf $f^{(r)}_*\omega_{X^{(r)}/\hat{Y}}^\nu$ remains the same if one replaces $X^{(r)}$ by a non-singular modification. Hence by abuse of notations one may assume that for some normal crossing divisor $\Pi$ on $\hat{X}$ the evaluation map induces a surjection
$$
{f^{(r)}}^*{\sF_{\hat{Y}}^{(\nu)}}^{\otimes r} \>>> \omega_{X^{(r)}/\hat{Y}}^\nu \otimes \sO_{X^{(r)}}(-\Pi).
$$
Let $\sH$ be an ample invertible sheaf on $\hat{Y}$ and define
$$
s(\nu) = {\rm Min} \big\{ \mu >0 ; {\sF_{\hat{Y}}^{(\nu)}}\otimes
\sH^{\nu\cdot\mu} \mbox{ \ is nef} \big\} .
$$
So $(\sH^{s(\nu)\cdot\nu} \otimes {\sF_{\hat{Y}}^{(\nu)}})^{\otimes r}=\sH^{s(\nu)\cdot\nu\cdot r} \otimes \hat{f}_*\omega_{X^{(r)}/\hat{Y}}^\nu$ is nef.
Let $\sA$ be a very ample invertible sheaf on $\hat{Y}$. By~\ref{pd1.3}
$$
\sA^{m+2}\otimes  f^{(r)}_*\Big(\omega_{X^{(r)}}\otimes
\big(\omega_{X^{(r)}/\hat{Y}} \otimes {f^{(r)}}^*\sH^{s(\nu)\cdot r}\big)^{\nu-1}
\otimes \sO_{\hat{Y}} \big(-\big[\frac{(\nu-1)\cdot \Pi}{\nu }\big]\big)\Big)
$$
is generated by global sections. It is a subsheaf of
$\sA^{m+2}\otimes \omega_{\hat{Y}} \otimes {\sF_{\hat{Y}}^{(\nu)}}^{\otimes r} \otimes \sH^{s(\nu)\cdot r \cdot (\nu-1)}$ and it contains the sheaf
\begin{multline*}
\sA^{m+2}\otimes \omega_{\hat{Y}}\otimes \sH^{s(\nu)\cdot r \cdot (\nu-1)}\otimes
f^{(r)}_*(\omega_{X^{(r)}/\hat{Y}}^\nu \otimes \sO_{\hat{Y}} (-\Pi) )=\\
\sA^{m+2}\otimes \omega_{\hat{Y}} \otimes \sH^{s(\nu)\cdot r\cdot (\nu-1)}\otimes
{\sF_{\hat{Y}}^{(\nu)}}^{\otimes r}.
\end{multline*}
So the three sheaves are equal, and the quotient sheaf
$$
\sA^{m+2}\otimes \omega_{\hat{Y}}\otimes S^r (
\sH^{s(\nu)(\nu-1)}\otimes \sF_{\hat{Y}}^{(\nu)}).
$$
is generated by global sections as well. Hence $\sH^{s(\nu)\cdot (\nu-1)}\otimes \sF_{\hat{Y}}^{(\nu)}$
is weakly positive over $\hat{Y}$. Since $\sH^{(s(\nu)-1)\cdot\nu} \otimes \sF_{\hat{Y}}^{(\nu)}$
does not have this property, one obtains
$$
s(\nu)\cdot(\nu-1) > (s(\nu)-1)\cdot \nu
\mbox{ \ \ or \ \ }
s(\nu)  < \nu.
$$
So $\sH^{\nu^2}\otimes \sF_{\hat{Y}}^{(\nu)}$ is weakly positive over $\hat{Y}$, hence nef.

Since the same exponent $\nu^2$ works for all $\hat{Y}$ mapping to $W$ and for all
ample invertible sheaves $\hat{\sH}$ on $\hat{Y}$, the nefness of $\sF_{W}^{(\nu)}$
follows from the equivalence of (1) and (5) in~\ref{ne.3}.
\end{proof}
\begin{variant}\label{pd1.4}
Assume in Theorem ~\ref{in.1} that the normalization of $W_0$ is non-singular and that for some $\eta>0$ the sheaf $\det(\sF_W^{(\eta)})$ is ample with respect to $W_0$. Then for $\nu\geq 2$ 
the sheaf $\sF_W^{(\nu)}$ is ample with respect to $W_0$ or zero.
\end{variant}
\begin{proof}
Let us sketch the argument. Knowing that $\sF_{W}^{(\nu)}$ is nef, Lemma~\ref{ne.9} allows to replace $W$ by a desingularization $\hat{Y}$. Consider for some $\gamma$ and for $r=\gamma\cdot \rk(\sF_{\hat{Y}}^{(\eta)})$ the tautological map
$$
\Xi:\det(\sF_{\hat{Y}}^{(\eta)})^\gamma \>>> \bigotimes^{r} \sF_{\hat{Y}}^{(\eta)}=f^{(r)}_*\omega_{X^{(r)}/\hat{Y}}^{\eta}.
$$
Assume that $\det(\sF_{\hat{Y}}^{(\eta)})=\sN^\alpha$ for some invertible sheaf $\sN$ and for some
$\alpha >0$. Then $\Xi$ induces a section of $\omega_{X^{(r)}/\hat{Y}}^{\eta}\otimes {f^{(r)}}^*\sN^{-\gamma\cdot \alpha}$ with zero divisor $\Gamma$.

We know already that $\sF_{\hat{Y}}^{(\rho)}$ is nef. So we can apply Lemma~\ref{pd1.3} to
the sheaf 
$$
\omega_{X^{(r)}/\hat{Y}}^{\rho} = \omega_{X^{(r)}/\hat{Y}}^{\rho + \eta}\otimes {f^{(r)}}^*\sN^{-\gamma\cdot \alpha}\otimes 
\sO_{X^{(r)}}(-\Gamma).
$$ 
If $\alpha$ is divisible by $\rho + \eta$, for $\sA$ very ample on $\hat{Y}$, the sheaf
$$
\sA^{\dim(\hat{Y})+2} \otimes \omega_{\hat{Y}} \otimes f^{(r)}_*\Big(\omega_{X^{(r)}/\hat{Y}}^{\nu} \otimes {f^{(r)}}^*\sN^{-\frac{(\nu-1)\gamma\cdot \alpha}{\rho+\eta}} \otimes 
\sJ\big(-\frac{\nu-1}{\rho+\eta}\Gamma\big)\Big)
$$
is globally generated over the preimage $\hat{Y}_0$ of $W_0$. For $\rho$ sufficiently large one may assume that
$\frac{\rho+\eta}{\nu-1} \geq e(\omega_F^\eta)$, for all smooth fibres of $f_0$. Then by \cite[Lemma 5.14 and Corollary 5.21]{Vie} one finds that $\sJ\big(-\frac{\nu-1}{\rho+\eta}\Gamma\big)$ is trivial over ${f^{(r)}}^{-1}(\hat{Y}_0)$. Then the inclusion of sheaves
\begin{equation}\label{eqpd1.3}
\sG(r):=f^{(r)}_*\Big(\omega_{X^{(r)}/\hat{Y}}^{\nu} \otimes\sJ\big(-\frac{\nu-1}{\rho+\eta}\Gamma\big)\Big)
\subset \bigotimes^r \sF_{\hat{Y}}^{(\nu)}
\end{equation}
is an isomorphism over $\hat{Y}_0$ and
\begin{multline*}
\sA^{\dim(\hat{Y})+2} \otimes \omega_{\hat{Y}} \otimes
\sN^{-\frac{(\nu-1)\gamma\cdot \alpha}{\rho+\eta}}\otimes 
\bigotimes^r \sF_{\hat{Y}}^{(\nu)}=\\
\sA^{\dim(\hat{Y})+2} \otimes \omega_{\hat{Y}} \otimes
\bigotimes^\gamma \Big( \sN^{-\frac{(\nu-1)\alpha}{\rho+\eta}}\otimes \bigotimes^{\rk(\sF_{\hat{Y}}^{(\nu)})} \sF_{\hat{Y}}^{(\nu)}\Big)
\end{multline*}
is globally generated over $\hat{Y}_0$. This being true for all $\gamma$, the Lemmata~\ref{ne.6} and~\ref{ne.5} imply that
$\sF_{\hat{Y}}^{(\nu)}$ is ample over $\hat{Y}_0$.

By \cite[Lemma 2.1]{Vie} the assumption that $\det(\sF_{\hat{Y}})=\sN^\alpha$,
with $\alpha$ divisible by $\rho+\eta$, will always be true over sufficiently large 
finite coverings of $\hat{Y}$, and by Lemma~\ref{ne.7} we are done. 
\end{proof}
\begin{comments}\label{pd1.5}

We repeated the well known proof of Variant~\ref{pd1.4} just to point out the difficulties we will encounter, trying to get rid of the additional assumption ``$W_0$ non-singular''. 
The notion ``ample with respect to $W_0$'' is not compatible with 
blowing ups of $W_0$, if the center meets $W_0$. We may assume that $\sN$ is the pullback of an invertible sheaf $\sN_W$ on $W$. In addition we have to construct a sheaf $\sG_W(r)$ whose pullback to a desingularization is the sheaf $\sG(r)$ considered in the proof of Variant~\ref{pd1.4}. 
In order to be allowed to use the functorial property~\ref{ne.3} the sheaf 
$$
\sG(r)\otimes \sN^{-\frac{(\nu-1)\gamma\cdot \alpha}{\rho+\eta}}
$$ 
must be nef, and not just weakly positive over $\hat{Y}_0$. This would hold, if the inclusion \eqref{eqpd1.3} is an isomorphism, but giving bounds for the threshold in bad fibres of a morphism does not seem to work.

So as a way out we will modify the construction of $W$ in such a way, that Theorem~\ref{in.1} remains true for the direct images $\sG_\bullet(r)$ of invertible sheaves tensorized by multiplier ideals.
\end{comments}
\section{On the construction of moduli schemes}\label{co}
The weak positivity and ampleness of the direct image sheaves in Corollary~\ref{in.3}, together with the stability criterion \cite[Theorem 4.25]{Vie}, allows the construction of a quasi-projective moduli scheme of canonically polarized manifolds. Following a suggestion of the referee, we sketch the argument before entering the quite technical 
details needed for the construction of $W$, hence for the proof of Corollary~\ref{in.3}, b).

Let $\fM_h$ be the moduli functor of canonically polarized manifolds with Hilbert polynomial $h$.
So as in \cite[Examples 1.4]{Vie} we consider for a scheme $Y_0$ the set
\begin{multline*}
\fM_h(Y_0)=\big\{f_0:X_0\to Y_0; \ f_0 \mbox{ smooth, projective, } \omega_{X_0/Y_0} \ f_0\mbox{-ample}\\
\mbox{and } h(\nu)=\rk(f_{0*}\omega_{X_0/Y_0}^\nu) ,\mbox{ for } \nu \geq 2 \big\}/_\sim.
\end{multline*}
In order to allow the canonical models of surfaces we could also consider 
\begin{multline*}
\fM'_h(Y_0)=\big\{f_0:X_0\to Y_0; \ f_0 \mbox{ flat, projective; all fibres }F\mbox{ normal}\\
\hspace*{1.5cm}\mbox{ with at most rational Gorenstein singularities, } \omega_{X_0/Y_0} \\
 f_0\mbox{-ample and } h(\nu)=\rk(f_{0*}\omega_{X_0/Y_0}^\nu) ,\mbox{ for } \nu \geq 2 \big\}/_\sim.
\end{multline*}
We leave the necessary changes of the arguments to the reader.\\[.1cm]
{\em Outline of the construction of a coarse quasi-projective moduli scheme $M_h$ for $\fM_h$.}\ \\
One first has to verify that the functor $\fM_h$ is a nice moduli functor, i.e. locally closed, separated and bounded (see \cite[Lemma 1.18]{Vie}). This implies that for some $\eta \gg 1$ 
one has the Hilbert scheme $H$ of $\eta$-canonically embedded manifolds in $\fM_h({\rm Spec}(\C))$, together with the universal family $g:\sX\to H$.

The universal property gives an action of $G=\BP{\rm Gl}(h(\eta))$ on $H$ and, as explained in \cite{Mum} or \cite[Lemma 7.6]{Vie}, the separatedness of the moduli functor implies that this action is proper and with finite stabilizers. The sheaves $\lambda_\eta=\det(g_*\omega_{\sX/H}^\eta)$ are all $G$-linearized.

The moduli scheme $M_h$, if it exists, should be a good quotient $H/G$. So one has to verify that all points in $H$ are stable for the group action and for a suitable ample sheaf. At this point one is allowed to
replace $H$ by $H_{\rm red}$; the set of stable points will not change. So by abuse of notations
we will assume that $H$ (and hence $M_h$) is reduced.

In order to apply the stability criterion \cite[Theorem 4.25]{Vie} one has to verify that the invertible sheaf
$\lambda_\eta$ on $H$ is ample on $H$, and that for a certain family $f_0:X_0\to Y_0$ in $\fM_h(Y_0)$
the sheaf $f_{0*}\omega_{X_0/Y_0}^\eta$ is weakly positive over $Y_0$.

The second statement follows from Corollary~\ref{in.3}, a). For the first one we start with the Pl\"ucker embedding showing that the invertible sheaves $\lambda_{\eta\cdot \mu}^{h(\eta)}\otimes \lambda_{\eta}^{-h(\eta\cdot\mu)\cdot\mu}$ are ample, for all $\mu$ sufficiently large. 
By Corollary~\ref{in.3} a) the sheaf $\lambda_{\eta}$ is weakly positive over $H$, hence by Lemma \ref{ne.3} (3) $\lambda_{\eta\cdot\mu}$ is ample. Using Corollary~\ref{in.3} b) one finds that the sheaf $g_*\omega_{\sX/H}^\eta$ is ample on $H$, hence its determinant $\lambda_{\eta}$ as well.
\qed\vspace{.1cm}

Let us express what we have shown in terms of stability of Hilbert points. On $H$ the sheaf
$\lambda_\eta$ is $G$ linearized and ample. The stability criterion says that 
all the points in $H$ are stable with respect to the polarization $\lambda_\eta$ of $H$. 
This in turn shows the ampleness of the sheaf $\lambda^{(p)}_{0,\eta}$ for $\eta \gg 1$.

One can consider the sheaf $\lambda_{\nu}$ on $H$ for all $\nu\geq 2$ with $h(\nu)>0$.
Those sheaves are $G$-linearized and for some $p>0$ the $p$-th power of $\lambda_\nu$ descends to
an invertible sheaf $\lambda_{0,\nu}^{(p)}$ on $M_h$. Using a slightly different stability criterion, stated in \cite[Addendum 4.26]{Vie}, one obtains:
\begin{addendum}\label{co.1}
For all $\nu\geq 2$ with $h(\nu)>0$ and for some $p>0$ there exists an ample invertible sheaf $\lambda_{0,\nu}^{(p)}$ on $M_h$ whose pullback to $H$ is $\lambda_\nu=\det(h_*\omega_{\sX/H}^\nu)$. In particular the sheaf $\lambda_{0,\nu}^{(p)}$ will satisfy the condition ($\divideontimes$) stated in the Notations~\ref{in.4}.
\end{addendum}
In Section~\ref{pr} we will even show that $\lambda_{0,\nu}^{(p)}$ extends to an invertible sheaf $\lambda_{\nu}^{(p)}$ on a suitable compactification of $M_h$ and that this sheaf is ample with respect to $M_h$.

For points of the Hilbert scheme of $\eta$-canonically embedded curves or surfaces of general type
the stability has been verified with respect to the Pl\"ucker embedding (see \cite{Mum} and \cite{Gie}). So one obtains
$M_h$ the ampleness of $\lambda_{0,\eta\cdot\mu}^{h(\eta)} \otimes \lambda_{0,\eta}^{-\mu\cdot h(\eta \cdot\mu)}$.\vspace{.1cm} 

Before turning our attention to moduli schemes of polarized minimal models, let us formulate the
generalization of Corollary~\ref{in.3}, needed for their construction. The proof
will be given at the end of Section~\ref{pd2}.
Here we use again the threshold $e(\sN)$ defined in~\ref{pd1.1}.
\begin{lemma}\label{co.2}
Let $f_0:X_0\to Y_0$ be a smooth family of minimal models, and let $\sL_0$ be an $f_0$-ample invertible sheaf. Assume that for some $\kappa >0$ the direct image $f_{0*}(\sL_0^\kappa)$ is non-zero, locally free 
and compatible with arbitrary base change. Choose some $\epsilon > e(\sL_0^\kappa|_F)$, for all fibres $F$ of $f_0$. Then:
\begin{enumerate}
\item For all positive integers $\eta$ the sheaf
$$
S^{\rk(f_{0*}(\sL_0^\kappa))}(f_{0*}(\omega_{X_0/Y_0}^{\epsilon\cdot\eta}\otimes \sL_0^{\kappa\cdot\eta}))\otimes \det(f_{0*}(\sL_0^\kappa))^{-\eta}
$$
is weakly positive over $W_0$ or zero. 
\item If for some $\eta'>0$ the sheaf
$$
\det(f_{0*}(\omega_{X_0/Y_0}^{\epsilon\cdot\eta'}\otimes \sL_0^{\kappa\cdot\eta'}))^{\rk(f_{0*}(\sL_0^\kappa))} \otimes \det(f_{0*}(\sL_0^\kappa))^{-\eta'\cdot \rk(f_{0*}(\omega_{X_0/Y_0}^{\epsilon\cdot\eta'}\otimes
\sL_0^{\kappa\cdot\eta'}))}
$$
is ample, then $S^{\rk(f_{0*}(\sL_0^\kappa))}(f_{0*}(\omega_{X_0/Y_0}^{\epsilon\cdot\eta}\otimes \sL_0^{\kappa\cdot\eta}))\otimes \det(f_{0*}(\sL_0^\kappa))^{-\eta}$ is ample, if not zero.
\end{enumerate}
\end{lemma}
The moduli functor $\fM_h$ of minimal polarized manifolds is given by 
\begin{multline*}
\fM_h(Y_0)= 
\big\{ (f_0:X_0\>>> Y_0, \sL_0) ; \ f_0 \mbox{ smooth, projective; } \omega_{X_0/Y_0} \ f_0\mbox{-semiample; }\\ 
\sL_0 \ f_0\mbox{-ample, with Hilbert polynomial } h\big\}/_{\sim}.
\end{multline*}
Recall that $(f_0:X_0\to Y_0, \sL_0)\sim (\widetilde{f}_0:\widetilde{X}_0\to Y_0, \widetilde{\sL}_0)$
if there is an $Y_0$-isomorphism $\iota:X_0\to \widetilde{X}_0$ and an invertible sheaf $\sA$ on $Y_0$ with
$\iota^*\widetilde{\sL}_0=\sL_0\otimes f_0^*\sA$. 

As we will see it is easier to study the moduli functor $\fM'_{h}$ with
\begin{multline*}
\fM'_{h}(Y_0)=
\big\{ (f_0:X_0\>>> Y_0, \sL_0)\in \fM_h; \sL_0 \ f_0\mbox{-very ample}\\
\mbox{with Hilbert polynomial } h; \ R^if_{0*}\sL_0^\mu=0 \mbox{ for } i>0, \mbox{ and } \mu>0,\big\}/_{\sim}.
\end{multline*}
For families of minimal varieties $F$ of Kodaira dimension zero the second condition will hold automatically.
In fact, if $\omega_F^\nu=\sO_F$ and if $\sA$ is ample, $\sA\otimes \omega_F^{\nu-1}$ is ample and Kodaira's Vanishing Theorem implies that $H^i(F,\sA)=H^i(F,\sA\otimes \Omega^\nu)=0$, for $i>0$. So here we should consider the functors $\fM_h^{(\nu)}$ with
\begin{multline*}
\fM^{(\nu)}_{h}(Y_0)=
\big\{ (f_0:X_0\>>> Y_0, \sL_0)\in \fM_h; \ f_0^*f_{0*}\omega_{X_0/Y_0}^\nu=\omega_{X_0/Y_0}; \\ 
\sL_0 \ f_0\mbox{-very ample with Hilbert polynomial } h\big\}/_{\sim}.
\end{multline*}
\begin{lemma}\label{co.3} \ 
\begin{enumerate}
\item Assume that for all $\hat{h}$ the moduli functor $\fM'_{\hat{h}}$ has a coarse quasi-projective moduli scheme 
$M'_{\hat{h}}$. Then the same holds true $\fM_{h}$.
\item To prove Theorem~\ref{in.6} it is sufficient to consider the moduli functors $\fM^{(\nu)}_{\hat{h}}$. 
\end{enumerate}
\end{lemma}
\begin{proof}
The boundedness of the moduli functor $\fM_h$ allows to find some $\gamma_0$ such that
for all $(F,\sA)\in \fM_h(\C)$ and for all $\gamma \geq \gamma_0$ the sheaf $\sA^\gamma$ is very ample and without higher cohomology. For suitable polynomials $h_1$ and $h_2$ one defines a map
\begin{gather*}
\fM_{h}\>>> \fM'_{h_1}\times \fM'_{h_2},
\mbox{ \ \ \ \ by}\\
(f_0:X_0\to Y_0,\sL_0) \mapsto [(f_0:X_0\to Y_0,\sL_0^{\gamma_0}),(f_0:X_0\to Y_0,\sL_0^{\gamma_0+1})].
\end{gather*}
It is easy to see that the image is locally closed. Hence if one is able to construct the corresponding moduli schemes $M'_{h_1}$ and $M'_{h_2}$ as quasi-projective schemes, $M_{h}$ is a locally closed subscheme. 
And if one finds a nice projective compactifications $\overline{M}'_{h_1}$ and $\overline{M}'_{h_2}$
of $M'_{h_1}$ and $M'_{h_2}$, one chooses $\overline{M}_{h}$ as the closure of $M_{h}$ 
in $\overline{M}'_{h_1}\times\overline{M}'_{h_2}$.
 
The additional condition $\omega_F^\nu=\sO_F$ considered in Theorem~\ref{in.6} just 
signalizes certain irreducible components of $M_h$. So by abuse of notations let
$M_h$ be one of those. Then the image of $M_h$ lies in the product
${M}^{(\nu)}_{h_1}\times {M}^{(\nu)}_{h_2}$. If one has constructed the compactifications
$\overline{M}^{(\nu)}_{h_1}$ and $\overline{M}^{(\nu)}_{h_2}$ according to Theorem~\ref{in.6}
one can choose $\overline{M}_h$ as the closure of $M_h$ and for $\lambda_\nu^{(2\cdot p)}$ 
the restriction of the exterior tensor product of the corresponding sheaves on $\overline{M}'_{h_1}$ and $\overline{M}'_{h_2}$ for $p$ instead of $2\cdot p$.
\end{proof}
\noindent
{\em Outline of the construction of coarse quasi-projective moduli scheme $M'_{h}$ for $\fM'_h$}\\
The construction is parallel to the one in the canonically polarized case. One constructs the Hilbert scheme $H$ parametrizing the elements $(F,\sA)$ of $\fM'_{h}(\C)$ together with an isomorphism
$\BP(H^0(F,\omega_F^\epsilon\otimes\sA))\cong \BP^{N}$. Here $\nu$ is chosen such that $\omega_F^\nu$ is globally generated and $\epsilon$ should be a multiple of $\nu$, larger than the threshold $e(\sA)$.

The Pl\"ucker embedding provides us with an ample invertible sheaf of the form
$\varpi_\mu^{r(1)}\otimes \varpi_1^{-\mu\cdot r(\mu)}$, where 
$$
\varpi_\nu=\det(g_*(\omega_{\sX/H}^{\epsilon\cdot\nu}\otimes\sL_{\sX}^\nu))\mbox{ \ \ and \ \ }
r(\nu)=\rk(g_*(\omega_{\sX/H}^{\epsilon\cdot\nu}\otimes\sL_{\sX}^\nu))
$$ 
for the universal family $(g:\sX\to H,\sL_{\sX})\in \fM'(H)$. 

By~\ref{co.2} (1) the sheaf $\varpi_1^{h(1)}\otimes \det(g_*\sL_\sX)^{-r(1)}$ is weakly positive over $H$, hence
$$
\varpi_\mu^{h(1)\cdot r(1)}\otimes \det(g_*\sL_\sX)^{-r(1)\cdot\mu\cdot r(\mu)}=
(\varpi_\mu^{h(1)}\otimes \det(g_*\sL_\sX)^{-\mu\cdot r(\mu)})^{r(1)}
$$ 
is ample. Using~\ref{co.2} (2) one finds that $\varpi_1^{h(1)}\otimes \det(g_*\sL_\sX)^{-r(1)}$ must be ample.

In order to apply the stability criterion \cite[Theorem 4.25]{Vie} to obtain the stability of all points of $H$ with respect
to the sheaf $\varpi_1^{h(1)}\otimes \det(g_*\sL_\sX)^{-r(1)}$, it remains to show that for a special family $(f_0:X_0\to Y_0,\sL_0)$ the rigidified direct image sheaf is weakly positive over $Y_0$. This is exactly the sheaf 
$$
S^{h(1)}(f_{0*}(\omega_{X_0/Y_0}^{\epsilon}\otimes \sL_0))\otimes \det(f_{0*}(\sL_0))^{-1}
$$
considered in Corollary~\ref{co.2} (1).
\qed
\begin{addendum}\label{co.4}
If $\omega_F^\nu=\sO_F$ for all $(F,\sA)\in \fM_h(\C)$, and if $\nu>0$ then for some $p\gg1$ there exists
an ample invertible sheaf $\lambda_{0,\nu}^{(p)}$ satisfying the property ($\divideontimes$) in the Notations~\ref{in.4} in the introduction.
\end{addendum}
\begin{proof}
The existence of the sheaf $\lambda_{0,\nu}^{(p)}$ satisfying the property ($\divideontimes$) 
follows from the construction of moduli schemes as a quotient of the Hilbert scheme.
In order to verify the ampleness, write for the universal family $G:\sX\to H$ over the Hilbert scheme
$\omega_{\sX/H}^\nu=g^*\lambda_\nu$. One has $r(1)=h(1)$ 
and the ample sheaf $\varpi_1^{h(1)}\otimes \det(g_*\sL_\sX)^{-r(1)}$ is 
\begin{equation}\label{co.5}
\det(g_*(\omega_{\sX/H}^{\epsilon}\otimes\sL_{\sX}))^{r(1)}
\otimes \det(g_*\sL_\sX)^{-r(1)}= \lambda_{\nu}^{\frac{\epsilon\cdot r(1)^2}{\nu}}. 
\end{equation}
\end{proof}
\begin{remark}\label{co.6}
Trying in the following sections to extend the polarization to degenerate fibres, we have to
keep the equality stated in \eqref{co.5}. As explained in Remark \ref{sa.1} this will force us to choose ``saturated extensions''
of the polarizations.
\end{remark}
\begin{remark}\label{co.7}
So for polarized minimal models we verified the stability of the points of $H$ for the polarization given by
$\det(h_*(\omega_{\sX/H}^\epsilon\otimes \sL))^{\hat{h}(1)}\otimes \det(h_*\sL)^{-r(1)}$. Let us assume for a moment, that $\omega_F^\eta$ is very ample for all $F\in \fM^1_{\hat{h}}$.
One can replace $H$ by the locally closed subscheme given by the condition that $\sL \sim \omega^\eta_{\sX/H}$. 
Of course this only can happen if for the Hilbert polynomial $h$ of $\omega_F$ one has
$\hat{h}(t)=h(\eta\cdot t)$. Let us assume that $\epsilon$ is divisible by $\eta$ and let us write 
$\mu=\frac{\epsilon}{\eta}+1$. Then
$$
\det(h_*(\omega_{\sX/H}^\epsilon\otimes \sL))^{\hat{h}(1)}\otimes \det(h_*\sL)^{-r(1)}=
\det(h_*\omega_{\sX/H}^{\mu\cdot\eta})^{h(\eta)}\otimes \det(h_*\omega_{\sX/H}^\eta)^{-h(\mu\cdot\eta)}.
$$
So we are still missing a factor $\mu$ on the right hand side, compared with the ample sheaf obtained by Mumford and Gieseker for moduli of curves or surfaces.
\end{remark}

\section{Weak semistable reduction}\label{ws}

Let us recall the Weak Semistable Reduction Theorem in \cite{AK} and some of the steps used in its proof.
The presentation is influenced by \cite{VZ1} and \cite{VZ2}, but all the concepts and results are due to D. Abramovich and K. Karu.

\begin{definition} \label{ws.1}
A projective morphism $\hat{g} : \hat{Z} \to \hat{Y}$ between quasi-projective
varieties is called {\em mild} if:
\begin{enumerate}
\item[(i)] $\hat{g}$ is flat, Gorenstein, and all fibres are reduced.
\item[(ii)] $\hat{Y}$ is non-singular and $\hat{Z}$ is normal with at most rational
singularities. There exists an open dense subscheme $\hat{Y}_g\subset \hat{Y}$
with $\hat{g}^{-1}(\hat{Y}_g)\to \hat{Y}_g$ smooth.
\item[(iii)] Given a dominant morphism $\hat{Y}_1 \to \hat{Y}$ from a normal quasi-projective variety $\hat{Y}_1$ with at
most rational Gorenstein singularities, $\hat{Z} \times_{\hat{Y}} \hat{Y}_1$ is normal with at most rational Gorenstein singularities.
\item[(iv)] Given a non-singular curve $\hat{C}$ and a morphism $\tau:\hat{C} \to \hat{Y}$ whose image meets $\hat{Y}_g$, the fibered product
$\hat{Z}\times_{\hat{Y}} \hat{C}$ is normal, Gorenstein and with at most rational singularities.
\end{enumerate}
\end{definition}
For a curve $\hat{Y}$ an example of a mild morphism is a semistable one, i.e. a morphism $\hat{g}:\hat{Z}\to \hat{Y}$ with $\hat{Z}$ a manifold and with all fibres reduced normal crossing divisors.  

Obviously the property iii) implies that for two mild morphisms $\hat{g}_i:\hat{Z}_i\to \hat{Y}$
the fibre product $\hat{Z}_1\times_{\hat{Y}}\hat{Z}_2 \to \hat{Y}$ is again mild. So one has:
\begin{lemma}\label{ws.2}
If $\hat{g}_i:\hat{Z}_i\to \hat{Y}$ are mild morphisms, for $i=1,\ldots,s$, then the fibre product
$\hat{Z}^r=\hat{Z}_1 \times_{\hat{Y}}\cdots \times_{\hat{Y}}\hat{Z}_s \to \hat{Y}$ is mild.
\end{lemma}
\begin{definition}\label{ws.3} Let $\hat{Y}$ be a projective manifold, $\hat{Y}_0\subset \hat{Y}$ open and dense,
and let $\hat{f}_0:\hat{X}_0\to \hat{Y}_0$ be a dominant morphism.
Then $\hat{f}_0$ has a {\em mild model} if there exists a mild morphism $\hat{g}:\hat{Z}\to \hat{Y}$,
with $\hat{Z}$ birational to some compactification of $\hat{X}$ over $\hat{Y}$.
\end{definition}
The Weak Semistable Reduction Theorem implies that after a non-singular alteration of
the base, every morphism $f_0:X_0\to Y_0$ has a mild model:
\begin{constr}\label{ws.4} \ \\[.1cm]
{\bf Start.} {\it Let $f_0:X_0\to Y_0$ be a flat surjective projective morphism
between quasi-projective varieties of pure dimension $n+m$ and $m$,respectively, 
and with a geometrically integral generic fibre.}\\[.1cm]
We will consider two cases. Either $f_0$ is smooth, or $Y_0$ is non-singular and $f_0$ a flat morphism.\\[.1cm]
{\bf Step I.} {\it Choose a flat projective model $f:X\to Y$ of $f_0$.
If $\widetilde{f}:\widetilde{X}\to \widetilde{Y}$ is any projective model of $f_0$
one may choose $Y$ and $X$ to be modifications of $\widetilde{Y}$ and $\widetilde{X}$, respectively.}\\[.1cm]
Start with any compactification $\widetilde{f}:\widetilde{X}\to \widetilde{Y}$ and with an embedding
$\widetilde{X} \to \BP^\ell$. Then $f_0$ defines a morphism
$\vartheta:Y_0\to \mathfrak{Hilb}$ to the Hilbert scheme of subvarieties of
$\BP^{\ell}$. We choose a modification $Y$ of $\widetilde{Y}$ such that
the morphism $\vartheta$ extends to $\vartheta:Y\to \mathfrak{Hilb}$. The family $f:X\to Y$
is defined as the pullback of the universal family.\\[.1cm]
{\bf Step II.} {\it There exist modifications $\sigma$ and $\sigma'$ and a diagram
\begin{equation}\label{eqws.1}
\xymatrix{X' \ar[r]^{\sigma'}\ar[d]_{f'}& X\ar[d]^{f} \\
Y' \ar[r]^\sigma&Y}
\end{equation}
with $Y'$ non-singular, such that for some open dense subschemes $U_Y\subset Y'$ and $U_X\subset X'$ the morphism
$$
f':(U_X\subset X')\to (U_Y\subset Y')
$$ 
is equidimensional, toroidal, and where $X'$ is without horizontal divisors, i.e.
where none of the irreducible components of $X'\setminus U_X$ is dominant over $Y'$.}\\[.1cm]
The construction is done in \cite{AK} in several steps. Replacing $Y$ by its normalization and $X$ by
the pullback family one may assume that $Y$ is integral. Theorem 2.1 (loc.cit.) allows to find
the diagram (\ref{eqws.1}) with $f'$ toroidal for suitable subsets
$U_X\subset X'$ and $U_Y\subset Y'$, and with $X'$ and $Y'$ non-singular.
Next Section 3 (loc.cit.) explains how to get rid of horizontal divisors in $X'$, without changing $f'$.  

In Proposition 4.4 (loc.cit.) the authors construct a non-singular projective modification of $Y'$
and a projective modification of $X'$ such that the induced rational map is in fact 
an equidimensional toroidal morphism. \\[.1cm]
{\bf Step III.} {\it For each component $D_i$ of $Y'\setminus U_Y$ there exists a positive integer
$m_i$ with the following property. 

For a ``Kawamata covering package'' $(D_i,m_i,H_{i,j})$ (defined on page 261 (loc.cit.)) consider the diagram
$$
\xymatrix{\hat{Z}\ar[r]^{\pi'}\ar[d]_{\hat{g}}&
X' \ar[r]^{\sigma'} \ar[d]_{f'}& X\ar[d]^{f}\\
\hat{Y}\ar[r]^\pi&Y' \ar[r]^\sigma &Y}
$$
where $\pi:\hat{Y}\to Y'$ is the covering given by $(D_i,m_i,H_{i,j})$, and where $\hat{Z}$ is the normalization of
$X'\times_{Y'}\hat{Y}$. Then $\hat{g}:\hat{Z}\to \hat{Y}$ is mild.}\\[.1cm]
The definition of the numbers $m_i$ is given in \cite[Page 264]{AK}, and the rest is contained in
Propositions 5.1 and 6.4 (loc.cit.).  There however the authors define a mild morphism as one satisfying
the condition~\ref{ws.1} (i)--(iii). As pointed out by K. Karu in \cite{Kar}, proof of 2.12, 
the arguments used to verify the property~\ref{ws.1} (iii) carry over ``word by word'' to show
the property (iv). So there is no harm in adding this condition.
\end{constr}
Summing up what we obtained in Construction~\ref{ws.4} one has:
\begin{proposition}\label{ws.5}
Starting with a flat projective morphism $f:X\to Y$ as in step I, one finds a commutative
diagram
\begin{equation}\label{eqws.2}
\xymatrix{X\ar[d]_{f}&\ar[l]_{\hat{\varphi}}\hat{Z}\ar[d]^{\hat{g}}\\
Y &\ar[l]_\varphi\hat{Y}}
\end{equation}
of projective morphisms with:
\begin{enumerate}
\item[(a)] $\hat{Y}$ is non-singular and $\varphi$ is an alteration. In particular if $f_0:X_0\to Y_0$
is smooth, then $X_0\times_{Y_0}\varphi^{-1}(Y_0)$ is non-singular.
\item[(a')] If $Y_0$ is non-singular and if $f_0:X_0\to Y_0$ is a mild morphism, then the variety
$X_0\times_{Y_0}\varphi^{-1}(Y_0)$ is normal with at most rational Gorenstein singularities.
\item[(b)] $\hat{g}:\hat{Z} \to \hat{Y}$ is mild.
\item[(c)] The induced morphism $\hat{Z}\to X\times_Y\hat{Y}$ is a modification.
\end{enumerate}
\end{proposition}
A more natural object to study is a desingularization $\hat{X}$ of the pullback family
${\rm pr}_2:X\times_Y\hat{Y} \to \hat{Y}$. Although the resulting morphism not necessarily flat
we will use both constructions joint by a non-singular modification $Z$:
\begin{setup}\label{ws.6} Assume that $f_0:X_0\to Y_0$ is smooth.
Starting with the diagram (\ref{eqws.2}), one can find projective morphisms
\begin{equation}\label{eqws.3}
\xymatrix{
X \ar[dr]_f & \ar[l]_{\hat{\varphi}}  \hat{Z}\ar[dr]_{\hat{g}} & \ar[l]_{\hat{\delta}} Z \ar[r]^\delta \ar[d]^g & \hat{X} \ar[dl]^{\hat{f}} \ar[r]^{\rho} & X \ar[dl]^{f} \\
& Y & \ar[l]^{\varphi} \hat{Y} \ar[r]_{\varphi} & Y,}
\end{equation}
\begin{enumerate}
\item[i.] such that $\rho:\hat{X}\to X$ factors through a desingularization 
$\rho':\hat{X}\to X\times_Y\hat{Y}$,
\item[ii.] $\hat{\delta}$ and $\delta$ are modifications, and $Z$ is non-singular.
\end{enumerate}
\end{setup}
\begin{notations}\label{ws.7}
We will denote in \ref{ws.6} by $\hat{Y}_0$, $\hat{Z}_0$ $\hat{X}_0$ (and so on) the preimages of the open subscheme $Y_0\subset Y$, and by $\hat{\varphi}_0$, $\hat{g}_0$, $\hat{f}_0$ (and so on) the restriction of the corresponding morphisms. The condition i) implies in particular that $\hat{X}$ contains $\hat{X}_0=X_0\times_{Y_0}\hat{Y}_0$ as an open dense subscheme. Later we will also consider a ``good'' dense open subscheme $Y_g\subset Y_0$ and correspondingly its preimages will be denoted by $\hat{Y}_g$, $\hat{Z}_g$, $\hat{X}_g$ (and so on).

Obviously the properties in \ref{ws.5} are compatible with replacing $\hat{Y}$ by any non-singular
alteration $\hat{Y}_1\to \hat{Y}$. We will do so several times, in order to
add additional conditions on the morphism $\hat{g}$. We will write
$\hat{Z}_1=\hat{Z}\times_{\hat{Y}}\hat{Y}_1$ and $\hat{g}_1$ for the second projection.
For $\hat{X}_1$ and $Z_1$ choose desingularizations of the main components of $\hat{X}\times_{\hat{Y}}\hat{Y}_1$ and $Z\times_{\hat{Y}}\hat{Y}_1$, respectively.
All the morphisms in the diagram corresponding to (\ref{eqws.3}) will keep their names,
decorated by a little ${}_1$. Once the additional property is verified, we usually will change back notations and drop the lower index ${}_1$.

We are also allowed to replace $Y$ by a modification with center in $Y\setminus Y_0$,
provided we modify the other schemes in the diagram (\ref{eqws.2}) accordingly.
\end{notations}
As said in the introduction, we are also interested in the polarized case, starting with a morphism $f_0:X_0\to Y_0$ and an $f_0$-ample invertible sheaf $\sL_0$. In order to 
have a reference sheaf one starts with the extension of $\sL_0$ to some projective compactification.
\begin{variant}\label{ws.8}
Assume in \ref{ws.4} that $\sL_0$ is an $f_0$-ample invertible sheaf.
Then one may choose $X$ such that the sheaf $\sL_0$ extends to an invertible sheaf $\sL$ on $X$.
Moreover, given $\widetilde{f}:\widetilde{X}\to \widetilde{Y}$ with $Y_0\subset \widetilde{Y}$ open and dense and with $\widetilde{f}^{-1}(Y_0)$ isomorphic to $X_0$ over $Y_0$, one may choose $Y$ and $X$ to be modifications of $\widetilde{Y}$ and $\widetilde{X}$, respectively.
\end{variant}
\begin{proof}[Proof of~\ref{ws.8}]
In fact, one just has to modify the first step in the construction~\ref{ws.4}.
Start with any compactification $\widetilde{f}:\widetilde{X}\to \widetilde{Y}$. Blowing up $\widetilde{X}$ one may assume that  $\sL_0$ extends to an invertible sheaf $\widetilde{\sL}$. Choose an invertible sheaf $\sA$ on $\widetilde{X}$ with $\sA$ and $\sA\otimes \widetilde{\sL}$ very ample.  
Those two sheaves define embeddings
$\iota: \widetilde{X} \to \BP^\ell$ and $\iota':\widetilde{X}\to \BP^{\ell'}$.
The restriction of the diagram
$$ 
\xymatrix{\widetilde{X} \ar[rr]^{ (\iota,\iota',\widetilde{f}) \hspace{.6cm}} \ar[dr]_{\widetilde{f}}&&
\BP^\ell \times \BP^{\ell'} \times \widetilde{Y}\ar[dl]^{{\rm pr}_3}\\
&\widetilde{Y}&}
$$
to $Y_0$ gives rise to a morphism $\vartheta:Y_0\to \mathfrak{Hilb}$ to the Hilbert scheme of subvarieties of $\BP^{\ell}\times \BP^{\ell'}$. We choose a projective compactification $Y$ of $Y_0$ such that the morphism $\vartheta$ extends to $\vartheta:Y\to \mathfrak{Hilb}$. The family $f:X\to Y$ 
is defined as the pullback of the universal family, and $\sL$ as the pullback of 
$\sO_{\BP^{\ell}\times \BP^{\ell'}}(-1,1)$.
\end{proof}
\section{Direct images and base change}\label{di}
We start by recalling some well known corollaries of ``Cohomology and Base Change'' for projective morphisms. 
\begin{lemma}\label{di.1}
Let $Y$ be quasiprojective, let $f:X\to Y$ be a projective morphism and let $\sN$ be a coherent sheaf on $X$, flat over $Y$.
\begin{enumerate}
\item[i.] There exists a maximal open dense subscheme $Y_m\subset Y$
such that the sheaf $f _*\sN|_{Y _m}$ is locally free and compatible with 
base change for morphisms $T\to Y $, factoring through $Y _m$.
\item[ii.]
If $f _*\sN$ is locally free and compatible with base change for all modifications
$\theta:Y '\to Y $, then it is compatible with base change for all morphisms 
$\varrho:T\to Y $ with $\varrho^{-1}(Y _m)$ dense in $T$. 
\item[iii.] There exists a modification $Y ' \to Y $ with center in $Y \setminus Y _m$
such that for
$$
\xymatrix{
X' =X \times_{Y }Y'  \ar[r]^{\ \ \ \ \ \ \ \ \theta'} \ar[d]_{f'} & X \ar[d]^f \\
Y'  \ar[r]_{\theta} & Y 
}
$$
the sheaf $f'_*(\theta'^*\sN)$ is locally free and compatible with base change for morphisms $\varrho:T\to Y' $ with $\varrho^{-1}\theta^{-1}(Y _m)$ dense in $T$.
\end{enumerate}
\end{lemma}
\begin{proof}
One can assume that $Y $ is affine. By ``Cohomology and Base Change'' there is a complex
\begin{equation}\label{eqdi.1}
E_0 \> \delta_0 >> E_1 \> \delta_1 >> \cdots \>\delta_{m-1} >> E_m
\end{equation}
of locally free sheaves, whose $i$-th cohomology calculates $R^if _*\sN$, as well as its base change. 
We choose $Y_m$ to be the open dense subscheme, where the image $\sC$ of $\delta_0$ locally splits in $E_1$. 
One has an exact sequence on $Y$
\begin{equation}\label{eqdi.2}
0\>>> \sK={\rm Ker}(\delta_0) \>>> E_0 \>>> \sC \>>> 0.
\end{equation}
Part ii) can be extended in the following way: 
\begin{claim}\label{di.2} The following conditions are equivalent.
\begin{enumerate}
\item[a.] $\sC$ is locally free.
\item[b.] $f_*\sN$ is locally free and compatible with base change for 
all modifications $\varrho:T\to Y$. 
\item[c.] $f_*\sN$ is locally free and compatible with base change for 
all morphisms $\varrho:T\to Y$ with $\varrho^{-1}(Y _m)$ dense in $T$. 
\end{enumerate}
\end{claim}
\begin{proof}
Of course c) implies b). If $\sC$ is locally free $\sK=f_*\sN$ is locally free,
and for all morphisms $\varrho:T\to Y$ the sequence
$$
0\>>> \varrho^* \sK \>>> \varrho^* E_0 \>>> \varrho^*\sC \>>> 0
$$
remains an exact sequence of locally free sheaves. If $\varrho^{-1}(Y_m)$ is dense in $T$
the morphism $\varrho^*\sC \to \varrho^* E_1$ is injective on some open dense subset,
hence injective. Recall that the complex
\begin{equation}\label{eqdi.3}
\varrho^* E_0 \> \delta'_0 >> \varrho^* E_1 \> \delta'_1 >> \cdots \>\delta'_{m-1} >> \varrho^* E_m
\end{equation}
calculates the higher direct images of ${\rm pr}_1^*\sN$ on the pullback family 
$X \times_{Y }T\to T$. As we have just seen, $\varrho^*\sK$ is the kernel of $\delta'_0$, hence equal to ${\rm pr}_{2*}{\rm pr}_1^*\sN$.

So a) implies c) and it remains to show that b) implies a). By assumption $\sK=f _*\sN$ is locally free,
so $\sC$ is the cokernel of a morphism between locally free sheaves of rank 
$\ell=\rk(\sK)$ and $e=\rk(E_0)$, and $r=e-\ell=\rk(\sC)$. So $\sC$ is not locally free    
if and only if the $r$-th Fitting ideal is non-trivial (see for example \cite[Proposition 20.6]{Eis}). 
Choose for $\varrho:T\to Y $ a blowing up, such that $\varrho^*\sC/_{\rm torsion}$ is locally free.
The fitting ideal is compatible with pullback (see \cite[Corollary 20.5]{Eis}), hence 
$\varrho^*\sC$ itself is not locally free. Then, using the notation from (\ref{eqdi.3}),
$\varrho^* \sK \varsubsetneq {\rm Ker}(\delta'_0) = {\rm pr}_{2*}{\rm pr}_1^*\sN$, 
violating b).
\end{proof}
The argument used at the end of the proof of~\ref{di.2} also implies that the subscheme $Y _m$ is maximal with the property asked for in ii). In fact, if the image $\sC$ does not split locally in a neighborhood of a general point
of $\varrho(T)$ the map $\varrho^*\sC \to \varrho^*E_1$ can not be injective and one finds again that
$\varrho^* \sK \varsubsetneq {\rm Ker}(\delta'_0)$.

By the choice of $Y_m$ the sequence (\ref{eqdi.2}) locally splits on $Y_m$,
and there is a blowing up $\theta:Y'\to Y $ with center in $Y \setminus Y_m$, such that
$\theta^*(\sC)/_{\rm torsion}$ is locally free. This sheaf is a subsheaf of $\theta^*(E_1)$, hence it is the
image of $\theta^*(\delta_0)$. So the latter is locally free, and by Claim~\ref{di.2} we found the modification we are looking for in iii).
\end{proof}
For relatively semiample sheaves on the total space of a mild morphism the modification of the base in~\ref{di.1} iii) is not needed. Let us recall the following base change criterion, essentially due to Koll\'ar:

\begin{lemma}\label{di.3}
Let $\hat{g}:\hat{Z}\to \hat{Y}$ be a mild morphism, and let $\hat{\sL}$ be a $\hat{g}$-semiample invertible sheaf on $\hat{Z}$. Then for all $i\geq 0$ the sheaves $R^i\hat{g}_*(\omega_{\hat{Z}/\hat{Y}}\otimes\hat{\sL})$
are locally free and compatible with arbitrary base change.
\end{lemma}
\begin{proof}
By ``Cohomology and Base Change'', i.e. using the complex $E_\bullet$ in (\ref{eqdi.1})
for $\hat{g}:\hat{Z}\to \hat{Y}$ instead of $f:X\to Y$, one finds that it is sufficient to show that the sheaves $R^i\hat{g}_*(\omega_{\hat{Z}/\hat{Y}}\otimes\hat{\sL})$ are locally free, or equivalently that the cohomology sheaves $\sH^i(E_\bullet)$ are all locally free.

As recalled in \cite[Corollary 6.12]{EV} Koll\'ar's vanishing theorem (loc.sit. Corollary 5.6) implies that
the sheaves $R^i\hat{g}_*(\omega_{\hat{Z}/\hat{Y}}\otimes\hat{\sL})$ are torsion free. In particular, if $\dim(\hat{Y})=1$
we are done.

In general consider the largest open subscheme $\hat{Y}_g$ of $\hat{Y}$ with $\hat{g}^{-1}(\hat{Y}_g)\to \hat{Y}_g$ smooth. Let $\iota:C\to \hat{Y}$ be a morphism 
from a projective curve to $\hat{Y}$ whose image meets $\hat{Y}_g$. Then $h:S=\hat{Z}\times_{\hat{Y}}C \to C$
is again mild, in particular $S$ is again normal with rational Gorenstein singularities.
Hence $R^ih_*(\omega_{S/C}\otimes{\rm pr}_1^*\hat{\sL})$ is locally free. This implies that
for all points $y\in \iota(C)$ the dimensions 
$$
h^i(y)=\dim H^i(\hat{g}^{-1}(y),\omega_{\hat{g}^{-1}(y)}\otimes\hat{\sL}|_{\hat{g}^{-1}(y)})
$$
are the same. Since $\hat{Y}$ is covered by such curves $h^i(y)$ is constant on $\hat{Y}$.
Hence $\sH^i(E_\bullet)$ is locally free.
\end{proof}
The proof of~\ref{di.3} gives a first indication why we need weak semistable models. In general even if
$X$ has at most rational Gorenstein singularities, and if $f:X\to Y$ is flat, the arguments used in
the proof of~\ref{di.3} fail. Given $\iota:C\to Y$ and $S=X\times_Y C$ we would not know that
$S$ again has rational Gorenstein singularities. 

Let us return to the notations introduced in the last section. 
So starting from a smooth morphism $f_0:X_0\to Y_0$ consider again morphisms $\varphi:\hat{Y}\to Y$ and $\hat{g}:\hat{Z}\to \hat{Y}$
satisfying the conditions (a)--(c) in~\ref{ws.5}. We choose the diagram
(\ref{eqws.3}) in~\ref{ws.6} in such a way that the conditions i) and ii) hold.

\begin{setup}\label{di.4}
Let $\sL_0$ be an invertible sheaf on $X_0$, either equal to $\sO_{X_0}$ or 
$f_0$-ample sheaf. In the first case we write $\sL=\sO_X$, in the second one we fix
an invertible extension $\sL$ of $\sL_0$ to $X$, as constructed in Variant~\ref{ws.8}.
Assume that $\sM_Z$, $\sM_{\hat{Z}}$ and $\sM_{\hat{X}}$ are invertible sheaves on
$Z$, $\hat{Z}$ and $\hat{X}$, respectively, with
\begin{gather*}
\hat{\delta}_*\sM_{Z}=\sM_{\hat{Z}}, \ \  \delta_*\sM_Z=\sM_{\hat{X}},  \ \ \hat{\varphi}^*\sL\subset \sM_{\hat{Z}}\\ \sM_{\hat{Z}_0}=\sM_{\hat{Z}}|_{\hat{Z}_0}=\hat{\varphi}^*_0\sL_0
\mbox{ \ \ and \ \ }\sM_{\hat{X}_0}=\sM_{\hat{X}}|_{\hat{X}_0}=\rho^*_0\sL_0 .
\end{gather*}
We fix some finite set $I$ of tuples $(\nu,\mu)$ of non-negative integers and define
$$
\sF_{\hat{Y}}^{(\nu,\mu)}=\hat{g}_*(\omega_{\hat{Z}/\hat{Y}}^\nu\otimes\sM_{\hat{Z}}^\mu).
$$
We choose for $\hat{Y}_g$ an open dense subscheme of $\hat{Y}_0$ such that
$\hat{g}^{-1}(\hat{Y}_g)\to \hat{Y}_g$ is smooth and such that the sheaves
$R^i\hat{g}_*(\omega_{\hat{Z}/\hat{Y}}^\nu\otimes\sM_{\hat{Z}}^\mu)$
are locally free and compatible with base change for morphisms
$\varrho:T\to \hat{Y}_g$, for all $(\nu,\mu)\in I$ and for all $i$.

If $\sL_0=\sO_{X_0}$ we choose $\sM_\bullet=\sO_\bullet$. In this case 
$I=I'\times \{0\}$ for some finite set of natural numbers $I'$. 
\end{setup}
Given an invertible sheaf $\sL$ on $X$ one could define $\sM_Z$, $\sM_{\hat{Z}}$ and $\sM_{\hat{X}}$ as the pullbacks of $\sL$. In particular this choice seems to be the most natural one if $\sL$ is $f$-ample. For families of polarized minimal models we will define in Section~\ref{sa} other extensions of $\sM_{\hat{Z}_0}=\hat{\varphi}^*_0\sL_0$ and
$\sM_{\hat{X}_0}=\rho^*_0\sL_0$, better suited for a generalization of the addendum~\ref{co.4}. 

If $\hat{Y}_1\to \hat{Y}$ is a non-singular alteration (see~\ref{ws.7} for our standard notations) 
the sheaves $\sM_{\hat{Z}_1}$, $\sM_{\hat{X}_1}$ and $\sM_{Z_1}$ are defined by pullback, and they obviously satisfy again the properties asked for in~\ref{di.4}, with $\hat{Y}_g$ replaced by its preimage in $\hat{Y}_1$. 

\begin{corollary}\label{di.5}
In~\ref{di.4} one may choose $\hat{Y}$ and $\hat{Z}$ in Proposition~\ref{ws.5} such that in addition to the conditions
(a)--(c) one has:
\begin{enumerate}
\item[(d)] For $(\nu,\mu)\in I$ the sheaves $\sF_{\hat{Y}}^{(\nu,\mu)}$
are locally free and compatible with base change for morphisms $\varrho:T\to \hat{Y}$ with
$\varrho^{-1}(\hat{Y}_g)$ dense in $T$. Hence writing
$$
\xymatrix{\hat{Z}_T \ar[r]^{\varrho'}\ar[d]_{\hat{g}_T}& \hat{Z}\ar[d]^{\hat{g}}\\
T \ar[r]^\varrho & \hat{Y}}
$$
for the fibre product, one has $\sF_{T}^{(\nu,\mu)}:=\varrho^* \sF_{\hat{Y}}^{(\nu,\mu)}=\hat{g}_{T*}(\omega^\nu_{\hat{Z}_T/T}\otimes \varrho'^*\sM_{\hat{Z}}^\mu)$.
\end{enumerate}
\end{corollary}
\begin{proof}
The properties (a)--(c) in~\ref{ws.5} are compatible with base change by non-singular alterations
$\hat{Y}_1\to \hat{Y}$. So using part iii) in~\ref{di.1}, we may assume that for a given tuple $(\nu,\mu)$
and $\sN=\omega_{\hat{Z}/\hat{Y}}^\nu\otimes \sM_{\hat{Z}}^\mu$ the condition ii) in~\ref{di.1},
holds true on $\hat{Y}$ itself. Again (d) is compatible with
base  change for alterations, and repeating the construction for the
other tuples in $I$ one obtains~\ref{di.5}.
\end{proof}
One important example are the $r$-fold fibre products. Recall that by Lemma~\ref{ws.2} the morphism
$$\hat{g}^r:\hat{Z}^r=\hat{Z}\times_{\hat{Y}}\cdots \times_{\hat{Y}}\hat{Z} \>>> \hat{Y}$$
is mild. One has $\omega_{\hat{Z}^r/\hat{Y}}={\rm pr}_1^*\omega_{\hat{Z}/\hat{Y}}\otimes \cdots \otimes {\rm pr}_r^*\omega_{\hat{Z}/\hat{Y}}$.
For 
$\sM_{\hat{Z}^r}={\rm pr}_1^*\sM_{\hat{Z}}\otimes \cdots \otimes {\rm pr}_r^*\sM_{\hat{Z}}$
flat base change and the projection formula give:
\begin{corollary}\label{di.6}
The condition d) in~\ref{di.5} implies that $\hat{g}^r_*(\omega_{\hat{Z}^r/\hat{Y}}^\nu\otimes\sM_{\hat{Z}^r}^\mu) =\bigotimes^r \sF_{\hat{Y}}^{(\nu,\mu)}$,
for $(\nu,\mu)\in I$. In particular those sheaves are again locally free and compatible with base change for morphisms $\varrho:T\to \hat{Y}$ with $\varrho^{-1}(\hat{Y}_g)$ dense in $T$.
\end{corollary}
In order to define the sheaves $\sF_{\hat{Y}}^{(\nu,\mu)}$ and to study their behavior under base change and products we used the mild model $\hat{g}:\hat{Z}\to \hat{Y}$. However since we might have blown up the smooth fibres of $X_0\to Y_0$ in order to find the mild model, this is not really the right object to study.
As a next step we will use the right hand side of the diagram (\ref{eqws.3}) in~\ref{ws.6}
to derive properties of the geometrically more meaningful morphism $\hat{f}:\hat{X}\to \hat{Y}$.
Jumping from one side of \eqref{eqws.3} to the other is possible by:
\begin{lemma}\label{di.7}
For all $\nu, \mu \geq 0$ the natural maps
\begin{gather*}
g_*(\omega_{Z/\hat{Y}}^\nu\otimes \sM_{Z}^\mu) \>>> \hat{f}_*(\omega_{\hat{X}/\hat{Y}}^\nu\otimes \sM_{\hat{X}}^\mu) 
\mbox{ \ \ and \ \ }\\
g_*(\omega_{Z/\hat{Y}}^\nu\otimes \sM_{Z}^\mu) \>>> \hat{g}_*(\omega_{\hat{Z}/\hat{Y}}^\nu\otimes \sM_{\hat{Z}}^\mu)= \sF_{\hat{Y}}^{(\nu,\mu)}
\end{gather*}
are both isomorphisms.
\end{lemma}
\begin{proof}
The morphisms $\delta$ and $\hat{\delta}$ are both birational. Since $\hat{X}$ is smooth and
$\hat{Z}$ Gorenstein with rational singularities one can find effective divisors $E_{\hat{Z}}$
and $E_{\hat{X}}$, contained in the exceptional loci of $\hat{\delta}$ and $\delta$, with
$$
\omega_{Z/\hat{Y}}=\hat{\delta}^* \omega_{\hat{Z}/\hat{Y}}\otimes\sO_Z(E_{\hat{Z}})=
\delta^*\omega_{\hat{X}/\hat{Y}}\otimes\sO_Z(E_{\hat{X}}).
$$
On the other hand, $\sM_{\hat{X}} = \delta_* \sM_{Z}$ and 
$\sM_{\hat{Z}} = \hat{\delta}_*\sM_Z$, hence for some
effective divisors $F_{\hat{Z}}$ and $F_{\hat{X}}$, contained again in the exceptional loci of 
$\hat{\delta}$ and $\delta$, one has
$$
\sM_{Z}=\hat{\delta}^*\sM_{\hat{Z}}\otimes \sO_Z(F_{\hat{Z}})=\delta^*\sM_{\hat{X}}\otimes \sO_Z(F_{\hat{X}}).
$$
The projection formula implies that 
\begin{gather*}
\hat{\delta}_*(\omega_{Z/Y}^\nu\otimes \sM_Z^\mu)=\omega_{\hat{Z}/\hat{Y}}^\nu\otimes \sM_{\hat{Z}}^\mu\otimes
\hat{\delta}_*\sO_Z(\nu\cdot E_{\hat{Z}}+\mu\cdot F_{\hat{Z}})=\omega_{\hat{Z}/\hat{Y}}^\nu\otimes \sM_{\hat{Z}}\\
\mbox{and \ \ } 
\delta_*(\omega_{Z/Y}^\nu\otimes \sM_Z^\mu)=\omega_{\hat{X}/\hat{Y}}^\nu\otimes \sM_{\hat{X}}^\mu\otimes
\delta_*\sO_Z(\nu\cdot E_{\hat{X}}+\mu\cdot F_{\hat{X}})=\omega_{\hat{X}/\hat{Y}}^\nu\otimes \sM_{\hat{X}},
\end{gather*}
hence~\ref{di.7}
\end{proof}
As we just have seen, the isomorphisms of sheaves in~\ref{di.7} are given over some
open dense subscheme by the birational maps $\hat{\delta}$ and $\delta$. We will write
in a sloppy way $=$ instead of $\cong$ for all such isomorphisms and for those induced by base change. 

Since $f:\hat{X}\to \hat{Y}$ is not necessarily flat, we can not apply ``Cohomology and Base Change''
to the right hand side of the diagram (\ref{eqws.3}), except if the (unnatural) assumptions of the next 
lemma hold true, for example for embedded semistable reductions over curves considered in Section~\ref{ec}.
\begin{lemma}\label{di.8}
Assume in~\ref{di.5} that for $(0,\mu)\in I$ the sheaves $f_{0*}\sM_{X_0}^\mu$
are locally free and compatible with arbitrary base change.
Let $U \subset \hat{Y}$ be an open subscheme, such that $V =\hat{f}^{-1}(U)\to U$ is flat.
Let $T\subset U$ be a curve, meeting $\hat{Y}_0$, and assume that for all coverings $T'\to T$ the variety $V \times_{U}T'$ is normal with at most rational Gorenstein singularities.
Then for $(\nu,\mu)\in I$ the direct image $\hat{f}_{*}(\omega_{\hat{X}/\hat{Y}}^\nu\otimes \sM_{\hat{X}}^\mu)$ is compatible with base change for all $T'\to T\subset U$.
\end{lemma}
\begin{proof} By Lemma~\ref{di.3} the sheaves $f_{0*}(\omega_{X_0/Y_0}^\nu\otimes \sM_{X_0}^\mu)$ are locally free and compatible with arbitrary base change for all $(\nu,\mu)\in I$ with $\nu>0$, and by  assumption the same holds for $\nu=0$. 

Let $\theta:\hat{Y}_1 \to \hat{Y}$ be a modification. By the choice of $I$ in Corollary
\ref {di.5} one knows that $\theta^* \sF^{(\nu,\mu)}_{\hat{Y}}=\sF^{(\nu,\mu)}_{\hat{Y}_1}$, and 
by Lemma~\ref {di.7} 
$$
\sF^{(\nu,\mu)}_{\hat{Y}}=\hat{f}_{*}(\omega_{\hat{X}/\hat{Y}}^\nu\otimes \sM_{\hat{X}}^\mu), \mbox{ \ \ and \ \ }
\sF^{(\nu,\mu)}_{\hat{Y}_1}=\hat{f}_{1*}(\omega_{\hat{X}_1/\hat{Y}_1}^\nu\otimes \sM_{\hat{X}_1}^\mu).
$$
Cutting $\hat{Y}_1$ with hyperplanes on finds through any point $p\in \theta^{-1}(T)$ a curve $T'$ mapping surjectively to $T$. Then, as we will see in Lemma~\ref{ec.3}, $\hat{X}\times_{\hat{Y}}\hat{Y}_1$ will be normal with at most rational Gorenstein singularities in a neighborhood of $\theta^{-1}(T)$.
Replacing $\hat{Y}$ by $U$, we may assume that this is the case for $\hat{X}\times_{\hat{Y}}\hat{Y}_1$ itself. So $\hat{f}_{*}(\omega_{\hat{X}/\hat{Y}}^\nu\otimes \sM_{\hat{X}}^\mu)$ is locally free and compatible with base change for modifications. On the other hand by assumption the sheaves
$\hat{f}_{0*}(\omega_{\hat{X}_0/\hat{Y}_0}^\nu\otimes \sM_{\hat{X}_0}^\mu)$
are locally free and compatible with arbitrary base change, hence the open subscheme
$\hat{Y}_m$ in Lemma~\ref{di.1} ii) contains $\hat{Y}_0$, and~\ref {di.8} follows from
Lemma~\ref {di.1} ii).
\end{proof}
Remark that Lemma~\ref{di.8} does not imply that $\hat{g}_{*}(\omega_{\hat{Z}/\hat{Y}}^\nu\otimes \sM_{\hat{Z}}^\mu)$ is compatible with base change for morphisms $\varrho:T\to \hat{Y}$ with $\varrho^{-1}(\hat{Y}_0)$ dense. If $\varrho^{-1}(\hat{Y}_g)$ is not dense, we do not know that $\hat{Z}\times_{\hat{Y}}T \to T$ is mild, hence we can not use Lemma~\ref{di.7}.  
\begin{comments}\label{di.9}
The proof of Theorem~\ref{in.1} could be finished at this stage.
Let us sketch the line of arguments, hoping that it will serve as an introduction to the remaining part of the article. 

The Extension Theorem of Gaber starts with a projective scheme $Y$, an open dense subscheme $Y_0$, a non-singular alteration $\varphi:\hat{Y} \to Y$. Write $\hat{Y}_0=\varphi^{-1}(Y_0)$ and $\varphi_0=\varphi|_{\hat{Y}_0}$. Next one considers locally free sheaves $\sF_{Y_0}$ and $\sF_{\hat{Y}}$ with $\varphi_0^*(\sF_{Y_0})=\sF_{\hat{Y}}|_{\hat{Y}_0}$. In addition one has a finite covering $\phi:W\to Y$ with a splitting trace map, whose normalization is the Stein factorization of $\varphi$.  
 
In addition one needs a sheaf $\sF_C$ for each curve $\pi:C \to W$ whose image meets $W_0=\phi^{-1}(Y_0)$.
If $\pi$ factors through $\chi:C\to \hat{Y}$ one requires that $\chi^*\sF_{\hat{Y}}=\sF_C$, and
$\sF_C$ must be compatible with replacing $C$ by a covering. The conclusion is the existence of 
the sheaf $\sF_W$, perhaps after replacing $W$ by a modification with center in $W\setminus W_0$.

Let us try to verify those conditions for $\sF_{Y_0}^{(\nu)}$. The the alteration $\hat{Y}$ and the sheaf $\sF_{\hat{Y}}^{(\nu)}$ have been constructed in Sections~\ref{ws} and~\ref{di}.
For $\sF_C^{(\nu)}$ there is little choice. It has to be the direct image $h_*\omega_{S/C}^\nu$ for a desingularization $h:S\to C$ of the pullback family. The compatibility with finite coverings 
enforces the assumption that $h: S\to C$ has a semistable or a mild model. 
So we will have to verify two conditions:
\begin{enumerate}
\item If $\pi_0:C_0\to W_0$ is a morphism from a non-singular curve, then the pullback family
$h_0:S_0\to C_0$ allows a mild model $f:S\to C$ over the smooth compactification $C$ of $C_0$.
\item If the morphism $\pi_0$ factors through the restriction of $\chi:C\to \hat{Y}$ to
$C_0$, then $\chi^*\sF_{\hat{Y}}^{(\nu)}=h_*\omega_{S/C}^\nu$.
\end{enumerate}
For (2) Lemma~\ref{di.8} will be of help. Its application is made possible by the embedded weak semistable reduction over curves, discussed in the first part of Section~\ref{ec} and stated in Proposition
\ref{ec.8}. There we first flatify the morphism in a neighborhood of a given curve. Replacing
this neighborhood by an alteration we may assume that the pullbackfamily is semistable
over $C$. Using this construction we will verify (2) in Proposition~\ref{mi.5}. 

For (1) remark, that it holds for curves $\pi_0:C_0\to W_0$
whose image meets the open dense subscheme $W_g$ where $\hat{Y} \to W$ is an isomorphism,
and whose lifting to $\hat{Y}$ meets the open set $\hat{Y}_g$, defined in~\ref{ws.1}.
This allows to verify (1) in Section~\ref{ex} on coverings of certain locally closed subschemes of
$W$. The necessary gluing is made possible by studying embeddings of $W$ into projective spaces,
which at the same time will take care of the splitting trace map. 

So both sides of the diagram \eqref{ws.3} play their role. The left hand side is needed for the definition of $\sF_{\hat{Y}}^{(\nu)}$, for its compatibility with alterations and for the verification of (1).
The right hand side is needed to control the compatibility of the sheaves to curves, as stated in (2). 

As indicated in the comments in \ref{pd1.5}, the properties of $W$ and $\sF_W^{(\nu)}$
stated in Theorem~\ref{in.1} are not sufficient to get the condition iv) in Theorem
\ref{in.2}. This forces us to allow multiplier ideals, which will 
also help to extend Corollary~\ref{in.3} to sheaves of the form $\sF_\bullet^{(\nu,\mu)}$ with 
$\mu>0$. Since we do not want to repeat the same construction in two slightly different set-ups, we will
first try to understand multiplier ideals in families in Section~\ref{fm}, or to be more precise, base change properties of ``multiplier sheaves'', i.e. of the tensor product of a multiplier ideal with an invertible sheaf. As we will see in \ref{fm.3} one does not even have a reasonable ``base change morphism'' to start with. And in general such multiplier sheaves will not be flat over the base. If they occur as in Lemma~\ref{pd1.1} as a direct factor of the direct image of an alteration $Z'\to Z$, and if $Z$ is the total space of a family $Z\to Y$ we will perform the weakly stable reduction for $Z'\to Y$, in order to find a nice model. 

The corresponding application of the embedded weak semistable reduction is given in the second half of Section~\ref{ec}. In Section~\ref{de} we define certain multiplier sheaves on fibre products, depending
on certain tautological maps, similar to the determinant $\Xi$ in the proof of
Variant~\ref{pd1.4}. Here we have to give the definitions for both sides of the diagram \eqref{ws.3}
and to verify certain compatibilities. In Section~\ref{va} we verify similar compatibilities for the restriction to curves,and we extends the proof of (2) to multiplier sheaves. 

We invite the reader to jump directly to Section~\ref{ex}, stopping shortly
at Sections~\ref{ec} and \ref{mi}, and to fill in the details needed for the ampleness property later. 
As we said already in Remark~\ref{co.5}, the Section~\ref{sa} is mainly needed for the identification of the natural ample sheaf on the compactification of the moduli scheme in Theorem~\ref{in.6}.
\end{comments}

\section{Flattening and pullbacks of multiplier ideals}\label{fm}

If multiplier ideals occur as in Lemma~\ref{pd1.1} as a direct factor of the direct image of an alteration $Z'\to Z$, and if $Z$ is the total space of a family $Z\to Y$, the weakly stable reduction for $Z'\to Y$ will allow to verify certain functorial properties. This will later be applied to the 
mild morphism $\hat{g}:\hat{Z}\to \hat{Y}$. The constructions will force us to replace the base by a new alteration, an excuse to drop the $\hat{ \ }$ and to start with:

\begin{ass-not}\label{fm.1} Let $g:Z\to Y$ be a flat projective surjective Gorenstein morphism over a non-singular variety $Y$. Assume that the $r$-fold fibre product $Z ^r=Z\times_Y\cdots\times_YZ$ is normal with at most rational singularities.

Let $\sN$ be an invertible sheaf on $Z $, let $\Delta$ be an effective Cartier divisor on $Z $ and let $N>1$ be a natural number. Assume that
there is a locally free sheaf $\sE$ together with a morphism $\sE \to g_*\sN^N$ on $Y $ with $g^*\sE \to \sN^N\otimes \sO_{Z }(-\Delta)$ surjective. 

Let $\fC$ be a set of morphisms from normal varieties $T$ with at most rational Gorenstein singularities to $Y$, such that for all $(\theta:T\to Y)\in \fC$ and for all $r>0$ the variety
$Z^r_T=Z^r\times_YT$ is again normal with at most rational Gorenstein singularities. We will need in addition that $g^{-1}(\theta(T))$ is not contained in the support of $\Delta$.

For $(\varrho:T\to Y)\in \fC$ we will write $\varrho_T:Z_T  \to Z $ and $g_T:Z_T  \to T$ for the induced morphisms. On the products the corresponding morphisms will be denoted by 
$$
\varrho_T^r:Z^r_T  \to Z^r\mbox{ \ \ and \ \ }g^r_T:Z^r_T  \to T.
$$
We consider $\Delta^{r}={\rm pr}_1^*\Delta + \cdots + {\rm pr}_r^*\Delta$ on $Z^r$
and $\Delta_T$ or $\Delta^{r}_T$ denote the pullbacks of those divisors to $Z_T$ or $Z^r_T$. 
We write
$$
\sN_{Z^{r}}={\rm pr}_1^*\sN \otimes \cdots \otimes {\rm pr}_r^*\sN\mbox{ \ \ and \ \ }
\sA_{Z^{r}}={\rm pr}_1^*\sA \otimes \cdots \otimes {\rm pr}_r^*\sA 
$$
for an invertible auxiliary sheaf $\sA$ on $\hat{Z}$, usually ample or semiample.
\end{ass-not}
If $g:Z\to Y$ is a mild morphism, smooth over a dense open subscheme $Y_g$, and if $\Delta$ does not contain
$g^{-1}(y)$ for $y\in Y_g$, then $\fC$ can be chosen as the set of morphisms $\varrho:T\to Y$ with $T$ a normal variety with at most rational Gorenstein singularities, where either $\varrho$ is dominant, or $\dim(T)=1$ and $\varrho^{-1}(Y_g)$ is dense in $T$. 
\begin{conditions}\label{fm.2} 
In \ref{fm.1} write $\varepsilon=\frac{1}{N}$ and consider for $\varrho\in \fC$ the following statements:
\begin{enumerate} 
\item[a.] $\sJ(-\varepsilon\cdot\Delta)$ is compatible with $r$-th products, i.e.
$$
\sJ(-\varepsilon\cdot\Delta^{r})\cong\big[{\rm pr}_1^*\sJ(-\varepsilon\cdot\Delta)\otimes \cdots \otimes {\rm pr}_r^*\sJ(-\varepsilon\cdot\Delta)\big]/_{\rm torsion}.
$$
\item[b.] For all $r\geq 1$ there is a natural isomorphism
$$
\varrho_T^{r*}\sJ(-\varepsilon\cdot\Delta^{r})/_{\rm torsion}\> \cong >> \sJ(-\varepsilon\cdot\Delta^{r}_T).
$$
\item[c.] For all $g $-semiample invertible sheaves $\sA$ on $Z $
the direct image 
$$
g^r_*(\omega_{Z^r /Y }\otimes \sA_{Z^r} \otimes\sN_{Z^r}\otimes \sJ(-\varepsilon\cdot\Delta^{r}))
$$ 
is locally free and the composite 
\begin{multline*}
\hspace*{1cm} {\varrho}^* g^r_*(\omega_{Z^r /Y }\otimes \sA_{Z^r} \otimes\sN_{Z^r}\otimes \sJ(-\varepsilon\cdot\Delta^{r}))\>\gamma >> \\
g^r_{T*}(\omega_{Z^r_T/T}\otimes\varrho_T^{r*}(\sA_{Z^r}\otimes \sN_{Z^r}\otimes \sJ(-\varepsilon\cdot\Delta^{r})))
\>\eta >> \\
g^r_{T*}(\omega_{Z^r_T/T}\otimes \varrho_T^{r*}(\sA_{Z^r}\otimes \sN_{Z^r}\otimes \sJ(-\varepsilon\cdot\Delta^{r}))/_{\rm torsion})
\> \cong >>\\
g^r_{T*}(\omega_{Z^r_T/T}\otimes \varrho_T^{r*}(\sA_{Z^r}\otimes \sN_{Z^r}) \otimes \sJ(-\varepsilon\cdot \Delta^{r}_T))
\end{multline*} 
of the base change morphism and the quotient map in b) is an isomorphism.
\item[d.] One has an isomorphism
$$ \ \ \ \ \ 
\bigotimes^r g_*(\omega_{Z/Y}\otimes\sA\otimes\sN\otimes \sJ(-\varepsilon\cdot\Delta))\cong
g^r_*(\omega_{Z^r/Y}\otimes\sA_{Z^r}\otimes\sN_{Z^r}\otimes \sJ(-\varepsilon\cdot\Delta^{r})).
$$
\end{enumerate}
\end{conditions}
\begin{example}\label{fm.3} In general multiplier ideals behave badly under base change.
Even if $T\subset Y$ is a smooth divisor curve $\sJ(-\varepsilon\cdot \Delta)|_{Z_T}$  
might be larger that $\sJ(-\varepsilon\cdot \Delta_T)$. So in general one can not even expect the existence
of a map
$$
\varrho_T^*\sJ(-\varepsilon\cdot \Delta)\>>> \sJ(-\varepsilon\cdot \Delta_T)
$$
in \ref{fm.2}, b).
\end{example}
\begin{lemmadef}\label{fm.4} Under the assumptions made in \ref{fm.1}
we say that $\sJ(-\varepsilon\cdot\Delta)$ is compatible with pullback, base change and products 
for $\varrho\in \fC$ if the conditions a)--d) in \ref{fm.2} hold true, or if equivalently:
\begin{enumerate} 
\item[i.] For all $r>0$ the sheaves 
$$
{\rm pr}_1^*\sJ(-\varepsilon\cdot\Delta)\otimes \cdots \otimes {\rm pr}_r^*\sJ(-\varepsilon\cdot\Delta)
\mbox{ \ \ and \ \ } \varrho_T^{r*}\sJ(-\varepsilon\cdot\Delta^r)
$$
are torsion free and isomorphic to $\sJ(-\varepsilon\cdot\Delta^r)$ and $\sJ(-\varepsilon\cdot\Delta^r_T)$,
respectively.
\item[ii.] For all $g $-semiample invertible sheaves $\sA$ on $Z $
the direct image 
$$
g^r_*(\omega_{Z^r /Y }\otimes \sA_{Z^r} \otimes\sN_{Z^r}\otimes \sJ(-\varepsilon\cdot\Delta^{r}))
$$ 
is locally free and compatible with base change for $\varrho\in \fC$.
\end{enumerate}
Moreover the conditions i) and ii) imply:
\begin{enumerate}
\item[iii.] The multiplier ideal $\sJ(-\varepsilon\cdot\Delta^r)$ is flat over $Y$.
\end{enumerate}
\end{lemmadef}
\begin{proof}
Let us remark first, that by Grothendieck's cohomological criterion for flatness \cite[Proposition 7.9.14]{EGA} the local freeness of the sheaves 
$$g^r_*(\omega_{Z^r /Y }\otimes \sA_{Z^r} \otimes\sN_{Z^r}\otimes \sJ(-\varepsilon\cdot\Delta^{r}))
$$ 
for all powers of a given ample sheaf implies that $\sJ(-\varepsilon\cdot\Delta^{r}))$ is flat over $Y$.
Hence ii) or the condition c) in \ref{fm.2} imply the condition iii).

Let us assume that i) and ii) hold true. Then a) and b) follow from i), and the condition
iii) allows to get d) in \ref{fm.2} by flat base change. 
By i) the morphism $\eta$ in c) is the identity, and $\gamma$ is the usual base change map, hence an isomorphism by ii).

For the other direction we have already seen that iii) in \ref{fm.4} holds.
So for $\sA$ sufficiently ample, the base change morphism $\gamma$ in~\ref{fm.2} c) is an isomorphism.
Since its composite with $\eta$ is assumed to be an isomorphism as well, $\eta$ must be an isomorphism.
This being true for all ample sheaves $\sA$ one finds that $\varrho_T^{r*}\sJ(-\varepsilon\cdot\Delta^{r})$
is torsion free, and b) implies 
$$
\varrho_T^{r*}\sJ(-\varepsilon\cdot\Delta^{r})\>\cong >> \varrho_T^{r*}\sJ(-\varepsilon\cdot\Delta^{r})/_{\rm torsion} \> \cong >> \sJ(-\varepsilon\cdot\Delta^r_T).
$$
Using iii) for $r=1$, one finds by flat base change and by the projection formula that 
\begin{multline}\label{eqfm.1}
g^r_*(\omega_{Z^r/Y }\otimes\sA_{Z^r}\otimes\sN_{Z^r}\otimes{\rm pr}_1^*\sJ(-\varepsilon\cdot\Delta)\otimes \cdots \otimes {\rm pr}_r^*\sJ(-\varepsilon\cdot\Delta) )\hspace*{3cm}\\ \hspace*{.5cm}\cong
\bigotimes^r g_*(\omega_{Z /Y }\otimes\sA\otimes\sN\otimes \sJ(-\varepsilon\cdot\Delta)).
\end{multline}
In particular both sheaves are locally free, and the cohomological criterion for flatness implies that 
${\rm pr}_1^*\sJ(-\varepsilon\cdot\Delta)\otimes \cdots \otimes {\rm pr}_r^*\sJ(-\varepsilon\cdot\Delta)$ 
is a flat over $Y$. The condition d) tells us that the direct images in (\ref{eqfm.1}) are isomorphic to 
$$
g^r_*(\omega_{Z^r/Y}\otimes \sA_{Z^r}\otimes \sN_{Z^r} \otimes \sJ(-\varepsilon\cdot \Delta^{r})),
$$
hence that ${\rm pr}_1^*\sJ(-\varepsilon\cdot\Delta)\otimes \cdots \otimes {\rm pr}_r^*\sJ(-\varepsilon\cdot\Delta)$ is isomorphic to
$\sJ(-\varepsilon\cdot\Delta^r)$ and torsion free.

The condition ii) now follows from i) and c).
\end{proof}

The main result of this section is a complement to the Weak Semistable Reduction Theorem.
\begin{theorem}\label{fm.5} Assume in~\ref{fm.1} that $g:Z\to Y$ is mild, that $Y_g\subset Y$ is open
with $g^{-1}(Y_g)\to Y_g$ smooth, and that $\Delta$ does not contain any fibres $g^{-1}(y)$ for $y\in Y_g$. Then there exists a fibre product diagram of morphisms
$$
\xymatrix{Z_1 \ar[r]^{\theta'}\ar[d]_{g_1}& Z\ar[d]^g\\
Y_1\ar[r]^\theta&Y,}
$$
with $\theta$ a non-singular alteration, and an open dense subscheme $Y_{1g}$ of $\theta^{-1}(Y_g)$, such that for 
\begin{multline*}
\fC_1=\big\{\varrho:T \to Y_1 \mbox{ with either } \varrho \mbox{ dominant and } T \mbox{ normal with at most rational}\\
\mbox{ Gorenstein singularities, or } T  \mbox{ a non-singular curve and } \varrho_T^{-1}(Y_{1g}) \mbox{ dense in }T \big\}
\end{multline*}
and for $\Delta_1=\theta'^* \Delta$ the sheaf $\sJ(-\varepsilon\cdot \Delta_1)$ is compatible with pullback, base change and products for all $(\varrho:T \to Y_1)\in \fC_1$.
\end{theorem}
\begin{proof} We will verify the conditions a)--d) stated in~\ref{fm.2}.\\[.2cm]
{\bf Step I.} \ As a first step, under the additional assumption
\begin{equation}\label{eqfm.2}
\sN^N\otimes \sO_{Z}(-\Delta)=\sO_{Z}
\end{equation}
we will construct a non-singular alteration $Y_1\to Y$ such that the pullback family $g_1:Z_1\to Y_1$ satisfies the condition~\ref{fm.2} b), for $r=1$.

Consider the cyclic covering $W\to Z$ obtained by taking the $N$-th root out of $\Delta$ and a log-resolution $\delta':\widetilde{Z} \to Z$ for $\Delta$. One has a diagram
\begin{equation}\label{eqfm.3}
\xymatrix{
\widetilde{W}\ar[d]_{\widetilde\pi}\ar[dr]^\pi \ar[r] & W \ar[d]\\
\widetilde{Z} \ar[r]_{\delta'}& Z
}
\end{equation}
where $\widetilde{W}$ is a desingularization of the fibre product. Since $\varepsilon=\frac{1}{N}$
Lemma~\ref{pd1.2} implies that
$$
\sN \otimes \delta'_*(\omega_{\widetilde{Z}/Y}\otimes \sO_{\widetilde{Z}}(-\big[\varepsilon\cdot \delta'^* \Delta\big]))=
\sN \otimes \omega_{Z/Y} \otimes \sJ(-\varepsilon\cdot \Delta)
$$
is a direct factor of $\pi_*\omega_{\widetilde{W}/Y}$. As we have seen there,
the assumption that $\widetilde{W}\to Z$ factors through $\widetilde{Z}$ is not needed. Similarly it is sufficient to require $\widetilde{W}$ to have rational Gorenstein singularities. 

Nevertheless let us start with $\widetilde{W}$ as in (\ref{eqfm.3}).
We choose $Y_1\to Y$ to be a non-singular alteration, such that ${\rm pr}_2:\widetilde{W}\times_{Y} Y_1 \to Y_1$ has a mild model
$h_1: W_1\to Y_1$. By construction, one has a morphism $W_1\to \widetilde{W}$ and hence $\pi_1:W_1 \to Z_1$.
Remark that the divisor $\pi_1^*\theta'^*\Delta$ is divisible by $N$.

Let us formulate what we know up to now and what will be used in the next step:
\begin{setup}\label{fm.6}
$Y_1\to Y$ is a non-singular alteration, $h_1: W_1\to Y_1$ is a flat Gorenstein morphism 
factoring through an alteration $\pi_1:W_1 \to Z_1$. The morphism $h_1$ has reduced fibres and it is smooth over an open dense subscheme $Y_{1g}$ of $Y_1$. 

Moreover $g_1:Z_1\to Y_1$ is mild and for all $g_1$-semiample sheaves $\sA$ on $Z_1$ by Lemma \ref{di.3} the sheaf
$h_{1*}(\pi_1^*\sA\otimes \omega_{W_1/Y_1})$ is locally free and compatible with arbitrary base change. 
\end{setup}
Given $\varrho:T\to Y_1$, as in Theorem~\ref{fm.5}, one has 
$$
\xymatrix{
W_T \ar[r]^{\varrho'}\ar[d]_{\pi_T}& W_1\ar[d]^{\pi_1}\\
Z_T \ar[r]^{\varrho_T}\ar[d]_{g_T}& Z_1 \ar[d]^{g_1}\\
T \ar[r]^{\varrho}&Y_1.}
$$
So $g_1$ and $h_1=g_1\circ \pi_1$, as well as $g_T$ and $h_T=g_T\circ \pi_T$ are flat. 

Let us write $\Delta_T=\varrho_T^*\Delta_1$, let $\sA$ be an invertible sheaf on $Z_1$ and $\sA_T=\varrho_T^*\sA$. One has compatible base change morphisms
\begin{gather*}
\varrho_T^*\pi_{1*}(\pi_1^* \sA\otimes \omega_{W_1/Y_1})=
\sA_T \otimes \varrho_T^*\pi_{1*}\omega_{W_1/Y_1}
\>\alpha >>  \sA_T \otimes \pi_{T*}\omega_{W_T/T},\\
\varrho^*h_{1*}(\pi_1^* \sA\otimes \omega_{W_1/Y_1})=\varrho^*g_{1*}(\sA\otimes \pi_{1*} \omega_{W_1/Y_1}) \>\gamma >>  
g_{T*}(\sA_T\otimes \varrho_T^*\pi_{1*}\omega_{W_1/Y_1}),\\
\mbox{and \ \ \ \ }
\varrho^*h_{1*}(\pi_1^* \sA\otimes \omega_{W_1/Y_1})\> \ \beta=(g_{T*}(\alpha))\circ \gamma \ >>  h_{T*}(\pi_T^* \sA_T \otimes \omega_{W_T/T}).
\end{gather*}
\begin{claim}\label{fm.7} In~\ref{fm.6} for all invertible sheaves $\sA$ on $Z_1$ the morphism $\alpha$ is surjective and induces an isomorphism
$$
[\sA_T \otimes \varrho_T^*\pi_{1*}\omega_{W_1/Y_1}]/_{\rm torsion}
\>>>  \sA_T \otimes \pi_{T*}\omega_{W_T/T}.
$$
\end{claim}
\begin{proof}
Remark that $h_T:W_T\to T$ is flat, Gorenstein, with reduced fibres and with a non-singular general fibre.
So the singular locus of $W_T$ lies in codimension at least two, and $W_T$ has to be normal, hence 
it is a disjoint union of irreducible schemes, each one flat over an irreducible components of $T$. So 
$\pi_{T*}\omega_{W_T/T}$ will be a torsion free $\sO_T$ module.

It is sufficient to prove Claim~\ref{fm.7} for one invertible sheaf $\sA$. So we may assume that $\sA$ is ample, hence $\pi_1^*\sA$ semiample. By assumption $\beta$ is an isomorphism, and 
$$
g_{T*}(\alpha):g_{T*}(\sA_T \otimes \varrho_T^*\pi_{1*}\omega_{W_1/Y_1})
\>>>  g_{T*}(\sA_T \otimes \pi_{T*}\omega_{W_T/Y_T})
$$ 
has to be surjective. For $\sA$ sufficiently ample, 
the evaluation map induces a surjection
$$
g^*_T g_{T*}(\sA_T \otimes \varrho_T^*\pi_{1*}\omega_{W_1/Y_1} ) \> 
g^*_T(g_{T*}(\alpha)) >>
g^*_T g_{T*}( \sA_T \otimes \pi_{T*}\omega_{W_T/Y_T}) \>>> \sA_T \otimes \pi_{T*}\omega_{W_T/Y_T}.
$$
Since it factors through $\alpha: \sA_T \otimes \varrho_T^*\pi_{1*}\omega_{W_1/Y_1}
\to  \sA_T \otimes \pi_{T*}\omega_{W_T/Y_T}$,
the latter must be surjective as well. By flat base change $\alpha$ is an isomorphism over some open dense subscheme of $Z_T$, hence its kernel is exactly the torsion subsheaf.
\end{proof}
Let us return to the notations used in the beginning, hence $\widetilde{W}$ is a desingularization of the cyclic covering obtained by taking the $N$-th root out of $\Delta$ and $Y_1$ is chosen, such that $\widetilde{W}\to Y$ has a mild reduction
$h_1:W_1\to Y_1$. So the conditions in~\ref{fm.6} hold true by the definition of a mild morphism and by Lemma~\ref{di.3}.

Since $W_T$ has at most rational Gorenstein singularities one obtains $\sJ(-\varepsilon\cdot \Delta_T)$ as a direct factor of 
$$
\varrho_T^*\theta'^*\sN^{-1} \otimes \pi_{T*}\omega_{W_T/Z_T}=\varrho_T^*\theta'^*\sN^{-1} \otimes \omega_{Z_T/T}^{-1}\otimes \pi_{T*}\omega_{W_T/Y_T}.
$$
By flat base change this factor coincides with $\varrho_T^*\sJ(-\varepsilon\cdot \Delta_1)$ on some open dense subscheme
of $Z_1$. Applying~\ref{fm.7} for $\sA=\theta'^*\sN^{-1} \otimes \omega_{Z_1/Y_1}^{-1}$, the morphism $\alpha$ induces an isomorphism
$$
\varrho_T^*\sJ(-\varepsilon\cdot\Delta)/_{\rm torsion}\> \cong >> 
\sJ(-\varepsilon\cdot\varrho_T^*\Delta).
$$
\noindent
{\bf Step II.} \ Next we will verify b) for $r=1$ without the additional assumption (\ref{eqfm.2}), by gluing
local models constructed in the first step.

To construct a non-singular alteration $Y_1$ such that the properties
b) in~\ref{fm.2} holds true for the family $g_1:Z_1 \to Y_1$,
one may replace $\sN$ by $\sN\otimes g^*\sH$ and correspondingly $\sE$ by $\sE\otimes \sH^N$. 

So choosing $\sH$ sufficiently ample, one may assume that $\sE$ is generated by global sections, as well
as $\sN^N\otimes \sO_{Z}(-\Delta)$. Next choose $H_1,\ldots, H_\ell$ to be zero divisors of general global sections 
of $\sN^N\otimes \sO_{Z}(-\Delta)$ and $U_i=Z\setminus H_i$, with
\begin{equation}\label{eqfm.4}
\bigcap_{i=1}^\ell H_i=\emptyset \mbox{ \ \ or \ \ }
\bigcup_{i=1}^\ell U_i=Z.
\end{equation}
By step I, for $H_i+\Delta$ instead of $\Delta$ and for each $i$, one has a non-singular alteration
$Y^{[i]}_1\to Y$ and a fibre product
$$
\xymatrix{Z^{[i]}_1 \ar[r]^{\theta^{[i]}}\ar[d]_{g^{[i]}_1}& Z\ar[d]^{g}\\
Y^{[i]}_1\ar[r] & Y}
$$
such that $\sJ(-\varepsilon\cdot {\theta^{[i]}}^* (H_i+\Delta))$ is compatible with pullback up to torsion.
Fix a non-singular alteration $\theta:Y_1\to Y$ dominating all the $Y^{[i]}_1$. For $Y_{1,g}$ choose the intersection of the preimages of the different good loci $Y^{[i]}_{1,g}$ and for $Z_1$ the pullback family.

By construction $\sJ(-\varepsilon\cdot (\Delta +H_i))|_{U_i}=\sJ(-\varepsilon\cdot \Delta)|_{U_i}$ and
$$
\sJ(-\varepsilon\cdot (\Delta_1 +\theta'^* H_i))|_{\theta^{-1}(U_i)}=\sJ(-\varepsilon\cdot \Delta_1)|_{\theta^{-1}(U_i)}.
$$
Since $\sJ(-\varepsilon\cdot (\Delta_1+\theta'^* H_i))$ is compatible with pullback up to torsion, 
the sheaf of ideals $\sJ(-\varepsilon\cdot \Delta_1)$ has the same property over $U_i$. Since $\{U_i; \ i=1,\ldots,\ell\}$
is an open covering of $Z$ the condition~\ref{fm.2} b) follows for $\Delta_1$ and for $r=1$. \\[.2cm]
{\bf Step III.}
For the model $Z_1\to Y_1$ constructed in step II we will verify the property b) for $r>1$ and the
compatibility with products, stated in~\ref{fm.2} a). Let us formulate again the set-up we will refere to at this point.
\begin{setup}\label{fm.8}
$\pi_1^{[i]}:W^{[i]}_1 \to Z_1$ are alterations such that the induced morphisms $h_1^{[i]}:W^{[i]}_1 \to Y_1$ 
satisfy the assumptions stated in~\ref{fm.6}. 
\end{setup}
Choose a tuple $\underline{i}$ consisting of $r$ elements $i_1,\ldots,i_r \in \{1,\ldots,\ell\}$ and the induced morphisms
$h^r_1:W^r=W^{[i_1]}_1\times_{Y_1} \cdots \times_{Y_1}W^{[i_r]}_1 \to Y_1$
and $\pi^r_1:W^r \to Z^r$. Let 
\begin{equation}\label{eqfm.5}
\sA_{Z^r}\otimes \bigotimes_{\iota=1}^r {\rm pr}_{i_\iota}^*\pi^{[i_\iota]}_*\omega_{W^{[i_\iota]}/Y_1}
\>\alpha^r >>  \sA_{Z^r} \otimes \pi^r_* \omega_{W^r/Y_1}=\sA_{Z^r} \otimes \pi^r_* 
\bigotimes^r {\rm pr}_{i_\iota}^*\omega_{W^{[i_\iota]}/Y_1}
\end{equation}
be induced by the tensor products of the base change maps
$$
{\rm pr}_{i_\iota}^*\pi^{[i_\iota]}_*\omega_{W^{[i_\iota]}/Y_1}\>>> \pi^r_* 
{\rm pr}_{i_\iota}^*\omega_{W^{[i_\iota]}/Y_1}.
$$
By assumption, for $\sA$ ample the sheaves $h^{[i]}_{1*} {\pi_1^{[i]}}^* \sA \otimes \omega_{W^{[i]}/Y_1}$
are locally free. By flat base change and the projection formula, one has an isomorphism
$$
\bigotimes_{\iota=1}^r h^{[i_\iota]}_{1*} ({\pi^{[i_\iota]}}^*\sA\otimes\omega_{W^{[i_\iota]}/Y_1})
\>\beta^r >> h^r_{1*} ({\pi^r_1}^* \sA_{Z^r} \otimes \omega_{W^r/Y_1}).
$$
\begin{claim}\label{fm.9}
There is a natural morphism
\begin{multline*}
\bigotimes_{\iota=1}^r h^{[i_\iota]}_{1*} ({\pi^{[i_\iota]}}^*\sA\otimes\omega_{W^{[i_\iota]}/Y_1})
=\bigotimes_{\iota=1}^r g_{1*} (\sA\otimes {\pi^{[i_\iota]}}_* \omega_{W^{[i_\iota]}/Y_1}) \\
\>\gamma^r >> g^r_{1*}\Big(\sA_{Z^r}\otimes \bigotimes_{\iota=1}^r {\rm pr}_{i_\iota}^*\pi^{[i_\iota]}_*\omega_{W^{[i_\iota]}/Y_1}\Big).
\end{multline*}
\end{claim}
\begin{proof}
Let $p_1:Z_1^r\to Z_1$ and $p_2:Z_1^r \to Z_1^{r-1}$ denote the projection to the first and the last $r-1$
factors of the fibre product. Assume one has constructed $\gamma^{r-1}$, hence the morphism
\begin{multline*}
\bigotimes_{\iota=1}^{r} g_{1*} (\sA\otimes {\pi^{[i_\iota]}}_* \omega_{W^{[i_\iota]}/Y_1})
\>\gamma^{r-1}\otimes{\rm id} >>\\
g_{1*}(\sA\otimes {\pi^{[i_r]}}_* \omega_{W^{[i_r]}/Y_1})\otimes
g^{r-1}_{1*}\Big(\sA_{Z^{(r-1)}}\otimes \bigotimes_{\iota=1}^{r-1} {\rm pr}_{i_\iota}^*\pi^{[i_\iota]}_*\omega_{W^{[i_\iota]}/Y_1}\Big).
\end{multline*}
The right hand side maps to
$$
g^{r-1}_{1*}{g^{r-1}_1}^*g_{1*}(\sA\otimes {\pi^{[i_r]}}_* \omega_{W^{[i_r]}/Y_1})\otimes
g_{1*}g^*_1 g^{r-1}_{1*}\Big(\sA_{Z^{(r-1)}}\otimes \bigotimes_{\iota=1}^{r-1} {\rm pr}_{i_\iota}^*\pi^{[i_\iota]}_*\omega_{W^{[i_\iota]}/Y_1}\Big)$$
The  tensor product of the two base change maps and the multiplication of sections
map this sheaf to 
\begin{multline*}
g^{r}_{1*}p_2^*\Big(\sA_{Z^{(r-1)}}\otimes \bigotimes_{\iota=1}^{r-1} {\rm pr}_{i_\iota}^*\pi^{[i_\iota]}_*\omega_{W^{[i_\iota]}/Y_1}\Big)\otimes
g^{r}_{1*}p_1^*(\sA\otimes {\pi^{[i_r]}}_* \omega_{W^{[i_r]}/Y_1})\\
\> {\rm mult} >> g^{r}_{1*} \Big(\sA_{Z^r} \otimes \bigotimes_{\iota=1}^{r} {\rm pr}_{i_\iota}^*\pi^{[i_\iota]}_*\omega_{W^{[i_\iota]}/Y_1}\Big).
\end{multline*}
\end{proof}
Again the isomorphism $\beta^r$ is equal to $g^r_{1*}(\alpha^r)\circ \gamma^r$, hence
$g^r_{1*}(\alpha^r)$ has to be surjective. As in the proof of~\ref{fm.7}, for $\sA$
sufficiently ample, one finds that $\alpha^r$ is surjective. Let us state what we obtained.
\begin{claim}\label{fm.10}
In~\ref{fm.8} the base change map $\beta^r$ in (\ref{eqfm.5}) is an isomorphism for all $g_1$-semiample sheaves $\sA$. The morphism $\alpha^r$ is surjective and its kernel is a torsion sheaf.
\end{claim}
Let us return to the situation considered in step II. So we have chosen 
alterations $\pi_1^{[i]}:W^{[i]}_1 \to Z_1$, dominating the cyclic covering obtained by taking the $N$-th root out of $\Delta_1+\theta'^* H_i$, such that the induced morphisms $h_1^{[i]}:W^{[i]}_1 \to Y_1$ are mild. 
By Lemma~\ref{ws.2} the morphism $h^r_1$ is again mild and $W^r$ has rational Gorenstein singularities.
$W^r$ dominates the cyclic covering obtained by taking the $N$-th root out of $\Delta_1+\theta'^* H_i$.
So $\pi^r_1:W^r \to Z^r$
is again an alteration, dominating the cyclic covering obtained by taking the
$N$-th root out of
$$
\Gamma={\rm pr}_{i_1}^*(\Delta_1+\theta'^* H_{i_1})+ \cdots +
{\rm pr}_{i_r}^*(\Delta_1+\theta'^* H_{i_r}).
$$
By step I, up to the tensor product with an invertible sheaf,
$\sJ(-\varepsilon\cdot \Gamma)$ is a direct factor of 
$$
\pi^{r}_{1*}\omega_{W^r/Y_1}={\rm pr}_{i_1}^*\pi^{[i_1]}_*\omega_{W^{[i_1]}/Y_1} \otimes \cdots \otimes
{\rm pr}_{i_r}^*\pi^{[i_r]}_*\omega_{W^{[i_r]}/Y_1} .
$$
On some open dense subscheme this factor is isomorphic to
$$
{\rm pr}_{i_1}^*\sJ(-\varepsilon\cdot (\Delta_1+\theta'^* H_{i_1}))\otimes \cdots \otimes
{\rm pr}_{i_r}^*\sJ(-\varepsilon\cdot (\Delta_1+\theta'^* H_{i_r})).
$$
So the first part of Claim~\ref{fm.10} implies that $\alpha^r$ induces
an isomorphism
$$
[{\rm pr}_{i_1}^*\sJ(-\varepsilon\cdot (\Delta_1+\theta'^* H_{i_1}))\otimes \cdots \otimes
{\rm pr}_{i_r}^*\sJ(-\varepsilon\cdot (\Delta_1+\theta'^* H_{i_r}))]/_{\rm torsion}
\> \cong >> \sJ(-\varepsilon\cdot \Gamma).
$$
For $U_i=Z\setminus H_i$ and $U_{\underline{i}}=U_{i_1}\times \cdots \times U_{i_r}$
one has 
$$
\sJ(-\varepsilon\cdot (\Delta_1+\theta'^* H_{i_\iota}))|_{U_{i_\iota}}=
\sJ(-\varepsilon\cdot \Delta_1)|_{U_{i_\iota}}
$$
and $\sJ(-\varepsilon\cdot \Gamma)|_{U_{\underline{i}}}=\sJ(-\varepsilon\cdot \Delta^{r})|_{U_{\underline{i}}}$.
Since by (\ref{eqfm.4}) each point of $Z^r$ lies in $U_{\underline{i}}$
for some choice of the tuple $\underline{i}$, one obtains the property~\ref{fm.2} a).\\[.2cm] 
The same construction gives the proof of property b) for $r>1$. One just has to remark that
$\pi^r_1:W^r_1\to Z^r_1$ and $\hat{h}^r_1:W^r_1\to Y$ satisfy again the assumption made in
\ref{fm.6}. So one just has to replace the ample sheaf $\sA$ by a sufficiently high power of $\sA_{Z^r}$,
and one obtains an isomorphism 
$$
\varrho_T^{r*}\sJ(-\varepsilon\cdot\Gamma)/_{\rm torsion}\> \cong >> 
\sJ(-\varepsilon\cdot\varrho_T^{r*}\Gamma),
$$ 
for the divisor $\Gamma$ introduced above. So~\ref{fm.2} b) holds for $\Delta^{r}$ on $U_{\underline{i}}$,
hence everywhere.\\[.2cm]
{\bf Step IV.} \
It remains to verify the properties~\ref{fm.2} c) and d). To simplify notations, let us drop the lower index ${}_1$ and assume that the properties a) and b) in~\ref{fm.2} hold true for $g:Z\to Y$ itself.

Let us first remark that we know c) and d) if $\sN^N\otimes\sO_{Z}(-\Delta)$ is the pullback of an invertible sheaf on $Y$. In fact, the base change morphisms in~\ref{fm.2} c) and d) are just direct factors of the base change morphisms $\beta$ in step I or $\beta^r$ in Claim~\ref{fm.10} in step III. So we will reduce everything to this case. 

As we have seen this can be done by adding the zero divisor $H$ of a general section of $\sN^N\otimes\sO_{Z}(-\Delta)$
to $\Delta$. There is a problem with the term ``general''. We can choose $H$ to be general for a  fibre of $g_1:Z_1\to Y_1$, hence (\ref{eqpd1.1}) holds for $F$ and for $F$ replaced by a small neighborhood. However we can not choose $H$ such that this remains true for neighborhoods of fibres of all $g_T:Z_T\to T$ and for the pullback of $H$. So we will argue in a different way.

Let us assume that the construction in step II was possible over $Y$. In particular $\sE$ and hence $\sN^N\otimes\sO_{Z}(-\Delta)$ are generated by global sections, and for some section of $\sN^N\otimes\sO_{Z}(-\Delta)$ with zero divisor $H_1$ the cyclic covering obtained by taking the $N$-th root out of $\Delta+H_1$ has a mild model $h^{[1]}:W^{[1]}\to Y$ factoring through $\pi^{[1]}:W^{[1]}\to Z$. 

As before $\delta':\widetilde{Z} \to Z$ denotes a log-resolution for $\Delta$. Fix a point
$y\in Y$. For the zero set $H$ of a general section of $\sN^N\otimes\sO_{Z}(-\Delta)$
the divisor $\delta'^*H$ will be smooth meeting $\delta'^*\Delta$ transversely. So $\sJ(-a\cdot \Delta)=\sJ(-a\cdot(\Delta+H))$ for $0 \leq a < 1$. Moreover, ${\pi^{[1]}}^*H$ will not contain any component of ${h^{[1]}}^{-1}(y)$. 

On $W^{[1]}$ the divisor ${\pi^{[1]}}^*\Delta$ is divisible by $N$. Hence the sheaf 
$$
{\pi^{[1]}}^*(\sN^N\otimes\sO_{Z}(-\Delta))=({\pi^{[1]}}^*\sN^N)\otimes\sO_{Z}(-{\pi^{[1]}}^*\Delta))
=\sO_{W^{[1]}}({\pi^{[1]}}^*H)
$$
is the $N$-th power of an invertible sheaf $\sL$. We choose
$\phi: W\to W^{[1]}$ to be the cyclic covering obtained by taking the $N$-th root out of ${\pi^{[1]}}^*H$ and $\pi=\pi^{[1]}\circ \phi$.
\begin{claim}\label{fm.11} For $H$ sufficiently general, replacing $Y$ by a neighborhood of $y$, one has:\vspace{-.2cm}
\begin{enumerate}
\item[i.] \ \hspace*{\fill} $\displaystyle\phi_*\omega_{W/Y}=\bigoplus_{\iota=0}^{N-1}\omega_{W^{[1]}/Y}\otimes \sL^{\iota}$.\hspace*{\fill} \ 
\item[ii.] The induced morphism $h:W\to Y$ is flat and Gorenstein. 
\item[iii.] The fibres of $h$ are reduced and the general fibre is non-singular.
\item[iv.] If $\sA$ is $g$-semiample the direct image sheaves $h_*(\pi^*\sA\otimes \omega_{W/Y})$ are locally free and compatible with arbitrary base change.
\item[v.] The sheaf $\sJ(-\varepsilon\cdot \Delta)\otimes \sN \otimes \omega_{Z/Y}$ is a direct factor of
$\pi_*\omega_{W/Y}$.   
\end{enumerate}
\end{claim}
\begin{proof}
The first part follows from \cite[Section 3]{EV}. However there we considered cyclic coverings over a non-singular base
and we have to explain, how to reduce the statement to this case.

Let $\tau:V\to W^{[1]}$ be a desingularization. For $H$ sufficiently general, $\tau^*H$ is non-singular.
The normalization $V'$ of $V$ in the function field of $W$ is non-singular and isomorphic to
$$
{\rm\bf Spec}(\sF) \mbox{ \ \ for \ \ } \sF=\bigoplus_{\iota=0}^{N-1}\tau^*\sL^{-\iota}.
$$
The canonical sheaf $\omega_{V'}$ is the invertible sheaf corresponding to 
$$
\sF\otimes \sL^{N-1}=\bigoplus_{\iota=0}^{N-1}\omega_{V}\otimes \tau^* \sL^{\iota}.
$$
Since $W^{[1]}$ is Gorenstein with rational singularities,
$\displaystyle \phi_*\sO_{W}=\tau_*\sF=\bigoplus_{\iota=0}^{N-1}\sL^{-\iota}$.

So $\phi_*\omega_{W}$ contains $\tau_*\sF \otimes \sL^{N-1}$ and both are isomorphic outside of a codimension two
subset. The second sheaf is a locally free $\tau_*\sF$ module of rank $1$, hence equal to $\phi_*\omega_{W}$.  
In particular $\phi:W\to W^{[1]}$ is flat, and $\omega_{W}$ is invertible. 

For iii) remark that $g^{[1]}$ is smooth over some open dense subset $Y_g$ of $Y$. The restriction of a general divisor $H$ to one fibre will be non-singular, and thereby $g$ has at least one non-singular fibre. Choosing $Y$ small enough,
we may assume that $H$ does not contain components of any fibre of $g^{[1]}$. Since the fibres of $g^{[1]}$ are reduced, the fibres of $h$ have the same property.

Part iv) follows from~\ref{di.3}, applied to the sheaves $g^{[1]}_*(\omega_{W^{[1]}/Y}\otimes \sL^{\iota}\otimes \sA)$,
and by the direct sum decomposition in i). So it remains to verify v).

Let $\pi':W' \to Z$ be the cyclic covering obtained by taking the $N$-th root out of $\Delta+H$. Then
$W$ is just the normalization of the fibre product $W'\times_{Z}W^{[1]}$. In fact, the latter is
the cyclic covering of $W^{[1]}$, obtained by taking the $N$-th root out of ${\pi^{[1]}}^*\Delta +{\pi^{[1]}}^*H$. However, ${\pi^{[1]}}^*\Delta$
is divisible by $N$, hence it is the same to take the $N$-th root out of ${\pi^{[1]}}^*H$.

So $\pi'_*\omega_{W'/Y}$ is a direct factor of $\pi_*\omega_{W/Y}$, and
$$
\sJ(-\varepsilon\cdot (\Delta+H))\otimes \sN \otimes \omega_{Z/Y}=\sJ(-\varepsilon\cdot \Delta)\otimes \sN \otimes \omega_{Z/Y}
$$
is a direct factor of both of them.
\end{proof}
Parts ii), iii) and iv) of~\ref{fm.11} imply that the assumptions stated in~\ref{fm.6} hold.
Hence by Claim~\ref{fm.7}  for all $\varrho:T\to Y$ considered in~\ref{fm.5}
the morphism 
$$
\alpha: \varrho_T^*\pi_{*}\omega_{W/Y} \to \pi_{T*}\omega_{W_T/T}
$$
is a surjection with torsion kernel. Moreover the composite
$$
\beta:\varrho^*g_{*}(\sA\otimes \pi_{*} \omega_{W/Y}) \>\gamma >>  
g_{T*}(\sA_T\otimes \varrho_T^*\pi_{*}\omega_{W/Y})\> g_* \alpha >> g_{T*}(\sA_T \otimes \pi_{T*}\omega_{W_T/T})
$$ 
is an isomorphism for all $g$-semiample sheaves $\sA$ on $Z$. By~\ref{fm.11} v) the
sheaf 
$$
\varrho^*g_{*}(\sA\otimes \sJ(-\varepsilon\cdot \Delta)\otimes \sN \otimes \omega_{Z/Y})
$$ 
is a direct factor of the left hand side, and by the property b), which we verified in steps I. and II. the corresponding
direct factor of the right hand side is 
$$
g_{T*}(\omega_{Z_T/T}\otimes \varrho_T^*(\sA\otimes \sN) \otimes \sJ(-\varepsilon\cdot \Delta_T)).
$$ 
So we obtained the property c) for $r=1$. 

For $r>1$ the argument is the same. Using the notations from step II for $\underline{i}=(1,\ldots,1)$ we just have to replace $Z$ by $Z^r$ and the divisor $H_1$ by ${\rm pr}_1^*H_1 + \cdots + {\rm pr}_r^*H_1$.

For d) we choose for the morphisms $h^{[i]}:W^{[i]} \to Z$ in step III the same morphism $h:W\to Y$. By~\ref{fm.11} ii), iii), and iv) the assumptions made in~\ref{fm.8} hold true, and by Claim~\ref{fm.10} the composite 
\begin{multline*}
\bigotimes_{\iota=1}^r g_*(\sA\otimes\pi_*\omega_{W/Y})=
\bigotimes_{\iota=1}^r h_*(\pi^*\sA\otimes\omega_{W/Y})
\>\beta^r >> h_* (\pi^* \sA_{Z^r} \otimes \omega_{W^r/Y})=\\
g_* (\sA_{Z^r} \otimes \pi^r_* \omega_{W^r/Y})\< \cong < \alpha^r <
g_* \big(\sA_{Z^r} \otimes \big[ \bigotimes_{\iota=1}^r {\rm pr}_\iota^* \pi_* \omega_{W/Y}\big]/_{\rm torsion}\big)
\end{multline*}
is an isomorphism. The left hand side contains 
$\displaystyle
\bigotimes_{\iota=1}^r g_*(\omega_{Z/Y}\otimes\sA\otimes\sN\otimes\sJ(-\varepsilon\cdot \Delta))$
as a direct factor, and the corresponding direct factor of the right hand side is
$$
g^r_*\Big(\omega_{Z^r/Y}\otimes\sA_{Z^r}\otimes\sN_{Z^r}\otimes\big[ \bigotimes_{\iota=1}^r {\rm pr}_\iota^* \sJ(-\varepsilon\cdot \Delta) \big]/_{\rm torsion}\Big).
$$
By part a) of~\ref{fm.2} this is 
$\displaystyle
g^r_*(\omega_{Z^r/Y}\otimes\sA_{Z^r}\otimes\sN_{Z^r}\otimes\sJ(-\varepsilon\cdot \Delta^{r}))$
and we obtain d).
\end{proof}
\begin{remark}\label{fm.12}
Even if one adds to~\ref{fm.1} the additional condition $N>e(\Delta)$ one can not expect in Theorem~\ref{fm.5} that $\sJ(-\varepsilon\cdot \Delta_1)$ is isomorphic to $\sO_{Z_1}$. At least we were not able to compare 
$e(\Delta_1)$ and $e(\Delta)$. So for the equality $\sJ(-\varepsilon_1\cdot \Delta_1)=\sO_{Z_1}$ with $\varepsilon_1=\frac{1}{N_1}$ we have to choose $N_1$ larger than $N$ and we loose the compatibility with base change and products.
\end{remark}
\section{Embedded weak semistable reduction over curves}\label{ec}

For a morphism to a curve with smooth general fibre, a semistable model is mild. The existence of such a model over some covering of the base has been shown by Kempf, Knudsen, Mumford, and Saint-Donat in \cite{KKMS}. Applying it to a family over a discrete valuation ring one obtains the semistable reduction theorem in codimension one: 
\begin{theorem}\label{ec.1}
Let $U$ and $V$ be a quasi-projective manifolds and let $E\subset U$ be a submanifold of codimension one. Let $f:V\to U$ be a surjective projective morphism with connected general fibre.
Then there exists a finite covering $\theta:U'\to U$, a desingularization $V'$ of the main component of $V\times_UU'$, and an open neighborhood $\widetilde{U}$ of the general points of $\theta^{-1}(E)$
such that for the induced morphism $f':V'\to U'$ the restriction $f'^{-1}(\widetilde{U})\to \widetilde{U}$
is flat and $f'^{-1}(\widetilde{U}\cap\theta^{-1}(E))$ a reduced relative normal crossing divisor over 
$\widetilde{U}\cap\theta^{-1}(E)$.
\end{theorem}
As indicated in~\ref{di.9} we will need some ``embedded version'' of the semistable reduction in a neighborhood of a given curve.
This will allow to apply the base change criterion in Lemma~\ref{di.8}.

Since as in Section~\ref{fm} we have to allow multiplier ideals the notations get again a bit complicated in the second half of the Section. The multiplier ideal of the restriction of a divisor
is contained in the restriction of the multiplier ideal. As we will see, for total spaces of families of varieties one can 
enforce an isomorphism after an alteration of the base.

\begin{lemma}\label{ec.2} 
Let $f:X\to Y$ be a projective surjective morphism between quasi-projective manifolds with smooth part
$f_0=f|_{X_0}:X_0\to Y_0$. Let $\pi:C\to Y$ be a morphism from a non-singular curve $C$ with $C_0=\pi^{-1}(Y_0)$ dense in $C$. Then
there exists a non-singular alteration $\theta:Y_1\to Y$ and a desingularization $\theta':X_1\to X\times_{Y}Y_1$ of the main component such that for the induced morphism $f_1:X_1\to Y_1$ the following holds:
\begin{enumerate}
\item[a.] $C\to Y$ lifts to an embedding $C\subset Y_1$.
\item[b.] There exists a neighborhood $U_1$ of $C$ in $Y_1$ with $f^{-1}_1(U) \to U$ flat.
\item[c.] $S=f^{-1}_1(C)$ is non-singular and $f^{-1}_1(C\setminus C_0)$ an normal crossing divisor in $S$.
\end{enumerate}
\end{lemma}
\begin{proof} Replacing $Y$ by a hyperplane in $C\times Y$, containing the 
graph of $\pi:C\to Y$, one may assume that $C\to Y$ is an embedding.
Next replace $X$ by an embedded log-resolution of the closure $S$ of $f^{-1}(C)\cap X_0$ for the divisor $f^{-1}(C\setminus C_0)$. 
So we may assume that the closure $S$ of $f^{-1}(C_0)$ is non-singular and that the singular fibres of $S\to C$ are normal crossing divisors. Consider for a very ample invertible sheaf $\sA$ on $X$, the induced embedding
$\iota:X\to \BP^M$, and the diagram
$$
\xymatrix{X \ar[rr]^{(\iota,f)\hspace*{.3cm}}\ar[dr]_{f}&&\BP^M\times Y\ar[dl]^{{\rm pr}_2}\\
&Y.&}
$$
Since $X_0\to Y_0$ is flat, it gives rise to a morphism $\vartheta_0:Y_0\to \mathfrak{Hilb}$ to the Hilbert
scheme of subvarieties of $\BP^M$. Since $S\to C$ is also flat the restriction of $\vartheta_0$ to $C\cap Y_0$ extends to a morphism
$\varrho:C\to \mathfrak{Hilb}$, and the pullback of the universal family over $\mathfrak{Hilb}$
to $C$ coincides with $S$. 

We choose a modification $\theta:Y_1\to Y$ with center outside of $Y_0$ such that $\vartheta_0$ extends to a morphism
$\vartheta:Y_1\to \mathfrak{Hilb}$. For $f_1:X_1\to Y_1$ we choose the pullback of the universal family.
Remark that $f_1$ satisfies the conditions a), b) and c), however $X_1$ might be singular.
Since we are allowed to modify $X_1$ outside of a neighborhood of $S$ it remains to verify
that $X_1$ is non-singular in such a neighborhood. This will be done in the next Lemma.
\end{proof}
\begin{lemma}\label{ec.3}
Let $f:V\to U$ be a flat morphism, with $U$ non-singular. Let $C\subset U$ be a non-singular curve
and $S=f^{-1}(C)$. Then one can find an open neighborhood $U_0$ of $C$ in $U$ with:
\begin{enumerate}
\item[i.] If $S$ is non-singular, $f^{-1}(U_0)$ is non-singular.
\item[ii.] If $S$ is reduced, normal, Gorenstein with at most rational singularities then $f^{-1}(U_0)$ is normal, Gorenstein with at most rational singularities.
\item[iii.] If $S$ is reduced, and Gorenstein, and if for some open subscheme $U_g$ of $U$, meeting $C$
the preimage $f^{-1}(U_g)$ is non-singular, then $V$ is normal and Gorenstein.
\end{enumerate}
\end{lemma}
\begin{proof} $C$ is a smooth curve in $U$. For a point $p\in C$ choose local parameter $t_1, \cdots , t_\ell$
such that $C$ is the zero-set of $(t_1,\cdots,t_{\ell-1})$. The parameters $(t_1,\cdots,t_{\ell-1})$
define a smooth morphism ${\rm Spec}\sO_{p,U}\to {\rm Spec}\sO_{0,\A^{\ell-1}}$. The composite flat morphism
$$
\Phi:V \times_{U}{\rm Spec}\sO_{p,U} \>>> {\rm Spec}\sO_{p,U}\>>> {\rm Spec}\sO_{0,\A^{\ell-1}}
$$
has $S_0=S\times_{C}{\rm Spec}\sO_{p,C}$ as closed fibre. If the latter is smooth, $\Phi$ is smooth and one obtains i). 

Assume that $S$ is Gorenstein. Then $S_0$ is Gorenstein, and $\Phi$ is a Gorenstein morphism.
 
If in addition $S$ is reduced and normal, it is smooth outside of a codimension 
one subset, hence $V \times_{U}{\rm Spec}\sO_{p,U}$ will be normal. And if $S$ has at most rational singularities,
the same holds true for $V \times_{U}{\rm Spec}\sO_{p,U}$. 

In iii) the assumptions imply that the singular locus $\Upsilon$ of $V\times_{U}{\rm Spec}\sO_{p,U}$ does not meet the general fibre 
of $\Phi$. On the other hand, since the special fibre $S_0$ is reduced, $\Upsilon$ contains no component of $S_0$.
So again $\Upsilon$ is of codimension two and since $V\times_{U}{\rm Spec}\sO_{p,U}$ is Gorenstein it is normal.
\end{proof}
\begin{variant}\label{ec.4} Under the assumptions made in~\ref{ec.2} 
one can find a finite covering $\hat{C}\to C$, a non-singular alteration $\theta:Y_1\to Y$ and a desingularization $\theta':X_1\to X\times_{Y}Y_1$ such that for the induced morphism $f_1:X_1\to Y_1$ in addition to the properties a), b) and c) (for $\hat{C}$ instead of $C$) in~\ref{ec.2} one has:
\begin{enumerate}
\item[d.] $f^{-1}_1(\hat{C}\setminus \hat{C}_0)$ is a reduced normal crossing divisor in $\hat{S}=f^{-1}_1(\hat{C})$.
\end{enumerate}
\end{variant}
\begin{proof}
We use the notations from the proof of~\ref{ec.2}, except that we assume that the conditions a)--c) hold true for $Y$ itself. So $C\subset Y$ and the morphism $f$ is flat in a neighborhood of $S=f^{-1}(C)$. The latter is non-singular and the fibres of $S\to C$ are normal crossing divisors. 

Choose $\hat{C}\to C$ to be a covering, such that $S\times_{C}\hat{C}\to \hat{C}$ has a semistable model $\hat{S}\to S$. In particular there is a morphism
$\hat{S}\to S$ inducing $\tau:\hat{S}\to S\times_{C}\hat{C}$. As in the proof of~\ref{ec.2} we can choose $Y_1$ such that
$\hat{C}\to C \to Y$ lifts to an embedding $\hat{C}\to Y_1$. Consider the fibre product $X\times_{Y}Y_1$. It contains 
$S\times_{C}\hat{C}$. Since $\tau$ is birational and projective, it is given by the blowing up of a sheaf
of ideals $\sI$ on $S\times_{C}\hat{C}$. Let $\sJ$ be a sheaf of ideals on $X\times_{Y}Y_1$, whose restriction
to $\hat{S}\to S\times_{C}\hat{C}$ is $\sI$, and let $\delta:X_1 \to X\times_{Y}Y_1$ be the blowing up of $\sJ$. 
Then one obtains a closed immersion $\hat{S}\to X_1$, whose image is contained in $f^{-1}_1(\hat{C})$.

Repeating the argument in the proof of~\ref{ec.2} we replace $X_1$ by some modification and 
$X_1\to Y_1$ by the pullback of a universal family over a Hilbert scheme, with $f^{-1}_1(\hat{C})=\hat{S}$.   
\end{proof}
\begin{definition}\label{ec.5} Let $U$ be a quasi-projective manifold, let $C$ be a smooth curve
and $\pi:C \to U$ a morphism.
We call $\theta: U_1\to U$ a {\em local alteration for $C$} if $\theta$ is the restriction of a 
non-singular alteration to some open subscheme, and if there is a smooth curve
$C_1\subset \theta^{-1}(C)$ with $C_1\to C$ finite. We call such a curve $C_1$ a lifting of $C$.
\end{definition}
\begin{lemma}\label{ec.6}
Let us assume that $C\subset Y$ is a smooth curve, that $S=f^{-1}(C)$ is a non-singular variety, semistable over $C$, that $f$ is flat over a neighborhood $U$ of $C$, and that $V=f^{-1}(U)$ is nonsingular. 
Let $\theta:U_1\to U$ be a local alteration for $C$, let $C_1\in U_1$ be a lifting of $C$
and $f_1:V_1=X\times_{U}U_1\to U_1$ the pullback family. Write
$f_1^r:V_1^r=V_1\times_{U_1}\cdots\times_{U_1}V_1 \to U_1$ for the $r$-fold fibre product. 
Then 
\begin{enumerate}
\item[($\diamond$)] For each $r>0$ there exists a neighborhood $\widetilde{U}$ of $C_1$ in $U_1$ such that
$\widetilde{V}^r=(f^r_1)^{-1}(\widetilde{U})$ is normal, Gorenstein with at most rational singularities and the induced morphism $\widetilde{f}^r:\widetilde{V}^r \to \widetilde{U}$ is flat and projective.\\
Moreover $S^r_1=(\widetilde{f}^r)^{-1}(C_1)$ is normal with at most rational Gorenstein singularities,
and $S^r_1\to C_1$ has reduced fibres. 
\end{enumerate}
\end{lemma}
\begin{proof}
As the pullback of a semistable family $S_1=F_1^{-1}(C_1)=S\times_{C}C_1$ is normal, Gorenstein with quotient singularities. The same holds true for the $r$-fold product $S^r_1=S_1\times_{C_1}\cdots \times_{C_1} S_1$. So one can apply Lemma~\ref{ec.3}.
\end{proof}
\begin{definition}\label{ec.7}
Let $f:X\to Y$ be a projective surjective morphism between quasi-projective manifolds with smooth part
$f_0=f|_{X_0}:X_0\to Y_0$, let $\pi:C\to Y$ be a morphism from a non-singular curve $C$ with $C_0=\pi^{-1}(Y_0)$ dense in $C$,
let $\theta:U_1\to Y$ be a morphism, and let $V_1\to X\times_{Y}U_1$ be a modification of the main component
with center outside of the preimage of $Y_0$. We call the induced family $f_1:V_1\to U_1$ an {\em embedded weak semistable reduction (of $X\to Y$) over $C$} if $\theta:U_1\to Y$ is a local alteration for $C$ and if for some lifting $C_1 \in U_1$ 
the condition ($\diamond$)~\ref{ec.6} hold true.

We call $f_1:V_1\to U_1$ an {\em embedded semistable reduction over $C$} if in addition
$S_1=f_1^{-1}(C_1)$ is non-singular and semistable over $C_1$.
\end{definition}
Usually we will replace $U_1$ by some neighborhood $\widetilde{U}$ and assume that
the condition in ($\diamond$) holds for $\widetilde{U}$.
Let us restate what we obtained:
\begin{proposition}\label{ec.8} 
Let $f:X\to Y$ be a projective surjective morphism between quasi-projective manifolds with smooth part
$f_0=f|_{X_0}:X_0\to Y_0$ and let $\pi:C\to Y$ be a morphism from a non-singular curve $C$ with $C_0=\pi^{-1}(Y_0)$ dense in $C$.
\begin{enumerate}
\item[a.] There exists an embedded semistable reduction $V_1 \to U_1$ over $C$.
\item[b.] Let $Y_1\to Y$ be a non-singular alteration. Then there exists a scheme $U_2$ and a morphism 
$U_2\to Y_1$ such that the image of the composed morphism $U_2\to Y$ is in $U_1$ and such that
$V_2=V_1\times_{U_1}U_2\to U_2$ is a weak semistable reduction over $C$.
\end{enumerate}
\end{proposition}
Proposition~\ref{ec.8} will allow to apply the base change criterion in Lemma~\ref{di.8}. As in Section~\ref{fm} we will need a similar criterion for multiplier sheaves. We start with a variant of Theorem~\ref{fm.5} 
replacing the mild morphism by an embedded weak semistable reduction over a curve.
\begin{assumptions}\label{ec.9}
$f:V\to U$ is an embedded weak semistable reduction for $C\subset U$, with smooth
part $f_0:V_0\to U_0$ for $U_0$ dense in $U$. There exists
a mild morphism $g:Z \to U$ factoring through a modification $\tau:Z\to V$.
Let $\sN$ be an invertible sheaf on $V$, and let $\Delta$ be an effective Cartier divisor on $V$ 
not containing fibres of $f_0$ and let $N>1$ be a natural number. 
There is a morphism $\sE\to f_*\sN^N$ on $U$ with $\sE$ locally free
and with $f^*\sE \to \sN^N\otimes \sO_{V}(-\Delta)$ surjective. As in the last Section we write $\varepsilon=\frac{1}{N}$.

Assume that $\sJ(-\varepsilon\cdot \tau^*\Delta)$ is compatible
with pullback, base change and products, for all alterations of $U$, as defined in
\ref{fm.4}, and (for simplicity) that on the general fibre of $S\to C$ the multiplier
sheaf $\sJ(-\varepsilon\cdot \Delta|_S)$ is isomorphic to $\sO_S$. 
\end{assumptions}
\begin{lemma}\label{ec.10} In~\ref{ec.9} let $\fC$ be the set of local alterations
$\theta: U_1\to U$ such that $f_1:V_1=V\times_UU_1\to U_1$ is an embedded 
weak semistable reduction for $f:V\to U$ over $C$. Then $\sJ(-\varepsilon\cdot \Delta)$ is 
flat over $U$ and compatible with pullback, base change and products for $(\varrho:U_1\to U)\in \fC$ in a neighborhood of each lifting $C_1$ of $C$, i.e. the conditions i) and ii) in Definition~\ref{fm.4} hold true over a neighborhood $\widetilde{U}\subset U_1$ of $C_1$, possibly depending on $r$. 
\end{lemma}
\begin{proof}
Choose a log-resolution $\widetilde{\delta}:\widetilde{Z}\to Z$. For $\delta=\tau\circ\widetilde{\delta}:\widetilde{Z}\to V$
one has 
\begin{multline*}
\sJ(-\varepsilon\cdot\Delta)=
\delta_*(\omega_{\widetilde{Z}/V}\otimes \sO_{\widetilde{Z}}(-[\varepsilon\cdot\delta^*\Delta]))=\\
\tau_*\widetilde{\delta}_*(\omega_{\widetilde{Z}/V}\otimes \sO_{\widetilde{Z}}(-[\varepsilon\cdot\widetilde{\delta}^*\tau^*\Delta]))=
\tau_*(\omega_{\hat{Z}/V}\otimes \sJ(-\varepsilon\cdot\tau^*\Delta)).
\end{multline*}
Then
$$
g_*(\tau^*\sA\otimes\omega_{Z/V}\otimes \tau^* \sN\otimes \sJ(-\varepsilon\cdot\tau^*\Delta))=
f_*(\sA\otimes\omega_{V/U}\otimes \sN \otimes \sJ(-\varepsilon\cdot\Delta)),
$$
and by~\ref{fm.4} ii) both are locally free, and the left hand side is compatible with pullbacks.
The cohomological criterion \cite[Proposition 7.9.14]{EGA} implies that $\sJ(-\varepsilon\cdot\Delta)$ 
is flat over $U$.

For the compatibility with base change for $\varrho:U_1\to U$ consider the induced
fibre products
$$
\xymatrix{Z_1 \ar[r]^{\hat{\varrho}}\ar[d]_{\tau_1}& Z\ar[d]^{\tau}\\
V_1 \ar[r]^{\varrho'}\ar[d]_{f_1}& V\ar[d]^{f}\\
U_1 \ar[r]^{\varrho}&U.}
$$
One has for $\sA$ ample on $Z$ the base change map
\begin{multline*}
\varrho'^*(\omega_{V/U}\otimes \sN \otimes \sA \otimes \sJ(-\varepsilon\cdot \Delta))=
\varrho'^*\tau_*(\omega_{Z/U}\otimes \tau^*(\sN \otimes \sA) \otimes \sJ(-\varepsilon\cdot\tau^*\Delta))\\
\>\alpha >> \tau_{1*}(\omega_{Z_1/U_1}\otimes \tau_1^*\varrho'^*(\sN \otimes \sA) \otimes \sJ(-\varepsilon\cdot\tau_1^*\varrho'^*\Delta)).
\end{multline*}
The base change map for $g_*(\tau^*\sA\otimes\omega_{Z/V}\otimes \tau^* \sN\otimes \sJ(-\varepsilon\cdot\tau^*\Delta))$ factors through $f_{1*}(\alpha)$, so
the latter must be surjective. This being true for all ample sheaves
$\sA$, as in the proof of~\ref{fm.7} one finds that $\alpha$ is surjective. By flat base change,
$\alpha$ is an isomorphism on some open dense subscheme. 

By assumption on the general fibre of $S\to C$ the multiplier sheaf $\sJ(-\varepsilon\cdot \Delta|_S)$ is trivial. By \cite[Section 5.4]{Vie} or \cite[Properties 7.5]{EV} this implies that 
$\sJ(-\varepsilon\cdot\Delta)$ is isomorphic to $\sO_V$ in a neighborhood of a general fibre of $f$. Since the latter is flat over $U$, the sheaf $\varrho'^*\sJ(-\varepsilon\cdot\Delta)$ is torsion free,
hence isomorphic to 
$$
\tau_{1*}\omega_{Z_1/V_1}\otimes\sJ(-\varepsilon\cdot\tau_1^*\varrho'^*\Delta)=
\sJ(-\varepsilon\cdot\varrho'^*\Delta).
$$ 
In addition $f_{1*}(\alpha)$ is an isomorphism, hence 
$f_*(\sA\otimes\omega_{V/U}\otimes \sN \otimes \sJ(-\varepsilon\cdot\Delta))$
is compatible with base change for $\varrho\in\fC$. 

A similar argument allows to identify the multiplier ideals on the $r$-fold fibre products, for $r>1$.
Let us write $\tau^r:Z^r\to V^r$ for the modification, ${\rm pr}_\iota:V^r\to V$ and ${\rm p}_\iota:Z^r\to Z$
for the projections. By flat base change one has a natural isomorphism
$$
{\rm pr}_\iota^* \sJ(-\varepsilon\cdot \Delta)\>>> \tau^r_*(\omega_{Z/V}\otimes \sJ(-\varepsilon\cdot
{\rm p}_\iota^*\tau^* \Delta)).
$$ 
Since the multiplier ideal on $Z$ is compatible with products, as formulated in~\ref{fm.4} i)
multiplication of sections induces a morphism $\alpha^r$ from $\displaystyle\bigotimes_{\iota=1}^r {\rm pr}_\iota^* \sJ(-\varepsilon\cdot \Delta)$ to
$$
\tau^r_*(\omega_{Z/V}\otimes \sJ(-\varepsilon\cdot ({\rm p}_1^*\tau^* \Delta + \cdots +
{\rm p}_r^*\tau^* \Delta))) = \sJ(-\varepsilon\cdot ({\rm pr}_1^*\Delta + \cdots +
{\rm pr}_r^* \Delta)).
$$
By flat base change 
$$
f^r_{*}\Big(\bigotimes_{\iota=1}^r{\rm pr}_\iota^*(\omega_{V/U}\otimes \sA\otimes \sN \otimes \sJ(-\varepsilon\cdot\Delta))\Big)=
\bigotimes_{\iota=1}^r f_{*}(\omega_{V/U}\otimes \sA\otimes \sN \otimes \sJ(-\varepsilon\cdot\Delta))
$$
is locally free, hence on $V^r$ the sheaf $\displaystyle\bigotimes_{\iota=1}^r{\rm pr}_\iota^*\sJ(-\varepsilon\cdot\Delta)$ is flat over $U$ and torsion free.
So 
$$
\bigotimes_{\iota=1}^r{\rm pr}_\iota^*\sJ(-\varepsilon\cdot\Delta)\> \alpha^r >>
\sJ(-\varepsilon\cdot ({\rm pr}_1^*\Delta + \cdots + {\rm pr}_r^*\Delta))
$$
is injective. Finally, writing again $\sA_{V^r}$ for the exterior tensor product
and $\sA_{Z^r}$ for the pullback to $Z^r$, the composite
\begin{multline*}
f^r_{*}\Big(\bigotimes_{\iota=1}^r{\rm pr}_\iota^*(\omega_{V/U}\otimes \sA\otimes \sN \otimes \sJ(-\varepsilon\cdot\Delta))\Big)\> f^r_*(\alpha^r) >> \\
f^r_{*}( \omega_{V^r/U}\otimes \sA_{V^r}\otimes \sN_{V^r} \otimes\sJ(-\varepsilon\cdot ({\rm pr}_1^*\Delta + \cdots + {\rm pr}_r^*\Delta)))=\hspace*{1cm}\\ \hspace*{1cm}
f^r_{*}\tau^r_*( \omega_{Z^r/U}\otimes \sA_{Z^r}\otimes \sN_{Z^r} \otimes\sJ(-\varepsilon\cdot ({\rm p}_1^*\tau^*\Delta + \cdots + {\rm p}_r^*\tau^*\Delta)))=\\
\bigotimes f_{*}\tau_*(\omega_{Z/U}\otimes \tau^*\sA\otimes \tau^*\sN \otimes \sJ(-\varepsilon\cdot\tau^*\Delta))
\end{multline*}
is an isomorphism. For $\sA$ sufficiently ample, as in the proof of~\ref{fm.7}, this implies that
$\alpha^r$ is an isomorphism.

Since $Z^r\to U$ is again mild, one may replace in the first part of the proof
$Z$ and $V$ by $Z^r$ and $V^r$, respectively, and obtains the compatibility with pullbacks, required in~\ref{fm.4} ii), for all $r$.
\end{proof}
As promised we can now formulate and prove the compatibility of multiplier ideal sheaves 
in total spaces of families with the restriction to subfamilies over curves.
This will lead for suitable models to the compatibility of certain direct image sheaves with restriction to curves.
\begin{proposition}\label{ec.11} 
Under the assumptions made in~\ref{ec.9} there exists a local alteration $\theta: U_1 \to U$ for
$C$ such that:
\begin{enumerate}
\item $f_1:V_1 =V\times_UU_1 \to U_1$ is an embedded weak semistable reduction of $f$ over $C$.  
\item For a lifting $C_1\subset U_1$ of $C$, for $S_1=f_1^{-1}(C_1)$ denote the induced morphisms by
$$
\xymatrix{S_1 \ar[r]^{\chi'}\ar[d]_{\zeta}&V\ar[d]^{f}\\
C_1\ar[r]^{\chi}& U.}
$$ 
Then there is an isomorphism $\sJ(-\varepsilon\cdot \chi'^* \Delta) \cong \chi'^* \sJ(-\varepsilon\cdot \Delta)$.
\item Let $\sA$ be an $f$-semiample sheaf on $V$. Then
$$
\chi^*f_*(\sA\otimes \sN \otimes \omega_{V/U}\otimes \sJ(-\varepsilon\cdot\Delta))=
\zeta_*(\chi'^*(\sA\otimes\sN)\otimes\omega_{S_1/C_1}\otimes \sJ(-\varepsilon\cdot\chi^*\Delta)).
$$
\end{enumerate}
\end{proposition}
\begin{proof}
Let us first show, that (1) and (2) imply (3). 

By Lemma~\ref{ec.10} the sheaf
$\sJ(-\varepsilon\cdot \Delta)$ is flat over $U$ and compatible with pullbacks and base change
for $\theta: U_1 \to U$. So by abuse of notations is is sufficient in (3) to consider the case
$U_1=U$, and to assume that $C\subset U$. On a general fibre of $S\to C$ the multiplier ideal sheaf is isomorphic to the structure sheaf, hence by \cite[Properties 7.5]{EV} the same holds over a neighborhood of the general point of 
$C$ in $U$. As in the proof of~\ref{di.3} Koll\'ar's Vanishing Theorem implies that
over this neighborhood the direct image of 
$\sA\otimes \sN \otimes \omega_{V/U}\otimes \sJ(-\varepsilon\cdot\Delta)$
is locally free and compatible with arbitrary base change. Hence applying~\ref{di.1} to this sheaf
the open dense subscheme $U_m$ in part i) contains a general point of $C$. Then
(3) follows from~\ref{di.1} ii).

To construct $U_1$ with the properties (1) and (2), we may assume that $\sE$, hence $\sN^N\otimes\sO_{V}(-\Delta)$ is globally generated. Since the question is local on $V$, as in the second step in the proof of~\ref{fm.5} we can cover $V$ by the complements of divisors of general sections of 
$\sN^N\otimes\sO_{V}(-\Delta)$. Hence we may replace $\Delta$ by $\Delta+H$ and assume that $\sN^N = \sO_{V}(\Delta)$.

Choose a desingularization of the cyclic covering, obtained by taking the $N$-th root out of $\Delta$.
Over some alteration, this desingularization will have a mild model. Since this property is compatible with pullbacks, we
may choose a local alteration for $C$, dominating the alteration, and we find some $U_1$ such that (1) holds and such that $V_1\to U_1$ has a mild model. The compatibility for local alterations, shown in Lemma~\ref{ec.10} 
allows to assume that $U_1=U$, hence that the mild model exists over $U$ itself. Let us call it $T \to U$, and the induced morphism $\psi:T\to V$. So $\psi^* \Delta$ is the $N$-th power of a Cartier divisor. 

Next we want to construct a desingularization $W$ of $T$, which is flat over a general point of the curve $C$. To this aim, let $\widetilde{U}\to U$ be the blowing up of $C$, or a finite covering 
of such a blowing up. Let $\widetilde{V}\to \widetilde{U}$ be the pullback family. The preimage of the exceptional divisor $E$ in $\widetilde{U}$ is covered by curves $\widetilde{C}$, finite over $C$. Lemma~\ref{ec.3} allows to shrink
$\widetilde{U}$ such that the total space $\widetilde{V}$ is still normal with at most rational Gorenstein singularities.

Let $\widetilde{\phi}:\widetilde{W}\to \widetilde{T}=T\times_U\widetilde{U}$ be a desingularization. It dominates 
the finite covering obtained by taking the $N$-th root out of $\widetilde{\Delta}={\rm pr}_1^*\Delta$. If $\widetilde{h}:\widetilde{W}\to\widetilde{U}$ denotes the induced map, we also assume that $\widetilde{h}^{-1}(E)$ is a normal crossing divisor. Over the complement $\widetilde{U}_g$ of a codimension two subset of $\widetilde{U}$ the morphism $\widetilde{h}$ will be flat and $\widetilde{h}^{-1}(E)\cap \widetilde{h}^{-1}(\widetilde{U}_g)$ will be equisingular over $E\cap\widetilde{U}_g$. 

The divisor $\widetilde{h}^{-1}(E\cap\widetilde{U}_g)$ might be non-reduced. 
If so we perform the semistable reduction in codimension one, described in Theorem~\ref{ec.1}. Replacing $\widetilde{U}$ by some alteration and choosing $\widetilde{U}_g$ sufficiently small, this allows to assume that $\widetilde{h}^{-1}(E\cap\widetilde{U}_g)$ is a reduced relative normal crossing divisor.

For a curve $\hat{C}\subset E$ meeting $\widetilde{U}_g$ choose a neighborhood $U'$ in $\widetilde{U}$. By construction
$\widetilde{h}^{-1}(\hat{C}\cap\widetilde{U}_g)$ has non-singular components, meeting transversely. For $W'$ choose an embedded desingularization of the components of $\widetilde{h}^{-1}(\hat{C})$, and assume that the closure $\Sigma$ of $\widetilde{h}^{-1}(\hat{C}\cap\widetilde{U}_g)$ is the union of manifolds, meeting transversely. Remark that the induced morphism $h':W'\to U'$ is still flat over some open subscheme
$U'_g$, meeting $\hat{C}$, and that there are morphisms 
$$
\psi':T'=T\times_UU'\to V'=V\times_{U}U'\mbox{ \ \ and \ \ }\phi:W'\to T'.
$$
For $\hat{C}$ sufficiently general, $\phi'$ is birational and $\psi'$ an alteration.

As in the proof of~\ref{ec.2} one obtains a morphism $\vartheta_0:U'_g \to \mathfrak{Hilb}$ to the Hilbert scheme of subvarieties of some $\BP^M$, parametrizing the fibres of $h'$. 

Since $\Sigma \to \hat{C}$ is flat the restriction of $\vartheta_0$ to $\hat{C} \cap U'_g$ extends to a 
morphism $\hat{C}\to \mathfrak{Hilb}$, and the pullback of the universal family over $\mathfrak{Hilb}$ to $\hat{C}$ coincides with $\Sigma$. 

Blowing up $U'$ with centers in $U'\setminus U'_g$ we obtain a new family, again denoted by $h':W'\to U'$, which is flat and such that $h'^{-1}(\hat{C})=\Sigma$. By~\ref{ec.3} ii), choosing the neighborhood $U'$ of $\hat{C}$ small enough, $W'$ will be normal and Gorenstein.

Let us drop again all the ${}'$ and assume that the morphisms we just constructed exists over $V$ itself. So we will assume that we have alterations
$$
W\>\phi>> T\>\psi >> V, \ \ \pi=\psi\circ\phi \mbox{ \ \ and \ \ } \gamma:\Sigma=\pi^*(S) \to S
$$
such that:
\begin{enumerate}
\item[i.] $T\to U$ is mild and $\psi^*\Delta$ is divisible by $N$.
\item[ii.] $W$ is normal and Gorenstein, flat over $U$ and $\phi$ is birational.
\item[iii.] $\Sigma$ is reduced, and the union of manifolds, meeting transversely.
\end{enumerate}
The multiplier ideal $\sJ(-\varepsilon\cdot\Delta)$ is a direct factor of $\psi_*\omega_{T/V}\otimes\sN^{-1}$. 
Let $\delta:\widetilde{W} \to W$ be a desingularization, Then one has 
$$
\delta_*\omega_{\widetilde{W}} \>\subset >> \omega_{W} \mbox{ \ \ and \ \ } \phi_*\delta_*\omega_{\widetilde{W}} \>\subset >> \phi_*\omega_{W}\>\subset>> \omega_{T}.
$$
Since $T$ has rational singularities, $\phi_*\delta_*\omega_{\widetilde{W}}=\omega_{T}$ and
$\sN\otimes \sJ(-\varepsilon\cdot\Delta)$ is a direct factor of $\pi_*\omega_{W/V}$. 

The base change map induces a morphism
$$
\eta: \sN\otimes\sJ(-\varepsilon\cdot\Delta)|_S \>>> \pi_*\omega_{W/V}|_S \>>> \gamma_* \omega_{\Sigma/S}.
$$
Recall that the sheaf $\sJ(-\varepsilon\cdot\Delta)|_S$ is flat over $C$. By \cite[Properties 7.5]{EV} it contains $\sJ(-\varepsilon\cdot\Delta|_S)$, and by assumption both coincide on the general fibre of $S\to C$. Hence $\sJ(-\varepsilon\cdot\Delta)|_S$ is torsion free and  
$\eta$ is injective.

Choose $\hat{\Sigma}$ as the union of all components of $\Sigma$ which
dominate the irreducible variety $S$, and $R$ the union of the other irreducible components
$R_1,\ldots,R_\ell$. By construction, the components of $\Sigma$  
are non-singular, and meet transversely. So one has an exact sequences
\begin{gather*}
\hspace*{.6cm} 0\>>> \omega_{\hat{\Sigma}} \>>> \omega_{\Sigma} \>>> \omega_{R}\otimes\sO_{R}(R\cap \hat{\Sigma})\>>> 0\mbox{ \ \ \ and}\\
0 \>>> \gamma_* \omega_{\hat{\Sigma}} \>>> \gamma_* \omega_{\Sigma} \>>> \gamma_*(\omega_{R}\otimes\sO_{R}(R\cap \hat{\Sigma})) 
\end{gather*}
The non-singular alteration $\hat{\Sigma}\to S$ dominates the covering obtained by taking the $N$-th root out of $\Delta|_S$. By Lemma~\ref{pd1.2} the multiplier ideal
$\sJ(-\varepsilon\cdot\Delta|_S)$ is a direct factor of $\sN^{-1}|_S\otimes \gamma_* \omega_{\hat{\Sigma}}$. On the other hand,
the sheaf $\gamma_*(\omega_{R}\otimes\sO_{R}(R\cap \hat{\Sigma}))$ is contained in 
$$
\bigoplus_{\iota=1}^\ell \gamma_*(\omega_{R_\iota}\otimes\sO_{R_\iota}(\Gamma_\iota))
$$
where $\Gamma_\iota$ is the intersection of $R_\iota$ with the other components. Each of the sheaves
$\gamma_*(\omega_{R_\iota}\otimes\sO_{R_\iota}(\Gamma_\iota))$ is torsion free over its support $\pi(R_\iota)$. By construction
$\pi(R_\iota)$ is dominant over $C$. By assumption the composite
$$
\eta: \omega_{V}\otimes\sN\otimes\sJ(-\varepsilon\cdot\Delta)|_S \>\eta >> \gamma_* \omega_{\Sigma}\>>> \bigoplus_{\iota=1}^\ell \gamma_*(\omega_{R_\iota}\otimes\sO_{R_\iota}(\Gamma_\iota))
$$
is zero along the general fibre of $S\to C$, hence it is zero. So 
$\sJ(-\varepsilon\cdot\Delta)|_S$ maps to $\sJ(-\varepsilon\cdot\Delta|_S)$, and both must be equal.
\end{proof}

\section{Saturated extensions of polarizations}\label{sa}
As in Section~\ref{ws} $f:X\to Y$ will denote a projective morphism with smooth part $f_0:X_0\to Y_0$ for $Y_0\subset Y$ dense and with $\omega_{X_0/Y_0}$ relative semiample over $Y_0$. 
If $\sL_0$ is $f_0$-ample, we choose $X$ and $\sL$ as in Variant~\ref{ws.8}. Recall the diagram \eqref{eqws.3}
\begin{equation}\label{eqsa.1}
\xymatrix{
X \ar[dr]_f & \ar[l]_{\hat{\varphi}}  \hat{Z}\ar[dr]_{\hat{g}} & \ar[l]_{\hat{\delta}} Z \ar[r]^\delta \ar[d]^g & \hat{X} \ar[dl]^{\hat{f}} \ar[r]^{\rho} & X \ar[dl]^{f} \\
& Y & \ar[l]^{\varphi} \hat{Y} \ar[r]_{\varphi} & Y,}
\end{equation}
with $\hat{g}$ mild, $\hat{X}$ non-singular, and such that the left and right hand diagrams are birational to fibre products.
$Z$ is a modification of both, $\hat{Z}$ and $\hat{X}$. Recall moreover, that starting from an invertible sheaf $\sL$ on $X$ with $\sL_0=\sL|_{X_0}$ we considered in~\ref{di.4} invertible sheaves $\sM_\bullet$ on $\bullet$, where $\bullet$ stands for $Z$, $\hat{Z}$ and $\hat{X}$. They should satisfy the compatibilities $\hat{\delta}_*\sM_{Z}=\sM_{\hat{Z}}$,   $\delta_*\sM_Z=\sM_{\hat{X}}$, $\hat{\varphi}^*\sL\subset \sM_{\hat{Z}}$,  $\sM_{\hat{Z}_0}=\sM_{\hat{Z}}|_{\hat{Z}_0}=\hat{\varphi}^*_0\sL_0$
and $\sM_{\hat{X}_0}=\sM_{\hat{X}}|_{\hat{X}_0}=\rho^*_0\sL_0$.

In this section we will pose more conditions on those sheaves. First of all, if the evaluation map
for $\omega_{X_0/Y_0}^\nu\otimes \sL_0^\mu$ is surjective, we may replace $\hat{X}$ by a modification with center in $\hat{X}\setminus \hat{X}_0$ and assume that the image of the evaluation map 
$$
\hat{f}^*\hat{f}_*\big(\omega_{\hat{X}/\hat{Y}}^\nu\otimes \sM_{\hat{X}}^\mu\big) \>>> 
\omega_{\hat{X}/\hat{Y}}^\nu\otimes \sM_{\hat{X}}^\mu
$$
is invertible. As a first step we have to show that the same is possible with $\hat{X}$ replaced by $\hat{Z}$ without loosing the mildness of $\hat{g}$. 
This will simplify some of the constructions in the next Sections, but mainly it will be needed
to define {\em saturated extensions}: 
\begin{remark}\label{sa.1} 
Assume for a moment that $\dim{Y}=1$, hence that
we can choose $\hat{Z}=Z=\hat{X} \to \hat{Y}$ to be the semistable model, and assume moreover that the smooth fibres of $f_0$ are
of Kodaira dimension zero. Then $\omega_{\hat{X}/\hat{Y}}^\nu=\sO_{\hat{X}}(\Pi)\otimes \hat{f}^*\hat{f}_*\omega_{\hat{X}/\hat{Y}}^\nu$.
for some $\nu$ and for some effective divisor contained in the singular fibres of $f$. For $\sM=\sM_{\hat{X}}$ the sheaf corresponding to the left hand side in \eqref{co.5} is the $r(1)$-th power of 
$$
\det(\hat{f}_*(\omega_{\hat{X}/\hat{Y}}^{\epsilon}\otimes\sM))
\otimes \det(\hat{f}_*\sM)^{-1}= \det(\hat{f}_*(\sO_{\hat{X}}(\frac{\epsilon}{\nu}\Pi)\otimes\sM))\otimes \det(\hat{f}_*\sM)^{-1}\otimes \lambda_{\nu}^{\frac{\epsilon}{\nu}}. 
$$
Roughly speaking $\sM$ will be a saturated extension of the polarization $\sM|_{\hat{X}_0}$
if 
$$\hat{f}_*(\sO_{\hat{X}}(\frac{\epsilon}{\nu}\cdot \Pi)\otimes\sM)=\hat{f}_*(\sO_{\hat{X}}(*\cdot \Pi)\otimes\sM)=\hat{f}_*\sM.$$ 
\end{remark}
\begin{lem-not}\label{sa.2} 
Consider in Corollary~\ref{di.5} for a given tuple $(\nu,\mu) \in I$ 
a locally free sheaf $\sE_{\hat{Y}}$ and a morphism $\sE_{\hat{Y}} \to \hat{f}_*(\omega_{\hat{X}/\hat{Y}}^\nu\otimes\sM_{\hat{X}}^\mu)$
such that the evaluation map $\hat{f}^*\sE_{\hat{Y}} \to \omega_{\hat{X}/\hat{Y}}^\nu\otimes\sM_{\hat{X}}^\mu$
is surjective over $\hat{X}_0$. Then, replacing $\hat{Y}$ by some non-singular alteration, $\hat{Z}$ by a modification of the pullback family and $\sE_{\hat{Y}}$ by its pullback, one can assume that beside the conditions (a)--(c) in~\ref{ws.5} and beside  the condition (d) in~\ref{di.5} one has:
\begin{enumerate}
\item[(e)] The images of the the evaluation maps 
$$
\hat{g}^*\sE_{\hat{Y}} \to \omega_{\hat{Z}/\hat{Y}}^\nu\otimes\sM_{\hat{Z}}^\mu
\mbox{ \ \ and \ \ }\hat{f}^*\sE_{\hat{Y}} \to \omega_{\hat{X}/\hat{Y}}^\nu\otimes\sM_{\hat{X}}^\mu
$$
are invertible sheaves. So for some divisors $\Sigma_{\hat{Z}}$ and $\Sigma_{\hat{X}}$ those images are of the form
$$ \hspace*{1cm}
\sB_{\hat{Z}}=\omega_{\hat{Z}/\hat{Y}}^\nu\otimes\sM_{\hat{Z}}^\mu\otimes\sO_{\hat{Z}}(-\Sigma_{\hat{Z}})
\mbox{ \ \ and \ \ }
\sB_{\hat{X}}=\omega_{\hat{X}/\hat{Y}}^\nu\otimes\sM_{\hat{X}}^\mu\otimes\sO_{\hat{X}}(-\Sigma_{\hat{X}}).
$$
\end{enumerate}
On the common modification $Z$ one has $\hat{\delta}^*\sB_{\hat{Z}}=\delta^*\sB_{\hat{X}}$. We denote this sheaf by $\sB_Z$.
\end{lem-not}
\begin{proof}
Consider a blowing up $\tau:Z'\to \hat{Z}$ such that the image $\sB_{Z'}$ of 
$$
\tau^*\hat{g}^*\sE_{\hat{Y}} \>>> \omega_{Z'/\hat{Y}}^\nu\otimes \tau^* \sM_{\hat{Z}}^\mu
$$
is invertible.

Let us perform the weak semistable reduction~\ref{ws.4} a second time, starting from a flattening of the morphism
$Z' \to \hat{Y}$ as explained in~\ref{ws.4} step I. By~\ref{ws.5} we obtain a mild morphism $\widetilde{g}:\widetilde{Z}_1\to \hat{Y}_1$ and a diagram
$$
\xymatrix{\hat{Z}\ar[dr]_{\hat{g}}&\ar[l]_\tau Z'\ar[d]_{g'}& \ar[l]_{\widetilde{\varphi}_1}\widetilde{Z}_1\ar[d]^{\widetilde{g}_1}\\
&\hat{Y}& \ar[l]_{\varphi_1} \hat{Y}_1.}
$$
So over $\hat{Y}_1$ we have two different mild models, $\widetilde{g}_1:\widetilde{Z}_1\to \hat{Y}_1$ and
$g_1:\hat{Z}_1\to \hat{Y}_1$, and a morphism $\tau':\widetilde{Z}_1\to \hat{Z}_1$. 
We define $\sM_{\widetilde{Z}_1}$ as the pullback of $\sM_{\hat{Z}_1}$.

The sheaf $\sF^{(\nu,\mu)}_{\hat{Y}_1}$ is independent of the mild model, and Lemma~\ref{di.5} implies that $\varphi_1^*\sF^{(\nu,\mu)}_{\hat{Y}}=\sF^{(\nu,\mu)}_{\hat{Y}_1}$. So for $\sE_{\hat{Y}_1}=\varphi^*\sE_{\hat{Y}}$ 
the pullback $\widetilde{g}_1^*\sE_{\hat{Y}_1}= \widetilde{\varphi}_1^*\tau^*\hat{g}^*\sE_{\hat{Y}}$
maps surjectively to the invertible sheaf $\sB_{\widetilde{Z}}=\widetilde{\varphi}_1^*\sB_{Z'}$. 

Since the evaluation map $\hat{f}^*_0 \sE_{\hat{Y}_0} \to \omega_{\hat{X}_0/\hat{Y}_0}^\nu\otimes\sM_{\hat{X}_0}^\mu$
is surjective, the same holds true for the pullback family, and the 
image sheaf $\sB_{\hat{X}_1}$ is locally free outside of the preimage of $\hat{Y}_0$. So replacing
$\hat{X}_0$ by a suitable non-singular modification, we may assume that it is invertible.

Replacing $\widetilde{Z}_1$ by $\hat{Z}$ and dropping the index ${}_1$ we found the invertible sheaf $\sB_{\hat{Z}}$ and $\sB_{\hat{X}}$.
Both, $\hat{\delta}^*\sB_{\hat{Z}}$ and $\delta^*\sB_{\hat{X}}$ are the images of the evaluation map
$g^*\sE_{\hat{Y}} \to  \omega_{Z/\hat{Y}}^\nu\otimes\sM_Z^\mu$, hence they coincide.
\end{proof}
Remark that the divisor $\Sigma_{\hat{X}}$ is supported in the boundary, whereas in general the divisor $\Sigma_{\hat{Z}}$ can meet $\hat{Z}_0$.

For dominant morphisms $\theta:\hat{Y}_1\to \hat{Y}$ or for morphisms from curves, whose images meet $\hat{Y}_g$, the sheaves $\sB_{\hat{Z}}$ and $\sB_{\hat{X}}$ are compatible with base change
in the following sense.

Consider $\hat{Z}_1=\hat{Z}\times_{\hat{Y}}\hat{Y}_1$ and a desingularization $\iota:\hat{X}_1 \to \hat{X}\times_{\hat{Y}}\hat{Y}_1$ 
of the main component. Writing $\sE_{\hat{Y}_1}=\theta^*\sE_{\hat{Y}}$,
the evaluation maps factor through surjections
\begin{equation}\label{eqsa.2}
\hat{g}^*_1\sE_{\hat{Y}_1} \>>> {\rm pr}_1^*\sB_{\hat{Z}} \mbox{ \ \ and \ \ }
\hat{f}^*_1\sE_{\hat{Y}_1}\>>> \iota^*{\rm pr}_1^*\sB_{\hat{X}}.
\end{equation}
On the other hand, $\sM_{\hat{Z}_1}={\rm pr}_1^*\sM_{\hat{Z}}$ and $\omega_{\hat{Z}_1/\hat{Y}_1}={\rm pr}_1^*\omega_{\hat{Z}/\hat{Y}}$. So
${\rm pr}_1^*\sB_{\hat{Z}}$ is a subsheaf of $\omega_{\hat{Z}_1/\hat{Y}_1}^\nu\otimes\sM_{\hat{Z}_1}^\mu$, and
we write $\sB_{\hat{Z}_1}={\rm pr}_1^*\sB_{\hat{Z}}$. By Corollary~\ref{di.5} 
$\sF_{\hat{Y}_1}^{(\nu,\mu)}=\theta^*\sF_{\hat{Y}}^{(\nu,\mu)}$
and Lemma~\ref{di.7} implies that the images of the second evaluation maps in (\ref{eqsa.2}) lies in
$\omega_{\hat{X}_1/\hat{Y}_1}^\nu\otimes\sM_{\hat{X}_1}^\mu$. Then $\sB_{\hat{Z}_1}$ and $\sB_{\hat{X}_1}=\iota^*{\rm pr}_1^*\sB_{\hat{X}}$
satisfy again the conditions stated in~\ref{sa.2}. 

However in~\ref{sa.2} we also changed the mild model. Using the notations from the proof
of~\ref{sa.2} we replaced $\hat{Z}_1\to \hat{Y}_1$ by a new mild model $\widetilde{Z}_1\to \hat{Y}_1$.
One is allowed to do so, if there is a birational morphism $\tau':\widetilde{Z}_1\to \hat{Z}_1$, as it
is the case in~\ref{sa.2}. One chooses $\sM_{\widetilde{Z}}$ as the pullback of $\sM_{\hat{Z}_1}$. Then
$\sB_{\widetilde{Z}}=\tau'^*\sB_{\hat{Z}_1}$ satisfies again the conditions stated in~\ref{sa.2}.
\begin{addendum}\label{sa.3} \ 
Assume that $\hat{Y}$ and $\hat{Z}$ are chosen such that the conclusion of~\ref{sa.2} holds true.
Then we may replace $\hat{Y}$ by a non-singular alteration $\hat{Y}_1$ and the pullback of the given mild model 
$\hat{Z}_1\to \hat{Y}_1$ by any mild morphism $\widetilde{Z}_1\to \hat{Y}_1$ provided there is a morphism $\tau':\widetilde{Z}_1\to \hat{Z}_1$, birational over $\hat{Y}_1$.
\end{addendum}
In particular, given a finite number of $(\nu,\mu)\in I$, and a finite number of sheaves $\sE_{\hat{Y}}$, one can apply~\ref{sa.2} successively.

Since we assumed that $\sF_{\hat{Y}}^{(\nu,\mu)}$ is locally free, one possible choice for $\sE_{\hat{Y}}$
is the sheaf $\sF_{\hat{Y}}^{(\nu,\mu)}$ itself.
\begin{notations}\label{sa.4}
Consider in~\ref{di.5} a subset $\widetilde{I}\subset I$ and assume that for $(\nu,\mu) \in \widetilde{I}$ the evaluation map
$$
f_0^*f_{0*}(\omega_{X_0/Y_0}^\nu\otimes\sL_0^\mu) \>>> \omega_{X_0/Y_0}^\nu\otimes\sL_0^\mu
$$
is surjective. If one chooses in~\ref{sa.2} $\sE_{\hat{Y}}=\sF_{\hat{Y}}^{(\nu,\mu)}$,
we will write $\Sigma^{(\nu,\mu)}_{\bullet}$ and $\sB_\bullet^{(\nu,\mu)}$ instead of
$\Sigma_{\bullet}$ and $\sB_\bullet$, where $\bullet$ stands for $\hat{Z}$, $\hat{X}$ or $Z$.
In particular 
$$
\sB^{(\nu,\mu)}_{\bullet}=\omega_{\bullet/\hat{Y}}^\nu\otimes\sM_{\bullet}^\mu\otimes\sO_{\bullet}(-\Sigma^{(\nu,\mu)}_{\bullet}).
$$
If $\mu=0$ we will write $\varpi_\bullet^{(\nu)}$ and $\Pi_\bullet^{(\nu)}$ instead
of $\sB_\bullet^{(\nu,0)}$ and $\Sigma^{(\nu,0)}_{\bullet}$.
\end{notations}
Let us collect the properties we can require for a well chosen non-singular alteration
$\hat{Y}\to Y$ and for the morphisms in the diagram (\ref{eqsa.1}). 
\begin{con-not}\label{sa.5} 
We start with a finite set $I$ of tuples $(\nu,\mu)$ of natural numbers, and 
with a subset $\widetilde{I}$ of $I$. We assume that for some $\eta_0>0$ with $(\eta_0,0) \in \widetilde{I}$ the evaluation map $f_0^*f_{0*}\omega_{X_0/Y_0}^{\eta_0} \to \omega_{X_0/Y_0}^{\eta_0}$
is surjective, and that for all other $(\eta,0)\in \widetilde{I}$ the natural number
$\eta$ is divisible by $\eta_0$.

Then we can find $\hat{Y}$ and the diagram \eqref{eqws.3} (recalled in \eqref{eqsa.1}) such that:
\begin{enumerate}
\item[i.] The conditions (a), (b) and (c) in Proposition~\ref{ws.5} hold true, as well as the conditions 
i) and ii) in~\ref{ws.6}.
\item[ii.] For $(\eta,0)\in \widetilde{I}$ there are invertible sheaves $\varpi_{\hat{Z}}^{(\eta)}$,
$\varpi_{Z}^{(\eta)}$, and $\varpi_{\hat{X}}^{(\eta)}$ on $\hat{Z}$, $Z$ and on $\hat{X}$, respectively,
with surjective evaluation maps, with
\begin{gather*}
\varpi_{Z}^{(\eta)}=\hat{\delta}^*\varpi_{\hat{Z}}^{(\eta)}=\delta^*\varpi_{\hat{X}}^{(\eta)}
\mbox{ \ \ and with}\\
\sF_{\hat{Y}}^{(\eta)}:=
\sF_{\hat{Y}}^{(\eta,0)}=\hat{g}_*\omega_{\hat{Z}/\hat{Y}}^{\eta}=
\hat{g}_*\varpi_{\hat{Z}}^{(\eta)}=\hat{f}_*\varpi_{\hat{X}}^{(\eta)}.
\end{gather*}
\item[iii.] For all $(\nu,0)\in I$ the sheaves
$\sF_{\hat{Y}}^{(\nu)}:=
\sF_{\hat{Y}}^{(\nu,0)}=\hat{g}_*\omega_{\hat{Z}/\hat{Y}}^{\nu}$
are locally free.
\item[iv.] There is an open dense subscheme $\hat{Y}_g$ with $\hat{g}^{-1}(\hat{Y}_g)\to \hat{Y}_g$ smooth
such that for all $(\nu,0)\in I$ the sheaves $\sF_{\hat{Y}}^{(\nu)}=\hat{g}_*\omega_{\hat{Z}/\hat{Y}}^{\nu}$
are compatible with base change for morphisms $\varrho:T\to \hat{Y}$ 
with $\varrho^{-1}(\hat{Y}_g)$ dense in $T$.
\item[v.] In ii) $\Pi^{(\eta)}_{Z}$, $\Pi^{(\eta)}_{\hat{Z}}$ and $\Pi^{(\eta)}_{\hat{X}}$ denote the divisors with
\begin{gather*}
\omega_{Z/\hat{Y}}^{\eta}=\varpi_{Z}^{(\eta)}\otimes \sO_{Z}(\Pi^{(\eta)}_{Z}) \ \ \ \ \ \ 
\omega_{\hat{Z}/\hat{Y}}^{\eta}=\varpi_{\hat{Z}}^{(\eta)}\otimes \sO_{\hat{Z}}(\Pi^{(\eta)}_{\hat{Z}})\\
\mbox{ \ \ and \ \ }
\omega_{\hat{X}/\hat{Y}}^{\eta}=\varpi_{\hat{X}}^{(\eta)}\otimes \sO_{\hat{X}}(\Pi^{(\eta)}_{\hat{X}}).
\end{gather*}
\end{enumerate}
\end{con-not}
\begin{con-not}\label{sa.6} {\bf (Canonical polarizations)}\ \\
All we need in this case is collected in \ref{sa.5}. We will of course choose
$\widetilde{I}$ and $I$ as subsets of $\N\times \{0\}$. 
\end{con-not}
If $\sL_0\neq \sO_{X_0}$, i.e. if we consider polarized manifolds, we will need
more:
\begin{con-not}\label{sa.7} {\bf (Polarizations)}\ \\
We consider in~\ref{sa.5} an invertible sheaf $\sL$ on $X$ with $\sL_0=\sL|_{X_0}$
$f_0$-ample, and we choose $\gamma_0>0$ such that the evaluation map
$f_0^*f_{0*}\sL_0^{\gamma_0} \to \sL_0^{\gamma_0}$
is surjective. We fix some subset $\widetilde{I}$ of $I$ consisting of tuples
$(\beta,\alpha)$ of natural numbers with $\alpha$ divisible by $\gamma_0$ and with $\beta$ divisible by $\eta_0$. By Lemma~\ref{di.3} the direct images $f_{0*}(\omega_{X_0/Y_0}^\nu\otimes \sL_0^\mu)$
are locally free and compatible with arbitrary base change, whenever $\nu >0$ and $\mu\geq 0$.
For $(0,\mu)\in I$ we have to add the corresponding statement to the list of assumptions.

Then we can find $\hat{Y}$ and the diagram \eqref{eqws.3} such that
the conditions i) -- v) in~\ref{sa.5} hold true and in addition:
\begin{enumerate}
\item[vi.] $\sM_{\hat{Z}}$, $\sM_Z$, and $\sM_{\hat{X}}$ are the pullback of $\sL$.
\item[vii.] For $(\beta,\alpha)\in \widetilde{I}$ there are invertible sheaves $\sB_{\hat{Z}}^{(\beta,\alpha)}$,
$\sB_{Z}^{(\beta,\alpha)}$, and $\sB_{\hat{X}}^{(\beta,\alpha)}$ on $\hat{Z}$, $Z$ and on $\hat{X}$, respectively, with surjective evaluation maps, with
\begin{gather*}
\sB_{Z}^{(\beta,\alpha)}=\hat{\delta}^*\sB_{\hat{Z}}^{(\beta,\alpha)}=\delta^*\sB_{\hat{X}}^{(\beta,\alpha)}
\mbox{ \ \ and with}\\
\sF_{\hat{Y}}^{(\beta,\alpha)}=\hat{g}_*(\omega_{\hat{Z}/\hat{Y}}^{\beta}\otimes\sM_{\hat{Z}}^{\alpha})=
\hat{g}_*\sB_{\hat{Z}}^{(\beta,\alpha)}=\hat{f}_*\sB_{\hat{X}}^{(\beta,\alpha)}.
\end{gather*}
\item[viii.] For all $(\nu,\mu) \in I$ the sheaves
$\sF_{\hat{Y}}^{(\nu,\mu)}=\hat{g}_*\big(\omega_{\hat{Z}/\hat{Y}}^{\nu}\otimes \sM_{\hat{Z}}^{\mu}\big)$
are locally free.
\item[ix.] There is an open dense subscheme $\hat{Y}_g$ with $\hat{g}^{-1}(\hat{Y}_g)\to \hat{Y}_g$ smooth
such that for all $(\nu,\mu)\in I$ the sheaves 
$\sF_{\hat{Y}}^{(\nu,\mu)}=\hat{g}_*\big(\omega_{\hat{Z}/\hat{Y}}^{\nu}\otimes \sM_{\hat{Z}}^{\mu}\big)$
are compatible with base change for morphisms $\varrho:T\to \hat{Y}$ 
with $\varrho^{-1}(\hat{Y}_g)$ dense in $T$.
\item[x.] $\Sigma^{(\beta,\alpha)}_{Z}$, $\Sigma^{(\beta,\alpha)}_{\hat{Z}}$ and $\Sigma^{(\beta,\alpha)}_{\hat{X}}$ denote the divisors with
\begin{gather*}
\omega_{Z/\hat{Y}}^\beta\otimes\sM_{Z}^{\alpha}=\sB_{Z}^{(\beta,\alpha)}\otimes \sO_{Z}( \Sigma^{(\beta,\alpha)}_{Z}) \ \ \ \ \ \ 
\omega_{\hat{Z}/\hat{Y}}^\beta\otimes\sM_{\hat{Z}}^{\alpha}=\sB_{\hat{Z}}^{(\beta,\alpha)}\otimes \sO_{\hat{Z}}(\Sigma^{(\beta,\alpha)}_{\hat{Z}})\\
\mbox{ \ \ and \ \ }
\omega_{\hat{X}/\hat{Y}}^\beta\otimes\sM_{\hat{X}}^{\alpha}=\sB_{\hat{X}}^{(\beta,\alpha)}\otimes \sO_{\hat{X}}( \Sigma^{(\beta,\alpha)}_{\hat{X}}).
\end{gather*}
\end{enumerate}
\end{con-not}
\begin{aconstr}\label{sa.8} The conditions stated in~\ref{sa.5} and~\ref{sa.7}
and the sheaves $\sF_\bullet^{(\nu,\mu)}$ for $(\nu,\mu)\in I$ are compatible with the following constructions:
\begin{enumerate}
\item[I.] Replace $\hat{Y}$ by a non-singular alteration, $\hat{Z}$ by its pullback, and $\hat{X}$
by a desingularization of the main component of its pullback.
\item[II.] Replace $\hat{Z}$ by a mild morphism $\widetilde{Z}\to \hat{Y}$, for which there is a birational $\hat{Y}$-morphism $\tau:\widetilde{Z}\to \hat{Z}$.
\end{enumerate}
In particular assume that for some open set $U\subset \hat{Y}$ containing $\hat{Y}_0$ the morphism $f^{-1}(U)\to \hat{Y}$
is flat. Then one can choose a mild morphism $\widetilde{Z}_1\to \hat{Y}_1$
factoring through $\tau_1:\widetilde{Z}_1\to \hat{X}_1$, and still assume that~\ref{sa.5} and~\ref{sa.7}
holds true.
\end{aconstr}
\begin{proof}
This has been shown in Addendum~\ref{sa.3}. For the last part, one performs the weak semistable
reduction, starting with $\hat{X}\to \hat{Y}$ instead of $\widetilde{X}\to \widetilde{Y}$ in step I of~\ref{ws.4}.
\end{proof}
Next we will start to construct the saturated extensions of the polarizations. Although this will only be applied for families of Kodaira dimension zero, we will allow $\omega_{X_0/Y_0}$ to be $f_0$-semiample. 
\begin{lemma}\label{sa.9} Let $\sM_{\hat{Z}}$, $\sM_{\hat{X}}$ and $\sM_{Z}$ be invertible sheaves on $\hat{Z}$, $\hat{X}$ and $Z$, respectively, satisfying the compatibility conditions in~\ref{di.4}.
Assume that $\kappa$ is a positive integer with $(0,\kappa)\in I$.
Using the notations and conditions in~\ref{sa.5} one has:
\begin{enumerate}
\item For all $\varepsilon\geq 0$ and for all alterations $\hat{Y}_1$ of $\hat{Y}$
\begin{gather*}
\hat{\delta}_*(\sM_{Z_1}^\kappa\otimes \sO_{Z_1}(\varepsilon\cdot \Pi^{(\eta_0)}_{Z_1}))=
\sM_{\hat{Z}_1}^\kappa\otimes \sO_{\hat{Z}_1}(\varepsilon\cdot \Pi^{(\eta_0)}_{\hat{Z}_1})\mbox{ \ \ and}\\
 \delta_*(\sM_{Z_1}^\kappa\otimes \sO_{Z_1}(\varepsilon\cdot \Pi^{(\eta_0)}_{Z_1}))=
\sM_{\hat{X}_1}^\kappa\otimes \sO_{\hat{X}_1}(\varepsilon\cdot \Pi^{(\eta_0)}_{\hat{X}_1}).
\end{gather*}
\item For each $\kappa>0$ there exists some $\varepsilon_0 \geq 0$ such that
$$
\iota:\hat{g}_{1*}\sM_{\hat{Z}_1}^\kappa\otimes \sO_{\hat{Z}_1}(\varepsilon_0\cdot \Pi^{(\eta_0)}_{\hat{Z}_1})\>>> \hat{g}_{1*}\sM_{\hat{Z}_1}^\kappa\otimes \sO_{\hat{Z}_1}(\varepsilon\cdot \Pi^{(\eta_0)}_{\hat{Z}_1})
$$
are isomorphisms for all $\varepsilon \geq \varepsilon_0$, and for all alterations $\hat{Y}_1$ of $\hat{Y}$.
\end{enumerate}
\end{lemma}
Remark that (1) and (2) imply that for all $\varepsilon \geq \varepsilon_0$ one also has
$$
\hat{f}_{1*}\sM_{\hat{X}_1}^\kappa\otimes \sO_{\hat{X}_1}(\varepsilon_0\cdot \Pi^{(\eta_0)}_{\hat{X}_1})\>\cong >> \hat{f}_{1*}\sM_{\hat{X}_1}^\kappa\otimes \sO_{\hat{X}_1}(\varepsilon\cdot \Pi^{(\eta_0)}_{\hat{X}_1}).
$$
\begin{proof}[Proof of~\ref{sa.9}] Let us replace $\sM^\kappa_\bullet$ by $\sM_\bullet$, hence assume that $\kappa=1$.
For (1) consider the common modification $Z$. By~\ref{sa.5} ii)
\begin{gather*}
\varpi^{(\eta_0)}_Z=\hat{\delta}^*\varpi_{\hat{Z}}^{(\eta_0)}=\delta^*\varpi^{(\eta_0)}_{\hat{X}},
\mbox{ \ \ and}\\
\Pi^{(\eta_0)}_Z = \hat{\delta}^* \Pi^{(\eta_0)}_{\hat{Z}}+\eta_0\cdot E_{\hat{Z}}=\delta^* \Pi^{(\eta_0)}_{\hat{X}}+\eta_0\cdot E_{\hat{X}},
\end{gather*}
where $E_\bullet$ are effective relative canonical divisors for $Z/\bullet$.
The assumptions $\delta_*\sM_{Z}=\sM_{\hat{X}}$ and $\delta_*\sM_{Z}=\sM_{\hat{X}}$
imply that 
$$
\sM_{Z}=\hat{\delta}^*\sM_{\hat{Z}}\otimes \sO_{Z}(F_{\hat{Z}})=\delta^*\sM_{\hat{X}}\otimes \sO_{Z}(F_{\hat{X}})
$$
for effective exceptional divisors $F_{\hat{Z}}$ and $F_{\hat{X}}$, and (1) for $\hat{Y}_1=\hat{Y}$ follows
from the projection formula. The same argument works over any alteration.

For (2) remark that one may replace $\hat{Y}_1$ by a modification $\theta:\hat{Y}_2$ and $\hat{Z}_1$ by the pullback
family $\hat{Z}_2=\hat{Z}_1\times_{\hat{Y}_1}\hat{Y}_2 \to \hat{Y}_2$. In fact, the divisor $\Pi_{\hat{Z}_1}$ is compatible with
pullback, and for all $\varepsilon \geq 0$ one has
$$
{\rm pr}_{1*}(\sM_{\hat{Z}_2}\otimes \sO_{\hat{Z}_2}(\epsilon\cdot\Pi_{\hat{Z}_2}))=
\sM_{\hat{Z}_1}\otimes \sO_{\hat{Z}_1}(\varepsilon\cdot\Pi_{\hat{Z}_1}).
$$
Hence $\theta_*\hat{g}_{2*}\sM_{\hat{Z}_2}^\kappa\otimes \sO_{\hat{Z}_2}(\varepsilon\cdot \Pi^{(\eta_0)}_{\hat{Z}_2})=
\hat{f}_{1*}\sM_{\hat{X}_1}^\kappa\otimes \sO_{\hat{X}_1}(\varepsilon\cdot \Pi^{(\eta_0)}_{\hat{X}_1})$
and if the first sheaf is independent of $\varepsilon$, for $\varepsilon$ sufficiently large,
the same holds for the second one.
 
The fibres of $\hat{Z}\to \hat{Y}$ are reduced. Then the compatibility of $\sF_{\hat{Y}}^{(\eta_0)}$ with
pullback under alterations and the surjectivity of the evaluation map for
$\omega_{\hat{Z}/\hat{Y}}^{\eta_0}\otimes\sO_{\hat{Z}}(- \Pi^{(\eta_0)}_{\hat{Z}})$ imply that 
$\Pi^{(\eta_0)}_{\hat{Z}}$ can not contain a whole fibre. Otherwise, for
some sheaf of ideals $\sJ$ on $\hat{Y}$ one would have $\varpi_{\hat{Z}}^{(\eta_0)}\subset \hat{g}^*\sJ\otimes\omega_{\hat{Z}/\hat{Y}}^{\eta_0}$. Blowing up $\hat{Y}$ one gets the same, with
$\sJ=\sO_{\hat{Y}}(-\Gamma)$, for an effective divisor $\Gamma$. Then the projection formula implies that
$\hat{g}_*\varpi_{\hat{Z}}^{(\eta_0)}\subset \sJ\otimes \hat{g}_*\omega_{\hat{Z}/\hat{Y}}^{\eta_0}$, contradicting~\ref{sa.5} ii).

By flat base change, the question whether $\iota$ is an isomorphism is local for the \'etale topology. So by abuse of notations we may replace $\hat{Y}$ by any \'etale neighborhood. Hence given $y\in \hat{Y}$ 
we may assume that $\hat{g}$ has a section $\sigma:\hat{Y}\to \hat{Z}$ whose image does not meet $\Pi^{(\eta_0)}_{\hat{Z}}$, but meets the open set $V_0$ where $\hat{\varphi}_0:\hat{Z}_0\to \hat{X}_0$ is an isomorphism.
Let $\sI$ be the ideal sheaf of $\sigma(\hat{Y})$. For a general fibre $F$ of $\hat{f}$
and for $\upsilon$ sufficiently large $H^0(F,(\varphi_*\sI^\upsilon)\otimes\sM_{\hat{X}}|_F)=0$.
Then
$$
\hat{g}_{0*}((\sI^\upsilon\otimes \sM_{\hat{Z}}\otimes\sO_{\hat{Z}}(\varepsilon\cdot \Pi^{(\eta_0)}_{\hat{Z}}))|_{\hat{Z}_0})=
\hat{f}_{0*}(\hat{\varphi}_{0*}(\sI^\upsilon\otimes\sO_{\hat{Z}}(\varepsilon\cdot \Pi^{(\eta_0)}_{\hat{Z}}))|_{\hat{Z}_0} \otimes \sM_{\hat{X}_0})=0,
$$
and $\hat{g}_{*}\sM_{\hat{Z}}\otimes \sO_{\hat{Z}}(\varepsilon\cdot \Pi^{(\eta_0)}_{\hat{Z}})$ is a subsheaf of
$\hat{g}_*\sM_{\hat{Z}}/\sI^\upsilon=\hat{g}_*(\sO_{\hat{Z}}(\varepsilon\cdot \Pi^{(\eta_0)}_{\hat{Z}})
\otimes\sM_{\hat{Z}}/\sI^\upsilon)$.
So $\sC=\hat{g}_{*}\sM_{\hat{Z}}\otimes \sO_{\hat{Z}}(*\cdot\Pi^{(\eta_0)}_{\hat{Z}})$ as a subsheaf of a fixed locally free sheaf is isomorphic to $\hat{g}_{*}\sM_{\hat{Z}}\otimes \sO_{\hat{Z}}(\varepsilon_1\cdot \Pi^{(\eta_0)}_{\hat{Z}})$
for some $\varepsilon_1$. 

Let $\theta:\hat{Y}_2\to \hat{Y}$ be a modification, such that $\sC_2=\theta^*\sC/_{\rm torsion}$ is locally free,
and contained in a locally free locally splitting subsheaf $\sC'$ of $\theta^*\hat{g}_{1*}(\sM_{\hat{Z}_1}/\sI^\upsilon_1)$ with $\rk(\sC')=\sC_2)$. Writing $\sI_2$ for the pullback of the sheaf of ideals $\sI$, the latter is of the form $\hat{g}_{2*}(\sM_{\hat{Z}_2}/\sI^\upsilon_2)$. For some effective divisor $D$ one has an inclusion $\sC'\subset \sC_2\otimes \sO_{\hat{Y}_2}(D)$. The base change morphism
$$
\theta^*\hat{g}_{*}\sM_{\hat{Z}}\otimes \sO_{\hat{Z}}(\varepsilon\cdot \Pi^{(\eta_0)}_{\hat{Z}})\>>>
\hat{g}_{2*}\sM_{\hat{Z}_2}\otimes \sO_{\hat{Z}_2}(\varepsilon\cdot \Pi^{(\eta_0)}_{\hat{Z}_2})
$$
implies that for all $\varepsilon \geq \varepsilon_1$
\begin{multline*}
\sC_2 \subset \hat{g}_{2*}\sM_{\hat{Z}_2}\otimes \sO_{\hat{Z}_2}(\varepsilon\cdot \Pi^{(\eta_0)}_{\hat{Z}_2})
\subset \sC'\subset \sC_2\otimes \sO_{\hat{Y}_2}(D)\\
\subset \hat{g}_{2*}\sM_{\hat{Z}_2}\otimes \sO_{\hat{Z}_2}(\varepsilon_1\cdot \Pi^{(\eta_0)}_{\hat{Z}_2}+\hat{g}^*_2 D).
\end{multline*}
Let us choose $\varepsilon_0\geq \varepsilon_1$ such that for an irreducible Weil divisors $\Pi$ the multiplicity in $(\varepsilon_0-\varepsilon_1)\cdot \Pi^{(\eta_0)}_{\hat{Z}_2}$ is either zero, or larger that its multiplicity
in $g_2^*D$. Remark already, that this choice of $\varepsilon_0$ is compatible with further pullback.

For $\varepsilon \geq \varepsilon_0$ the image of the evaluation map 
$$
\hat{g}^*_2\hat{g}_{2*}\sM_{\hat{Z}_2}\otimes \sO_{\hat{Z}_2}(\varepsilon\cdot \Pi^{(\eta_0)}_{\hat{Z}_2})\>>> \sM_{\hat{Z}_2}\otimes \sO_{\hat{Z}_2}(*\cdot \Pi^{(\eta_0)}_{\hat{Z}_2})
$$
is contained in the image of $\hat{g}^*_2\sC'\to \sM_{\hat{Z}_2}\otimes \sO_{\hat{Z}_2}(*\cdot \Pi^{(\eta_0)}_{\hat{Z}_2})$
hence in 
$$
\sM_{\hat{Z}_2}\otimes \sO_{\hat{Z}_2}(\varepsilon_1\cdot \Pi^{(\eta_0)}_{\hat{Z}_2}+\hat{g}^*_2 D)\cap
\sM_{\hat{Z}_2}\otimes \sO_{\hat{Z}_2}(*\cdot \Pi^{(\eta_0)}_{\hat{Z}_2})
\subset
\sM_{\hat{Z}_2}\otimes \sO_{\hat{Z}_2}(\varepsilon_0\cdot \Pi^{(\eta_0)}_{\hat{Z}_2}).
$$
We found $\varepsilon_0$ after replacing $\hat{Y}$ by some non-singular modification $\hat{Y}_2$, hence as remarked above the same $\varepsilon_0$ works for $\hat{Y}$ itself. Moreover, the same $\varepsilon_0$ works for
all alterations dominating $\hat{Y}_2$. Since for any alteration $\hat{Y}_1$ of $\hat{Y}$ one can find a non-singular modification, dominating $\hat{Y}_2$, one obtains the same for $\hat{Y}_1$.
\end{proof}
\begin{definition}\label{sa.10}
Assume that $\sL$ is an invertible sheaf on $X$, and let $\kappa$ be a positive integer. Assume
that $f_{0*}\sL_0^\kappa$ is locally free and compatible with arbitrary base change.
\begin{enumerate}
\item An invertible sheaf $\sM_{\hat{Z}}$ on $\hat{Z}$ is a {\em $\kappa$-saturated extension} of $\sL$ if
\begin{gather}\label{eqsa.3}
\hat{\varphi}^*\sL \subset \sM_{\hat{Z}}\subset \hat{\varphi}^*\sL\otimes \big( \sO_{\hat{Z}}(*\cdot\Pi_{\hat{Z}}^{(\eta_0)})\cap
\sO_{\hat{Z}}(*\cdot\hat{g}^{-1}(\hat{Y}\setminus \hat{Y}_0)\big),\\
\notag \mbox{and if \ \ \ \ }  
\hat{g}_{1*}\sM_{\hat{Z}_1}^\kappa = \hat{g}_{1*}(\sM_{\hat{Z}_1}^\kappa\otimes \sO_{\hat{Z}_1}(\varepsilon\cdot\Pi^{(\eta_0)}_{\hat{Z}_1}))
\end{gather}
for all $\varepsilon\geq 0$ and for all alterations $\hat{Y}_1\to \hat{Y}$.
Moreover we require~\ref{di.5} d) to hold for $(\nu,\mu)=(0,\kappa)$, i.e. that there exists an open dense subscheme $\hat{Y}_g$ of $\hat{Y}$ such that $\hat{g}_*\sM_{\hat{Z}}^\kappa$ is locally free and compatible with pullback for morphisms $\theta:T\to \hat{Y}$ with $\theta^{-1}(\hat{Y}_g)$ dense in $T$.
\item  
We call a tuple of invertible sheaves $\sM_{\hat{Z}}$, $\sM_{\hat{X}}$ and $\sM_{Z}$ on $\hat{Z}$, $\hat{X}$ and $Z$ a {\em $\kappa$-saturated extension of the polarization $\sL$}, if $\sM_{\hat{Z}}$ is $\kappa$ saturated and if
(as in~\ref{di.4}) $\hat{\delta}_*\sM_{Z}=\sM_{\hat{Z}}$, $\delta_*\sM_Z=\sM_{\hat{X}}$, $\sM_{\hat{Z}_0}=\hat{\varphi}^*_0\sL_0$ and $\sM_{\hat{X}_0}=\rho^*_0\sL_0$.
\end{enumerate}
\end{definition}
\begin{lemma}\label{sa.11} Assume that the conditions in~\ref{sa.5} hold true.
\begin{enumerate}
\item[a.]
If $\sM_{\hat{Z}}$ is a $\kappa$-saturated extension of $\sL$,
one can always find $\sM_{\hat{X}}$ and $\sM_{Z}$ such that
$(\sM_{\hat{Z}},\sM_{\hat{X}},\sM_{Z})$ is $\kappa$-saturated.
\item[b.] The condition (\ref{eqsa.3}) in (1) is equivalent to the existence of an effective Cartier
divisor $\hat{\Pi}$, supported in $\hat{g}^{-1}(\hat{Y}\setminus \hat{Y}_0)\cap (\Pi_{\hat{Z}}^{(\eta_0)})_{\rm red}$, 
and with
$$
\sM_{\hat{Z}}=\hat{\varphi}^*\sL\otimes \sO_{\hat{Z}}(\hat{\Pi}).
$$ 
\item[c.] If $(\sM_{\hat{Z}},\sM_{\hat{X}},\sM_{Z})$ is $\kappa$-saturated
$$
\hat{f}_{1*}\sM_{\hat{X}_1}^\kappa = \hat{f}_{1*}(\sM_{\hat{X}_1}^\kappa\otimes \sO_{\hat{X}_1}(\varepsilon\cdot\Pi^{(\eta_0)}_{\hat{X}_1}))=
\hat{f}_{*}(\rho^*\sL^\kappa \otimes \sO_{\hat{X}_1}(*\cdot\Pi^{(\eta_0)}_{\hat{X}}))
$$
for all $\varepsilon\geq 0$ and for all alterations $\hat{Y}_1\to \hat{Y}$.
\item[d.] Let $\widetilde{g}:\widetilde{Z}\to \hat{Y}$ be a second mild morphism and $\tau':\widetilde{Z}_1 \to \hat{Z}_1$
a birational morphism over $\hat{Y}$. If $\sM_{\hat{Z}}$ is $\kappa$-saturated the same holds for 
$\sM_{\widetilde{Z}}=\tau'^*\sM_{\hat{Z}}$.
\item[e.] If $\sM_{\hat{Z}}$ (or $(\sM_{\hat{Z}},\sM_{\hat{X}},\sM_{Z})$) is $\kappa$-saturated, and if $\kappa'$ divides $\kappa$ then $\sM_{\hat{Z}}$ (or $(\sM_{\hat{Z}},\sM_{\hat{X}},\sM_{Z})$) is also $\kappa'$-saturated, provided that
$\hat{g}_*\sM_{\hat{Z}}^{\kappa'}$ is locally free and compatible with base change
for morphisms $\theta:T\to \hat{Y}$ with $\theta^{-1}(\hat{Y}_g)$ dense in $T$.
\end{enumerate}
\end{lemma}
\begin{proof}
b) is just a translation and the first equality in c) follows directly from~\ref{sa.9}. 
For the second one, apply~\ref{sa.9} first to the pullback of $\sL$ and then to $\sM_\bullet$. One finds that
$\hat{f}_{1*}\sM_{\hat{X}_1}^\kappa$ is given by
$$
\hat{f}_{1*}(\rho_1^*\sL^\kappa \otimes \sO_{\hat{X}_1}(*\cdot\Pi^{(\eta_0)}_{\hat{X}_1}))=
\hat{g}_{1*}\hat{\varphi}^*_1\sL^\kappa \otimes \sO_{\hat{Z}_1}(*\cdot\Pi_{\hat{Z}_1}^{(\eta_0)}) = \hat{g}_{1*}\sM_{\hat{Z}_1}^\kappa \otimes \sO_{\hat{Z}_1}(*\cdot\Pi_{\hat{Z}_1}^{(\eta_0)}).
$$
For a) consider $\Pi=\hat{\delta}^*\hat{\Pi}$ and the divisor $\delta_*\Pi$ on $\hat{X}$. 
Define $\sM_{\hat{X}}=\rho^*\sL\otimes \sO_{\hat{X}}(\delta_*\Pi)$.
 
Since $\delta$ is a modification of a manifold, $\Pi-\delta^*\delta_*\Pi$ is supported in exceptional divisors for $\delta$, and 
$$
\delta^*\sM_{\hat{X}}\subset \delta^*\rho^*\sL \otimes \sO_{Z}(\Pi)=
\hat{\delta}^*\varphi^*\sL \otimes \sO_{Z}(\hat{\delta}^*\hat{\Pi})= \hat{\delta}^*\sM_{\hat{Z}},
$$
and $\sM_{\hat{X}} = \delta_* \hat{\delta}^*\sM_{\hat{Z}}$. So we can choose $\sM_{Z}=\hat{\delta}^*\sM_{\hat{Z}}$.

In d) remark that $\varpi_{\hat{Z}}^{(\eta_0)}$ is invertible and its pullback is $\varpi_{\widetilde{Z}}^{(\eta_0)}$.
So $\Pi_{\widetilde{Z}}^{(\eta_0)}-\tau'^*\Pi_{\hat{Z}}^{(\eta_0)}$ is an effective divisor, supported in the exceptional locus of $\tau$. By the projection formula, for all $\varepsilon \geq 0$,
$$
\tau'_*\sM_{\widetilde{Z}}^{\kappa}\otimes \sO_{\widetilde{Z}}(\varepsilon\cdot \Pi_{\widetilde{Z}}^{(\eta_0)})=
\sM_{\hat{Z}}^\kappa\otimes \sO_{\hat{Z}}(\varepsilon\cdot \Pi_{\hat{Z}}^{(\eta_0)}),
$$
hence $\widetilde{g}_{*}(\sM_{\widetilde{Z}}^\kappa\otimes \sO_{\widetilde{Z}}(\varepsilon\cdot \Pi_{\widetilde{Z}}^{(\eta_0)}))= \hat{g}_{*}(\sM_{\hat{Z}}^\kappa\otimes \sO_{\hat{Z}}(\varepsilon\cdot \Pi_{\hat{Z}}^{(\eta_0)}))$.
Since the right hand side is independent of $\varepsilon$ the polarization $\sM_{\widetilde{Z}}$ is again $\kappa$-saturated.

For e) remark first, that the condition (2) in Definition~\ref{sa.10} is independent of $\kappa$, as well
as (\ref{eqsa.3}) in (1). If for some $\hat{Y}_1\to \hat{Y}$ and some $\varepsilon >0$ the sheaf
$$
\hat{g}_{1*}\sM_{\hat{Z}_1}^{\kappa'} \neq \hat{g}_{1*}(\sM_{\hat{Z}_1}^{\kappa'}\otimes \sO_{\hat{Z}_1}(\varepsilon\cdot\Pi^{(\eta_0)}_{\hat{Z}_1})),
$$
then the multiplication map shows, that the same holds for all multiples of $\kappa'$, in particular
for $\kappa$.
\end{proof}
\begin{lemma}\label{sa.12} Given a natural number $\kappa$ one may choose $\hat{Y}$ and $\hat{Z}$ in~\ref{di.4} and the sheaf $\sM_{\hat{Z}}$ such that $\sM_{\hat{Z}}$ is a $\kappa$-saturated extension of $\sL$.
\end{lemma}
\begin{proof}
Start with any $\hat{Y}$ as in~\ref{sa.5} and with $\sM_{\hat{Z}}$ the pullback of the invertible sheaf $\sL$ in~\ref{ws.8}. Apply~\ref{sa.9} to the polarization $\sM_{\hat{Z}}^\kappa$, and replace $\varepsilon_0$
by some larger natural number, divisible by $\kappa$. 

Define $\hat{\Pi}$ to be the sum over all components of $\Pi^{(\eta_0)}_{\hat{Z}}$ 
whose image in $\hat{Y}$ does not meet $\varphi_1^{-1}(Y_0)$, and choose 
$$
\widetilde{\sM}_{\hat{Z}}=\sM_{\hat{Z}}\otimes \sO_{\hat{Z}}\big(\frac{\varepsilon_0}{\kappa}\cdot \hat{\Pi}\big).
$$
Remark that $\hat{\Pi}$ might be just a Weil divisor, hence $\widetilde{\sM}_{\hat{Z}}$ is reflexive, but not necessarily invertible. So choose a modification $\sigma: W \to \hat{Z}$, such that $\sM_W=\sigma^*\widetilde{\sM}_{\hat{Z}}/_{\rm torsion}$ is invertible. By Proposition~\ref{ws.5} there exists a non-singular alteration $\theta:\hat{Y}_1\to \hat{Y}$
such that $W\otimes_{\hat{Y}}\hat{Y}_1$ has a mild model $W'\to \hat{Y}_1$. Again we may assume that the conditions in~\ref{sa.5} hold for $W'\to \hat{Y}_1$. One has a factorization $W'\to W \to \hat{Z}$ of $\sigma$, inducing a birational morphism
$\sigma':W'\to \hat{Z}_1=\hat{Z}\times_{\hat{Y}}\hat{Y}_1$.
By~\ref{sa.9} (2) we know that the evaluation map
$$
\hat{g}^*_1\hat{g}_{1*}\sM^\kappa_{\hat{Z}_1}(*\cdot \Pi_{\hat{Z}_1}^{(\eta_0)})\>>> \sM^\kappa_{\hat{Z}_1}(*\cdot \Pi_{\hat{Z}_1}^{(\eta_0)})
$$
has image $\sC$ in $\sM^\kappa_{\hat{Z}_1}(\varepsilon_0\cdot \Pi_{\hat{Z}_1}^{(\eta_0)})$. On the other hand,
on $\hat{g}^{-1}(\theta^{-1}(\hat{Y}_0))$ the sheaf $\sC$ is equal to $\sM^\kappa_{\hat{Z}_1}$ and
$\sC$ lies in the reflexive hull $\widetilde{\sM}_{\hat{Z}_1}^{(\kappa)}$ of ${\rm pr}_1^*\widetilde{\sM}_{\hat{Z}}^\kappa$.
By construction $\sM_{W'}=\sigma^*\widetilde{\sM}_{\hat{Z}_1}/_{\rm torsion}$ is invertible and
$\sigma^*\widetilde{\sM}_{\hat{Z}_1}^{(\kappa)}/_{\rm torsion}=\sM_{W'}^\kappa$.

Writing again $\Pi_{W'}^{(\eta_0)}$ for the relative fix locus of $\omega_{W'/\hat{Y}_1}^{\eta_0}$ one has
$$
\varpi_{W'}^{(\eta_0)}=\omega_{W'/\hat{Y}_1}^{\eta_0}\otimes\sO_{W'}(-\Pi_{W'}^{(\eta_0)})=
\sigma'^*\varpi_{\hat{Z}_1}^{(\eta_0)}.
$$
For all $\varepsilon\geq 0$ one obtains
$$
\sigma_*(\sM_{W'}^\kappa\otimes \sO_{W'}(\varepsilon\cdot\Pi_{W'}^{(\eta_0)}))=
\widetilde{\sM}_{\hat{Z}_1}^{(\kappa)}\otimes \sO_{\hat{Z}_1}(\varepsilon\cdot\Pi_{\hat{Z}_1}^{(\eta_0)}),
$$
and 
\begin{multline}\label{eqsa.4}
\hat{g}_{1*}\sigma_*\sM_{W'}^\kappa= \hat{g}_{1*}\widetilde{\sM}_{\hat{Z}_1}^{(\kappa)}=
\hat{g}_{1*}(\widetilde{\sM}_{\hat{Z}_1}^{(\kappa)}\otimes \sO_{\hat{Z}_1}(\varepsilon\cdot\Pi_{\hat{Z}_1}^{(\eta_0)}))=\\
\hat{g}_{1*}\sigma_*(\sM_{W'}^\kappa\otimes \sO_{W'}(\varepsilon\cdot\Pi_{W'}^{(\eta_0)})). 
\end{multline}
So on $W'$ we found the sheaf we are looking for. Finally, Corollary~\ref{di.5} allows to
replace $\hat{Y}_1$ by some modification, and to assume that the condition d) in~\ref{sa.7} holds
for $(0,\kappa)$. 
\end{proof}
By~\ref{sa.11} a) one can construct $\widetilde{\sM}_{\hat{Z}}$, $\widetilde{\sM}_{Z}$ and $\widetilde{\sM}_{\hat{X}}$ such that this tuple forms an $\kappa$-saturated extension of $\sL_0$. 
Perhaps some of the sheaves $\sB_\bullet^{(\nu,\mu)}$ or the sheaves $\sB_\bullet$, depending on
$\sE_{\hat{Y}}$ in~\ref{sa.2} are no longer invertible. If so, for $\widetilde{\sM}_\bullet$ and for the given set $I$ we have to perform again the alterations needed to get the invertible sheaves in~\ref{sa.4}. Lemma~\ref{sa.11} d) allows to do so, without loosing the $\kappa$-saturatedness. 
So one is allowed to modify the condition vi) in~\ref{sa.7}, keeping all the other ones:

\begin{con-not}\label{sa.13} {\bf (Saturated polarizations)} \\
We consider an invertible sheaf $\sL$ on $X$, with $\sL_0=\sL|_{X_0}$ relative ample over $Y_0$, and 
we start again with a finite set $I$ of tuples $(\nu,\mu)$ of natural numbers.
We choose $\eta_0>0$ and $\gamma_0>0$ such that the evaluation maps
$$
f_0^*f_{0*}\omega_{X_0/Y_0}^{\eta_0} \to \omega_{X_0/Y_0}^{\eta_0} \mbox{ \ \ and \ \ }
f_0^*f_{0*}\sL_0^{\gamma_0} \to \sL_0^{\gamma_0}
$$
are surjective. 

We fix some subset $\widetilde{I}$ of $I$ consisting of tuples
$(\beta,\alpha)$ with $\alpha$ divisible by $\gamma_0$ and with $\beta$ divisible by $\eta_0$.
We also fix a positive number $\kappa$ with $(0,\kappa)\in \widetilde{I}$.

Then we can find $\hat{Y}$ and the diagram \eqref{eqws.3} (or in \eqref{eqsa.1}) such that
the conditions i) -- v) in~\ref{sa.5} hold true and
and the conditions vii) -- x) in~\ref{sa.7} with $\sM_\bullet$ given by:
\begin{enumerate}
\item[vi.] There exists a tuple of $\kappa$-saturated extensions 
$(\sM_{\hat{Z}}, \sM_Z ,\sM_{\hat{X}})$ of $\sL$.
\end{enumerate}
\end{con-not}
Remark that by Lemma~\ref{sa.11} d) the ``Allowed Constructions'' in~\ref{sa.8} remain allowed, i.e.
they respect the condition vi) in~\ref{sa.13}.
\begin{corollary}\label{sa.14}
The conditions in~\ref{sa.13} imply that for all $\varepsilon\geq 0$ the direct images  
$$
\hat{g}_*\sB^{(0,\kappa)}_{\hat{Z}}, \ \ \ \hat{g}_{*}\sM_{\hat{Z}}^{\kappa} \mbox{ \ \ and \ \ }\hat{g}_*(\sM_{\hat{Z}}^{\kappa}\otimes \sO_{\hat{Z}}(\varepsilon\cdot\Pi_{\hat{Z}}^{(\eta_0)}))
$$ 
coincide, and that they are locally free and compatible with base change for 
morphisms $\varrho:T\to \hat{Y}$ with $\varrho^{-1}(\hat{Y}_g)$ dense in $T$.
\end{corollary}
\begin{proof} By definition of ``saturated'' and by the choice of $\sB^{(0,\kappa)}_{\hat{Z}}$
$$
\hat{g}_*\sB^{(0,\kappa)}_{\hat{Z}}=\hat{g}_{*}\sM_{\hat{Z}}^{\kappa}=\hat{g}_{*}(\sM_{\hat{Z}}^{\kappa}\otimes \sO_{\hat{Z}}(\varepsilon\cdot\Pi_{\hat{Z}}^{(\eta_0)})).
$$
Since we assumed that $(0,\kappa)\in I$ the direct image
$\hat{g}_{*}\sM_{\hat{Z}}^{\kappa}$ is compatible with base change for alterations.
By Addendum~\ref{sa.3} the same holds true for $\hat{g}_*\sB^{(0,\kappa)}_{\hat{Z}}$
and by~\ref{sa.9} for $\hat{g}_{*}(\sM_{\hat{Z}}^{\kappa}\otimes \sO_{\hat{Z}}(\varepsilon\cdot\Pi_{\hat{Z}}^{(\eta_0)}))$.
So~\ref{sa.14} follows from Lemma~\ref{di.1} ii).
\end{proof}
So for $\kappa=1$ we could choose $\sM_{\hat{Z}}$ to be equal to 
$\sB^{(0,1)}_{\hat{Z}}$, but we will allow other choices. Anyway, it is easy to see that the direct image sheaves 
are independent of the choices. 

\section{The definition of certain multiplier ideals}\label{de}
The alterations, sheaves and divisors as described in the Conclusions and
Notations~\ref{sa.6},~\ref{sa.7} or~\ref{sa.13} depend on the choice of certain numbers and data.
Each time we add some numbers, we have to reemploy the constructions of Section~\ref{sa}.
In order not to run into an infinite circle of constructions we have
to give a complete list of data at some point, and this is done in the first part of this section.

However, we still have to extend the base change property stated in~\ref{sa.7} ix) to certain multiplier ideals
$$
\hat{g}_*(\omega_{\hat{Z}/\hat{Y}}^\nu\otimes\sM_{\hat{Z}}^\mu\otimes \sJ(-e\cdot D)),
$$
using Theorem~\ref{fm.5}. As in the proof of the Variant~\ref{pd1.4} the multiplier ideals we want to consider depend on the tautological map $\Xi$ and on a large number of integers. So in this Section we will include those maps in our bookkeeping. In order to get the locally freeness and the
compatibility with base change for certain morphisms, we will use again the left hand side of the diagram \eqref{eqws.3}. Then the compatibility conditions for the sheaves $\sM_\bullet$ will allow, as in Lemma~\ref{di.7}, to pass to the right hand side. 

\begin{conventions}\label{de.1}
Consider for a smooth fibre $F$ a finite tuple $\Xi$ of determinants and their natural inclusion in the tensor products, i.e. $\Xi=(\Xi_1,\ldots,\Xi_s)$ and
$$
\Xi_i: \bigwedge^{r_i}H^0(F,\omega_F^{\eta_i}\otimes \sL_0^{\gamma_i}|_F)\>>> \bigotimes^{r_i}
H^0(F,\omega_F^{\eta_i}\otimes \sL_0^{\gamma_i}|_F),
$$
where $r_i=\dim(H^0(F,\omega_F^{\eta_i}\otimes \sL_0^{\gamma_i}|_F))$. So for any $r$,
divisible by $r_1, \ldots , r_s$ and for each $i$ one obtains a map
$$
\Big(\bigwedge^{r_i}H^0(F,\omega_F^{\eta_i}\otimes \sL_0^{\gamma_i}|_F)\Big)^{\otimes \frac{r}{r_i}}
\>>> \bigotimes^{r} H^0(F,\omega_F^{\eta_i}\otimes \sL_0^{\gamma_i}|_F)
$$
and finally, for $\gamma=\gamma_1+\cdots+\gamma_s$ and for $\eta=\eta_1+\cdots+\eta_s$ one has the product
$$
\bigotimes_{i=1}^s \Big(\bigwedge^{r_i}H^0(F, \omega_F^{\eta_i}\otimes\sL_0^{\gamma_i}|_F)\Big)^{\otimes \frac{r}{r_i}}
\>\Xi^{(r)} >> \bigotimes^{r} H^0(F,\omega_F^{\eta}\otimes \sL_0^{\gamma}|_F).
$$
\end{conventions}
We will later require certain divisibilities. For example, we will need
that the integers $\eta_0$ and $\gamma_0$ in \ref{sa.5} or \ref{sa.7} divide $\eta$ and $\gamma$.
This can be achieved by replacing $\Xi$ by $s'$ copies $(\Xi,\ldots,\Xi)$
for a suitable $s'$. 

Remark that those conventions carry over to the smooth part of our families,
provided the direct image sheaves $f_{0*}(\omega_{X_0/Y_0}^{\eta_i} \otimes \sL_0^{\gamma_i})$
are all locally free and compatible with base change. This holds by~\ref{di.3} for $\eta_i>0$.
If $\eta_i=0$ we listed this already in~\ref{sa.7} as an additional condition.

Next we have to explain how to choose the finite sets $\widetilde{I}$ and $I$ in~\ref{sa.5} or~\ref{sa.7}.
\begin{setup}\label{de.2} \ \\[.1cm]
{\bf Canonically polarized case.} Here we start by choosing an integer $\eta_0>0$ 
such that the evaluation map $f_0^*f_{0*}\omega_{X_0/Y_0}^{\eta_0} \to \omega_{X_0/Y_0}^{\eta_0}$
is surjectiv and we choose $\ell>0$, divisible by $\eta_0$.  
We will assume in $\ref{de.1}$ that $\gamma_i=0$ for all $i$, and that $\eta$ is divisible by $\ell$, hence by $\eta_0$. We choose $\widetilde{I}=\{(\eta_0,0),(\eta,0)\}$. The set $I\subset \N\times \{0\}$ should contain $\widetilde{I}$, the tuples $(\eta_i,0)$ for $i=1,\ldots,s$ and for some $\beta\geq 1$
the tuple $(\beta+\frac{\eta}{\ell},0)$. For compatibility of the notations we write $\alpha=\kappa=0$ and $b$ will denote any positive integer. We choose $\hat{Y}$ and the different sheaves and divisors according to~\ref{sa.5}.\\[.1cm] 
{\bf Polarized case.} If $\sL_0$ is $f_0$-ample we start with integers $\eta_0>0$ and $\gamma_0>0$ such that the evaluation maps
$$
f_0^*f_{0*}\omega_{X_0/Y_0}^{\eta_0} \to \omega_{X_0/Y_0}^{\eta_0}\mbox{ \ \ and \ \ }
f_0^*f_{0*}\sL_0^{\gamma_0} \to \sL_0^{\gamma_0}
$$ 
are both surjective. In addition we will require that for $N\geq \gamma_0$ and for all fibres
$F$ of $f_0$ the sheaves $\sL_0^N|_F$ have no higher cohomology.
We choose $\ell>0$, divisible by $\eta_0$ and $\gamma_0$.   
In \ref{de.1}, replacing $\Xi$ by $(\Xi,..., \Xi)$, and correspondingly $s$ by some multiple, one may assume that $\ell$ divides $\gamma$ and $\eta$.
Fix in addition some tuple $(\beta,\alpha)$ of natural numbers with $\beta\geq 1$ (or a finite set of such tuples), and some positive integer $b$, with $b\cdot(\beta -1 ,\alpha)\in \eta_0\cdot\N\times \gamma_0\cdot\N$. The finite set of tuples $\widetilde{I}$ should contain 
$$
\{(\eta_0,0), \ (0,\gamma_0), \ (\eta,\gamma) \},
$$
and $I$ should contain $\widetilde{I}$,
$$
(\beta+\frac{\eta}{\ell},\alpha+\frac{\gamma}{\ell})
\mbox{ \ \ and \ \ } (\eta_i,\gamma_i) \mbox{ \ \ for \ \ } i=1,\cdots,s.
$$
For compatibility reasons we choose $\kappa=0$ in this case, and we choose $\hat{Y}$ and the different sheaves and divisors according to~\ref{sa.7}.\\[.1cm]
{\bf Saturated polarized case.} Everything is as in the polarized case, except that we also choose some
positive multiple $\kappa$ of $\gamma_0$ and assume that $(0,\kappa)\in I'$, and apply~\ref{sa.13} instead of~\ref{sa.7}.\\[.1cm]
In all three cases we fix a natural number $e$ with
\begin{equation}\label{eqde.1}
e\geq \frac{e(\omega_F^{\eta}\otimes \sL_0^{\gamma}|_F)}{\ell}
\end{equation}
for all fibers $F$ of $f_0$, where $e$ denotes the threshold introduced in~\ref{pd1.1}.
The sheaves $\sF_{\hat{Y}}^{(\eta_i,\gamma_i)}$ are locally free.
Replacing $\hat{Y}$ by a non-singular alteration one finds an invertible sheaf $\sV$ on $\hat{Y}$ with 
$$
\bigotimes_{i=1}^s \det(\sF_{\hat{Y}}^{(\eta_i,\gamma_i)})^{\frac{r}{r_i}}
=\bigotimes_{i=1}^s \det (\hat{g}_*(\omega_{\hat{Z}/\hat{Y}}^{\eta_i}\otimes\sM^{\gamma_i}_{\hat{Z}} ))^{\frac{r}{r_i}}=\sV^{r\cdot e\cdot \ell}.
$$
Remark that all assumptions remains true when we replace $r$ by some multiple or $\hat{Y}$ by an alteration.
\end{setup}
We are not yet done. We will need another auxiliary sheaf.
\begin{ass-not}\label{de.3} Consider for a locally free sheaf $\sE_{\hat{Y}}$
a morphism $\sE_{\hat{Y}}\to \sF_{\hat{Y}}^{(\beta_0,\alpha_0)}$ where
$$
\beta_0=b\cdot(\beta-1) \cdot e \cdot \ell+\eta\cdot b \cdot (e-1) \mbox{ \ and \ }
\alpha_0=b\cdot\alpha\cdot e \cdot \ell+\gamma\cdot b \cdot (e-1).
$$ 
For $b$ sufficiently large, the evaluation map 
\begin{equation}\label{eqde.2}
\hat{f}^*\sE_{\hat{Y}} \>>> \omega_{\hat{X}/\hat{Y}}^{\beta_0}\otimes 
\sM^{\alpha_0}_{\hat{X}}
\end{equation}
is surjective over $\hat{X}_0$. We will choose 
\begin{equation}\label{eqde.3}
(\beta_0,\alpha_0)=(b\cdot(\beta-1) \cdot e \cdot \ell+\eta\cdot b \cdot (e-1),
b\cdot\alpha\cdot e \cdot \ell+\gamma\cdot b \cdot (e-1)) \in \widetilde{I},
\end{equation}
where of course $\alpha_0=0$ in the canonically polarized case.

Lemma~\ref{sa.2} allows to assume (replacing $\hat{Y}$ by an alteration) that
the image of the evaluation map (\ref{eqde.2}) is an invertible sheaf $\sB_{\hat{X}}$, and that the image of
$\hat{g}^*\sE_{\hat{Y}} \>>> \omega_{\hat{Z}/\hat{Y}}^{\beta_0}\otimes 
\sM^{\alpha_0}_{\hat{Z}}$
is an invertible sheaf $\sB_{\hat{Z}}$. We write again $\Sigma_{\hat{Z}}$ for the effective divisor
with $\sB_{\hat{Z}}=\omega_{{\hat{Z}}/\hat{Y}}^{\beta_0}\otimes 
\sM^{\alpha_0}_{{\hat{Z}}}\otimes \sO_{\hat{Z}}(-\Sigma_{\hat{Z}}).$
\end{ass-not}

\begin{Variant}\label{de.4}
In the application we have in mind $\sE_{\hat{Y}}$ will be a subsheaf of
$$
\sF_{\hat{Y}}^{(\beta_1,\alpha_1)} \otimes \cdots \otimes \sF_{\hat{Y}}^{(\beta_s,\alpha_s)},
$$
with cokernel supported in $\hat{Y}\setminus \hat{Y}_0$. Here we have to assume that for all $\iota\in\{1,\ldots,s\}$ the evaluation map for $\omega_{\hat{X}/\hat{Y}}^{\beta_\iota}\otimes \sM_{\hat{X}}^{\alpha_\iota}$ is surjective over $\hat{X}_0$.

The morphism $\sE_{\hat{Y}} \to \sF_{\hat{Y}}^{(\beta_0,\alpha_0)}$
will be induced by the multiplication map
$$
\sF_{\hat{Y}}^{(\beta_1,\alpha_1)} \otimes \cdots \otimes \sF_{\hat{Y}}^{(\beta_s,\alpha_s)} \>{\rm m} >>
\sF_{\hat{Y}}^{(\beta_0,\alpha_0)}.
$$
Of course one needs that $ \beta_1 + \cdots + \beta_s = \beta_0$ and 
$\alpha_1+\cdots+\alpha_s=\alpha_0$.
In this case one can replace the condition (\ref{eqde.3}) by
\begin{equation}\label{eqde.4}
(\beta_1,\alpha_1), \ldots, (\beta_s,\alpha_s) \in \widetilde{I}.
\end{equation}
Finally remark that here $\sB_{\hat{Z}}$ is contained in the tensor product of the sheaves $\sB_{\hat{Z}}^{(\beta_\iota,\alpha_\iota)}$ that both coincide on $\hat{Z}_0$.
\end{Variant}
We need a long list of different sheaves and divisors on certain products.
\begin{notations}\label{de.5} \ \\[.1cm] 
{\bf (Saturated) polarized case.}
Let $\hat{g}:\hat{Z}\to \hat{Y}$ be the mild morphism constructed in \ref{sa.5}, \ref{sa.7} and \ref{sa.13} using the data given in \ref{de.1}--\ref{de.4}
(or by abuse of notations, its pullback under a morphism from a curve to $\hat{Y}$, assuming that it is mild). Consider the $r$-fold product 
\begin{gather*}
\hat{g}^r:\hat{Z}^r=\hat{Z}\times_{\hat{Y}}\cdots\times_{\hat{Y}}\hat{Z}\> >> \hat{Y},\mbox{ \ \ and \ \ }
\sM_{\hat{Z}^r}={\rm pr}_1^*\sM_{\hat{Z}} \otimes \cdots \otimes {\rm pr}_r^*\sM_{\hat{Z}}.
\end{gather*}
For $(\nu,\mu)\in I$ one obtains by flat base change 
\begin{equation}\label{eqde.5}
\hat{g}^r_*(\omega^{\nu}_{\hat{Z}^r/\hat{Y}}\otimes\sM_{\hat{Z}}^{\mu})=\bigotimes^r \hat{g}_*(\omega_{\hat{Z}/\hat{Y}}^{\nu}
\otimes\sM_{\hat{Z}}^\mu)= \bigotimes^r\sF_{\hat{Y}}^{(\nu,\mu)}.
\end{equation}
For $(\nu,\mu)=(\eta,\gamma)$ the equality (\ref{eqde.5}) implies that the image of the evaluation map
$$
{\hat{g}^r}{}^*\hat{g}^r_*\omega^{\eta}_{\hat{Z}^r/\hat{Y}}\otimes\sM_{\hat{Z}^r}^{\gamma} \>>> \omega^{\eta}_{\hat{Z}^r/\hat{Y}}\otimes\sM_{\hat{Z}^r}^{\gamma}
$$
is the invertible sheaf
$\sB^{(\eta,\gamma)}_{\hat{Z}^r}:=
{\rm pr}_1^*\sB^{(\eta,\gamma)}_{\hat{Z}} \otimes \cdots \otimes {\rm pr}_r^*\sB^{(\eta,\gamma)}_{\hat{Z}}$.
So the definition of $\sB^{(\eta,\gamma)}_{\hat{Z}^r}$ is compatible with the one in~\ref{sa.7}, and
$\sB^{(\eta,\gamma)}_{\hat{Z}^r}$ can be written as  
$$
\omega^{\eta}_{\hat{Z}^r/\hat{Y}}\otimes\sM_{\hat{Z}^r}^{\gamma}\otimes \sO_{\hat{Z}^r}(-\Sigma^{(\eta,\gamma)}_{\hat{Z}^r})
\mbox{ \ \ for \ \ }
\Sigma^{(\eta,\gamma)}_{\hat{Z}^r}=\sum_{i=1}^r {\rm pr}_i^* \Sigma^{(\eta,\gamma)}_{\hat{Z}}.
$$
Since $\hat{g}^r_*(\omega_{\hat{Z}^r/\hat{Y}}^\eta\otimes\sM_{\hat{Z}^r}^\gamma)=\hat{g}^r_*\sB_{\hat{Z}^r}^{(\eta,\gamma)}$,      one has an inclusion
$$
\sV^{r\cdot e\cdot\ell}=\bigotimes_{i=1}^s \det (\hat{g}_*\sM^{\gamma_i}_{\hat{Z}}\otimes \varpi_{\hat{Z}}^{\eta_i}))^{\otimes \frac{r}{r_i}}\\ 
\>\Xi^{(r)} >>
\bigotimes^r \hat{g}_*(\omega_{\hat{Z}/\hat{Y}}^\eta\otimes\sM_{\hat{Z}}^\gamma)=
\hat{g}^r_*\sB_{\hat{Z}^r}^{(\eta,\gamma)}
$$
which splits locally, hence a section of $\sB_{\hat{Z}^r}^{(\eta,\gamma)}\otimes {\hat{g}^r}{}^* \sV^{-r\cdot e\cdot\ell}$
whose zero divisor $\Gamma_{\hat{Z}^r}$ does not contain any fibre (but perhaps components of fibres). 

In~\ref{de.3} one can apply (\ref{eqde.5}) to see that 
the invertible sheaf 
$$
\sB_{\hat{Z}^r}={\rm pr}_1^*\sB_{\hat{Z}} \otimes \cdots \otimes {\rm pr}_r^*\sB_{\hat{Z}}
$$
is again the image of the evaluation map ${\hat{g}^r}{}^*\sE_{\hat{Y}}^{\otimes r} \to 
\omega_{\hat{Z}^r/\hat{Y}}^{\beta_0}\otimes \sM^{\alpha_0}_{\hat{Z}^r}$.

In Variant~\ref{de.4} the same holds true for the sheaves $\sB_{\hat{Z}^r}^{(\beta_\iota,\alpha_\iota)}$, hence for their tensor product and for the image $\sB_{\hat{Z}^r}$ of ${\hat{g}^r}{}^*\sE_{\hat{Y}}^{\otimes r}$.
In both cases one finds
$$
\sB_{\hat{Z}^r}=\omega_{{\hat{Z}^r}/\hat{Y}}^{\beta_0}\otimes 
\sM^{\alpha_0}_{{\hat{Z}^r}}\otimes \sO_{\hat{Z}^r}(-\Sigma_{\hat{Z}^r}), \mbox{ \ \ for \ \ }
\Sigma_{\hat{Z}^r}=\sum_{i=1}^r {\rm pr}_i^* \Sigma_{\hat{Z}}.
$$
To shorten the expressions, we put
\begin{gather*}
\Delta_{\hat{Z}^r}=b \cdot (\Gamma_{\hat{Z}^r}+\Sigma_{\hat{Z}^r}^{( \eta,\gamma)}) + \Sigma_{\hat{Z}^r}, \ \ \ \ N=b\cdot e\cdot\ell\\
\mbox{and \ \ \ }
\sG_{\hat{Y}}^{(\Xi^{(r)},\sE;\beta+\frac{\eta}{\ell},\alpha+\frac{\gamma}{\ell})}=\hat{g}^r_*\Big(\omega_{\hat{Z}^r/\hat{Y}}^{\beta+\frac{\eta}{\ell}} \otimes\sM_{\hat{Z}^r}^{\alpha+\frac{\gamma}{\ell}}\otimes\sJ\big(-\frac{1}{N}\cdot \Delta_{\hat{Z}^r}\big)\Big).
\end{gather*}
We will usually write $\sG_{\hat{Y}}^{(\beta+\frac{\eta}{\ell},\alpha+\frac{\gamma}{\ell})}$ instead of
$\sG_{\hat{Y}}^{(\Xi^{(r)},\sE;\beta+\frac{\eta}{\ell},\alpha+\frac{\gamma}{\ell})}$, keeping however in mind that this sheaf depends on the choice of $r$, of the tautological maps $\Xi_i$ and on $\sE$.\\[.1cm]
{\bf Canonically polarized case.} We will use the same notation. This is a bit queer, but it allows to
handle both cases in the same way. So in this case $\alpha=\gamma=0$ and $\sM_\bullet=\sO_\bullet$.
\end{notations}
\begin{lemma}\label{de.6}
Under the assumptions made in~\ref{de.1}--\ref{de.4} one may choose $\hat{Y}$ and $\hat{Z}$ in~\ref{sa.5}, \ref{sa.7} or~\ref{sa.13} and an open dense subscheme $\hat{Y}_g\subset \hat{Y}_0$ such that in addition to the conditions i)--v) in \ref{sa.5} or to i)--x) in \ref{sa.7} or \ref{sa.13} on has:
\begin{enumerate}
\item[xi.] The multiplier ideal sheaves $\sJ\big(-\frac{1}{b\cdot e\cdot \ell}\cdot \Delta_{\hat{Z}^r}\big)$ 
are compatible with pullback, base change and products with respect to $\hat{Y}_g$, as defined in~\ref{fm.4}.
In particular they are flat over $\hat{Y}$ and the direct image sheaves
$\sG_{\hat{Y}}^{(\beta+\frac{\eta}{\ell},\alpha+\frac{\gamma}{\ell})}$
are compatible with pullback for morphisms $\varrho:T\to \hat{Y}$ where $\varrho$ is either dominant and $T$ a normal variety with at most rational Gorenstein singularities, or where $T$ is a non-singular curve and $\varrho^{-1}(\hat{Y}_g)$ dense in $T$. Moreover for $r'>0$ 
$$
{\sG_{\hat{Y}}^{(\Xi^{(r)},\sE;\beta+\frac{\eta}{\ell},\alpha+\frac{\gamma}{\ell})}}^{\otimes r'}=
\sG_{\hat{Y}}^{(\Xi^{(r\cdot r')},\sE;\beta+\frac{\eta}{\ell},\alpha+\frac{\gamma}{\ell})}.
$$
\end{enumerate}
\end{lemma}
\begin{proof}
Choose $\sN=\omega_{\hat{Z}^r/\hat{Y}}^{\beta-1+\frac{\eta}{\ell}} \otimes\sM_{\hat{Z}^r}^{\alpha+\frac{\gamma}{\ell}}$.
Then $\sN^N\otimes\sO_{\hat{Z}^r}(-\Delta_{\hat{Z}^r})$ is equal to 
$$
\big[\omega_{\hat{Z}^r/\hat{Y}}^{\beta_0} \otimes\sM_{\hat{Z}^r}^{\alpha_0}\otimes\sO_{\hat{Z}^r}(-\Sigma_{\hat{Z}^r})\big]\otimes
\big[\omega_{\hat{Z}^r/\hat{Y}}^{\eta\cdot b} \otimes\sM_{\hat{Z}^r}^{\gamma\cdot b}\otimes\sO_{\hat{Z}^r}(-b(\Sigma^{(\eta,\gamma)}_{\hat{Z}^r}
+\Gamma_{\hat{Z}^r}))\big],
$$
where the first factor is the image of ${\hat{g}^r}{}^*\sE_{\hat{Y}}^{\otimes r}$ whereas the second one is the $b$-th power of 
$\sB_{\hat{Z}^r}^{(\eta,\gamma)}\otimes\sO_{\hat{Z}^r}(-\Gamma_{\hat{Z}^r})={\hat{g}^r}{}^* \sV^{r\cdot e\cdot\ell}$. So we obtain:
\begin{claim}\label{de.7}
For $\sN$, for $\Delta=\Delta_{\hat{Z}^r}$, and for $\sE=\sE_{\hat{Y}}^{\otimes r}\otimes \sV^{b\cdot r\cdot e\cdot\ell}$ the assumptions made in~\ref{fm.1} hold true (for $Z$ replaced by $\hat{Z}^r$). 
\end{claim}
So we are allowed to apply Theorem~\ref{fm.5}. Dropping the index ${}_1$, assume that $\hat{Y}=\hat{Y}_1$, hence that
$\sJ\big(-\frac{1}{N}\cdot \Delta_{\hat{Z}^r}\big)$ is compatible with pullback, base change and products with respect to $\hat{Y}_g$.

For $\sA=\sO_{\hat{Z}^r}$ in Definition~\ref{fm.4} the properties i) and  ii)
give the compatibility with pullback under $\varrho$, and by flat base change also the compatibility with products.
\end{proof}

Before proving an analog of Lemma~\ref{di.7} for the sheaves $\sG_{\hat{Y}}^{(\beta+\frac{\eta}{\ell},\alpha+\frac{\gamma}{\ell})}$
we have to extend the definition of the sheaves and divisors to desingularizations of compactifications 
of $\hat{X}^r_0\to \hat{Y}_0$ (or again of the pullback of this morphism to a curve, meeting $\hat{Y}_0$).
\begin{notations}\label{de.8}
Consider the $r$-fold product
$
\hat{f}^r: \hat{X}^r=\hat{X}\times_{\hat{Y}}\cdots\times_{\hat{Y}}\hat{X}\to \hat{Y}$. The morphism
$\rho':X^{(r)}\to \hat{X}^r$ is obtained by desingularizing the main component of $\hat{X}^r$.
By~\ref{ws.2} the morphism $\hat{g}^r:\hat{Z}^r\to \hat{Y}$ in~\ref{de.5} and~\ref{de.6} is again mild, hence it is a mild model of the induced morphism
$f^{(r)}:X^{(r)}\to \hat{Y}$. Let us write 
$$
\sM_{X^{(r)}}=\rho'^*({\rm pr}_1^*\sM_{\hat{X}} \otimes \cdots \otimes {\rm pr}_r^*\sM_{\hat{X}}).
$$
Recall that for $\nu$ divisible by $\eta_0$ and for $\mu$ divisible by $\gamma_0$ the evaluation map
$$
{\hat{f}^r_0}{}^* 
\hat{f}^r_{0*}(\omega^{\nu}_{\hat{X}^r_0/\hat{Y}_0}\otimes \sM_{\hat{X}^r_0}^\mu) \>>> 
\omega^{\nu}_{\hat{X}^r_0/\hat{Y}_0}\otimes \sM_{\hat{X}^r_0}^\mu
$$
is surjective, where again the index ${}_0$ refers to the preimages of $Y_0$ or for sheaves
to their restriction. 
Consider a non-singular modification $\hat{\delta}{}^r:Z^{(r)}\to \hat{Z}^r$ which allows
a morphism $\delta^{(r)}:Z^{(r)}\to X^{(r)}$, and which dominates the main component
of $Z\times_{\hat{Y}}\cdots \times_{\hat{Y}}Z$. Defining $\sM_{Z^{(r)}}$ as the pullback of
${\rm pr}_1^*\sM_{Z} \otimes \cdots \otimes {\rm pr}_r^*\sM_{Z}$,
one has ${\hat{\delta}{}^r}^*\sM_{\hat{Z}^r}\subset \sM_{Z^{(r)}}$ and ${\delta^{(r)}}^*\sM_{X^{(r)}}\subset \sM_{Z^{(r)}}$.
\end{notations}
\begin{lemma}\label{de.9}
The sheaves $\sM_{Z^{(r)}}$, $\sM_{\hat{Z}^r}$ and $\sM_{X^{(r)}}$ satisfy again the
Assumptions asked for in~\ref{di.4}.
\end{lemma}
\begin{proof}
Since $\hat{Z}^r$ is normal the assumption $\hat{\delta}_*\sM_Z=\sM_{\hat{Z}}$ in~\ref{di.4} implies that
$$
\hat{\delta}{}^r_*\sM_{Z^{(r)}}= \sM_{\hat{Z}^r}.
$$
For $\sM_{X^{(r)}}$ remark first, that 
${\delta}^*\sM_{\hat{X}} \otimes \sO_Z(F)=\sM_{Z}$, for some $\delta$-exceptional effective divisor $F$.
Consider the diagram 
$$
\xymatrix{Z^r\times_{\hat{X}^r}X^{(r)} \ar[r]^{\hspace{.6cm}\theta}\ar[d]_{p_1}& X^{(r)}\ar[d]^{\rho}\\
Z^r \ar[r]^{\delta^r}& X^r.}
$$ 
Then ${\delta^r}^*({\rm pr}_1^*\sM_{\hat{X}} \otimes \cdots \otimes {\rm pr}_r^*\sM_{\hat{X}})$ 
is a subsheaf of ${\rm pr}_1^*\sM_{Z} \otimes \cdots \otimes {\rm pr}_r^*\sM_{Z}$ and both
coincide outside of a divisor $F'$ with ${\rm codim}(\hat{\delta}{}^r(F'))\geq 2$. So the same holds true
for the subsheaf
$$
p_1^*{\delta^r}^*({\rm pr}_1^*\sM_{\hat{X}} \otimes \cdots \otimes {\rm pr}_r^*\sM_{\hat{X}})=
\theta^*\sM_{X^{(r)}}
$$ 
of $p_1^*({\rm pr}_1^*\sM_{Z} \otimes \cdots \otimes {\rm pr}_r^*\sM_{Z}).$
The statement is independent of the desingularization. Hence we may assume that 
$Z^{(r)}$ dominates the main component of $Z^r\times_{\hat{X}^r}X^{(r)}$. So
${\delta^{(r)}}^*\sM_{X^{(r)}} \otimes \sO_{Z^{(r)}}(F'')=\sM_{Z^{(r)}}$ for some effective ${\delta^{(r)}}$
exceptional divisor $F''$.  
\end{proof}
Lemma~\ref{de.9} allows to apply Lemma~\ref{di.7} and
\begin{equation}\label{eqde.6}
f^{(r)}_{*}(\omega^{\nu}_{X^{(r)}/\hat{Y}}\otimes \sM_{X^{(r)}}^\mu) =
\hat{g}^r_*(\omega^{\nu}_{\hat{Z}^r/\hat{Y}}\otimes \sM_{\hat{Z}^r}^\mu).
\end{equation}
For $(\nu,\mu)\in I$ one can use flat base change and the projection formula to identify
the right hand side as 
$$\bigotimes^r \hat{g}_*(\omega^{\nu}_{\hat{Z}/\hat{Y}}\otimes \sM_{\hat{Z}}^\mu).$$
Using~\ref{di.7} again, one finds
$$
f^{(r)}_{*}(\omega^{\nu}_{X^{(r)}/\hat{Y}}\otimes \sM_{X^{(r)}}^\mu) =
\bigotimes^r \hat{f}_*(\omega^{\nu}_{\hat{X}/\hat{Y}}\otimes \sM_{\hat{X}}^\mu).
$$
In particular those sheaves are locally free and compatible with base change for morphisms $\varrho:T\to \hat{Y}$ with $\varrho^{-1}(\hat{Y}_g)$ dense in $T$.\vspace{.1cm} 

Next we move to the right hand side of the diagram \eqref{eqws.3} and redefine all the sheaves and divisors 
from~\ref{de.5} with $\hat{Z}^r$ replaced by $\hat{X}^{(r)}$ or $Z^{(r)}$. As in~\ref{de.5} we will give the
definitions in the polarized case. For the canonically polarized case the last lines of~\ref{de.5}
apply.
\begin{notations}\label{de.10} As at the end of the proof of Lemma~\ref{sa.2}, blowing up $X^{(r)}$ with centers outside $\hat{X}^r_0$, one may assume that the image of
$$
{f^{(r)}}^* f^{(r)}_{*}(\omega^{\eta}_{X^{(r)}/\hat{Y}}\otimes \sM_{X^{(r)}}^\gamma) \>>> 
\omega^{\eta}_{X^{(r)}/\hat{Y}}\otimes \sM_{X^{(r)}}^\gamma
$$
is invertible and we denote it by $\sB^{(\eta,\gamma)}_{X^{(r)}}$. 
The effective divisor $\Sigma^{(\eta,\gamma)}_{X^{(r)}}$ is chosen such that 
$$
\sB^{(\eta,\gamma)}_{X^{(r)}}\otimes\sO_{X^{(r)}}(\Sigma^{(\eta,\gamma)}_{X^{(r)}})=
\omega^{\eta}_{X^{(r)}/\hat{Y}}\otimes \sM_{X^{(r)}}^\gamma,
$$
hence supported outside of $X_0^{(r)}$. If the condition (\ref{eqde.3}) holds, we can apply
(\ref{de.8}) for the tuple 
$$
(\beta_0,\alpha_0)=(b\cdot(\beta-1) \cdot e \cdot \ell+\eta\cdot b \cdot (e-1),
b\cdot\alpha\cdot e \cdot \ell+\gamma\cdot b \cdot (e-1)) 
$$ 
and obtain an inclusion $\sE_{\hat{Y}}^r \to f^{(r)}_{*}(\omega^{\beta_0}_{X^{(r)}/\hat{Y}}\otimes 
\sM_{X^{(r)}}^{\alpha_0})$. The image of ${f^{(r)}}^*\sE_{\hat{Y}}^r$ under the evaluation map will be denoted by $\sB_{X^{(r)}}$.

In Variant~\ref{de.4}, i.e. if (\ref{eqde.4}) holds, one applies (\ref{eqde.6})
for the tuples $(\beta_\iota,\alpha_\iota)$. So one has Morphisms
$$
\bigotimes_{\iota=1}^s
(\hat{g}_*\sB_{\hat{Z}}^{(\beta_\iota,\alpha_\iota)}) \>>> f^{(r)}_{*}(\omega^{\beta_\iota}_{X^{(r)}/\hat{Y}}\otimes 
\sM_{X^{(r)}}^{\alpha_\iota}).
$$
The image of ${f^{(r)}}^*(\hat{g}_*\sB_{\hat{Z}}^{(\beta_\iota,\alpha_\iota)})^{\otimes r}$ is
an invertible sheaf $\sB_{X^{(r)}}^{(\beta_\iota,\alpha_\iota)}$, and the image of
${f^{(r)}}^*\displaystyle\bigotimes_{\iota=1}^s(\hat{g}_*\sB_{\hat{Z}}^{(\beta_\iota,\alpha_\iota)})^{\otimes r}$
under the product map is 
$$
\bigotimes_{\iota=1}^s\sB_{X^{(r)}}^{(\beta_\iota,\alpha_\iota)}
\subset \omega^{\beta_0}_{X^{(r)}/\hat{Y}}\otimes 
\sM_{X^{(r)}}^{\alpha_0}.
$$
So the image of ${f^{(r)}}^*\sE_{\hat{Y}}$ is a subsheaf $\sB_{X^{(r)}}$. 

In both cases $\sB_{X^{(r)}}$ is isomorphic to 
$\omega^{\beta_0}_{X^{(r)}/\hat{Y}}\otimes \sM_{X^{(r)}}^{\alpha_0}$
on $\hat{X}^r_0={f^{(r)}}^{-1}(\hat{Y}_0)$. Blowing up $X^{(r)}$ we find a divisor $\Sigma_{X^{(r)}}$ with
$\omega^{\beta_0}_{X^{(r)}/\hat{Y}}\otimes 
\sM_{X^{(r)}}^{\alpha_0}=\sB_{X^{(r)}}\otimes \sO_{X^{(r)}}(\Sigma_{X^{(r)}})$.

Finally the equation (\ref{eqde.6}) implies that
$$
f^{(r)}_{*}\sB^{(\eta,\gamma)}_{X^{(r)}}=
f^{(r)}_{*}(\omega^{\eta}_{X^{(r)}/\hat{Y}}\otimes \sM_{X^{(r)}}^\gamma)=
\hat{g}^r_*(\omega^{\eta}_{\hat{Z}^r/\hat{Y}}\otimes \sM_{\hat{Z}^r}^\gamma).
$$
Hence $\Xi^{(r)}:\sV^{r\cdot e\cdot\ell}\to  \hat{g}^r_*(\omega^{\eta}_{\hat{Z}^r/\hat{Y}}\otimes \sM_{\hat{Z}^r}^\gamma)$ induces 
a section of $\sB^{(\eta,\gamma)}_{X^{(r)}}\otimes {f^{(r)}}^{*}\sV^{-r\cdot e\cdot\ell}$ whose zero divisor
will be denoted by $\Gamma_{X^{(r)}}$. We write again
$$
\Delta_{X^{(r)}}=b \cdot (\Gamma_{X^{(r)}}+\Sigma_{X^{(r)}}^{( \eta,\gamma)}) + \Sigma_{X^{(r)}},
$$
and recall that
$X^{(r)}_0=\hat{X}^r_0$, $\Sigma_{X^{(r)}_0}^{(\beta_0,\alpha_0)}=\Sigma_{X^{(r)}_0}^{(\eta,\gamma)}=0$  and ${\hat{\delta}{}^r}^* \Gamma_{\hat{Z}^r}={\delta^{(r)}}^*\Gamma_{X^{(r)}}$.
\end{notations}

\begin{lemma}\label{de.11}
The sheaf $\sG_{\hat{Y}}^{(\beta+\frac{\eta}{\ell},\alpha+\frac{\gamma}{\ell})}$ in~\ref{de.6} is equal to
$$
f^{(r)}_*\Big(\omega_{X^{(r)}/\hat{Y}}^{\beta+\frac{\eta}{\ell}}\otimes\sM_{X^{(r)}}^{\alpha+\frac{\gamma}{\ell}} \otimes 
\sJ\big(-\frac{1}{N}\cdot \Delta_{X^{(r)}}\big)\Big).
$$
On $\hat{X}^r_0={f^{(r)}}^{-1}(\hat{Y}_0)$ one has
$$
\sJ\big(-\frac{1}{N}\cdot \Delta_{X^{(r)}}\big)|_{\hat{X}^r_0}=\sJ\big(-\frac{1}{e\cdot\ell}\cdot\Gamma_{\hat{X}^r_0}\big)=\sO_{\hat{X}^r_0},
$$ 
and the inclusion $\displaystyle \sG_{\hat{Y}}^{(\beta+\frac{\eta}{\ell},\alpha+\frac{\gamma}{\ell})}\to \bigotimes^r\sF_{\hat{Y}}^{(\beta+\frac{\eta}{\ell},\alpha+\frac{\gamma}{\ell})}$
is an isomorphism on $\hat{Y}_0$. 
\end{lemma}
\begin{proof}
We keep the notations from~\ref{de.8} and assume in addition that the pullbacks of
$\Delta_{\hat{Z}^r}$ and of $\Delta_{X^{(r)}}$ to $Z^{(r)}$ are normal crossing divisors.

Since $\hat{\delta}{}^r_*\omega_{Z^{(r)}/\hat{Y}}=\omega_{\hat{Z}^r/\hat{Y}}$ and
$\delta^{(r)}_*\omega_{Z^{(r)}/\hat{Y}}=\omega_{X^{(r)}/\hat{Y}}$, and since by Lemma~\ref{de.9} 
the same holds for the sheaves $\sM_\bullet$ 
one can find for all $(\nu,\mu)$ effective $\hat{\delta}{}^r$-exceptional divisors $E_{Z^{(r)}/\hat{Z}^r}$ and $F_{Z^{(r)}/\hat{Z}^r}$
and $\delta^{(r)}$-exceptional divisors $E_{Z^{(r)}/X^{(r)}}$ and $F_{Z^{(r)}/X^{(r)}}$
with
\begin{gather*}
\omega_{Z^{(r)}/\hat{Y}}^\nu\otimes \sM_{Z^{(r)}}^\mu
={\hat{\delta}{}^r}^*(\omega_{\hat{Z}^r/\hat{Y}}^\nu\otimes
\sM_{\hat{Z}^r}^\mu)\otimes\sO_{Z^{(r)}}(\nu\cdot E_{Z^{(r)}/\hat{Z}^r}+\mu\cdot F_{Z^{(r)}/\hat{Z}^r})
\hspace{1cm}\\ \hspace{2.5cm}
={\delta^{(r)}}^*(\omega_{X^{(r)}/\hat{Y}}^\nu\otimes
\sM_{X^{(r)}}^\mu)\otimes\sO_{Z^{(r)}}(\nu\cdot E_{Z^{(r)}/X^{(r)}}+\mu\cdot F_{Z^{(r)}/X^{(r)}}).
\end{gather*}
By Lemma~\ref{sa.2} one has  
${\hat{\delta}{}^r}^*\sB^{(\eta,\gamma)}_{\hat{Z}^r}={\delta^{(r)}}^*\sB^{(\eta,\gamma)}_{X^{(r)}}$ and
${\hat{\delta}{}^r}^*\sB_{\hat{Z}^r}={\delta^{(r)}}^*\sB_{X^{(r)}}$.
This implies that  
$$
{\hat{\delta}{}^r}^*\Sigma^{(\eta,\gamma)}_{\hat{Z}^r}+\eta\cdot E_{Z^{(r)}/\hat{Z}^r}+\gamma\cdot F_{Z^{(r)}/\hat{Z}^r}=
{\delta^{(r)}}^*\Sigma^{(\eta,\gamma)}_{X^{(r)}}+\eta\cdot E_{Z^{(r)}/X^{(r)}}+\gamma\cdot F_{Z^{(r)}/X^{(r)}},
$$
and that
$$
{\hat{\delta}{}^r}^*\Sigma_{\hat{Z}^r}+\beta_0\cdot E_{Z^{(r)}/\hat{Z}^r}+\alpha_0\cdot F_{Z^{(r)}/\hat{Z}^r}=
{\delta^{(r)}}^*\Sigma_{X^{(r)}}+\beta_0\cdot E_{Z^{(r)}/X^{(r)}}+\alpha_0\cdot F_{Z^{(r)}/X^{(r)}}.
$$
Moreover ${\hat{\delta}{}^r}^*\Gamma_{\hat{Z}^r}={\delta^{(r)}}^*\Gamma_{X^{(r)}}$, and putting everything together one finds
\begin{multline*}
{\hat{\delta}{}^r}^*\Delta_{\hat{Z}^r}+(b\cdot(\beta-1)\cdot e \cdot \ell+\eta\cdot b \cdot e)\cdot E_{Z^{(r)}/\hat{Z}^r}+(b\cdot\alpha\cdot e \cdot \ell+\gamma\cdot b \cdot e)\cdot F_{Z^{(r)}/\hat{Z}^r}=\\
{\delta^{(r)}}^*\Delta_{X^{(r)}}+(b\cdot(\beta-1)\cdot e \cdot \ell+\eta\cdot b \cdot e)\cdot E_{Z^{(r)}/X^{(r)}}+(b\cdot\alpha\cdot e \cdot \ell+\gamma\cdot b \cdot e)\cdot F_{Z^{(r)}/X^{(r)}}
\end{multline*}
and  
\begin{multline*}
{\hat{\delta}{}^r}^*(\omega_{\hat{Z}^r/\hat{Y}}^{\beta+\frac{\eta}{\ell}-1}\otimes\sM_{\hat{Z}^r}^{\alpha+\frac{\gamma}{\ell}})
\otimes \sO_{Z^{(r)}} \big(-\big[\frac{1}{N}\cdot{\hat{\delta}{}^r}^*\Delta_{\hat{Z}^r}\big]\big)=\\
{\delta^{(r)}}^*(\omega_{X^{(r)}/\hat{Y}}^{\beta+\frac{\eta}{\ell}-1}\otimes\sM_{X^{(r)}}^{\alpha+\frac{\gamma}{\ell}})\otimes
\sO_{Z^{(r)}}\big(-\big[\frac{1}{N}\cdot{\delta^{(r)}}^*
\Delta_{X^{(r)}}\big]\big).
\end{multline*}
By the definition of multiplier ideals  this implies
\begin{multline*}
\sG_{\hat{Y}}^{(\beta+\frac{\eta}{\ell},\alpha+\frac{\gamma}{\ell})}= \hat{g}^r_*\hat{\delta}{}^r_*\Big(\omega_{Z^{(r)}/\hat{Y}}\otimes 
{\hat{\delta}{}^r}^*(\omega_{\hat{Z}^r/\hat{Y}}^{\beta+\frac{\eta}{\ell}-1}\otimes\sM_{\hat{Z}^r}^{\alpha+\frac{\gamma}{\ell}})
\otimes \sO_{Z^{(r)}} \big(-\big[\frac{1}{N}\cdot{\hat{\delta}{}^r}^*\Delta_{\hat{Z}^r}\big]\big)\Big)\\
= f^{(r)}_*\delta^{(r)}_*\Big(\omega_{Z^{(r)}/\hat{Y}}\otimes 
{\delta^{(r)}}^*(\omega_{X^{(r)}/\hat{Y}}^{\beta+\frac{\eta}{\ell}-1}\otimes\sM_{X^{(r)}}^{\alpha+\frac{\gamma}{\ell}})
\otimes \sO_{Z^{(r)}} \big(-\big[\frac{1}{N}\cdot{\delta^{(r)}}^*\Delta_{X^{(r)}}\big]\big)\Big)=\\
f^{(r)}_*\Big(\omega_{X^{(r)}/\hat{Y}}^{\beta+\frac{\eta}{\ell}}\otimes\sM_{X^{(r)}}^{\alpha+\frac{\gamma}{\ell}} \otimes 
\sJ\big(-\frac{1}{N}\cdot \Delta_{X^{(r)}}\big)\Big) .
\end{multline*}
as claimed in~\ref{de.11}. In particular one has a natural inclusion 
$$
\sG_{\hat{Y}}^{(\beta+\frac{\eta}{\ell},\alpha+\frac{\gamma}{\ell})}\to f^{(r)}_*(\omega_{X^{(r)}/\hat{Y}}^{\beta+\frac{\eta}{\ell}}\otimes\sM_{X^{(r)}}^{\alpha+\frac{\gamma}{\ell}})= \bigotimes^r\sF_{\hat{Y}}^{(\beta+\frac{\eta}{\ell},\alpha+\frac{\gamma}{\ell})},
$$
induced by $\sJ\big(-\frac{1}{N}\cdot\Delta_{X^{(r)}}\big)\subset \sO_{X^{(r)}}$.
It remains to show that the latter is an isomorphism over $X^{(r)}_0=\hat{X}^r_0$.

Since $\Sigma_{X^{(r)}}|_{X_0^{(r)}}=\Sigma^{(\eta,\gamma)}_{X^{(r)}}|_{X_0^{(r)}}=0$, 
$$
\sG_{\hat{Y}}^{(\beta+\frac{\eta}{\ell},\alpha+\frac{\gamma}{\ell})}|_{\hat{Y}_0}=f^{(r)}_{0*}\Big(\omega_{X^{(r)}_0/\hat{Y}_0}^{\beta+\frac{\eta}{\ell}} \otimes\sM_{X^{(r)}_0}^{\alpha+\frac{\gamma}{\ell}}\otimes \sJ\big(-\frac{1}{e\cdot\ell}\cdot\Gamma_{X^{(r)}_0}\big)\Big).
$$
By definition, $X^{(r)}_0=\hat{X}^r_0$ and by \cite[Proposition 5.19]{Vie}
$$
e(\Gamma_{\hat{X}^r_0}) \leq {\rm Max}\big\{ e(\omega_{F^r}^\eta\otimes \sM_{X^{(r)}}^\gamma|_{F^r}); \ F \mbox{ a fibre of } \hat{f}_0\big\}
$$
By \cite[Corollary 5.21]{Vie} the right hand side is equal to
$$
{\rm Max}\big\{ e(\omega_{F}^\eta\otimes \sM_{\hat{X}}^\gamma|_{F}); \ F \mbox{ a fibre of } \hat{f}_0\big\}.
$$
So the choice of $e$ in~\eqref{eqde.1} implies that $\sJ\big(-\frac{1}{e\cdot\ell}\cdot\Gamma_{\hat{X}^r_0}\big)=\sO_{\hat{X}^r_0}$. 
\end{proof}
\begin{remark}\label{de.12}
Replacing $e$ by some larger number one can force the multiplier ideal $\sJ\big(-\frac{1}{b\cdot e \cdot \ell}\cdot\Delta_{\hat{Z}^r}\big)$ to be equal to $\sO_{\hat{Z}^r}$ and 
$$
\sG_{\hat{Y}}^{(\beta+\frac{\eta}{\ell},\alpha+\frac{\gamma}{\ell})}\hookrightarrow \bigotimes^r\sF_{\hat{Y}}^{(\beta+\frac{\eta}{\ell},\alpha+\frac{\gamma}{\ell})}
$$ 
in~\ref{de.11} to be an isomorphism on $\hat{Y}$. However, changing $e$ one looses the compatibility of the multiplier ideals with pullbacks and, as remarked already in~\ref{fm.12}, one can not expect the same $e$ to work over the alterations needed to enforce this condition.
\end{remark}

\section{Mild reduction over curves}\label{mi}

The sheaves $\sF_{\hat{Y}}^{(\nu,\mu)}$ and $\sG^{(\beta+\frac{\eta}{\ell},\alpha+\frac{\gamma}{\ell})}=\sG_{\hat{Y}}^{(\Xi^{(r)},\sE;\beta+\frac{\eta}{\ell},\alpha+\frac{\gamma}{\ell})}$ are only compatible with base change for dominant morphisms, and for morphisms from curves whose image meets a certain open subscheme $\hat{Y}_g$ of $\hat{Y}_0$. We will extend the latter in Proposition~\ref{mi.5} to morphisms whose image meets $\hat{Y}_0$.
We will need in addition that the pullback family over $C$ has a semistable or mild model, as it will be defined in this section. 

First we consider the sheaves $\sF_{\hat{Y}}^{(\nu,\mu)}$. The necessary changes for $\sG_{\hat{Y}}^{(\beta+\frac{\eta}{\ell},\alpha+\frac{\gamma}{\ell})}$ will be discussed in
the next Section. For the canonically polarized case the last lines of~\ref{de.5} apply. One just has to
choose $\sM_\bullet=\sO_\bullet$ and choose $\alpha=\gamma=\kappa=0$.

We will need that the sheaves $\sM_{\bullet}$ are also well defined
for the restrictions of our families to curves. This is evidently true
for the dualizing sheaves, and for the pullback of the invertible sheaf $\sL$ on $X$.
For the saturated extensions of the polarization, we will need some additional
arguments. So at some points we will handle the two cases separately.

We keep in this section the setup and the assumptions from \ref{de.1}--\ref{de.4}, and we choose 
the morphisms in the diagram \eqref{eqws.3} according to \ref{sa.5}, \ref{sa.7} or \ref{sa.13}.
\begin{definition}\label{mi.1}
Consider a non-singular curve $\hat{C}$, an open dense subscheme $\hat{C}_0$ and a morphism
$\chi': \hat{C} \to Y$ with $\chi'(\hat{C}_0)\subset Y_0$.

We say that $\chi:\hat{C}\to \hat{Y}$ has a {\em mild reduction}, if there exists a commutative diagram
\begin{equation}\label{eqmi.1}
\xymatrix{\hat{S} \ar[r]^{\zeta \ \ \ } \ar[d]_{ } & X\times_Y\hat{C} \ar[dl]^{{\rm pr}_2}\\
\hat{C}}
\end{equation}
of morphisms of normal projective varieties with
\begin{enumerate}
\item[i.] $\hat{h}$ is mild.
\item[ii.] $\zeta:\hat{S} \to X\times_Y\hat{C}$ is a modification of $X\times_Y\hat{C}$.
\end{enumerate}
{\bf The canonically polarized case.}\\
We call $\hat{h}:\hat{S}\to \hat{C}$ a {\em mild reduction of $\chi':\hat{C}\to Y$}.\\[.1cm]
{\bf The polarized case.}
We call $(\hat{h}:\hat{S}\to \hat{C}, \sM_{\hat{S}})$ a {\em mild reduction of $\chi':\hat{C}\to Y$ (for $\sL$)}, 
if in addition to i) and ii) one has
\begin{enumerate}
\item[iii.] $\sM_{\hat{S}}=\zeta^*{\rm pr}_1^* \sL$.
\end{enumerate}
\end{definition}
It is easy to find a mild reduction over $\hat{C}$ whenever $\hat{C} \to \chi'(\hat{C})$ is sufficiently
ramified. As in Section~\ref{ec} one can desingularize $X\times_Y \hat{C}$ such that all the fibres
become normal crossing divisors, and then one can replace $\hat{C}$ by a larger covering, to
get rid of multiple fibre components.

In the saturated case we have to be more careful. We can not choose $\sM_{\hat{S}}$ as the pullback, since we do not
want to require the existence of a morphism from $\hat{S}$ to $\hat{X}$. 
\begin{definition}\label{mi.2}\ \\
{\bf The saturated case.} 
We call in \ref{mi.1} $(\hat{h}:\hat{S}\to \hat{C}, \sM_{\hat{S}})$ a {\em mild reduction of $\chi':\hat{C}\to Y$ (for $\sL$ 
or for $\sL$ and $\eta_0$)}, if in addition to i) and ii) in~\ref{mi.1} one has:
\begin{enumerate}
\item[iii.] There exists a Cartier divisor $\Pi^{(\eta_0)}_{\hat{S}}$ on $\hat{S}$ with
$$
\hat{h}^*\hat{h}_*\omega_{\hat{S}/\hat{C}}^{\eta_0} \>>> \varpi_{\hat{S}}^{(\eta_0)}=\omega_{\hat{S}/\hat{C}}^{\eta_0}\otimes \sO_{\hat{S}}(-\Pi_{\hat{S}}^{(\eta_0)})
$$
surjective. Moreover $\sM_{\hat{S}}$ is a $\kappa$-saturated extension of $\zeta^*{\rm pr}_1^* \sL$, i.e.
it satisfies the condition required for $\sM_{\hat{Z}}$ in~\ref{sa.10}:
$$
\zeta^*{\rm pr}_1^*\sL \subset
\sM_{\hat{S}}\subset \zeta^*{\rm pr}_1^*\sL\otimes ( \sO_{\hat{S}}(*\cdot\Pi_{\hat{S}}^{(\eta_0)})\cap
\sO_{\hat{S}}(*\hat{h}^{-1}(\hat{C}\setminus \chi'^{-1}(Y_0))),
$$
and $\hat{h}_{*}\sM_{\hat{S}}^\kappa = \hat{h}_{*}(\sM_{\hat{S}}^\kappa\otimes \sO_{\hat{S}}(\varepsilon\cdot\Pi^{(\eta_0)}_{\hat{S}}))
$ for all $\varepsilon\geq 0$. 
\end{enumerate}
\end{definition}
In all cases, if $(\hat{h}:\hat{S}\to \hat{C}, \sM_{\hat{S}})$ is a mild reduction of $\chi':\hat{C}\to Y$ 
for $\sL$, we define $\sF^{(\nu,\mu)}_{\hat{C}}=\hat{h}_*(\omega_{\hat{S}/\hat{C}}^\nu\otimes \sM_{\hat{S}}^\mu)$.
We will need the compatibility of this sheaf with pullback:
\begin{lemma}\label{mi.3} Let $(\hat{h}:\hat{S}\to \hat{C}, \sM_{\hat{S}})$ be a mild reduction
for $\chi':\hat{C}\to Y$ and for $\sL$.
\begin{enumerate}
\item If $\theta:\hat{C}_1\to \hat{C}$ is a finite morphism between non-singular curves, then
$$
(\hat{S}\times_{\hat{C}}\hat{C}_1 \to \hat{C}_1, {\rm pr}_1^*\sM_{\hat{S}})
$$
is a mild reduction for $\chi'\circ\theta$.
\item In (1) base change induces an isomorphism $\theta^*\sF^{(\nu,\mu)}_{\hat{C}} \> \cong >> \sF_{\hat{C}_1}^{(\nu,\mu)}$
(which we will write again as an equality of sheaves).
\item[(3)] Let $\sigma: S \to X\times_Y\hat{C}$ be a modification 
of $X\times_Y\hat{C}$ with $S$ non-singular, and $h={\rm pr}_2\circ\sigma$.
In the (canonically) polarized case choose $\sM_S=\sigma^*{\rm pr}_1^*\sL$. In the saturated case choose
$\sM_S$ according to Lemma~\ref{sa.11} a). Then
$$
\sF^{(\nu,\mu)}_{\hat{C}}=h_*(\omega_{S/\hat{C}}^\nu\otimes \sM_{S}^\mu).
$$
In particular, the sheaf $\sF^{(\nu,\mu)}_{\hat{C}}$ is independent of the mild model.
\end{enumerate}
\end{lemma}
\begin{proof} 
Since $\hat{S}\times_CC_1 \to C_1$ is again mild (1) is obvious and
(2) follows by flat base change. (3) is a special case of Lemma~\ref{di.7}, using
in the saturated case for a smooth model dominating both, $\hat{S}$ and $S$, Lemma~\ref{sa.11} a).
\end{proof}
\begin{lemma}\label{mi.4}
In~\ref{sa.5} and~\ref{sa.7}, or in~\ref{sa.13} one may choose an open dense subscheme $Y_g\subset Y_0$ such that
for all morphisms
$$
\chi':\hat{C}\>\pi' >> \hat{Y}\>\varphi>> Y
$$ 
with $\hat{C}_g=\chi'^{-1}(Y_g) \neq \emptyset$ the tuple  $(\hat{S}:=\hat{Z}\times_{\hat{Y}}\hat{C}\to \hat{C},\sM_{\hat{S}}:={\rm pr}_1^*\sM_{\hat{Z}})$ is a mild  reduction for $\chi'$
and such that
\begin{equation}\label{eqmi.2}
\sF^{(\nu,\mu)}_{\hat{C}}=\pi'^*\sF^{(\nu,\mu)}_{\hat{Y}} \mbox{ \ \ for \ \ } (\nu,\mu)\in I.
\end{equation}
\end{lemma}
\begin{proof}
Choose $Y_g$, such that $\varphi^{-1}(Y_g)$ is contained in the open set $\hat{Y}_g$ in~\ref{sa.5} iv) or~\ref{sa.7} ix)
and such that $\hat{Z}$ is smooth over $\varphi^{-1}(Y_g)$. Then the definition of a mild morphism in~\ref{ws.1}
implies that $\hat{h}={\rm pr}_2:\hat{S}=\hat{Z}\times_{\hat{Y}}\hat{C}\to \hat{C}$ is mild. In the diagram (\ref{eqws.3}) in~\ref{ws.6} we require
the existence of a morphism $\hat{\varphi}:\hat{Z}\to X$ lifting $\varphi:\hat{Y}\to Y$, hence 
there is a modification $\varphi':\hat{Z}\to X\times_{Y}\hat{Y}$. The fibres of $\hat{Z}$ and $X\times_Y\hat{Y}$ over $\hat{Y}_g$ are smooth,
and $\varphi'$ restricts to a modification of those fibres. This implies that the induced morphism  
$\hat{Z}\times_{\hat{Y}}\hat{C}\to \hat{X}\times_{\hat{Y}}\hat{C}$ is birational. The equality in (\ref{eqmi.2}) follows from~\ref{sa.7} ix) and from the choice of $Y_g$.

It remains to verify the condition iii) in the saturated case, as stated in~\ref{mi.2}.
By assumption~\ref{sa.5} iv) the direct image $\hat{g}_*\omega_{\hat{Z}/\hat{Y}}^{\eta_0}=\hat{g}_*\varpi^{(\eta_0)}_{\hat{Z}}$ is locally free and compatible with base change for $\pi'$. Then the evaluation map for $\varpi_{\hat{S}}^{(\eta_0)}:={\rm pr}_1^*\varpi^{(\eta_0)}_{\hat{Z}}$
is surjective, and the first part of the condition iii) in~\ref{mi.2} holds true. The second condition just says that
the pullback of $\sL$ to $\hat{S}$ coincides with $\sM$ over some open subscheme of $\hat{C}$. This follows, since the same holds for $\sM_{\hat{Z}}$ over $\hat{Y}_0$. The last condition follows from Corollary~\ref{sa.14}. 
\end{proof}

\begin{proposition}\label{mi.5}
Let $\hat{C}$ be an irreducible curve, and let $\pi':\hat{C}\to \hat{Y}$ be a morphism.
If $\hat{C}_0=\pi'^{-1}(\hat{Y}_0)\neq \emptyset$ and if
$\chi'=\varphi\circ\pi'$ admits a mild reduction $(\hat{h}:\hat{S}\to \hat{C},\sM_{\hat{S}})$, then
$\sF^{(\nu,\mu)}_{\hat{C}}=\pi'^*\sF^{(\nu,\mu)}_{\hat{Y}}$
for $(\nu,\mu)\in I$.
\end{proposition}

\begin{proof}
Remark that one may replace $\hat{Y}$ in~\ref{sa.5},~\ref{sa.7} or~\ref{sa.13} by any modification, without loosing the properties i)---x). In particular the sheaves $\sF^{(\nu,\mu)}_{\hat{Y}}$ are compatible with pullback by dominant
morphisms for $(\nu,\mu)\in I$. Part (1) of Lemma~\ref{mi.3} allows to replace $\hat{C}$ by any covering,
hence dropping as usual the lower index ${}_1$ one can assume
that $\hat{Y}=\hat{Y}_1$ in~\ref{ec.2} and use the three properties stated there.
Let us write $h:S\to \hat{C}$ for the induced morphism and $\sM_S=\sM_{\hat{X}}|_S$. 

In the (canonically polarized case $\sM_{S}$ is the pullback of $\sL$ to $S$. By assumption $\hat{C}\to Y$ has a mild reduction $(\hat{h}:\hat{S}\to \hat{C},\sM_{\hat{S}})$. By~\ref{mi.3} (3)
$\sF^{(\nu,\mu)}_{\hat{C}}= h_*(\omega_{S/\hat{C}}^\nu\otimes \sM_{S}^\mu)$ and by Lemma~\ref{di.8}
this is the pullback of $\sF^{(\nu,\mu)}_{\hat{Y}}$.\vspace{.1cm}

For the saturated case, we have to argue in a slightly different way.
Recall that we defined in~\ref{sa.13} the invertible sheaves $\sB_{\hat{X}}^{(0,\kappa)}$ 
and $\varpi_{\hat{X}}^{(\eta_0)}$ as the images of the evaluation maps
$$
\hat{f}^*\hat{f}_*\sM_{\hat{X}}^\kappa \>>> \sM_{\hat{X}}^\kappa \mbox{ \ \ and \ \ }
\hat{f}^*\hat{f}_* \omega_{\hat{X}/\hat{Y}}^{\eta_0}\>>> \omega_{\hat{X}/\hat{Y}}^{\eta_0}.
$$  
Lemma~\ref{di.8} implies that the direct images $\hat{f}_*\sM_{\hat{X}}^\kappa$ and $\hat{f}_* \omega_{\hat{X}/\hat{Y}}^{\eta_0}$ are compatible with pullback. The sheaves 
$$
\sB_{S}^{(0,\kappa)}= \sB_{\hat{X}}^{(0,\kappa)}|_S \mbox{ \ \  and \ \ }
\varpi_{S}^{(\eta_0)}=\varpi_{\hat{X}}^{(\eta_0)}|_S,
$$
are again invertible and the images of the evaluation maps for 
$\sM_S^\kappa$ and $\omega_{S/\hat{C}}^{\eta_0}$, respectively. The latter implies that the divisor $\Pi_S^{(\eta_0)}$
is the pullback of $\Pi_{\hat{X}}^{(\eta_0)}$. By the definition of $\kappa$-saturated in~\ref{sa.10} and by Lemma~\ref{sa.11} c) one knows that
$$
\hat{f}_*\sB_{\hat{X}}^{(0,\kappa)}=
\hat{f}_*\sM_{\hat{X}}^\kappa = \hat{f}_*(\sM_{\hat{X}}^\kappa\otimes \sO_{\hat{X}}(*\cdot\Pi_{\hat{X}}^{(\eta_0)}))=
\hat{f}_*(\rho^*\sL^\kappa\otimes \sO_{\hat{X}}(*\cdot\Pi_{\hat{X}}^{(\eta_0)})).
$$
Lemma~\ref{di.8} implies that the corresponding property holds true for $S$ instead of $\hat{X}$.

By assumption $\hat{C}\to Y$ has a mild $\kappa$-saturated reduction $(\hat{h}:\hat{S}\to \hat{C},\sM_{\hat{S}})$. 
Let $\Psi:W\to S$ and $\Psi':W\to \hat{S}$ be modifications, with $W$ smooth. By~\ref{di.7} 
$$
\hat{h}_*\varpi_{\hat{S}}^{(\eta_0)}=\hat{h}_*\omega_{\hat{S}/\hat{C}}^{\eta_0}=h_*\omega_{S/\hat{C}}^{\eta_0}=h_*\varpi_{S}^{(\eta_0)},
$$ 
hence
$\Psi'^*\varpi_{\hat{S}}^{(\eta_0)}=\Psi^*\varpi_{S}^{(\eta_0)}$. Call this sheaf
$\varpi_{W}^{(\eta_0)}$. The divisor $\Pi_{W}^{(\eta_0)}$ with
$$
\omega_{W/\hat{C}}^{\eta_0}=\varpi_{W}^{(\eta_0)}\otimes\sO_{W}(-\Pi_{W}^{(\eta_0)})
$$
is of the form $\Psi'^*\Pi_{\hat{S}}^{(\eta_0)}+\eta_0\cdot E_{W/\hat{S}}=\Psi'^*\Pi_{S}^{(\eta_0)}+\eta_0\cdot E_{W/S}$, 
where $E_{W/\hat{S}}$ and $E_{W/S}$ are relative canonical divisors. If $\sL_\bullet$ denotes the pullback of $\sL$,
as in~\ref{sa.9} one finds that for all $\varepsilon\geq 0$ 
$$
\hat{h}_*(\sL^\kappa_{\hat{S}}\otimes\sO_{\hat{S}}(\varepsilon\cdot\Pi_{\hat{S}}^{(\eta_0)}))=
h_*(\sL^\kappa_{S} \otimes\sO_{S}(\varepsilon\cdot\Pi_{S}^{(\eta_0)})),
$$
and that for  some $\varepsilon_0$ and all $\varepsilon\geq \varepsilon_0$, both sheaves are independent of $\varepsilon$.
Since for those $\varepsilon$ the left hand side is $\hat{h}_*\sB_{\hat{S}}^{(0,\kappa)}$ and the right hand side $h_*\sB_{S}^{(0,\kappa)}$
the two sheaves are equal. This implies that $\Psi'^*\sB_{\hat{S}}^{(0,\kappa)}=\Psi^*\sB_{S}^{(0,\kappa)}$.

The divisor $\Sigma_{\hat{S}}^{(0,\kappa)}$ and $\Sigma_{S}^{(0,\kappa)}$ have the same support as 
$\Pi_{\hat{S}}^{(\eta_0)}\cap \hat{h}^{-1}(\hat{C}\setminus \hat{C}_0)$ and $\Pi_{S}^{(\eta_0)}$, respectively. Define
$\Sigma$ to be the smallest divisor on $W$, larger than $\Psi'^*\Sigma_{\hat{S}}^{(0,\kappa)}
\cap \hat{h}^{-1}(\hat{C}\setminus \hat{C}_0)$ and $\Psi^*\Sigma_{S}^{(0,\kappa)}$. Adding components of
$\Pi_{W}^{(\eta_0)}$ one finds some $\Sigma_{W}^{(0,\kappa)}$ such that
$$
\Psi'^*\sB_{\hat{S}}^{(0,\kappa)}\otimes \sO_{W}(\Sigma_{W}^{(0,\kappa)})
$$
is the $\kappa$-th power of an invertible subsheaf $\sM_{W}$ of
$\sL_W\otimes\sO_{W}(*\cdot\Pi_{W}^{(\eta_0)})$. Obviously
$\Psi'_*\sM_{W}=\sM_{\hat{S}}$ and $\Psi_*\sM_{W}=\sM_{S}$, hence we are allowed to apply
\ref{di.7} and find 
$$
\hat{h}_*(\omega_{\hat{S}/C}^\nu\otimes\sM^\mu_{\hat{S}})=h_*(\omega_{S/C}^\nu\otimes\sM^\mu_{S})=
\sF_{\hat{Y}}^{(\nu,\mu)}|_{\hat{C}}.
$$
\end{proof}
\section{A variant for multiplier ideals}\label{va}
Let us return to the set-up in~\ref{de.1} and to the assumptions introduced in~\ref{de.3}
or in Variant~\ref{de.4}. 
As in Section \ref{mi} we assume that $\sM_{\hat{Z}}$ and $\sM_{\hat{X}}$ are either the structure sheaves, or the pullback of an invertible sheaf $\sL$ on $X$, or $\kappa$-saturated extensions of $\sL$.

Consider again a non-singular curve $\hat{C}$ and a morphism $\chi':\hat{C}\to Y$ whose image meets
$Y_0$, and a mild reduction $(\hat{h}:\hat{S}\to \hat{C},\sM_{\hat{S}})$ for $\sL$, as defined in~\ref{mi.2}. In particular one has a morphism $\zeta:\hat{S}\to X$, and the sheaves
$$
\sF^{(\nu,\mu)}_{\hat{C}}=\hat{h}_*(\omega_{\hat{S}/\hat{C}}^\nu\otimes \sM_{\hat{S}}^\mu)
$$
are defined. Lemma~\ref{mi.3} and Proposition~\ref{mi.5} imply that
$\chi'^*\sF^{(\nu,\mu)}_{\hat{Y}}=\tau^*\sF^{(\nu,\mu)}_{\hat{C}}$, whenever one has a lifting
\begin{equation}\label{eqva.1}
\xymatrix{C' \ar[r]^{\chi'}\ar[d]_{\tau}&\hat{Y}\ar[d]^{\varphi}\\
\hat{C} \ar[r]^{\chi}& Y}
\end{equation}
with $C'$ a non-singular curve.

We will need that the different invertible sheaves and divisors introduced in~\ref{sa.7},~\ref{sa.9} or~\ref{sa.13}, and in
Section~\ref{de} are defined for the morphism $\hat{h}:\hat{S}\to \hat{C}$. 

\begin{assumption}\label{va.1} Assume that the assumptions made in~\ref{de.1} and~\ref{de.3}
hold true, and that $\hat{Y}$, $\hat{Z}$ and $\hat{X}$ is chosen according to Lemma~\ref{de.6}.
\begin{enumerate}
\item[1.] $(\hat{h}:\hat{S}\to \hat{C},\sM_{\hat{S}})$ is a mild reduction for $\chi':\hat{C}\to Y$ and for $\sL$. For $\eta_0$ the image $\varpi_{\hat{S}}^{(\eta_0)}$ of $\hat{h}^*\hat{h}_*\omega_{\hat{S}/\hat{C}}^\nu$ in $\omega_{\hat{S}/\hat{C}}^{\eta_0}$ is locally free, as well as for $(\beta,\alpha)\in \widetilde{I}$ the images $\sB_{\hat{S}}^{(\beta,\alpha)}$ of the evaluation maps of $\omega_{\hat{S}/\hat{C}}^\beta\otimes \sM^\alpha_{\hat{S}}$.
\item[2.] There exists a subsheaf $\sE_{\hat{C}}$ of $\sF_{\hat{C}}^{(\beta_0,
\alpha_0)}$, with $\chi'^*\sE_{\hat{Y}}=\tau^*\sE_{\hat{C}}$, for all liftings $\hat{\varphi}$ as in (\ref{eqva.1}). Moreover the image $\sB_{\hat{S}}$ of the evaluation map
$$
\hat{h}^*\sE_{\hat{C}} \>>> \omega_{\hat{S}/\hat{C}}^{\beta_0}\otimes 
\sM^{\alpha_0}_{\hat{S}}
$$
is invertible.
\end{enumerate}
\end{assumption}
\begin{remark}\label{va.2}
If in~\ref{de.4} one has  
$\displaystyle \sE_{\hat{Z}}=  \bigotimes_{\iota=1}^s\hat{g}_*\sB_{\hat{Z}}^{(\beta_\iota,\alpha_\iota)}$ 
the condition 2) in~\ref{va.1} follows from the assumption 
$(\beta_1,\alpha_1), \ldots (\beta_s,\alpha_s) \in \widetilde{I}$ for $\iota=1,\ldots,s$.
 
In fact, the latter implies that the pullback of the sheaves
$\sF_{\hat{S}}^{(\beta_\iota,\alpha_\iota)}$ and $\sF_{\hat{Y}}^{(\beta_\iota,\alpha_\iota)}$
coincide on $C'$, and so does their image under the multiplication map.

If $\sE_{\hat{Z}}$ is smaller, we will need that it is defined on a compactification
of $Y$, in order to enforce the compatibility condition 2) in~\ref{va.1}.
\end{remark}
We will write again $\Pi_{\hat{S}}^{(\eta_0)}$, $\Sigma_{\hat{S}}^{(\beta,\alpha)}$ 
and $\Sigma_{\hat{S}}$ for the divisors given by the inclusions $\varpi_{\hat{S}}^{(\eta_0)}\subset \omega_{\hat{S}/\hat{C}}^{\eta_0}$, $\sB_{\hat{S}}^{(\beta,\alpha)}\subset \omega_{\hat{S}/\hat{C}}^\beta\otimes \sM^\alpha_{\hat{S}}$ and $\sB_{\hat{S}}\subset \omega_{\hat{S}/\hat{C}}^{\beta_0}\otimes 
\sM^{\alpha_0}_{\hat{S}}$.

As in~\ref{de.5} one defines the different products, models, sheaves and divisors,
with $\hat{g}:\hat{Z}\to \hat{Y}$ replaced by $\hat{h}:\hat{S}\to \hat{C}$. In particular we have again the divisor
$$
\Delta_{\hat{S}^r}=b \cdot (\Gamma_{\hat{S}^r}+\Sigma_{\hat{S}^r}^{( \eta,\gamma)}) + \Sigma_{\hat{S}^r},
$$
on the $r$-fold fibre product $\hat{h}^r:\hat{S}^r\to \hat{C}$, and we define
\begin{gather*}
\sG_{\hat{C}}^{(\beta+\frac{\eta}{\ell},\alpha+\frac{\gamma}{\ell})}=\hat{h}^r_*\Big(\omega_{\hat{S}^r/\hat{C}}^{\beta+\frac{\eta}{\ell}} \otimes\sM_{\hat{S}^r}^{\alpha+\frac{\gamma}{\ell}}\otimes\sJ\big(-\frac{1}{N}\cdot \Delta_{\hat{S}^r}\big)\Big)=\\ 
\hat{h}^r_*\Big(\omega_{\hat{S}^r/\hat{C}}^{\beta+\frac{\eta}{\ell}} \otimes\sM_{\hat{S}^r}^{\alpha+\frac{\gamma}{\ell}}
\otimes\sJ\big(-\frac{1}{e\cdot\ell}\cdot (\Gamma_{\hat{S}^r}+\Sigma_{\hat{S}^r}^{( \eta,\gamma)})-\frac{1}{N}\cdot \Sigma_{\hat{S}^r}\big)\Big),
\end{gather*}
where $N=b\cdot e\cdot \ell$ and where $\Gamma_{\hat{S}^r}$ is the zero divisor induced
by the natural inclusion
$$
\bigotimes_{i=1}^s \det (\hat{h}_*\sM^{\gamma_i}_{\hat{S}}\otimes \varpi_{\hat{S}}^{\eta_i}))^{\otimes \frac{r}{r_i}}\\ 
\>\Xi^{(r)} >>
\bigotimes^r \hat{h}_*(\omega_{\hat{S}/\hat{Y}}^\eta\otimes\sM_{\hat{Z}}^\gamma)=
\hat{h}^r_*\sB_{\hat{S}^r}^{(\eta,\gamma)}.
$$
Again we should have written $\sG_{\hat{C}}^{(\Xi^{(r)},\sE;\beta+\frac{\eta}{\ell},\alpha+\frac{\gamma}{\ell})}$
since the sheaf depends on $\Xi^{(r)}$ and $\sE$, but we hope that the reader will not forget.
\begin{lemma}\label{va.3}
Let $\theta:\hat{C}_1\to \hat{C}$ be a finite non-singular covering, and let 
$$
\xymatrix{\hat{S}^r_1 \ar[r]^{\theta'}\ar[d]_{\hat{h}^r_1}&\hat{S}^r\ar[d]^{\hat{h}^r}\\
\hat{C}_1 \ar[r]^{\theta}&\hat{C}}
$$
be the induced morphism. Then:
\begin{enumerate}
\item[a.] If $\hat{h}:\hat{S}\to \hat{C}$ satisfies the assumption~\ref{va.1} then $\hat{h}_1:\hat{S}_1\to \hat{C}_1$
satisfies the same assumption.
\item[b.] $\sJ\big(-\frac{1}{N}\cdot \Delta_{\hat{S}^r_1}\big)$ 
is a subsheaf of $\theta'^*\sJ\big(-\frac{1}{N}\cdot \Delta_{\hat{S}^r}\big)$
\item[c.] There is a natural inclusion
$$
\sG_{\hat{C}_1}^{(\beta+\frac{\eta}{\ell},\alpha+\frac{\gamma}{\ell})} \>>> \theta^*\sG_{\hat{C}}^{(\beta+\frac{\eta}{\ell},\alpha+\frac{\gamma}{\ell})}.
$$
\end{enumerate}
\end{lemma}
\begin{proof}
As in the proof of Lemma~\ref{sa.2} part a) of~\ref{va.3} follows from 
Lemma~\ref{mi.3} (1) and (2).

For b) remark that ${\rm pr}_1:\hat{S}^r_1\to \hat{S}^r$ is flat, hence 
$\theta'^*\sJ\big(-\frac{1}{N}\cdot \Delta_{\hat{S}^r}\big)$
has no torsion. Consider a desingularization $\tau:S\to \hat{S}^r$ such that all fibres
are normal crossing divisors, and such that $\tau^*\Gamma_{\hat{S}^r}$ is a relative normal crossing divisor. So $\tau^*(\Delta_{\hat{S}^r})$ is a normal crossing divisor, as well.

Let $\tau_1:S_1\to \hat{S}^r_1$ be the normalization of the pullback family,
$$
\xymatrix{S_1\ar[r]^{\sigma\hspace*{.5cm}}\ar[dr]_{\tau_1}
&S\times_{\hat{S}^r}\hat{S}^r_1 \ar[r]^{\hspace*{.5cm}{\rm pr}_1}\ar[d]_{{\rm pr}_2}& S\ar[d]^{\tau}\\
& \hat{S}^r_1 \ar[r]^{\theta'}&\hat{S}^r}
$$
and $\theta''={\rm pr}_1\circ\sigma$ the induced morphisms. By flat base change
$$
\theta'^*\tau_*\omega_{S/\hat{C}}\big(-\big[\frac{1}{N}\cdot \tau^*\Delta_{\hat{S}^r}\big]\big)=
{\rm pr}_{2*} {\rm pr}_1^*\Big( \omega_{S/\hat{C}_1}\big(- \big[\frac{1}{N}\cdot \tau^*\Delta_{\hat{S}^r})\big]\big)\Big).
$$
Dualizing sheaves become smaller under normalizations, and this sheaf contains 
$$
\tau_{1*} \omega_{S_1/\hat{C}_1}\otimes {\theta''}^* \sO_{S}\big(- \big[\frac{1}{N}\cdot \tau^*\Delta_{\hat{S}^r}\big]\big).
$$
Since $S_1$ has at most rational Gorenstein singularities, this sheaf remains the same
if we replace $S_1$ by a desingularization. Hence by abuse of notations we may assume that $S_1$ is non-singular, that the fibres of $S_1\to \hat{C}_1$ are normal crossing divisors, and that ${\theta''}^*\tau^*\Gamma_{\hat{S}^r}$ is a relative normal crossing divisor.

Obviously one has an inclusion
\begin{gather*}
\sO_{S_1}\big(- \big[\frac{1}{N}\cdot {\theta''}^*\tau^*\Delta_{\hat{S}^r})\big]\big) \subset
{\theta''}^* \sO_{S}\big(- \big[\frac{1}{N}\cdot \tau^*\Delta_{\hat{S}^r}\big]\big)\\
\mbox{and hence \ \ }
\tau_{1*} \omega_{S_1/\hat{C}_1}\big(- \big[\frac{1}{N}\cdot \tau^*_1 \Delta_{\hat{S}^r_1}\big]\big)
\subset
\theta'^*\tau_*\omega_{S/\hat{C}}\big(-\big[\frac{1}{N}\cdot \tau^*\Delta_{\hat{S}^r}\big]\big)
\end{gather*}
as claimed in b). By flat base change c) follows from b).
\end{proof}
Let $\tau:S\to X \times_Y\hat{C}$ be any desingularization of the main component, and
let $h:S\to \hat{C}$ denote the induced morphism.
Recall that we assumed $\chi':\hat{C}\to Y$ to have a mild reduction. So we may choose $\sM_S$ as the pullback of $\sL$ in the (canonically) polarized case or by Lemma~\ref{sa.11} a) in the saturated case. 
Blowing up, we may assume that for $(\nu,\mu)\in \widetilde{I}$ the images $\sB_S^{(\nu,\mu)}$ of the evaluation maps
are invertible, in particular the image $\varpi^{(\eta_0)}_S$ of $h^*h_*\omega^{\eta_0}_{S/\hat{C}}\to \omega^{\eta_0}_{S/\hat{C}}$.
We write $h^{(r)}:S^{(r)}\to \hat{C}$ for the family obtained by desingularizing
the $r$-fold product $S^r=S\times_{\hat{C}}\cdots\times_{\hat{C}}S$, where again we assume that 
$\varpi^{(\eta_0)}_{S^{(r)}}$ is invertible.  

As above, or in Section~\ref{de} one chooses the sheaf $\sM_{S^r}$ as the exterior tensor product.
$\sM_{S^{(r)}}$ will denotes its pullback to $S^{(r)}$. 
Since $\chi':\hat{C}\to Y$ has a mild reduction,~\ref{di.5} implies that
one has again the inclusions
$$
{h^{(r)}}^*\bigotimes_{i=1}^s \det (\hat{g}_*\sM^{\gamma_i}_{\hat{S}}\otimes \varpi_{\hat{S}}^{\eta_i}))^{\otimes \frac{r}{r_i}}
\>>> \sB_{S^{(r)}}^{(\eta,\gamma)}
$$
with zero locus $\Gamma_{S^{(r)}}$. Writing $S^{(r)}_0={h^{(r)}}^{-1}(\chi'^{-1}(Y_0))$ 
for the smooth part of $h^{(r)}$ one obtains by~\ref{de.11}:

\begin{lemma}\label{va.4}
\begin{gather*}
\sG_{\hat{C}}^{(\beta+\frac{\eta}{\ell},\alpha+\frac{\gamma}{\ell})}=
h^{(r)}_*\Big(\omega_{S^{(r)}/\hat{C}}^{\beta+\frac{\eta}{\ell}} \otimes\sM_{S^{(r)}}^{\alpha+\frac{\gamma}{\ell}}
\otimes\sJ\big(-\frac{1}{e\cdot\ell}\cdot (\Gamma_{S^{(r)}}+\Sigma_{S^{(r)}}^{( \eta,\gamma)})-\frac{1}{N}\cdot \Sigma_{S^{(r)}}\big)\Big),\\
\mbox{and \ \ }\sJ\big(-\frac{1}{e\cdot\ell}\cdot (\Gamma_{S^{(r)}}+\Sigma_{S^{(r)}}^{( \eta,\gamma)})-\frac{1}{N}\cdot \Sigma_{S^{(r)}}\big)|_{S^{(r)}_0}=\sO_{{S^{(r)}_0}}.
\end{gather*}
In particular the inclusion 
$$
\sG_{\hat{C}}^{(\beta+\frac{\eta}{\ell},\alpha+\frac{\gamma}{\ell})}\subset
\bigotimes^r \sF_{\hat{C}}^{(\beta+\frac{\eta}{\ell},\alpha+\frac{\gamma}{\ell})}
$$ 
is an isomorphism over $\chi'^{-1}(Y_0)$.
\end{lemma}

\begin{definition}\label{va.5}
The mild reduction $(\hat{h}:\hat{S}\to \hat{C},\sM_{\hat{S}})$ is {\em exhausting} (or {\em exhausting for $(\Xi^{(r)},\sE;\beta+\frac{\eta}{\ell},\alpha+\frac{\gamma}{\ell})$}) 
if the properties 1) and 2) in~\ref{va.1} hold true and if:
\begin{enumerate}
\item[3.] For all finite coverings of non-singular curves $\hat{C}_1 \to \hat{C}$ the inclusion
$$
\sG_{\hat{C}_1}^{(\beta+\frac{\eta}{\ell},\alpha+\frac{\gamma}{\ell})} \>>> \theta^*\sG_{\hat{C}}^{(\beta+\frac{\eta}{\ell},\alpha+\frac{\gamma}{\ell})}.
$$
in~\ref{va.3} c) is an isomorphism.
\end{enumerate}
\end{definition}
The Lemmata~\ref{sa.2} and~\ref{de.6} imply that given $\chi':\hat{C}\to Y$ one can always find a finite
covering $\hat{C}_1\to \hat{C}$ and a mild reduction for the induced morphism $\hat{C}_1\to Y$ which is exhausting. Repeating the argument used to prove~\ref{mi.4} one obtains in addition:
\begin{lemma}\label{va.6}
There exists in~\ref{de.6} an open dense subscheme $Y_g\subset Y_0$
such that for all morphisms 
$$
\chi':\hat{C}\>\pi' >> \hat{Y} \> \varphi >> Y
$$ 
from a non-singular curve $\hat{C}$, with $\chi'^{-1}(Y_g)$ dense, $\chi'$ 
admits a mild exhausting reduction $(\hat{h}:\hat{S}\to \hat{C},\sM_{\hat{S}})$.
Moreover
$$
\sG_{\hat{C}}^{(\beta+\frac{\eta}{\ell},\alpha+\frac{\gamma}{\ell})}=\pi'^*\sG_{\hat{Y}}^{(\beta+\frac{\eta}{\ell},\alpha+\frac{\gamma}{\ell})}.
$$
\end{lemma}

\begin{proposition}\label{va.7} \
Consider in~\ref{de.6} a morphism $\pi':\hat{C}\to \hat{Y}$ from a non-singular curve $\hat{C}$ with $\pi'^{-1}(\hat{Y}_0)\neq \emptyset$.

If $\chi'=\varphi\circ\pi'$ admits a mild exhausting reduction $(\hat{h}:\hat{S}\to \hat{C},\sM_{\hat{S}})$, then
$$
\sG_{\hat{C}}^{(\beta+\frac{\eta}{\ell},\alpha+\frac{\gamma}{\ell})}=\pi'^*\sG_{\hat{Y}}^{(\beta+\frac{\eta}{\ell},\alpha+\frac{\gamma}{\ell})}.
$$
\end{proposition}
\begin{proof}
By~\ref{mi.5} 
$$
\sF_{\hat{C}}^{(\beta+\frac{\eta}{\ell},\alpha+\frac{\gamma}{\ell})}=\pi'^*\sF_{\hat{Y}}^{(\beta+\frac{\eta}{\ell},\alpha+\frac{\gamma}{\ell})},
$$
and the both sheaves remain unchanged if one replaces $\hat{C}$ by some finite covering or
$\hat{Y}$ by some alteration. The same holds for the subsheaves
$\sG_{\hat{C}}^{(\beta+\frac{\eta}{\ell},\alpha+\frac{\gamma}{\ell})}$ and
$\pi'^*\sG_{\hat{Y}}^{(\beta+\frac{\eta}{\ell},\alpha+\frac{\gamma}{\ell})}$. Hence they
coincide, if and only if they coincide on some $C'$ finite over $\hat{C}$.

By assumption the multiplier ideal $\sJ\big(-\frac{1}{N}\cdot \Delta_{\hat{Z}}\big)$ is compatible with pullback, base change, and products for all alterations. In particular, $\sJ\big(-\frac{1}{N}\cdot \Delta_{\hat{Z}}\big)$ is flat over $\hat{Y}$. 

We are allowed to replace $\hat{Y}$ by an alteration or by an open neighborhood of the image of $C$, 
hence by an local alteration for $\hat{C}$.
So by Proposition~\ref{ec.8} we may assume that $\pi'$ is an embedding, that $\hat{f}$ is flat and that $S=\hat{f}^{-1}(\hat{C})$ is non-singular and semistable over $\hat{C}$. By abuse of notations, we will allow $\hat{X}$ to be normal with rational Gorenstein singularities. By Lemma~\ref{ec.6} this holds for the total space of pullbacks under local alterations for $\hat{C}$, and for the fibre products. So
we will work with the condition that $\hat{f}:\hat{X}\to \hat{Y}$ is a weak semistable reduction for $\hat{C}$,
a condition which is compatible with pullbacks and products.
In particular $S^r$ is normal with at most rational Gorenstein singularities and
$h^r:S^r\to \hat{C}$ has reduced fibres. 

As stated in~\ref{sa.8} one is allowed to replace the mild family
$\hat{g}^r:\hat{Z}^r\to \hat{Y}$ by some mild model, dominating the flat part of the weak semistable reduction $\hat{f}^r:\hat{X}^r\to \hat{Y}$. 
Here we might loose the compatibility of $\sJ\big(-\frac{1}{N}\cdot \Delta_{\hat{Z}}\big)$ with pullback, base change, and products for all alterations. Theorem~\ref{fm.5} allows to repair this defect, by replacing $\hat{Y}$ by some larger local alteration.

The morphism $\hat{f}$ is smooth over $\hat{Y}_0$, and $\sJ\big(-\frac{1}{N}\cdot \Delta_{\hat{X}^r_1}\big)|_{S_1^r}$
is trivial over $\hat{X}_0$. So we may apply Proposition~\ref{ec.11}
\end{proof}

\section{Uniform mild reduction and the extension theorem}\label{ex}
Constructing the locally free sheaves $\sF^{(\nu,\mu)}_{\hat{Y}}$ and $\sG_{\hat{Y}}^{(\beta+\frac{\eta}{\ell},\alpha+\frac{\gamma}{\ell})}$ we used the
Weak Semistable Reduction Theorem several times and we have no control 
on the alteration $\hat{Y}$ of $Y$. As already indicated in \ref{di.9} we will
show how to use Gabber's Extension Theorem, recalled in~\ref{ex.6}, to obtain
those sheaves on a finite covering of a projective compactification of $Y_0$. Again the latter will be denoted by $Y$ and the covering will be written as $\phi:W\to Y$. In Proposition~\ref{ex.8} the corresponding result is formulated for the sheaves
$\sF_\bullet^{(\nu,\mu)}$ whereas the extension to the direct images $\sG_\bullet^{(\beta+\frac{\eta}{\ell}, \alpha + \frac{\gamma}{\ell}}$ of multiplier sheaves is
given in Variant~\ref{ex.11}. The Variant~\ref{ex.10} handles the case of the determinants of
$\sF_\bullet^{(\nu,\mu)}$, starting from a covering of a coarse moduli scheme. As an application we will
state and prove a generalization of Theorem~\ref{in.1} allowing arbitrary polarizations in Theorem~\ref{ex.12} and a variant for saturated polarizations in~\ref{ex.13}.

We will need in all those cases that the trace map of $\phi:W\to Y$ splits, i.e. that $\sO_Y$ is a direct factor of $\phi_*\sO_W$. By \cite[Lemma 2.2]{Vie} each finite surjective morphism $\widetilde W\to Y$ of reduced schemes, with $\widetilde W$ normal, factors through a finite covering $\phi:W\to Y$ with a splitting trace map and with $\widetilde W \to W$ birational. We will give here a different construction, starting with a fixed embedding $Y\hookrightarrow \BP^N$, or more generally with any embedding $Y\hookrightarrow \BP$ for $\BP$ irreducible, normal and projective. The latter will have the advantage to allow the gluing, needed to show the uniform mild extension over curves, required by the extension theorem.

\begin{lemma}\label{ex.1}
Let $\Psi:\BP'\to \BP$ be a finite normal covering and let $Y\subset \BP$ be a closed subvariety. Then $\phi:W=\Psi^{-1}(Y)\to Y $ has a splitting trace map. 
\end{lemma}
\begin{proof}
Since $\BP'$ is normal $\sO_{\BP}$ is a direct factor of $\Psi_*\sO_{\BP'}$,  
hence there is a surjection $\Psi_*\sO_{\BP'} \to \sO_{Y}$. Obviously this factors through 
$\phi_*\sO_{W} \to \sO_{Y}$.
\end{proof}
\begin{definition}\label{ex.2} Let $\widetilde\phi:\widetilde W\to Y$ be a surjective finite map and $Y\to \BP$ a closed embedding with $\BP$ irreducible, normal and projective.
Then $\Psi:\BP'\to \BP$ {\em dominates} $\widetilde{W}$, if $\BP'$ is normal and irreducible, 
if $\Psi$ is a finite covering and if the normalization $V \to \Psi^{-1}(Y) \to Y$ factors like
$V \to \widetilde{W} \to Y$.
\end{definition}
\begin{lemma}\label{ex.3} Given $\widetilde\phi:\widetilde W\to Y$ and $\BP$ as in Definition~\ref{ex.2},
there exists a finite normal covering $\Psi:\BP'\to \BP$ dominating $\widetilde{W}$. 
Moreover one may choose $\Psi:\BP'\to \BP$ to be a Galois covering.
\end{lemma}
\begin{proof}
In order to prove Lemma~\ref{ex.3} we may assume that $\widetilde{W}$ is normal.

Let us first assume that $Y$ and $\widetilde{W}$ are both irreducible and write $K$ for the function field of $Y$. The function field of $\widetilde{W}$ can be written as $K[T]/f$ for a monic irreducible polynomial $f\in K[T]$. For some open subscheme $U\subset \BP$ the polynomial $f$ lies in $\sO_{Y}(U\cap Y)[T]$
and lifts to a monic irreducible polynomial $F\in \sO_{\BP}(U)[T]$. Choose $\BP'$
as the normalization of $\BP$ in $L[X]/ F$, where $L$ denotes the function field of $\BP$. The preimage of $Y$ in $\BP'$ is birational to $\widetilde{W}$ and since $\widetilde{W}\to Y$ is finite, the normalization $V$ of $\Psi^{-1}(Y)$ dominates $\widetilde{W}$.

Next assume that $Y$ is irreducible, and that $W_1,\cdots,W_s$ are the components of $\widetilde{W}$.
We already know how to construct $\Psi_\iota:\BP'_\iota \to \BP$, dominating $W_\iota \to Y$. We choose $\Psi:\BP'\to \BP$ as the normalization of $\BP_1\times_{\BP} \cdots \times_{\BP}\BP_s$. Then $\Psi^{-1}(Y)$ is finite over $W_1\times_Y \cdots \times_Y W_s$, and we obtain the
factorization of the normalization. This remains true if we replace $\BP'$ by a larger covering, hence
we can glue the construction for different components of $Y$ in the same way.

Finally, if $\Psi:\BP'\to \BP$ is a finite covering,
dominating $\widetilde{W}$, the normalization of $\BP'$ in the Galois hull
of the function field $\C(\BP')$ over $\C(\BP)$ will again dominate $\widetilde{W}$.
So we can add the property ``Galois'' as well.
\end{proof}
\begin{lemma}\label{ex.4}
Let $\Psi:\BP'\to \BP$ be a finite morphism between normal schemes, let $Y\subset \BP$ a closed subscheme and $Y_0\subset Y$ an open set. Let $\overline{W}$ be a modification of $W=\Psi^{-1}(Y)$ with centers outside $W_0=\Psi^{-1}(Y_0)$. Then there exist normal modifications 
$\hat{\BP}\to \BP$ and $\hat{\BP}'\to \BP'$ with centers in $Y\setminus Y_0$ and $W\setminus W_0$, respectively, such that the induced rational map $\Psi_1:\hat{\BP}' \to \hat{\BP}$ is  a finite morphism and
such that the proper transform of $W$ is $\overline{W}$.
\end{lemma}
\begin{proof} It is sufficient to consider irreducible varieties $\BP$ and $\BP'$.
Assume first that $\BP'$ is Galois over $\BP$, say with Galois group $\Gamma$. One can extend the modification $\overline{W}\to W$ to a modification $\BP'_{\rm id} \to \BP'$ by blowing up an ideal $\sJ$, with $\sO_{\BP'}/\sJ$ supported in $W\setminus W_0$.
Blowing up the conjugates of $\sJ$ under $\sigma\in \Gamma$ the action of $\sigma\in \Gamma$, given by
$\sigma:\BP'\to \BP'$ extends to a morphism $\sigma:\BP'_{\rm id} \to \BP'_\sigma$.
So $\Gamma$ acts on the fibre product $\bigtimes{\hspace{-.2cm}}_\BP \ \BP'_\sigma,$ taken over all $\sigma\in \Gamma$. 

Obviously $\bigtimes{\hspace{-.2cm}}_\BP \ \BP'_\sigma$ contains an open dense subset of $\BP'$, embedded diagonally, and we choose $\hat{\BP}'$ to be the normalization of its closure.
The projection to $\BP'_{\rm id}$ induces the morphism $\hat{\BP}' \to \BP'$.
The group $\Gamma$ acts on $\hat{\BP}'$, and we can choose for $\hat{\BP}$ the quotient.

If $\BP'$ is not Galois, we replace $\BP'$ by its normalization $\BP''$ in the Galois hull of the function field extension for $\BP'\to \BP$. So $\BP'$ is the quotient of $\BP''$ by some subgroup
$\Gamma'\subset \Gamma$. Having constructed $\hat{\BP}''$, we choose $\hat{\BP}'=\hat{\BP}''/{\Gamma'}$.
\end{proof}

Let us recall Gabber's Extension Theorem. We start with the following set-up. 
\begin{setup}\label{ex.5}
Let $\BP$ be a normal projective scheme, $\widetilde{Y}\subset \BP$ a closed reduced subscheme, and
let $\widetilde{Y}_0\subset \widetilde{Y}$ be open and dense. Let $\Psi:\BP' \to \BP$ be a finite covering,
with $\BP'$ normal, and write $W=\Psi^{-1}(\widetilde{Y})$, $W_0=\Psi^{-1}(\widetilde{Y}_0)$, $\widetilde{\phi}=\Psi|_{W}$ and
$\widetilde{\phi}_0=\Psi|_{W_0}$. Consider a modification
$\xi_0:\hat{Y}_0\to W_0$ with $\hat{Y}_0$ non-singular, and a projective manifold $\hat{Y}$
containing $\hat{Y}_0$ as an open dense subscheme. 

Let $\sC_{\widetilde{Y}_0}$ and $\sC_{\hat{Y}}$ be locally free sheaves on $\widetilde{Y}_0$ and $\hat{Y}$ respectively, such that
for $\sC_{W_0}=\widetilde{\phi}_0^*\sC_{\widetilde{Y}_0}$ one has:
\begin{enumerate}
\item[i.] $\xi_0^*\sC_{W_0}=\sC_{\hat{Y}}|_{\hat{Y}_0}$.
\item[ii.] For each morphism $\pi:C\to \BP'$ from a non-singular projective curve $C$
with $C_0= \pi^{-1}(W_0)$ dense in $C$, the sheaf $(\pi|_{C_0})^*\sC_{W_0}$
extends to a locally free sheaf $\sC_C$ such that:
\begin{enumerate}
\item[a.] If $\pi':\hat{C}\to \BP'$ factors through $\iota:\hat{C}\to C$ then $\sC_{\hat{C}}=\iota^*\sC_C$.
\item[b.] If $\pi:C\to \BP'$ lifts to a morphism $\chi:C\to \hat{Y}$ then
$\sC_C=\chi^*\sC_{\hat{Y}}$.
\end{enumerate}
\end{enumerate}
\end{setup}
\begin{theorem}\label{ex.6}
In~\ref{ex.5}, blowing up $\BP$ with centers not meeting $\widetilde{Y}_0$ and replacing 
$\BP'$ by the normalization of $\BP$ in its function field, one finds
an extension of $\sC_{W_0}$ to a locally free sheaf $\sC_W$ on $W=\Psi^{-1}(\widetilde{Y})$ such that for all commutative diagrams
\begin{equation*}
\xymatrix{\hat{Y}_0 \ar[r]^{\subset}\ar[d]_{\xi_0}&\hat{Y}&\ar[l]_{\psi}\ar[dl]^{\rho}\Lambda \\
W_0 \ar[r]^{\subset}& W}
\end{equation*}
with $\psi$ either a dominant morphism, or a morphism from a non-singular curve $\Lambda$
with $\psi^{-1}(W_0)\neq \emptyset$ one has $ \psi^*\sC_{W}= \rho^*\sC_{\hat{Y}}$.
\end{theorem}
\begin{proof}
This is more or less what is shown in \cite[Theorem 5.1]{Vie}. There
we constructed a compactification $\overline{W}$ of $W_0$ and the sheaf $\sC_{\overline{W}}$.
Of course, we may replace $\overline{W}$ by a modification of $W$, and by Lemma~\ref{ex.4}
we can embed $\overline{W}$ in a modification of $\BP'$, finite over a modification of $\BP$. 
\end{proof}
Applying Theorem~\ref{ex.6} to prove Theorem~\ref{in.1} we will start with
$\widetilde{Y}_0=Y_0$ and with any compactification $\widetilde{Y}$.
Theorem~\ref{ex.6} will force us to choose for $Y$ a modification of $\widetilde{Y}$ 
in the diagram (\ref{eqws.3}). This one is a closed subscheme of $\BP$. 
If we have to reapply Theorem~\ref{ex.6}, as it will happen
in the proof of Theorem~\ref{in.2} iv), we may have to replace
$Y$ by some modification. Since it is obtained by blowing up $\BP$ we can choose $W$ as the 
preimage in the corresponding normal modification of $\BP'$.

The statement of~\ref{ex.6} is compatible with further blowing ups of $\hat{Y}$. So by abuse of notations, we may assume that there is a morphism $\hat{Y}\to \widetilde{Y}$, as required in the diagram (\ref{eqws.3}) in case $\widetilde{Y}=Y$.
We will denote the induced morphisms by 
\begin{equation}\label{eqex.1}
\xymatrix{\hat{Y}_0 \ar[r]^{\iota'}\ar[d]_{\xi_0}&\hat{Y}\ar[d]_{\xi} \\
W_0 \ar[r]^{\iota}\ar[d]_{\widetilde{\phi}_0}& W \ar[r]^{\subset}\ar[d]_{\widetilde{\phi}}&\BP'\ar[d]^{\Psi}\\
\widetilde{Y}_0 \ar[r]^{\subset}&\widetilde{Y}\ar[r]^{\subset}&\BP}  
\end{equation}
and $\varphi=\widetilde{\phi}\circ\xi$. If $\widetilde{Y}_0=Y_0$ we will drop all the $\widetilde{ \ \ }$.
\begin{addendum}\label{ex.7} \ 
\begin{enumerate}
\item If we consider a finite set of sheaves $\sC_\bullet$, we can choose the
same compactification $W$ for all of them. Assume for example that $\sC_\bullet$ and
$\sC'_\bullet$ are two systems of locally free sheaves satisfying the conditions i) and ii)
in~\ref{ex.6}. Then one may choose $W$ such that both locally free sheaves, $\sC_W$ and $\sC'_W$ exist, as well as the morphism $\xi$ in \eqref{eqex.1}.
\item Let $\sR$ be the sheaf of $\sO_W$ algebras $\sR=\xi_*\sO_{\hat{Y}}\cap \iota_*\sO_{W_0}$. 
The scheme ${\rm\bf Spec}(\sR)$ is finite and birational over $W$, the inclusion $\iota$
lifts to an open embedding of $W_0$ in ${\rm\bf Spec}(\sR)$. Lemma~\ref{ex.4} allows to replace $W$ by this covering, hence we may assume that $\xi_*\sO_{\hat{Y}}\cap \iota_*\sO_{W_0}=\sO_W$.
\item If in (1) one has morphisms $\iota: \sC'_{\widetilde{Y}_0}\to \sC_{\widetilde{Y}_0}$ and $\iota':\sC'_{\hat{Y}}\to\sC_{\hat{Y}}$, compatible with the pullback in~\ref{ex.6} i) one has a natural map
$$
\sC'_W \>>> \xi_*\xi^*\sC'_W =\xi_*\sC'_{\hat{Y}} \>\iota'>> \xi_*\sC_{\hat{Y}}=\sC_W\otimes \xi_* \sO_{\hat{Y}}.
$$  
So $\sC'_W$ maps to $\sC_W \otimes \sR$, for the coherent sheaf $\sR$ considered in (2).
Replacing $W$ by ${\rm\bf Spec}(\sR)$ one obtains $\iota'':\sC'_W \to \sC_W$, and this morphism is compatible with all further pullbacks.
\end{enumerate}
\end{addendum}
Let us state three slightly different applications of Theorem~\ref{ex.6} which will be proved
together after Variant \ref{ex.11}.
\begin{proposition}\label{ex.8} 
One may choose $Y$, $\hat{Y}$ and $\hat{Z}$ in~\ref{sa.5} and~\ref{sa.7} (or~\ref{sa.13} in the saturated case) such that in addition to the conditions i)--x) one has a diagram (\ref{eqex.1}) with $\widetilde{Y}=Y$  and $\widetilde{\phi}=\phi$ such that:
\begin{enumerate}
\item[I.] $\Psi$ is a finite covering, $\BP$ and $\BP'$ are normal and projective, $\Psi^{-1}(Y)=W$, and $\xi$ is birational.
\item[II.] Let $C$ be a smooth curve and $\chi:C\to Y$ a morphism. Assume that
$\chi$ factors through $C \>\pi >> W \> \phi >> Y$,
and that $C_0=\chi^{-1}(Y_0)$ is dense in $C$. Then $\chi$ admits a mild reduction.
\item[III.] For $(\nu,\mu)\in I$ there exists a locally free sheaf $\sF_W^{(\nu,\mu)}$ on $W$ with
$\xi^*\sF_W^{(\nu,\mu)}=\sF_{\hat{Y}}^{(\nu,\mu)}$, and such that $\sF_W^{(\nu,\mu)}|_{W_0}=\phi_0^*f_*(\omega_{X/Y}^\nu\otimes \sL^\mu)$.
\item[IV.] For all curves considered in II one has $\pi^*\sF_W^{(\nu,\mu)}=\sF_{C}^{(\nu,\mu)}$. 
\end{enumerate}
\end{proposition}
Assume for a moment, that a coarse moduli scheme $M_h$ exist for families of polarized manifolds with Hilbert polynomial $h$, and that the family $f_0:X_0\to Y_0$ lies in $\fM_h(Y_0)$
for the corresponding moduli functor. Assume the induced morphism $Y_0\to M_h$ is finite. Then we want to factor $\hat{Y}_0\to \widetilde{Y}_0=M_h$ through some $W_0$, birational to $\hat{Y}_0$, and finite over $M_h$ with a splitting trace map. In general the different direct image sheaves do not descend to the moduli schemes, just their determinants. So we only can expect that certain ``natural'' invertible sheaves descend to the compactification of the moduli scheme. In the canonically polarized case, those sheaves will be of the form $\det(\sF_{\hat{Y}}^{(\nu)})$. If one allows arbitrary polarizations, one has to choose certain rigidifications.
\begin{definition}\label{ex.9}
Let $\iota$ and $\iota'$ be integers. We call the sheaf 
$$
\det(\sF_\bullet^{(\nu,\mu)})^\iota \otimes \det(\sF_\bullet^{(\nu',\mu')})^{\iota'}
$$ 
a {\em rigidified determinant sheaf}, if
$\iota\cdot\mu\cdot\rk(\sF_\bullet^{(\nu,\mu)})+\iota'\cdot\mu'\cdot\rk(\sF_\bullet^{(\nu',\mu')})=0$.
\end{definition}
Recall that for the moduli problem of polarized manifolds one does not distinguish between families 
$$
(f_0:X_0\to Y_0,\sL) \mbox{ \ \  and \ \ }(f_0:X_0\to Y_0,\sL\otimes f_0^*\sN),
$$ 
where $\sN$ is an invertible sheaf on $\hat{Y}$. The Definition~\ref{ex.9} is made in such a way, that the
rigidified determinant sheaves are invariant under this relation.
It follows from the construction of moduli schemes that some power of a rigidified determinant sheaf
descends to $M_h$ (see \cite[Proposition 7.9]{Vie}, for example). 
\begin{variant}\label{ex.10} Assume that $Y_0$ is normal and that 
the family $f_0:X_0\to Y_0$ (or $(f_0:X_0\to Y_0,\sL_0)$) induces a finite morphism $Y_0\to M_h$. Then one can find for a compactification $Y$ of $Y_0$, the schemes $\hat{Y}$ and $\hat{Z}$ in~\ref{sa.5} and~\ref{sa.7} (or~\ref{sa.13} in the saturated case), such that in addition to the conditions i)--x) one has
for $\widetilde{Y}_0=M_h$ the diagram (\ref{eqex.1}) and: 
\begin{enumerate}
\item[I.] $\Psi$ is a finite covering, $\BP$ and $\BP'$ are normal and projective, $\Psi^{-1}(\overline{M}_h)=W$, and $\xi$ is birational.
\item[II.] Let $C$ be a smooth curve and $\chi:C\to \overline{M_h}$ a morphism. Assume that
$\chi$ factors through $C \>\pi >> W \> \widetilde{\phi} >> \overline{M}_h$,
and that $C_0=\chi^{-1}(M_h)$ is dense in $C$. Then the induced morphism
$C\to Y$ admits a mild reduction.
\item[III.] For $(\nu,\mu), \ (\nu',\mu') \in I$ and $\iota, \ \iota'\in \Z$ let
$\det(\sF_\bullet^{(\nu,\mu)})^\iota \otimes \det(\sF_\bullet^{(\nu',\mu')})^{\iota'}$ 
be a rigidified determinant. Then there exists some $p\gg 1$ and an invertible sheaf $\sC_{\overline{M}_h}$ on $\overline{M}_h$ with 
$$
\sC_{\hat{Y}}:=(\det(\sF_{\hat{Y}}^{(\nu,\mu)})^\iota \otimes \det(\sF_{\hat{Y}}^{(\nu',\mu')})^{\iota'} )^p=\xi^*\widetilde{\phi}^*\sC_{\overline{M}_h}.
$$
\item[IV.] Under the assumption made in III, for all curves as in II 
$$
\sC_C:= (\det(\sF_{C}^{(\nu,\mu)})^\iota \otimes \det(\sF_{C}^{(\nu',\mu')})^{\iota'} )^p=\pi^*\widetilde{\phi}^*\sC_{\overline{M}_h}.
$$
\end{enumerate}
\end{variant}

\begin{variant}\label{ex.11}
Assume again that $\widetilde{Y}=Y$ and that the assumptions made in~\ref{sa.5} and~\ref{sa.7} (or~\ref{sa.13} in the saturated case) hold true, as well as those made in~\ref{de.1}. 

Assume there exist for $(\nu,\mu)\in I$
locally free sheaves $\sF_{Y}^{(\nu,\mu)}$ on $Y$ whose pullback to $\hat{Y}$ coincides
with $\sF_{\hat{Y}}^{(\nu,\mu)}$ and whose restriction to $Y_0$ is $f_*(\omega_{X/Y}^\nu\otimes \sL^\mu)$. 
Assume moreover, that there is a locally free sheaf $\sE_Y$ and a morphism $\sE_{Y}\to \sF_{Y}^{(\beta_0,\alpha_0)}$ satisfying the assumptions~\ref{de.3} or~\ref{de.4}.

Then, replacing $Y$ by a modification with center in $Y\setminus Y_0$, 
one can find $\hat{Y}$ and $\hat{Z}$ such that~\ref{sa.5},~\ref{sa.7} and~\ref{de.6} hold,
and such that one has a diagram (\ref{eqex.1}) with:
\begin{enumerate}
\item[I.] $\Psi$ is a finite covering, $\BP$ and $\BP'$ are normal and projective, $\Psi^{-1}(Y)=W$, and $\xi$ is birational.
\item[II.] Let $C$ be a smooth curve and $\chi:C\to Y$ a morphism. Assume that
$\chi$ factors through $C \>\pi >> W \> \phi >> Y$, 
and that $C_0=\chi^{-1}(Y_0)$ is dense in $C$. Then $\chi$ admits a mild exhausting reduction for $(\Xi^{(r)},\sE;\beta+\frac{\eta}{\ell},\alpha+\frac{\gamma}{\ell})$.
\item[III.] There exists a locally free sheaf $\sG_{W}^{(\beta+\frac{\eta}{\ell},\alpha+\frac{\gamma}{\ell})}$
on $W$ whose pullback to $\hat{Y}$ is the sheaf $\sG_{\hat{Y}}^{(\beta+\frac{\eta}{\ell},\alpha+\frac{\gamma}{\ell})}$,
defined in~\ref{de.5}. One has an inclusion 
$$
\sG_{W}^{(\beta+\frac{\eta}{\ell},\alpha+\frac{\gamma}{\ell})}\> \subset >>
\bigotimes^r\sF_{W}^{(\beta+\frac{\eta}{\ell},\alpha+\frac{\gamma}{\ell})}, 
$$
and over $W_0$ both sheaves are isomorphic.
\end{enumerate}
\end{variant}

\begin{proof}[Proof of~\ref{ex.8},~\ref{ex.10} and~\ref{ex.11}]
The verification of the properties I and II in each of the cases goes along the same line.

In~\ref{ex.8} and~\ref{ex.11} one starts with $Y=\widetilde{Y}$, where $X_0\to Y_0$ extends to a flat
morphism $f:X\to Y$, as required in step I of~\ref{ws.4} or in variant~\ref{ws.8}. 
We choose an embedding $Y\to \BP=\BP^M$. 

In~\ref{ex.10} we start with an embedding $M_h\to \BP=\BP^M$. In order to be able to use induction on certain strata, we will allow  $\widetilde{Y}_0$ to be a subscheme of $M_h$, and 
we choose $\widetilde{Y}$ as the closure of $\widetilde{Y}$ in $\BP$. Correspondingly
we will replace $Y_0$ by the preimage of $\widetilde{Y}$, and $Y$ will be
a compactification $Y_0$ such that $\tau_0:Y_0\to \widetilde{Y}_0$ extends to 
a morphism $\tau:Y\to \widetilde{Y}$, and such that $f_0:X_0\to Y_0$ extends to a flat projective morphism $f:X\to Y$. 

In all cases our starting point are morphisms 
$$
Y \> \tau >> \widetilde{Y} \> \subset >> \BP \mbox{ \ \ and \ \ } f:X\>>> Y,
$$
and $f_0:X_0\to Y_0$ is the smooth part of $f$. So in \ref{ex.8} and \ref{ex.10}  
the data $\widetilde{I}$ and $I$ allow by \ref{sa.5},~\ref{sa.7} or~\ref{sa.13}
to choose a diagram as in (\ref{eqws.3}):
$$
\xymatrix{
&X \ar[dr]_f & \ar[l]_{\hat{\varphi}}  \hat{Z}\ar[dr]_{\hat{g}} & \ar[l]_{\hat{\delta}} Z \ar[r]^\delta \ar[d]^g & \hat{X} \ar[dl]^{\hat{f}} \ar[r]^{\rho} & X \ar[dl]^{f} \\
\BP&\ar[l]_{\subset}\widetilde{Y} &\ar[l]_{\tau} Y & \ar[l]^{\varphi} \hat{Y} \ar[r]_{\varphi} & Y \ar[r]^{\tau}&\widetilde{Y}\ar[r]^{\subset}&\BP.}
$$
In Variant~\ref{ex.11} we use the Lemma \ref{de.6} to get the same diagram, starting with
the data collected in~\ref{de.1}--\ref{de.4}. 
Recall that all those conditions are compatible with pullback under alterations of $\hat{Y}$.

Consider the Stein factorization $\widetilde{\eta} : \widetilde {V} \to \widetilde{Y}$ of
$\varphi:\hat{Y}\to \widetilde{Y}$. By~\ref{ex.3} we can find an irreducible normal covering
$\Psi:\BP'\to \BP$ dominating $\widetilde {V} \to \widetilde{Y}$. So the normalization of
$W:=\Psi^{-1}(\widetilde{Y})$ dominates $\widetilde {V}$. 
The compatibility of our constructions with further pullback, allows to assume that $\hat{Y}$ is a non-singular modification of $V$, and we obtain all the morphisms in \eqref{eqex.1}, except that they are not yet coming by an application of the extension theorem. In the course of the verification of II we will have to replace $\BP'$ by finite coverings, and by some modification with center in $W\setminus W_0$. The Lemma~\ref{ex.4} allows to replace $Y$ by a modification with center in $Y\setminus Y_0$, and to keep the conditions in I.

By Lemma~\ref{mi.4} there exists an open dense subscheme $Y_g\subset Y_0$ such that $\chi:C\to Y$ admits a mild reduction if $\chi^{-1}(Y_g)\neq \emptyset$ and if $\chi$ lifts to a morphism $C\to \hat{Y}$. In~\ref{ex.11} we assumed that the sheaf $\sE_{\hat{Y}}$ is the pullback of a sheaf on $Y$. So
as remarked in \ref{mi.2} this allows to apply~\ref{va.6}, and the same holds
for mild exhausting reductions.

Replacing $Y_g$ by some open dense subscheme, we may assume in addition that:
\begin{enumerate}
\item In~\ref{ex.10} one has $Y_g=\tau^{-1}(\widetilde{Y}_g)$ for some open dense subscheme $\widetilde{Y}_g$ of $\widetilde{Y}$.
\item $W_g=\widetilde{\phi}^{-1}(\widetilde{Y}_g)$ is normal and the restriction of $\xi$ to
$\hat{Y}_g=\xi^{-1}(W_g)$ is an isomorphism $\hat{Y}_g\to W_g$.
\end{enumerate}
(2) implies that a morphism $\pi:C\to W$ from a non-singular projective curve $C$ whose image meets $W_g$ lifts to a morphism $C\to \hat{Y}$. So the conditions II in \ref{ex.8},~\ref{ex.10} or~\ref{ex.11} hold for morphisms $\pi:C\to W$ with $\pi^{-1}(W_g)$ dense in $C$.

The open set $W_g$ will be the large stratum, and next we will construct a similar
open subset of the complement. Let us write $\widetilde{Y}_b$ for the closure of $\widetilde{Y}_{0b}=\widetilde{Y}_0\setminus \widetilde{Y}_g$ in $\widetilde{Y}$. Correspondingly $Y_b$ will be equal to $\widetilde{Y}_b$ in 
\ref{ex.8} and~\ref{ex.11}, and equal to $\tau^{-1}(\widetilde{Y}_b)$ in~\ref{ex.10}.

The dimension of $\widetilde{Y}_b$ is strictly smaller than $\dim(\widetilde{Y})$.
By induction on the dimension we assume that we have found a non-singular alteration $\hat{Y}_b\to Y_b$ and the covering $\Psi_b:\BP'_b \to \BP$, satisfying the conditions i)--v) in~\ref{sa.5} and vi)--x) in~\ref{sa.7} (or~\ref{sa.13}) and the assumptions made in~\ref{de.1}
and~\ref{de.3} or~\ref{de.4}, such that the conditions II hold for $\widetilde{Y}_b$ instead of
$\widetilde{Y}$. 

We choose $\BP'_1$ to be one of the irreducible components of the normalization of $\BP'\times_{\BP}\BP'^{(b)}$. Writing $\Psi_1:\BP'_1\to \BP$ for the induced map, we choose 
$\hat{Y}_1$ to be a desingularization of $W_1=\Psi_1^{-1}(\widetilde{Y})$, which maps to $\hat{Y}$. So all the conditions needed in~\ref{sa.5},~\ref{sa.7},~\ref{sa.13} and stated in \ref{de.1},~\ref{de.3} and~\ref{de.4} remain true.

Let $\chi:C\to \widetilde{Y}$ be a morphism with $\chi^{-1}(\widetilde{Y}_0)\neq\emptyset$,
and factoring through $W_1$. If $\chi^{-1}(\widetilde{Y}_g) \neq \emptyset$, we are done. Otherwise
$\chi(C_0)$ is contained in $\widetilde{Y}_b$. By the choice of $\BP'_1$, the morphism $\chi$
factors through $\BP'_b$, hence $C\to Y$ allows again a mild (exhausting) reduction.

So in each of the three cases considered, we found a non-singular alteration satisfying I and II. 
Dropping as usual the lower index ${}_1$ we will use the notations from the diagram \eqref{eqex.1}.

The conditions III and IV will follow from the Extension Theorem, so again we will have to
modify all the morphisms in \eqref{eqex.1}. In order to apply it, we have to
define the sheaves $\sC_C$ in the Set-up~\ref{ex.5} and to verify the properties
i) and ii) stated there. This will be done in each case separately.\vspace{.1cm}

Let us start with~\ref{ex.8}. Recall that by~\ref{sa.5} and~\ref{sa.7} on $\widetilde{Y}_0=Y_0$ the sheaves $\sC_{Y_0}=f_{0*}(\omega_{X_0/Y_0}^\nu\otimes \sL_0^\mu)$ are locally free and compatible with base change for $(\nu,\mu)\in I$. Correspondingly we choose $\sC_{\hat{Y}}=\sF_{\hat{Y}}^{(\nu,\mu)}$, and
$\sC_{C}=\sF_{C}^{(\nu,\mu)}$, as defined in~\ref{mi.2}. Then i) is obviously true, and ii) follows from II, using Proposition~\ref{mi.5}.\vspace{.1cm}  

The same argument works for~\ref{ex.10}. However here we have to choose for $\sC_{Y_0}$ the rigidified determinant
$$
(\det(f_{0*}(\omega_{X_0/Y_0}^\nu\otimes \sL_0^\mu))^\iota \otimes \det(f_{0*}(\omega_{X_0/Y_0}^{\nu'}\otimes \sL_0^{\mu')})^{\iota'})^p.
$$
As mentioned already, by \cite[Proposition 7.9]{Vie} for $p$ sufficiently large, this sheaf is the pullback of an invertible sheaf $\sC_{M_h}$. Then for $\sC_{\hat{Y}}$ and $\sC_C$, as defined in~\ref{ex.10} III and IV, the property i) follows from the compatibility
of $\sC_{Y_0}$ with pullback, and ii) follows again from II, using Proposition~\ref{mi.5}. So
the Extension Theorem~\ref{ex.6} gives the existence of the sheaf $\sC_W$. It remains to show, that
$\sC_W$, or some tensor power of $\sC_W$ descends to $\overline{M}_h$. 

To this aim, we can replace $\BP'$ by a finite covering, and assume that $\C(\BP')$ is Galois over $\C(\BP)$.
So the Galois group $\Gamma$ acts on $W$ and the quotient is $\overline{M}_h$. For $\sigma\in \Gamma$ one has
$\sigma^*\sC_W=\sC_W$. In fact, this holds true on the open dense subscheme $W_0$, and on every curve mapping to $W$ and meeting $W_0$. Replacing $p$ by some multiple, one finds the sheaf $\sC_{\overline{M}_h}$.\vspace{.1cm} 

In~\ref{ex.11} we start with 
$$
\sC_{Y_0}=\bigotimes^r f_{0*}(\omega_{X_0/Y_0}^{\beta+\frac{\eta}{\ell}}\otimes \sL_0^{\alpha+\frac{\gamma}{\ell}})
$$
and with $\sC_{\hat{Y}}=\sG_{\hat{Y}}^{(\beta+\frac{\eta}{\ell},\alpha+\frac{\gamma}{\ell})}$.
Again, those sheaves are compatible with pullback, and i) follows from
Lemma~\ref{de.11}. Since $\sE_{\hat{Y}}$ is the pullback of a sheaf on $Y$, we are allowed 
to use the constructions in Section~\ref{va}. We choose for $\sC_C$ the sheaf $\sG_{C}^{(\beta+\frac{\eta}{\ell},\alpha+\frac{\gamma}{\ell})}$, defined just before Lemma 
\ref{va.3}. The condition ii) in the set-up~\ref{ex.5} follows again from II and from~\ref{va.7}.
So the Extension Theorem gives the existence of the locally free sheaf $\sG_{W}^{(\beta+\frac{\eta}{\ell},\alpha+\frac{\gamma}{\ell})}$
and as remarked in~\ref{ex.7} (5) we can assume that it is a subsheaf of $\bigotimes^r\sF_{W}^{(\beta+\frac{\eta}{\ell},\alpha+\frac{\gamma}{\ell})}$. By~\ref{de.11} the pullback of 
both to $\hat{Y}_0$ are equal, hence their restrictions to $W_0$ as well.
\end{proof}
Let us formulate what we obtained up to now for the sheaves $\sF_\bullet^{(\nu,\mu)}$.
\begin{theorem}\label{ex.12}
Let $f:X\to Y$ be a flat projective morphism of quasi-projective reduced schemes,
and let $\sL$ be an invertible sheaf on $X$. Let $Y_0\subset Y$ be a dense open set,
with $f_0:X_0=f^{-1}(Y_0)\to Y_0$ smooth. Assume that $\omega_{X_0/Y_0}$ and 
$\sL_0=\sL|_{X_0}$ are both $f_0$ semiample. 

Let $I$ be a finite set of tuples $(\nu,\mu)$ of natural numbers. Assume that
for all $(0,\mu')\in I$ the direct image $f_{0*}\sL_0^{\mu'}$ is locally free and compatible with arbitrary base change. Then, replacing $Y$ by a modification with centers in $Y\setminus Y_0$, there exists a finite covering $\phi:W\to Y$ with a splitting trace map and for $(\nu,\mu) \in I$ a locally free sheaf $\sF^{(\nu,\mu)}_{W}$ on $W$ with:
\begin{enumerate}
\item[i.] For $W_0=\phi^{-1}(Y_0)$ and $\phi_0=\phi|_{W_0}$ one has 
$\phi_0^*f_{0*}(\omega_{X_0/Y_0}^\nu\otimes \sL_0^\mu)=\sF^{(\nu,\mu)}_{W}|_{W_0}$.
\item[ii.] Let $\theta:T\to W$ be a morphism from a non-singular variety $T$.
Assume that either $T\to W$ is dominant or that $T$ is a curve and $T_0=\theta^{-1}(W_0)$ dense in $T$. For some $r\geq 1$ let $X^{(r)}_T$ be a desingularization of
$$
X^r_T=(X\times_{Y}\cdots \times_{Y}X)\times_{Y}T.
$$
Let $\hat{\varphi}_T:X^{(r)}_T\to X^r_T$ and $f^{(r)}_T:X^{(r)}_T\to T$ be the induced morphisms and
$$
\sM=\hat{\varphi}_T^*({\rm pr}_1^*\sL\otimes \cdots \otimes {\rm pr}_r^*\sL)
$$
Then $\displaystyle f^{(r)}_{T*}(\omega_{X^{(r)}_T/T}^\nu\otimes \sM^\mu)
=\bigotimes^r\theta^*\sF^{(\nu,\mu)}_W$.
\end{enumerate}
\end{theorem}
For $\mu=0$ one obtains in particular parts i) and ii) of Theorem~\ref{in.1}, and
we have seen in Section~\ref{pd1} already that those two conditions imply iii) in Theorem~\ref{in.2}, saying that the sheaf $\sF_W^{(\nu,0)}=\sF^{(\nu)}_{W}$ is nef. So it remains to prove the 
``weak stability'' condition iv). This will be done in Section~\ref{pd2}.
Let us formulate first a variant of the last Theorem allowing saturated extensions of polarizations.
\begin{variant}\label{ex.13}
In~\ref{ex.12} fix some $\eta_0$ such that the evaluation map for $\omega_{X_0/Y_0}^{\eta_0}$ is surjective, and some $\kappa>0$, with $(\eta_0,0), \ (0,\kappa)\in I$. Then there exists a finite covering $\phi:W\to Y$ with a splitting trace map, and for $(\nu,\mu) \in I$
a locally free sheaf $\sF^{(\nu,\mu)}_{W}$ on $W$ with the property i) and 
\begin{enumerate}
\item[ii.] Let $\theta:T\to W$ be a morphism from a non-singular variety $T$.
Assume that either $T\to W$ is dominant or that $T$ is a curve and $T_0=\theta^{-1}(W_0)$ dense in $T$. For some $r\geq 1$ let $X^{(r)}_T$ be a desingularization of
$(X\times_{Y}\cdots \times_{Y}X)\times_{Y}T$.
Let $\hat{\varphi}_T:X^{(r)}_T\to X^r$ and $f^{(r)}_T:X^{(r)}_T\to T$ be the induced morphisms.
Assume that $X^{(r)}_T$ is chosen such that the image of the evaluation map for $\omega^{\eta_0}_{X^{(r)}_T/T}$ is invertible, hence equal to $\omega^{\eta_0}_{X^{(r)}_T/T}\otimes\sO_{X^{(r)}_T}(\Pi_{X^{(r)}_T})$ for an effective Cartier divisor $\Pi_{X^{(r)}_T}$. 
Then for $\sM=\hat{\varphi}_T^*({\rm pr}_1^*\sL\otimes \cdots \otimes {\rm pr}_r^*\sL)$
one has 
$$
\displaystyle f^{(r)}_{T*}(\omega_{X^{(r)}_T/T}^\nu\otimes \sM^\mu\otimes \sO_{X^{(r)}_T}(*\cdot \Pi_{X^{(r)}_T}))
=\bigotimes^r\theta^*\sF^{(\nu,\mu)}_W.
$$
\end{enumerate}
\end{variant}
\begin{proof}[Proof of~\ref{ex.12} and~\ref{ex.13}]
Start with $\hat{Y}$, $\hat{Z}$ and $\hat{X}$ according to~\ref{sa.5} and~\ref{sa.7} (or~\ref{sa.13} in~\ref{ex.13}).
Choose the compactification $Y$, and $W$ using Proposition~\ref{ex.8}. 

So there are locally free sheaves $\sF_W^{(\nu)}$ (or $\sF_W^{(\nu,\mu)}$), whose
pullbacks under $\xi$ are the sheaves $\sF_{\hat{Y}}^{(\nu)}$ (or $\sF_{\hat{Y}}^{(\nu,\mu)}$).
It remains to verify the condition ii) in all cases.

Recall that, $\hat{X}\to \hat{Y}$ has a mild model $\hat{Z}\to \hat{Y}$, hence $X^{(r)}_T\to \hat{Y}$
has $\hat{Z}^r\to \hat{Y}$ as a mild model. If $T$ dominates $\hat{Y}$ the property ii) in~\ref{ex.12} follows for $r=1$ from~\ref{sa.5} and~\ref{sa.7}, and for $r>1$ by flat base change.
In~\ref{ex.13} the same argument works for a $\kappa$ saturated extension $\sM_{X^{(r)}_T}$, and one finds that
$$ 
f^{(r)}_{T*}(\omega_{X^{(r)}_T/T}^\nu\otimes \sM_{X^{(r)}_T}^\mu)=
\bigotimes^r\theta^*\sF^{(\nu,\mu)}_W.
$$
In general there is some non-singular modification $\theta':T'\to T$ such that ii) holds on $T'$.
The sheaf $f^{(r)}_{T*}(\omega_{X^{(r)}_T/T}^\nu\otimes \sM^\mu)$ is independent of the
desingularization $X^{(r)}_T$, and we may assume that $f^{(r)}_T$ factors through
$h:X^{(r)}_T\to T'$. Then 
$$
h_*(\omega_{X^{(r)}_T/T}^\nu\otimes \sM^\mu)=
\bigotimes^r\theta'^*\theta^*\sF^{(\nu,\mu)}_W\otimes \omega_{T'/T}^\nu,
$$
and the projection formula implies that
$$
f^{(r)}_{T*}(\omega_{X^{(r)}_T/T}^\nu\otimes \sM^\mu)=
\bigotimes^r\theta^*\sF^{(\nu,\mu)}_W\otimes \theta'_* \omega_{T'/T}^\nu=
\bigotimes^r\theta^*\sF^{(\nu,\mu)}_W,
$$
as claimed in~\ref{ex.12}. In the situation considered in~\ref{ex.13}
the same equality holds with $\sM$ replaced by the $\kappa$ saturated extension
$\sM_{X^{(r)}_T}$. However both differ by some positive multiple of $\Pi_{X^{(r)}_T}$ and 
\begin{multline*}
f^{(r)}_{T*}(\omega_{X^{(r)}_T/T}^\nu\otimes \sM_{X^{(r)}_T}^\mu)=
f^{(r)}_{T*}(\omega_{X^{(r)}_T/T}^\nu\otimes \sM_{X^{(r)}_T}^\mu\otimes \sO_{X^{(r)}_T}(*\cdot \Pi_{X^{(r)}_T}))=\\
f^{(r)}_{T*}(\omega_{X^{(r)}_T/T}^\nu\otimes \sM^\mu\otimes \sO_{X^{(r)}_T}(*\cdot \Pi_{X^{(r)}_T})).
\end{multline*}
If $T$ is a curve, then by Proposition~\ref{ex.8} II we know that $T\to W \to Y$ admits
a mild reduction, and by part IV the pullback of $\sF_W^{(\nu,\mu)}$ is the sheaf
$\sF_C^{(\nu,\mu)}$ defined in Section~\ref{mi}. So it is equal to
$h_*(\omega_{S/T}^\nu\otimes \sM_S^\mu)$ for a mild model $h:S\to T$ of the pullback family.

The $r$-fold fibre product $h^r:S^r\to T$ is again mild, and for the exterior tensor product
$\sM_{S^r}$ one has by flat base change  
$\displaystyle h^r_*(\omega_{S^r/T}^\nu\otimes \sM_{S^r}^\mu)
=\bigotimes^r\theta^*\sF^{(\nu,\mu)}_W$. So the property ii) in~\ref{ex.12} or~\ref{ex.13}
for $T$ a curve follows from~\ref{di.7}.
\end{proof}

\section{Positivity of direct images II}\label{pd2}
Now we are all set to finish the proof of Theorem~\ref{in.2} iv), allowing $W_0$ to be singular, contrary to Variant~\ref{pd1.4}.
The lemma~\ref{ne.9} allows to replace $W$ by a larger covering, and we will do so several times. 
We will also formulate and prove a generalization to the polarized case, which in particular will imply the Lemma~\ref{co.2}. 

As explained in~\ref{pd1.5} we will construct to this aim a locally free subsheaf $\sG$ of $\sF_W^{(\nu,\mu)}$, isomorphic to $\sF_W^{(\nu,\mu)}$ over $W_0$, whose pullback to $\hat{Y}$ remains nef after tensoring with a ``negative'' invertible sheaf. 
The sheaf $\sG$ will depend on the data defined in Section~\ref{de}. 
\begin{setup}\label{pd2.1}
We will specify in each case the tautological maps in~\ref{de.1}, and  the sheaf
$\sE_{\hat{Y}}$ according to~\ref{de.3} or~\ref{de.4}. We choose the sets $\widetilde{I}$ 
and $I$ as in the Set-up~\ref{de.2} and enlarge them such that~\ref{de.3} or~\ref{de.4} applies. 
By abuse of notations we will assume that the alteration $\hat{Y}$ and $W$ are chosen according to Theorem
\ref{ex.12} and Variant~\ref{ex.13}, and moreover we will assume that Variant~\ref{ex.11} applies,
i.e. that the locally free subsheaf
$\sG_{W}^{(\beta+\frac{\eta}{\ell},\alpha+\frac{\gamma}{\ell})}$
of $\sF_{W}^{(\beta+\frac{\eta}{\ell},\alpha+\frac{\gamma}{\ell})}$ exist.

In all situations the locally free sheaf $\sE_{\hat{Y}}$ in~\ref{de.3} or~\ref{de.4}
will be the pullback of a locally free sheaf $\sE_W$, and the invertible sheaf $\sV$ in~\ref{de.1} will be the pullback of an invertible sheaf $\sV_W$; hence the $r\cdot e \cdot \ell$-th root out of 
$\det(\sF_{W}^{(\eta_1,\gamma_1)})^{\frac{r}{r_1}}\otimes\cdots\otimes\det(\sF_{W}^{(\eta_s,\gamma_s)})^{\frac{r}{r_s}}$
will exist on $W$. 
\end{setup}
Remark already that all those conditions can be realized, after blowing up $Y$ with centers in
$Y\setminus Y_0$ for some finite covering $W\to Y$ with a splitting trace map.
The conclusion stated in the sequel remain true over any model, where
the different sheaves are defined on $W$, locally free and compatible with pullback.
\begin{proposition}\label{pd2.2} In~\ref{pd2.1} one has:
\begin{enumerate}
\item[a.] If $\sE_{W}$ is nef, the sheaf
$\displaystyle \sG_{W}^{(\beta+\frac{\eta}{\ell},\alpha+\frac{\gamma}{\ell})}\otimes \sV_W^{-r}$
is nef and the sheaf
$\displaystyle \sF_W^{(\beta+\frac{\eta}{\ell},\alpha+\frac{\gamma}{\ell})}\otimes \sV_W^{-1}$
is weakly positive over $W_0$.
\item[b.] If for some invertible sheaf $\sH$ on $W$ the sheaf
$\sE_{W}\otimes \sH^{b\cdot e \cdot \ell}$ is nef, the sheaf
$\displaystyle\sG_{W}^{(\beta+\frac{\eta}{\ell},\alpha+\frac{\gamma}{\ell})}\otimes \sH^r \otimes \sV_W^{-r}$ is nef and the sheaf
$\displaystyle\sF_W^{(\beta+\frac{\eta}{\ell},\alpha+\frac{\gamma}{\ell})}\otimes\sH\otimes \sV_W^{-1}$
is weakly positive over $W_0$.
\end{enumerate}
\end{proposition}
\begin{proof} Writing $\sH=\sO_W$ in a) we will handle both cases at once. By~\ref{ex.11} one has an inclusion
$$
\sG_{W}^{(\beta+\frac{\eta}{\ell},\alpha+\frac{\gamma}{\ell})}\subset
\bigotimes^r \sF_{W}^{(\beta+\frac{\eta}{\ell},\alpha+\frac{\gamma}{\ell})},
$$ 
and both sheaves are isomorphic on $W_0$. Hence using the equivalence
of (1) and (2) in Lemma~\ref{ne.6}, it is sufficient to verify
that the first sheaf, tensorized by $\sH^r\otimes \sV_W^{-r}$ is nef.
By Lemma~\ref{ne.3} this follows if for $\hat{\sH}=\xi^*\sH$ the sheaf
$$
\sG_{\hat{Y}}^{(\beta+\frac{\eta}{\ell},\alpha+\frac{\gamma}{\ell})}\otimes
\hat{\sH}^r \otimes \sV^{-r}
$$ 
is nef. We work with the mild model, and we use the notations from Claim~\ref{de.7}. There we verified that the sheaf $\sN^N \otimes {\hat{g}^r}{}^*\sV^{-N\cdot r}\otimes\sO_{\hat{Z}^r}(-\Delta_{\hat{Z}^r})$
is the image of ${\hat{g}^r}{}^*\sE_{\hat{Y}}^{\otimes r}$, for $N=b\cdot e \cdot \ell$. So Lemma~\ref{pd1.3} implies for $\sN$ replaced by $\sN\otimes {\hat{g}^r}{}^*\hat{\sH}^r$ that for a very ample sheaf $\sA$ on $\hat{Y}$ the sheaf
\begin{multline*}
\omega_{\hat{Y}}\otimes \sA^{m+2}\otimes\sG_{W}^{(\beta+\frac{\eta}{\ell},\alpha+\frac{\gamma}{\ell})}\otimes\sH^r
\otimes\sV_W^{-r}=\\
\omega_{\hat{Y}}\otimes \sA^{m+2}\otimes \hat{g}^r_*\big(\omega_{\hat{Z}^r/\hat{Y}}\otimes \sN \otimes
\sJ\big(-\frac{1}{N}\cdot \Delta_{\hat{Z}^r}\big)\big)\otimes \hat{\sH}^r \otimes \sV^{-r} 
\end{multline*}
is globally generated. This remains true for $r\cdot r'$ instead of $r$. Since
$$
\sG_{W}^{(\Xi^{(r\cdot r')},\sE;\beta+\frac{\eta}{\ell},\alpha+\frac{\gamma}{\ell})}=
\bigotimes^{r'}\sG_{W}^{(\Xi^{(r)},\sE;\beta+\frac{\eta}{\ell},\alpha+\frac{\gamma}{\ell})}=
\bigotimes^{r'}\sG_{W}^{(\beta+\frac{\eta}{\ell},\alpha+\frac{\gamma}{\ell})},
$$
\ref{pd2.2} follows from the equivalence of (1) and (3) in~\ref{ne.6}.
\end{proof}
\begin{proof}[Proof of Theorem~\ref{in.2}] 
Remark that we already obtained Theorem~\ref{in.1} in Theorem~\ref{ex.12}, and we
keep the choice of $\phi:W\to Y$ we made there. The part iii) of Theorem~\ref{in.2}
has been verified in Section~\ref{pd1}.

For $r_\iota=\dim(H^0(F,\omega_F^{\eta_\iota}))$ 
choose $\Xi=(\Xi_1, \ldots , \Xi_s)$ in Proposition~\ref{pd2.4} as the tuple of tautological maps
$$
\Xi_\iota:\bigwedge^{r_\iota}H^0(F,\omega_F^{\eta_\iota})\>>> \bigotimes^{r_\iota} H^0(F,\omega_F^{\eta_\iota}).
$$
For some $\eta_0$ the evaluation map for $\omega_{X_0/Y_0}^{\eta_0}$ is surjective.
Replacing $\Xi$ by $\xi,\ldots,\xi$ we may assume that 
$\eta_0$ divides $\eta=\eta_1+\cdots +\eta_s$. We choose 
$\ell=\eta$, for $r$ we choose some positive common multiple of $r_1, \ldots , r_s$, for $e$ any integer larger that $\frac{1}{\ell} e(\omega_F^\eta)$, and for $b$ we choose any natural number with
$b\cdot(\nu-2)$ divisible by $\eta_0$. We choose $\widetilde{I}$ and $I$ such that the numerical conditions in~\ref{de.1} hold true.

So $\beta=\nu-1$, and $\beta_0=b\cdot\beta\cdot e\cdot \ell + \eta\cdot b \cdot(e-1)$. As in~\ref{de.3}
we assume that $\beta_0\in \widetilde{I}$ and for $\sE_{\hat{Y}}$ we choose $\sF_{\hat{Y}}^{(\beta_0)}$,
hence $\sE_{W}=\sF_{W}^{(\beta_0)}$ in~\ref{pd2.1}. 

Lemma~\ref{de.6} and Proposition~\ref{ex.8} allow to replace $W$ by some larger covering with
a splitting trace map, and to assume that the conditions in the Set-up~\ref{pd2.1} hold.
Doing so we are allowed to apply Proposition~\ref{pd2.2} a) and we obtain the weak positivity
of
$$
\Big(\bigotimes^\alpha\sF_W^{(\nu)}\Big)\otimes \bigotimes_{\iota=1}^s \det(\sF_W^{(\eta_\iota)})^{-\frac{r}{r_\iota}}
$$
over $W_0$ for some $\alpha>0$. We know by part iii) of Theorem~\ref{in.2} that the sheaves $\det(\sF_W^{(\eta_i)})$
are all nef. Hence we can enlarge the $a_\iota$ and assume that $a_\iota \cdot r_\iota$ is independent of $\iota$, hence that $\bigotimes_{\iota=1}^s \det(\sF_W^{(\eta_\iota)})^{\frac{r}{r_\iota}}$
is ample with respect to $W_0$.
\end{proof}
\begin{remark}\label{pd2.3}
If one wants to avoid using $s$ copies of $\Xi$, one can also argue in the following way:\\[.1cm]
$\Xi^{(r)}$ defines an embedding of some linear combination $\sN$ of the sheaves $\det(\sF_W^{(\eta_\iota)})$ in the sheaf $\bigotimes \sF_W^{(\eta)}$, and it is easy to see that this
inclusion locally splits. By part iii) of the Theorem~\ref{in.1} the second sheaf is nef. 
Then the quotient is locally free and nef, hence the determinant of $\bigotimes \sF_W^{(\eta)}$ must be ample. So we can replace in the assumptions of theorem~\ref{in.2} iv) $s$ by $1$ and $\eta_i$ by some large number $\eta$. In particular, we may assume that the evaluation map for $\omega_{X_0/Y_0}^{\eta_0}$ is surjective.
\end{remark}
\begin{proof}[Proof of Corollary~\ref{in.3}] As already stated in the introduction, a) follows from
Theorem~\ref{in.1} and~\ref{in.2}. Moreover in order to prove b) the Lemma~\ref{ne.7} (4)
allows to apply Theorem~\ref{in.1}, and to replace $\sF_{Y_0}^{(\bullet)}$ by $\sF_{W_0}^{(\bullet)}$.
In fact, we will choose $\widetilde{I}$ and $I$ as in the proof of Theorem~\ref{in.2}, given above, 
and we will choose $W$ and $\sE_W$ as we did it there. There is however a subtle point:\\[.1cm]
Even if $\bigotimes \det(\sF_{W_0}^{(\eta_i)})^{a_i}$ is ample and 
even if $\bigotimes \det(\sF_{W}^{(\eta_i)})^{a_i}$ is nef, the latter does not have to be ample with respect to $W_0$. 

We will not refere to $Y$ anymore, so we may blow up the boundary $W\setminus W_0$, and choosing the $a_i$
large enough, we may assume that for some divisor $B$ supported in $W\setminus W_0$ the sheaf 
$$
\sO_W(B)\otimes \bigotimes\det(\sF_{W}^{(\eta_i)})^{a_i}
$$
is semiample and ample with respect to $W_0$. Moreover, replacing $W$ by a finite covering with a splitting trace map, we can assume that the multiplicities of $B$ are as divisible as needed.
So applying~\ref{pd2.2} b) instead of a) one finds a divisor $B'$, still supported in $W\setminus W_0$
with 
$$
\Big(\bigotimes^\alpha\sF_W^{(\nu)}\Big)\otimes \sO_W(-B')\otimes \bigotimes_{\iota=1}^s \det(\sF_W^{(\eta_\iota)})^{-\frac{r}{r_\iota}}
$$
nef and with $\sO_W(B')\otimes \bigotimes_{\iota=1}^s \det(\sF_W^{(\eta_\iota)})^{\frac{r}{r_\iota}}$
ample with respect to $W_0$, which implies part b) of the Corollary~\ref{in.3}. 
\end{proof}
Next we will show analogs of Theorem~\ref{in.2} for the
sheaves $\sF_W^{(\nu,\mu)}$. As in the proof of Theorem~\ref{in.2} we will rely on Proposition~\ref{pd2.2}, however it will be a bit more complicated to choose the right data to start with.
\begin{proposition}\label{pd2.4}
Assume in Theorem~\ref{ex.12} or Variant~\ref{ex.13} that for some $\kappa > 0$ with $(0,\kappa) \in \widetilde{I}$ one has $\det(\sF_{W}^{(0,\kappa)})=\sO_W$. In~\ref{ex.13} assume in addition, that
the sheaves $\sM_\bullet$ are $\kappa$-saturated.

Choose some $\eta_0>0$ such that the evaluation map for $\omega^{\eta_0}_{X_0/Y_0}$ is surjective, and let $\epsilon$ be a positive multiple of $\eta_0$, with $\epsilon \geq e(\sL^{\kappa\cdot\eta_0}|_F)$ for all fibres $F$ of $f_0:X_0\to Y_0$. 
\begin{enumerate}
\item[i.] Assume that $(\epsilon\cdot\nu,\kappa\cdot\nu)$ and $(\eta_0,0)$ are in $I$. Then
the sheaf $\sF_W^{(\epsilon\cdot\nu,\kappa\cdot\nu)}$ is weakly positive over $W_0$.
\item[ii.] Assume that for some $\nu'>0$, divisible by $\eta_0$ and $\nu$
$$
((\epsilon+1)\cdot\nu,\kappa\cdot\nu), \
(\epsilon\cdot\nu,\kappa\cdot\nu), \ ((\epsilon+1)\cdot\nu',\kappa\cdot\nu'), \ (\eta_0,0) \in I.
$$
Then for some positive integer $c$ the sheaf
$$
S^c(\sF_W^{((\epsilon+1)\cdot \nu,\kappa\cdot\nu)})\otimes \det(\sF_W^{((\epsilon+1)\cdot \nu',\kappa\cdot\nu')})^{-1}
$$ 
is weakly positive over $W_0$. 
\end{enumerate}
\end{proposition}
\begin{proof}
For simplicity we will replace $\sL$ by $\sL^\kappa$ and assume that $\kappa=1$.
Choose an ample invertible sheaf $\sH$ on $W$ and define 
$$
\rho={\rm Min}\big\{\mu >0; \ \sF_{W}^{(\epsilon\cdot\nu,\nu)}\otimes \sH^{\epsilon\cdot\nu\cdot \mu-1}\mbox{ weakly positive over }
W_0\big\}.
$$
\begin{claim}\label{pd2.5}
The sheaf $\sF_{W}^{(\epsilon\cdot\nu,\nu)}\otimes \sH^{a}$ is weakly positive over $W_0$ for
$a=\nu\cdot\rho\cdot(\epsilon - \frac{\ell}{\nu})$.
\end{claim}
Part i) follows directly from~\ref{pd2.5}. In fact, by the choice of $\rho$  
$$
\nu\cdot\rho\cdot (\epsilon - \frac{\ell}{\nu}) > \epsilon\cdot\nu\cdot(\rho-1),
\mbox{ \ \ or \ \ }\rho < \dfrac{\epsilon\cdot\nu}{\ell}.
$$ 
Then $\displaystyle \sF_{W}^{(\epsilon\cdot\nu,\nu)}\otimes \sH^{\dfrac{\epsilon^2\cdot\nu^2}{\ell}}$ 
is weakly positive over $W_0$. The exponent $\dfrac{\epsilon^2\cdot\nu^2}{\ell}$ is independent
of $W$ and $\sH$. So the same holds true for any ample invertible sheaf $\sH'$ on any finite covering
$W'$ of $W$, and the weak positivity of $\sF_{W}^{(\epsilon\cdot\nu,\nu)}$ over $W_0$ follows from
\ref{ne.6}.
\begin{proof}[Proof of Claim~\ref{pd2.5}.] In the proof we will blow up $W$ with centers in $W\setminus W_0$,
so we will not use the ampleness of $\sH$, just the condition that
$\sF_{W}^{(\epsilon\cdot\nu,\nu)}\otimes \sH^{\epsilon\cdot\nu\cdot \rho}$
is ample with respect to $W_0$.

For ${r'}=\rk(\sF_{W}^{(0,1)})$ one has the natural locally splitting inclusion
$$
\sO_W=\det(\sF_{W}^{(0,1)}) \>>> \bigotimes^{{r'}} \sF_{W}^{(0,1)},
$$ 
whose pullback to $\hat{Y}$ is 
$\displaystyle
\Xi_1: \sO_{\hat{Y}}=\det(\hat{g}_*\sM_{\hat{Z}}) \>>> \bigotimes^{{r'}} \hat{g}_*\sM_{\hat{Z}}.$ 

Choose in~\ref{de.1} $\ell=\eta_0$ and for $\Xi$ the tuple consisting of $\ell$ copies of $\Xi_1$.
Hence 
$$
\gamma_1=\cdots =\gamma_{\ell}=1, \ \ \gamma=\ell \mbox{ \ \ and \ \ }\eta_1=\cdots = \eta_{\ell}=\eta=0.
$$
By assumption $\ell\cdot e=\epsilon \geq e(\sL^{\gamma}|_F)$,
as required in~\ref{de.1}. We choose $\beta=\epsilon\cdot\nu=e\cdot\ell\cdot\nu$ and $\alpha=\nu-1$, and for $b'$ any positive integer satisfying $b'\cdot(\beta-1,\alpha)\in \ell\cdot \N\times \N$. We may assume that $\nu$ and $\ell=\eta_0$ divide $b'$.

By the choice of $\rho$ the sheaf 
$$
S^{b'\cdot \epsilon - \frac{b'}{\nu}\cdot \ell}(\sF_{W}^{(\epsilon\cdot\nu,\nu)})\otimes
S^{\frac{b'}{\ell}\cdot\epsilon\cdot(\ell-1)}(\sF_{W}^{(\eta_0)})\otimes \sH^{\epsilon\cdot\nu\cdot\rho\cdot
(b'\cdot \epsilon - \frac{b'}{\nu}\cdot \ell)}
$$
is ample with respect to $W_0$. We can find some $d \gg 1$, a very ample sheaf $\sA$ on $W$ and a morphism 
$$
\bigoplus \sA \>>>
S^d\Big(S^{b'\cdot \epsilon - \frac{b'}{\nu}\cdot \ell}(\sF_{W}^{(\epsilon\cdot\nu,\nu)})\otimes
S^{\frac{b'}{\ell}\cdot\epsilon\cdot(\ell-1)}(\sF_{W}^{(\eta_0)})\otimes \sH^{\epsilon\cdot\nu\cdot\rho\cdot
(b'\cdot \epsilon - \frac{b'}{\nu}\cdot \ell)}\Big)
$$
surjective over $W_0$. Blowing up $W$ with centers in $W\setminus W_0$ we can assume that
the image of this map is locally free, hence nef. We write this image as $\sE_{W}\otimes\sH^{\epsilon\cdot d\cdot b' \cdot a}$,
and its pullback to $\hat{Y}$ as $\sE_{\hat{Y}}\otimes\tau^*\sH^{\epsilon\cdot d\cdot b' \cdot a}$.
Let us choose $b=d\cdot b'$. Multiplication of sections gives a map to 
$\sF_{\hat{Y}}^{(\beta_0,\alpha_0)}\otimes \tau^* \sH^{\epsilon \cdot b \cdot a}$ 
for 
$$
\beta_0=b \cdot \epsilon^2 \cdot \nu  - b \cdot \ell\cdot \epsilon + b\cdot\epsilon\cdot(\ell-1)
\mbox{ \ \ and \ \ } \alpha_0= b\cdot \epsilon \cdot \nu - b \cdot \ell.
$$
Since $\epsilon=e\cdot\ell$, $\beta=\epsilon\cdot\nu$ and $\alpha=\nu-1$ one has
$$
\beta_0= b \cdot (\beta-1) \cdot e \cdot \ell \mbox{ \ \ and \ \ }
\alpha_0 = b\cdot \alpha \cdot e\cdot \ell + \ell \cdot b \cdot (e - 1).  
$$
Since $\eta=0$ and $\gamma=\ell$ this is just what we required in~\ref{de.1},
and for a suitable choice of $I$ the assumptions in~\ref{de.1} and~\ref{de.4} hold true.

Since the sheaf $\sE_{\hat{Y}}$ is the pullback of a locally free sheaf $\sE_W$ on $W$ we can use~\ref{ex.11} for $W$ instead of $Y$, and obtain $\hat{Y}_1\to W_1$ and a finite covering $\tau:W_1 \to W$ with a splitting trace map, such that the sheaf 
$\sG_{W_1}^{(\beta+\frac{\eta}{\ell},\alpha+\frac{\gamma}{\ell})}=\sG_{W_1}^{(\epsilon\cdot\nu,\nu)}$
exists on $W_1$. The conditions in the Set-up~\ref{pd2.1} hold on $W_1$, and for $\sH_1=\tau^*\sH$
the sheaf $\sE_{W_1}\otimes \sH_1^{e\cdot \ell \cdot b \cdot a}$ is globally generated, hence nef.
Proposition~\ref{pd2.2} b) implies that $\sF_{W_1}^{(\epsilon\cdot\nu,\nu)}\otimes \sH_1^{a}$ is weakly positive over $\tau^{-1}(W_0)$. By~\ref{ne.6} the sheaf
$\sF_{W}^{(\epsilon\cdot\nu,\nu)}\otimes \sH^{a}$ is weakly positive over $W_0$. 
\end{proof}
So we finished the proof of part one and we can use in ii) 
that the sheaf $\sF_{W}^{(\epsilon\cdot\nu,\nu)}$ is weakly positive over $W_0$. In particular in the first part we can choose $\rho=1$ and $\sF_{W}^{(\epsilon\cdot\nu,\nu)}\otimes \sH^{\epsilon\cdot\nu}$
is ample with respect to $W_0$. In the proof of Claim~\ref{pd2.5} we obtains a bit more. 
\begin{addendum}\label{pd2.6} Under the assumptions made in~\ref{pd2.4}, there exists a projective morphism $\tau:W_1\to W$ such that its restriction $\tau^{-1}(W_0)\to W_0$ is finite with a splitting trace map, and there exists an inclusion 
$$
\sG_{W_1}=\sG_{W_1}^{(\epsilon\cdot\nu,\nu)}\subset
\bigotimes^{\rk(\sF_{W}^{(0,\kappa)})}\sF_{W_1}^{(\epsilon\cdot\nu,\nu)},
$$ 
surjective over $\tau^{-1}(W_0)$ with $\sG \otimes \tau^*(\sH)^{\nu\cdot(\epsilon-1)\cdot \rk(\sF_{W}^{(0,\kappa)})}$ nef. 
\end{addendum}
Replacing $W$ by $W_1$ we will assume that the subsheaf $\sG_W$ of $\bigotimes^{r'} \sF_W^{(\epsilon\cdot \nu,\nu)}$ exists on $W$,
for $r'=\rk(\sF_{W}^{(0,\kappa)})$. We will use~\ref{pd2.2} a second time, so we have to choose again data as in Section
\ref{de}. For $r=\rk(\sF_W^{((\epsilon+1)\cdot \nu',\nu')})$, we start with the tautological morphism
$$
\Xi: \det(\sF_W^{((\epsilon+1)\cdot \nu',\nu')})^{r'} \>>> \bigotimes^{r\cdot r'} \sF_W^{((\epsilon+1)\cdot \nu',\nu')}.
$$
So $\eta=\eta_1=(\epsilon+1)\cdot \nu'$ and $\ell=\gamma=\gamma_1=\nu'$. Necessarily one needs
$\beta=(\epsilon+1)\cdot (\nu-1)$ and $\alpha=\nu-1$. For $e$ we choose a natural number with
$\ell\cdot e \geq e(\omega_F^{(\epsilon+1)\cdot \nu'}\otimes \sL^{\nu'})$, for all fibres $F$ of $f_0$. For $b$ we choose any positive integer with 
$$
b\cdot(\beta-1,\alpha) \in \eta_0 \cdot\N \times \N,
$$ 
such that $r'\cdot \epsilon \cdot \nu$ divides $\alpha_0=b\cdot(\nu-1)\cdot e \cdot \ell + \gamma\cdot b \cdot(e-1)$. Comparing the different constants one finds 
\begin{multline*}
\beta_0=b\cdot((\epsilon+1)\cdot (\nu-1)-1)\cdot e \cdot \ell + \eta\cdot b \cdot(e-1)=\\ 
b\cdot \epsilon\cdot(\nu-1) \cdot e \cdot \ell + \epsilon \cdot \ell \cdot b \cdot (e-1)+
b\cdot\ell\cdot((\nu-1)\cdot e -1)=
\epsilon\cdot\alpha_0 + b\cdot\ell\cdot((\nu-1)\cdot e -1).
\end{multline*}
We choose 
\begin{multline*}
\hspace*{1cm}\sE_W=\Big(\bigotimes^r \sG_W^\frac{\alpha_0}{\nu}\Big)\otimes \bigotimes^{r\cdot r'}{\sF_{W}^{(\eta_0)}}^{(\frac{b\cdot\ell\cdot((\nu-1)\cdot e -1)}{\eta_0})}\\
\subset
\Big(\bigotimes^{r\cdot r'} {\sF_W^{(\epsilon\cdot\nu,\nu)}}^\frac{\alpha_0}{\nu}\Big)\otimes \bigotimes^{r\cdot r'}{\sF_{W}^{(\eta_0)}}^{(\frac{b\cdot\ell\cdot((\nu-1)\cdot e -1)}{\eta_0})}\hspace*{1cm}
\end{multline*}
and $\sE_{\hat{Y}}$ will denote its pullback to $\hat{Y}$. The $r\cdot r'$-tensor product of the multiplication map
gives
$$
\sE_{\hat{Y}} \>>> \bigotimes^{r\cdot r'}\sF_{\hat{Y}}^{(\beta_0,\alpha_0)}.
$$
Since $\sF_{W}^{(\eta_0)}$ is nef, the choice of $\sG$ in~\ref{pd2.6} implies that $\sE_{W}\otimes \sH^{\alpha_0\cdot(\epsilon-1)\cdot r'}$ is nef. Replacing $W$ by a larger covering, we may also assume that 
$\det(\sF_W^{((\epsilon+1)\cdot \nu',\nu')})$ is the $r\cdot e \cdot \ell$-th power of an invertible
sheaf $\sV_W$, and that $\sH^{\alpha_0\cdot(\epsilon-1)}$ is the
$b\cdot e \cdot \ell$-th power of an invertible sheaf. 

So all the conditions made in~\ref{pd2.1} hold, and we can apply
Proposition~\ref{pd2.2}. One obtains the weak positivity 
over $W_0$ of 
$$
\sF_W^{(\epsilon\cdot\nu,\nu)} \otimes \sH^{\frac{\alpha_0\cdot(\epsilon-1)}{b\cdot e \cdot \ell}} \otimes \sV_W^{-1}.
$$
The exponent ${\frac{\alpha_0\cdot(\epsilon-1)}{b\cdot e \cdot \ell}}$ 
is independent of $W$ and of the ample invertible sheaf $\sH$. So~\ref{ne.6}
implies that $\sF_W^{(\epsilon\cdot\nu,\nu)}\otimes \sV_W^{-1}$
is already weakly positive over $W_0$, hence
$$
S^{r\cdot e\cdot\ell}(\sF_W^{(\epsilon\cdot\nu,\nu)}) \otimes \sV_W^{-r\cdot e\cdot\ell }=
S^{r\cdot e\cdot\ell}(\sF_W^{(\epsilon\cdot\nu,\nu)}) \otimes \det(\sF_W^{((\epsilon+1)\cdot \nu',\nu')})^{-1}
$$
as well, as claimed in Proposition~\ref{pd2.4} ii).
\end{proof}
\begin{proof}[Proof of Lemma~\ref{co.2}]
We start with the models of our morphisms, constructed in Theorem~\ref{ex.12} or its Variant
\ref{ex.13}, with $\epsilon+1$ replaced by $\epsilon$. We may add the condition $\det(\sF_W^{(\kappa)})=\sO_W$.
In fact, replacing $W$ by a larger covering with a splitting trace map,
\cite[Lemma 2.1]{Vie} allows to assume that $\det(\sF_W^{(\kappa)})=\sW^{\rk(\sF_W^{(\kappa)})}$ for an invertible sheaf $\sW$.
Then one can replace the polarization on $\hat{X}\to \hat{Y}$ and on $X_0 \times_{Y_0} W_{0} \to W_{0}$ by $\sM\otimes \hat{f}^*\sW^{-1}$ and ${\rm pr}_1^*\sL_0 \otimes {\rm pr}_2^*\sW^{-1}$. 

So we may assume that the assumptions in~\ref{pd2.2} and~\ref{pd2.4} hold and the Lemma~\ref{co.2} follows from Proposition~\ref{pd2.4}, using Lemma~\ref{ne.6}.
\end{proof}
\section{The Proof of Theorems~\ref{in.5} and~\ref{in.6}}\label{pr}
In the construction of the compactifications $\overline{M}_h$ and the sheaf $\lambda_\nu^{(p)}$
we will use the negativity of the kernels of the multiplication map, stated and proved in \cite[Theorem 4.33]{Vie}.
Unfortunately there we did not keep track on what happens along the boundary, so we have to indicate 
the necessary modifications of the statements and proofs.  
\begin{theorem}\label{pr.1} Let
$W$ be a reduced projective scheme, let $W_0\subset W$ be open and dense, let $\sP$ and $\sQ$
be locally free sheaves on $W$. For a morphism
${\rm m} : S^\mu (\sP) \to \sQ$, surjective over $W_0$, assume
that the kernel of ${\rm m}$ has maximal variation in all points $w\in W_0$.

If $\sP$ is weakly positive over $W_0$ then for $b \gg a \gg 0$ the sheaf $\det (\sQ)^a \otimes \det (\sP)^b$ is ample with respect to $W_0$.
\end{theorem}
We will not recall the definition of ``maximal variation'' given in \cite[Definition 4.32]{Vie}. 
Instead we will just explain this notion in the special situation where the theorem will be applied.
\begin{example}\label{pr.2}
Assume that over $W_0$ there exists a flat family $f_0:X_0\to W_0$ and an $f_0$-ample invertible sheaf $\sL_0$ on $X_0$. Assume that $\sL$ is fibrewise very ample, and without higher cohomology.
So for all fibres $F$ one has an embedding 
$$
F\hookrightarrow \BP=\BP(H^0(F,\sL_0|_F).
$$
Choose $\beta \gg 1$ such that the homogeneous ideal of $F$ is generated in degree $\beta$, for all fibres. 
Assume that $\sP|_{W_0}=S^\beta(f_{0*}\sL_0)$, that $\sQ|_{W_0}=f_{0*}\sL_0^\beta$ and that
${\rm m}$ is the multiplication map. Then the kernel of ${\rm m}$ has maximal variation in all points $w\in W_0$ if and only if for each fibre $F$ the set
$$
\big\{ w'\in W_0; \ \mbox{ for } F'=f_0^{-1}(w') \mbox{ there is an isomorphism }
(F,\sL_0|_F)\cong (F',\sL_0|_{F'})\big\}
$$
is finite. Moreover this condition is compatible with base change under finite morphisms.
\end{example}
\begin{proof}[Sketch of the proof of~\ref{pr.1}] 
We will just recall the main steps of the proof of \cite[Theorem 4.33]{Vie}, to convince the reader
that one controls the sections along the boundary, and explain
where the condition ``maximal variation'' enters the scene.
 
Writing $r=\rk(\sP)$ we consider the projective bundle $\BP=\BP(\bigoplus^r\sP^\vee)$
with $\pi:\BP\to W$. On $\BP$ one has the ``universal basis'' 
$$
\underline{s}:\bigoplus^r\sO_{\BP}(-1) \>>> \pi^*\sP,
$$
and $\underline{s}$ is an isomorphism outside of an effective divisor $\Delta$ on $\BP$
with 
$$
\sO_{\BP}(\Delta)=\sO_{\BP}(r)\otimes \pi^*\det(\sP).
$$
The universal basis is induced by the tautological map $\bigoplus^r\pi^*\sP^\vee \to \sO_{\BP}(1)$.
The latter gives a surjection
$$
\bigoplus^r\pi^*\Big(\bigwedge^{r-1}\sP\Big)\cong\bigoplus^r\pi^*(\sP^\vee\otimes \det(\sP))\>>>
 \sO_{\BP}(1)\otimes \pi^*\det(\sP).
$$
Hence $\sO_{\BP}(1)\otimes \pi^*\det(\sP)=\sO_{\BP}(r-1)\otimes\sO_{\BP}(\Delta)$ is weakly positive over $\pi^{-1}(W_0)$.

The sheaf $\sB$ denotes the image of the composite
$$
S^\mu\Big(\bigoplus^r\sO_{\BP}(-1)\Big)=\sO_{\BP}(-\mu)\otimes S^\mu\Big(\bigoplus^r\sO_{\BP}\Big)
\> S^\mu(\underline{s})>> S^\mu(\sP)\> \pi^*({\rm m}) >> \pi^*\sQ.
$$
Remark that $\sB\to \sQ$ is an isomorphism outside $\Delta\cup \pi^{-1}(W\setminus W_0)$.
So there is a modification $\tau:\BP'\to \BP$ with center in this set, such that
$\sB'=\sB/_{\rm torsion}$ is locally free. Writing $\sO_{\BP'}(-\eta)$ for the pullback of
$\sO_\BP(-\eta)$, the surjection
$$
S^\mu\Big(\bigoplus^r\sO_{\BP'}\Big)\>>> \sB'\otimes \sO_{\BP'}(\mu)
$$
defines a morphism to a Grassmann variety $\rho':\BP'\to \G r$. 

The condition on the ``maximal variation'' is used here. One needs that $\rho'$ is quasi-finite on $(\pi\circ\tau)^{-1}(W_0)\setminus \tau^{-1}\Delta$. In the situation considered in Example~\ref{pr.2} this is obviously true. The kernel of ${\rm m}$ determines the fibre $F$ as a subscheme of $\BP(H^0(F,\sL_0|_F))$. So by assumption there are only finitely many $\BP{\rm Gl}(r-1,\C)$ orbits, hence fibres of $\pi|_{\BP\setminus \Delta}$, whose image in $\G r$ can meet.
And obviously $\rho'$ is injective on those fibres.

The Pl\"ucker embedding gives an ample invertible sheaf on $\G r$, and its pullback to
$\BP'$ is $\det(\sB')\otimes \sO_{\BP'}(\gamma)$ with $\gamma=\mu\cdot\rk(\sQ)$. So this sheaf is ample with respect to $(\pi\circ\tau)^{-1}(W_0)\setminus \tau^{-1}\Delta$.

Next, blowing up $\BP'$ a bit more, one can also assume that for some $\nu>0$ and for some divisor
$E$, supported in $\tau^{-1}(\Delta)$ the sheaf 
$$
\det(\tau^*\pi^*\sQ)^\nu\otimes \sO_{\BP'}(\gamma\cdot\nu)\otimes \sO_{\BP'}(-E)
$$
is ample with respect to $(\pi\circ\tau)^{-1}(W_0)$. As the pullback of a weakly positive sheaf
$$
\tau^*\pi^*\det(\sP)^{r-1}\otimes \sO_{\BP'}(\tau^*\Delta)
$$
is weakly positive over $(\pi\circ\tau)^{-1}(W_0)$.

Using the equality $\sO_{\BP'}(r)=\tau^*\pi^*\det(\sP)^{-1}\otimes \sO_{\BP'}(\tau^*\Delta')$, one finds that for all $\eta>0$ the sheaf
\begin{multline*}
\tau^*\pi^*(\det(\sQ)^{\nu\cdot\gamma}\otimes \det(\sP)^{\eta\cdot r -\eta}) \otimes
\sO_{\BP'}(\nu\cdot r\cdot \gamma)\otimes \sO_{\BP'}(-r\cdot E + \eta\cdot \tau^*\Delta)=\\
\tau^*\pi^*(\det(\sQ)^{\nu\cdot\gamma}\otimes \det(\sP)^{\eta\cdot r -\eta-\nu\cdot\gamma}) \otimes
\sO_{\BP'}(-r\cdot E + (\eta+\nu\cdot\gamma)\cdot \tau^*\Delta)
\end{multline*}
is still ample with respect to $(\pi\circ\tau)^{-1}(W_0)$. For $\eta$ sufficiently large the correction
divisor $-r\cdot E + (\eta+\nu\cdot\gamma)\cdot \tau^*\Delta$ will be effective. So we found some effective divisor
$\Delta''$, supported in $\tau^{-1}(\Delta)$ and $a,b >0$ such that
$$
\tau^*\pi^*(\det(\sQ)^{a}\otimes \det(\sP)^{b}) \otimes \sO_{\BP'}(\Delta'')
$$
is ample with respect to $(\pi\circ\tau)^{-1}(W_0)$. 

Next, by \cite[Lemma 4.29]{Vie}, for all $c>0$  one has a natural splitting 
\begin{equation}\label{eqpr.1}
\sO_W \>>> (\pi\circ\tau)_*\sO_{\BP'}(c\cdot\Delta'')\>>> \sO_W ,
\end{equation}
compatible with pullbacks. As in \cite[Proposition 4.30]{Vie} this implies that ``ampleness with respect to $(\pi\circ\tau)^{-1}W_0$ descends from $\BP'$ to $W$:\\[.2cm]
Let us write $\sN=\det(\sQ)^{a}\otimes \det(\sP)^{b}$. Consider two points
$w$ and $w'$ in $W_0$ and $T = w \cup w'$. Let $\BP'_T$ be the proper transform of
$\pi^{-1} (T)$ in $\BP'$. The splitting (\ref{eqpr.1}) gives a commutative diagram
$$
\xymatrix{H^0 (\BP' ,\tau^*\pi^*\sN^{\nu}\otimes \sO_{\BP'}(\nu\cdot\Delta'') ) \ar[r]\ar[d]_{\chi'}&
H^0 (W, \sN^\nu)\ar[d]^{\chi}\\
H^0 (\BP'_T ,\tau^*\pi^*(\sN^{\nu}\otimes \sO_{\BP'}(\nu\cdot\Delta''))|_{\BP'_T} ) \ar[r]&
H^0 (T,\sN^{\nu} |_T )}
$$
with surjective horizontal maps. For some $\nu \geq \nu (w,w')$
the map $\chi'$ and hence $\chi$ will be surjective. For
those $\nu$ the sheaf $\sN^{\nu}$ is generated in a neighborhood of $w'$
by global sections $t$, with $t(w)=0$. By Noetherian
induction one finds some $\nu_0 >0$ such that, for $\nu \geq \nu_0$,
the sheaf $\sN^{\nu}$ is generated by global sections $t_1, \ldots ,t_r$,
on $W_0\setminus \{w\}$ with $t_1(w) = \cdots = t_r(w)=0$, and moreover there is a global section
$t_0$ with $t_0(w) \neq 0$. For the subspace $V_{\nu}$ of $H^0 (W, \sN^{\nu})$, generated by $t_0, \ldots ,t_r$, the morphism $g_{\nu} : W \to \BP (V_{\nu})$ is quasi-finite in a neighborhood of $g^{-1}_{\nu} (g_{\nu} (w))$. In fact, $g^{-1}_{\nu} (g_{\nu} (w))\cap W_0$ is equal to $w$.

Again by Noetherian induction one finds some $\nu_1$ and for
$\nu \geq \nu_1$ some subspace $V_{\nu}$ such that the restriction of $g_{\nu}$ to $W_0$ is
quasi-finite. Then $g^{*}_{\nu} \sO_{\BP (V_{\nu})} (1) =\sN^{\nu}$ is ample with respect to $W_0$.
\end{proof}
We keep the notations introduced in Section~\ref{co}. For Theorem~\ref{in.5} we consider the moduli functor
$\fM_h$ of canonically polarized manifolds (Case CP). As shown in Lemma~\ref{co.3} (2), for Theorem~\ref{in.6} 
it is sufficient to consider the moduli functor $\fM^{(\nu)}_{h}$ of minimal manifolds $F$ with $\omega_F^{\nu}=\sO_F$, and with a very ample polarization $\sL_F$ without higher cohomology (Case PO). 
As in the construction of $M_h$ or $M_h^{(\nu)}$ in Section~\ref{co} we will by abuse of notations consider
$M_h$ and $M^{(\nu)}_h$ with their reduced structure. 

In general $M_h$ is not a fine moduli space, hence there is no universal family.
However Seshadri's Theorem on the elimination of finite isotropies, recalled in \cite[Theorem 3.49]{Vie},
provides us with a finite normal covering $\phi_0:Y_0\to M_h$ which factors over the moduli stack, i.e. which is induced by a family $f_0:X_0\to Y_0$ (or by $(f_0:X_0\to Y_0,\sL_0)$). 
So we are in the situation considered in Variant~\ref{ex.10}, and for each rigidified determinant sheaf, as defined in Definition~\ref{ex.9}, we can find $\overline{M}_h$ and $\phi:W \to \overline{M}_h$ such that $\sC_{\overline{M}_h}$ exists. Recall that its pullback is the $p$-th tensor power of the given rigidified determinant.

We apply~\ref{ex.10} to $\det(\sF_\bullet^{(\nu)})$ and we
obtain a morphism $\phi:Y\to \overline{M}_h$. The corresponding sheaf $\sC_{\overline{M}_h}$ is just the sheaf $\lambda_\nu^{(p)}$ in Theorem~\ref{in.5} (or $\lambda_\nu^{(p)}$ in Theorem~\ref{in.6}). 
So in order to prove both Theorems, it remains to show:
\begin{enumerate}
\item[($\star$)] The sheaf $\lambda_\nu^{(p)}$ is nef and ample with respect to $M_h$.
\end{enumerate}
To do so, Lemma~\ref{ne.9} allows to replace $\overline{M}_h$ by any finite covering, for example by the normalization of $W$ or by a modification $Y$ of the latter with centers outside the preimage of $M_h$. 

The preimage of $M_h$ in $Y$ maps to $Y_0$, and we may assume that both are equal. So we are exactly in the situation considered in Section~\ref{ws}. 
Replacing $Y$ by some alteration, finite over $Y_0$, we can assume that the mild morphism
$\hat{Z}\to \hat{Y}$ in Proposition~\ref{ws.5} exists over a desingularization $\varphi:\hat{Y}\to Y$ of $Y$,
hence all the morphisms in the diagram (\ref{ws.2}). Moreover we can assume that
the locally free sheaf $\sF_{\hat{Y}}^{(\nu)}$ (invertible for $\fM^{(\nu)}_{h}$) in Theorem~\ref{in.1} exists and that it is the pullback of a locally free sheaf $\sF_{Y}^{(\nu)}$ on $Y$.
So ($\star$) and hence the Theorems~\ref{in.5} and~\ref{in.6} follow from:
\begin{claim}\label{pr.3}
The locally free sheaf $\sF_{Y}^{(\nu)}$ is nef and ample with respect to $Y_0$.
\end{claim}
\begin{proof}[Proof of~\ref{pr.3} in Case CP] Let us fix besides of $\nu$ some $\eta_0$ such that for all $F\in \fM_h({\rm Spec}(\C))$ the sheaf $\omega_F^{\eta_0}$ is very ample. Choose $\eta_1=\beta\cdot\eta_0$ such that the multiplication map
$$
{\rm m}:S^{\beta}(H^0(F,\omega_F^{\eta_0}))\> >> H^0(F,\omega_F^{\eta_1})
$$
is surjective and such that its kernel generates the homogeneous ideal, defining
$F\subset \BP(H^0(F,\omega_F^{\eta_0}))$. By Theorem~\ref{in.1} the sheaves
$\sF_W^{(\eta_0)}$, $\sF_W^{(\eta_1)}$ and $\sF_W^{(\nu)}$ exist on some alteration of $Y$, finite over $Y_0$. So we can replace $Y$ by the normalization of this alteration, and assume that they exist on $Y$ itself. The multiplication of sections defines a morphism
$S^{\beta}(\sF_{\hat{Y}}^{(\eta_0)}) \to \sF_{\hat{Y}}^{(\eta_1)}$,
hence as in Addendum~\ref{ex.7} (5) this is the pullback of 
${\rm m}:S^{\beta}(\sF_{Y}^{(\eta_0)}) \to \sF_{Y}^{(\eta_1)}$.

Both sheaves are locally free and by Theorem~\ref{in.5} iii) they are nef. The kernel of ${\rm m}$ is of maximal variation, as explained in Example~\ref{pr.2}.
By Theorem~\ref{pr.2} iii) one finds that for some positive integers the sheaf
$\det(\sF_{Y}^{(\eta_1)})^a \otimes\det(\sF_{Y}^{(\eta_0)})^b$ is ample with respect to $Y_0$ and by part iv) the same holds for $\sF_Y^{(\nu)}$.
\end{proof}
\begin{proof}[Proof of~\ref{pr.3} in Case PO]
The proof of Theorem~\ref{in.6} is similar. We choose a positive integer $\beta$, divisible by
$s=h(1)$ such that the multiplication map
$$
{\rm m}:S^\beta(H^0(F,\sL|_F)) \>>> H^0(F,\sL^\beta|_F)
$$
is surjective for all $F\in \fM_h^{(\beta)}(\C)$, and such that its kernel defines the homogeneous ideal of the image of
$F$ in $\BP(H^0(F,\sL|_F))$. We choose a natural number $\epsilon$ divisible by $\beta$ and with
$\epsilon > e(\sL^\beta|_F)$. 
 
Since we are allowed to replace $Y$ by some finite covering, we can apply~\ref{sa.13}, Proposition~\ref{ex.8}
and \cite[Lemma 2.1]{Vie} and assume:
\begin{enumerate}
\item The sheaves $(\sM_{\hat{Z}},\sM_Z,\sM_{\hat{X}})$ are $\beta$-saturated. 
\item The invertible sheaf $\lambda=\sF_Y^{(\nu)}$, and the locally free sheaves $\sF_Y^{(0,1)}$ and $\sF_Y^{(0,\beta)}$ exist on $Y$.
\item For $s=\rk(\sF_Y^{(0,1)})$ the sheaf $\det(\sF_Y^{(0,1)})$ is the $s$-th tensor power of an invertible sheaf $\sN$. 
\end{enumerate}
Replacing $(\sM_{\hat{Z}},\sM_Z,\sM_{\hat{X}})$ and $\sF_Y^{(\beta,\mu)}$ by 
$$
(\sM_{\hat{Z}}\otimes \hat{g}^*\varphi^*\sN^{-1},\sM_Z\otimes g^*\varphi^*\sN^{-1},\sM_{\hat{X}}\otimes \hat{f}^*\varphi^*\sN^{-1})\mbox{ \ \ and \ \ }
\sF_Y^{(\beta,\mu)}\otimes \sN^{-\mu}
$$
we can add:
\begin{enumerate}
\item[(4)] $\det(\sF_Y^{(0,1)})=\sO_Y$ and hence $\det(\sF_{\hat{Y}}^{(0,1)})=\sO_{\hat{Y}}$.
\end{enumerate}
\begin{claim}\label{pr.4} The assumptions (1)--(4) imply for all $\epsilon'$ divisible by $\nu$
that:
\begin{enumerate}
\item[(5)] \ \hspace*{2.8cm} \ \hspace*{\fill}$\displaystyle \sF_Y^{(\epsilon'\cdot \beta,\beta)}=\lambda^\frac{\epsilon'\cdot\beta}{\nu}\otimes
\sF_{Y}^{(0,\beta)}.$\hspace*{\fill} \ \ \vspace{.1cm}
\item[(6)] For $r=\rk(\sF_{Y}^{(0,\beta)}))$ \hspace*{\fill}$\displaystyle \det(\sF_Y^{(\epsilon'\cdot \beta,\beta)})=\lambda^\frac{\epsilon'\cdot\beta\cdot r}{\nu} \otimes \det(\sF_{Y}^{(0,\beta)})$\hspace*{\fill} \ \ \ \vspace{.1cm}
\item[(7)] \ \hspace*{2.8cm} \ \hspace*{\fill}$\displaystyle \sF_Y^{(\epsilon',1)}=\lambda^\frac{\epsilon'}{\nu}\otimes
\sF_{Y}^{(0,1)}.$\hspace*{\fill} \ \ \vspace{.1cm}
\item[(8)] For $s=\rk(\sF_{\hat{Y}}^{(0,1)}))$ \hspace*{\fill}$\displaystyle \det(\sF_Y^{(\epsilon',1)})=\lambda^\frac{s\cdot\epsilon'}{\nu}.$\hspace*{2.5cm}
\hspace*{\fill} \ \ 
\end{enumerate}
\end{claim}
\begin{proof} It is sufficient to verify those four equations on $\hat{Y}$.
Let $\Pi_{\hat{X}}^{(\nu)}$ be the divisor with 
$$
\hat{f}^*\sF_{\hat{Y}}^{(\nu)}=\hat{f}^*\hat{f}_*\omega_{\hat{X}/\hat{Y}}^\nu=\omega_{\hat{X}/\hat{Y}}^\nu\otimes \sO_{\hat{X}}(-\Pi_{\hat{X}}^{(\nu)}).
$$
By Lemma~\ref{sa.11} c)
\begin{equation}\label{eqpr.2}
\hat{f}_*(\omega_{\hat{X}/\hat{Y}}^{\epsilon'\cdot\beta} \otimes\sM_{\hat{X}}^\beta)=\lambda^\frac{\epsilon'\cdot\beta}{\nu} \otimes \hat{f}_*\Big(\sM_{\hat{X}}^\beta
\otimes \sO_{\hat{X}}\big(\frac{\epsilon'\cdot\beta}{\nu}\cdot \Pi_{\hat{X}}^{(\beta)}\big)\Big)=
\lambda^\frac{\epsilon'\cdot\beta}{\nu} \otimes \hat{f}_*(\sM_{\hat{X}}^\beta).
\end{equation}
So (5) holds true, and (6) as well. For (7) we apply Lemma~\ref{sa.11} e) saying that the sheaves 
$(\sM_{\hat{Z}},\sM_Z,\sM_{\hat{X}})$ are also $1$-saturated. Then the equality (\ref{eqpr.2}) holds for $\beta$ replaced by $1$. Since $\det(\hat{f}_*\sM_{\hat{X}})=\sO_{\hat{Y}})$ one obtains (8).
\end{proof}
Remark that Claim~\ref{pr.4} implies in particular, that the sheaves $\sF_{Y}^{(\epsilon'\cdot\beta,\beta)}$, and
$\sF_{Y}^{(\epsilon',1)}$ automatically exist, with all the properties asked for in~\ref{ex.8}. 

By Proposition~\ref{pd2.4} we may assume that the sheaves $\sF_{Y}^{(\epsilon,1)}$ and $\sF_{Y}^{(\epsilon\cdot\beta,\beta)}$ are both weakly positive over $Y_0$. Since $Y$ is normal, the multiplication of sections on $\hat{Y}$ is the pullback of a morphism
${m}:S^\beta(\sF_{Y}^{(\epsilon,1)}) \to \sF_{Y}^{(\epsilon\cdot\beta,\beta)}$
It is surjective over $Y_0$ with kernel of maximal variation, as explained in Example~\ref{pr.2}. 
By Theorem~\ref{pr.1}, for some positive integers $a$ and $b$ the sheaf
\begin{equation}\label{eqpr.3}
\det(\sF_{Y}^{(\epsilon,1)})^a\otimes \det(\sF_{Y}^{(\epsilon\cdot\beta,\beta)})^b=
\lambda^\frac{a\cdot s\cdot\epsilon+b\cdot \epsilon\cdot\beta\cdot r}{\nu}\otimes 
\det(\sF_{Y}^{(0,\beta)})^b
\end{equation}
is ample with respect to $Y_0$. Since $\sF_{Y}^{(\epsilon,1)}$ is nef, we can replace $a$ by a larger integer, and assume that $a\cdot s $ is divisible by $b\cdot \beta\cdot r$.
So for $\epsilon'=\epsilon\cdot (\frac{a\cdot s}{b\cdot \beta\cdot r}+1)$  the sheaf in (\ref{eqpr.3}) is of the form $\det(\sF_{Y}^{(\epsilon'\cdot\beta,\beta)})^b$ and~\ref{pd2.4} ii) implies that $\sF_{Y}^{(\epsilon',1)}$ is ample with respect to $Y_0$, hence $\sF_Y^{(\nu)}$ as well.
\end{proof}
%%%%%%%%%%%%%%%%%%%%%%%% References %%%%%%%%%%%%%%%%%
%\bibliographystyle{plain}

\end{document}